\documentstyle[11pt]{article}

\input amssym.def
\input amssym.tex

\newcount\skewfactor
\def\mathunderaccent#1#2 {\let\theaccent#1\skewfactor#2
\mathpalette\putaccentunder}
\def\putaccentunder#1#2{\oalign{$#1#2$\crcr\hidewidth
\vbox to.2ex{\hbox{$#1\skew\skewfactor\theaccent{}$}\vss}\hidewidth}}
\def\name{\mathunderaccent\tilde-3 }
\def\Name{\name}

\def\lcirc{{<\!\!\!\circ}}
\def\cstok#1{\leavevmode\thinspace\mbox{\vrule\vtop{\vbox{\hrule\kern1pt
  \mbox{\vphantom{\tt/}\thinspace{\tt#1}\thinspace}}
    \kern1pt\hrule}\vrule}\thinspace}
\def\deq{{\stackrel{\rm def}{=}}} 
\def\e{{\varepsilon}} 
\def\lan{{\langle}}
\def\ran{{\rangle}}
\def\n{{\noindent}}
\def\ZZ{{\cal Z}}
\def\YY{{\cal Y}}

\newcommand{\Rightleftarrow}{\Leftrightarrow}
\newcommand{\bbr}{{\Bbb R}}

\renewcommand{\lg}{{\rm \ell g\/}}

\renewcommand{\H}{{\bf H}}
\newcommand{\bfcc}{{\bf c}}
\newcommand{\E}{{\bf E}}
\newcommand{\bt}{{\bf t}}
\newcommand{\PP}{{\cal P}}
\newcommand{\DD}{{\cal D}}

\newcommand{\AAA}{{\cal A}}
\newcommand{\cBB}{{\cal B}}

\newcommand{\TT}{{\cal T}}

\newcommand{\II}{{\cal I}}
\newcommand{\JJ}{{\cal J}}
\newcommand{\B}{{\bf B}}

\newcommand{\HH}{{\cal H}}
\newcommand{\C}{{\bf C}}
\newcommand{\bfdd}{{\bf d}}
\newcommand{\aaa}{{\frak a}}

\newcommand{\bbb}{{\frak b}}
\newcommand{\bb}{{\frak b}}

\newcommand{\Ord}{{\rm Ord}}

\newcommand{\Dom}{{\rm dom}}

\newcommand{\lesdot}{\mathrel{\mathord{<}\!\!\raise 0.8
pt\hbox{$\scriptstyle\circ$}}}

\newcommand{\whp}{{{}^\frown\!}}

\title{Special Subsets of ${}^{{\rm cf}(\mu)}\mu$,
Boolean Algebras and Maharam measure Algebras}
\author{{\bf Saharon Shelah}\thanks{\ \ \ Publication no.\ 620.
The research partially supported by  NSF under grant \#144-EF67 and by 
``Basic Research Foundation'' administered by The Israel Academy of Sciences
and Humanities.}\\
Institute of Mathematics, The Hebrew University of Jerusalem\\
Department of Mathematics, Rutgers University\\
Department of Mathematics, University of WI, Madison
}
\date{done: June 12, 1996\\
printed: \today}

\begin{document}
\maketitle

\begin{abstract} 
The original theme of the paper is the existence proof of ``there is
$\bar{\eta}=\langle \eta_\alpha:\,\alpha<\lambda\ran$ which is a
$(\lambda,J)$-sequence for $\bar{I}=\lan I_i:\,i<\delta\ran$, a sequence of
ideals. This can be thought of as in a generalization to Luzin sets and
Sierpinski sets, but for the product $\prod_{i<\delta}\Dom(I_i)$, the
existence proofs are related to pcf.

The second theme is when does a Boolean algebra $\B$ has free caliber
$\lambda$ (i.e. if $X\subseteq\B$ and $|X|=\lambda$, then for some $Y\subseteq
X$ with $|Y|=\lambda$ and $Y$ is independent). We consider it for $\B$ being a
Maharam measure algebra, or $\B$ a (small) product of free Boolean algebras,
and $\kappa$-cc Boolean algebras. A central case $\lambda=(\beth_\omega)^+$ or
more generally, $\lambda=\mu^+$ for $\mu$ strong limit singular of ``small"
cofinality. A second one is $\mu=\mu^{<\kappa}<\lambda<2^\mu$; the main case
is $\lambda$ regular but we also have things to say on the singular case.
Lastly, we deal with ultraproducts of Boolean algebras in relation to irr(-)
and s(-) etc.
\end{abstract}
\vfill
\eject

\centerline{\Large\bf Annotated Content}
\bigskip
\n{\bf \S0. Introduction.}
\smallskip

\n{\bf \S1. The framework and an illustration}
\smallskip

We define when ``$\bar\eta=\lan\eta_\alpha:\alpha <\lambda\ran$ is a $(\lambda
,I,J)$-sequence for $\bar I=\lan I_i:i<\delta\ran$, which means ($I=J^{\rm
bd}_\lambda$ for simplicity) that each $\eta_\alpha\in \prod_{i<\delta}{\rm
Dom}(I_i)$ and that for $\bar A= \lan A_i:\,i<\delta\ran\in\prod_{i<\delta}
({\rm Dom} I_i)$ for all large enough $\alpha$, $\eta_\alpha$ ``run away''
from $\bar A$ i.e.\ for the $J$-majority of $i<\delta$, $\eta_\alpha
(i)\not\in A_i$.  We give the easy existence (if $I_i$ is
$\kappa_i$-complete), and $\le|{\rm Dom} I_i|$ for all $i<\delta$ and
$\lan\kappa_i:i<\delta\ran$ are strictly increasing converging to a strong
limit (singular) $\mu$, $\mu^+=2^\mu=\lambda$ (1.18). We define normality,
explain how by the existence of such $\bar\eta$, colouring properties can be
lifted (1.6). As an illustration we prove that (the well known result that)
e.g.\ if $\lambda=2^{\beth_\omega}=\beth^+_\omega$, then $\beth^+_\omega$ is
not a free caliber of the Maharam measure algebra (i.e.\ some set $X$ of
$\lambda$ elements, is non-independent (in fact in a more specific way).  For
this we use ideals related to the Erd\"os-Rado Theorem.
\medskip

\n {\bf \S2. There are large free subsets}

Why does the application in \S1 involve $\lambda$ "near'' a strong limit
singular $\mu$ of cofinality $\aleph_0$?  We show that this was necessary: if
$\mu^{\aleph_0}<\lambda\le 2^\lambda$ an cf$(\lambda )$ large enough
($>\beth_2$ OK, $>2^{\aleph_0}$ almost OK, but involves more pcf
considerations), then $\lambda$ is a free caliber of the Maharam measure
algebra.]
\medskip

\n {\bf \S3. Strong independence in Maharam measure algebras}

[We define when ``$\bar\eta$ is a super $(\lambda ,\bar I, \bar J)$-sequence
for $\bar I$''.  The strengthening is that we now can dealing with $n$-tuples
(any $n <\omega$) and prove the easy existence (see 3.1, 3.2). We define for a
set of $\lambda$ intervals in a Boolean algebra variants of independence and
strong negation of it (3.4) and apply it to prove existence of strongly
$\lambda$-anti-independent set in Maharam Measure algebra (3.6), which (by
3.7) suffices for having a subalgebra of dimension $\lambda$ with no
independent set of cardinality $\lambda$.  This completes the consistency part
of the solution of a problem, which was to characterize all cardinals
$\lambda$ which can have this property.  The question was asked for $\lambda
=\aleph_1$ by Haydon and appeared in Fremlin's book \cite{Fre}. Haydon
\cite{Ha1}, \cite{Ha2} and Kunen \cite{Ku81} independently prove it consistent
for $\lambda=\aleph_1$ assuming CH.  The question from \cite{Ha1} and
\cite{Fre} was what happens with $\aleph_1$ under $MA$.  Recently, Plebanek
\cite{Pl1}, \cite{Pl2} proved that under $MA$ all regular cardinals
$\ge\aleph_2$ fail the property, and finally Fremlin \cite{Fre} gave the
negative answer to the original question of Haydon by showing that under $MA$
the property fails for $\aleph_1$.  D\v zamonja and Kunen \cite{DK1},
\cite{DK2} considered the general case (any $\lambda$) and topological
variants.

We prove here e.g.\ if $\lambda =\beth_{\omega +1}=\beth^+_\omega$, then there
is a Hausdorff compact zero dimensional topological space with measure on the
family of the Borel subsets such that it has dimension $\lambda$, so as a
measure space is isomorphic to the Maharam measure space $\cBB (\lambda )$, but
there is no homomorphism from $X$ onto ${}^\mu 2$ (see 3.8).  We finish by
some easy examples.]
\medskip

\n {\bf \S4. The interesting ideals and the direct pcf application}

[We return to our original aim: existence of $\lambda$-sequences for $\bar I$.
In 4.1 we consider some ideals ($J^{\rm bd}_A$, $\prod_{\ell <n} J_\ell$,
$J^{\rm bd}_{\lan\lambda_\ell:\ell<n\ran}=\prod_\ell J^{\rm
bd}_{\lambda_\ell}$, each $\lambda_\ell$ regular, in the cases
$\lambda_\ell<\lambda_{\ell +1}$, $\lambda_\ell>\lambda_{\ell +1}$,
$\lambda_\ell >2^{\lambda_{\ell +1}}$). We point out (4.9) that for
$\bar I=\lan J^{\rm bd}_{\lambda_i}:i<\delta\ran$, if $\lambda=\hbox{tcf}(
\prod_{i\le\delta} \lambda_i/J^{\rm bd}_\delta )$ we get existence directly
from the pcf theory.  We then turn to the case $I_i=\prod_{\ell <n_i} J^{\rm
bd}_{\lambda_{i,\ell}}$, give a sufficient pcf condition for the existence
when $\lan\lambda_{i,\ell}: \ell <n\ran$ is increasing (4.11) and then prove
that this condition occurs not rarely (in 4.13, so if $\lambda=\prod_{i<
\delta}\lambda_i/ J^{\rm bd}_\delta$, $\lambda_i$ increasing, we can ``group
together'' intervals of $\lambda_i$; and the existence of such $\lan\lambda_i:
i<\delta \ran$ is an important theme of pcf theory.]
\medskip

\n {\bf \S5. $\lambda$-sequences for decreasing $\bar\lambda^i$ by pcf}

[We consider cases with $I_i=J^{\rm bd}_{\lan\lambda_{i,\ell}:i<n_i\ran}$,
$\lan\lambda_{i,\ell}:i<n_i\ran$ a decreasing sequence of regulars.  We prove
the existence by using twice cases of true cofinalities, and show that if the
pcf structure is not so simple then there are such cases (e.g.\ $\beth_{
\omega_i+1}>(\beth_{\omega_i}^{+\omega}$).  We concentrate on the case $i<
\delta\Rightarrow n_i=n$, and then indicate the changes needed in the general
case.] 
\medskip

\n {\bf \S6. Products of Boolean Algebras}

[Monk asks about the free caliber of products of $B_i=FBA(\chi_i)=$ the free
algebra with $\chi_i$ generators, for $i<\delta$.  In fact he asks whether
$\lambda =\beth^+_\omega$ is a free caliber of the product or the
$FBA(\beth_n)$ for $n<\omega$.  But we think that the intention must had been
to ask if $\lambda =\hbox{cf}(\lambda )>2^{|\delta |}$ is a free caliber of
$\prod_{i<\delta} B_i$.  Note that this product satisfies the $(2^{|\delta
|})^+$-c.c.  In fact it has cellularity $2^{|\delta |}$, so ``tends to have
free calibers''.  We show that if there is a normal super $(\lambda
,J)$-sequence $\bar\eta$ for appropriate $\bar I=\lan I_n:n<\omega\ran$, then
$\lambda$ is not a free caliber of $\prod_{n<\omega}FBA(|\hbox{Dom }I_n|)$
(see 6.3, 6.3A), so a negative answer is possible.  Now being ``near a strong
limit singular of cofinality $\aleph_0$'' is necessary as a result parallel to
that of \S2 holds (see 6.4).

Though the choice of $\beth_\omega$ was probably just natural as the first
case to consider, actually the product of uncountably many $FBA(\chi_i)$'s
behave differently e.g.\ $\prod_{i<\omega_i}FBA(\beth_i)$ has free caliber
$(\beth_{\omega_1})^+$! (see 6.5).  The proof involves pcf considerations
dealt with in \S7.  We turn to another problem of Monk (\cite[Problem
34]{M2}), this time giving unambivalent solution. If $\kappa$ is weakly
inaccessible with $\lan 2^\mu :\mu <\kappa\ran$ not eventually constant, then
there is a $\kappa$-c.c.\ Boolean algebra of cardinality $2^{<\kappa}$ and no
independent subsets of cardinality $\kappa^+$ (see 6.8, using the existence of
suitable trees). We note that results similar to countable products holds for
the completion of $FBA(\chi )$.

We end by deducing from Gitik Shelah \cite{GiSh:597} complementary consistency
results (so e.g.\ the first question is not answerable in ZFC) and phrase
principles involved, so slightly sharpening the previous results. (See
6.11--6.14.) So together with the earlier part of the section we have answered
\cite[Problems 35, 36]{M2} and \cite[Problems 32, 33]{M2} in the case we are
near a strongly limit singular cardinal.]
\medskip

\n {\bf \S7. A nice subfamily of function exists}

[For completeness we deal with the following: $f_\alpha\in{}^\theta\hbox{Ord}$
for $\alpha <\lambda$ are given, $2^\theta <\lambda =\hbox{cf}(\lambda )$ and
we would like to get approximation to ``for some $X\subseteq\lambda$,
$|X|=\lambda$, $\lan f_\alpha :\alpha\in X\ran$ is a $\Delta$-system,
continuing \cite[Claim 6.6D]{Sh:430}. We phrase a special case (7.1) and deal
with some variants.]
\medskip

\n {\bf \S8. Consistency of ``$\PP(\omega_1)$ has a free caliber and
discussion of pcf''}

[We deal with another of Monk's problems, \cite[Problem 37]{M2}, proving the
consistency of ``there is no complete Boolean algebra $\B$ of cardinality
$2^{\aleph_1}$ with empty free caliber (in fact $\aleph_{\omega_1+1}=
2^{\aleph_1}$ is always a free caliber of $\B$''). The universe is obtained by
adding $\aleph_{\omega_1+1}$ Cohens to a model of ZFC $+$ GCH, and the proof
uses \S7.  We finish by discussing some pcf problems.]
\medskip

\n {\bf \S9. Having a $\lambda$-sequence for a sequence of non-stationary
ideals}

[We return to the original theme, for a more restricted case.  We assume
$\lambda =\hbox{cf}(2^\mu )$ where $\mu$ is strong limit singular, and in this
section $\lambda =2^\mu$ i.e.\ $2^\mu$ is regular (for the singular case see
\S10).  We get quite strong results: (fix $n(*)<\omega$ for simplicity) for
some ideal $J$ on $\hbox{cf}(\mu )$ (usually $J^{\rm bd}_{{\rm cf}(\mu)}$,
always close to it) we can find $\lan\bar\lambda^i:i<\hbox{cf}(\mu )\ran$,
$i<j\Rightarrow\max (\bar\lambda^i)<\min (\bar\lambda^j)$,
$\bar\lambda^i=\lan\lambda_{i,\ell}:\ell <n(*)\ran$, $\lambda_{i,\ell
+1}>2^{\lambda_{i,\ell}}$ ($\lambda_{i,\ell}$ regular of course,
$\mu=\sup_{i<{\rm cf}(\mu )}\lambda_{i,0}$), such that there is a $(\lambda
,J)$-sequence for $\bar I=\lan J^{\rm bd}_{\prod\bar\lambda^i}:i<{\rm cf}(\mu
)\ran$.  This is nice (compare with \S5) but we get much more: $\bar I$ is a
sequence of nonstationary ideals and even $\lan\prod_{\ell <n(*)} J^{{\rm
nst},\sigma}_{\lambda_{i+\ell}}:i<{\rm cf }(\mu)\ran$ where $J^{{\rm 
nst},\sigma}_\chi=\{A:A\cap\{\delta <\chi :{\rm cf}(\delta )=\sigma\}$ is not
stationary$\}$ and $\sigma ={\rm cf}(\sigma )\in ({\rm cf}(\mu ),\mu)$.

We then work more and get versions with club guessing ideals.  We deal at
length with the version we get for the case ${\rm cf}(\mu )=\aleph_0$. (So it
is less clear which ideals $J$ can be used.)
\medskip

\n {\bf \S10. The power of a strong limit singular is itself singular:
existence} 

[We do the parallel of the first theorem of \S9 in the case $2^\mu$ is
singular.]
\medskip

\n {\bf \S11. Preliminaries to the construction of ccc Boolean algebras with
no large independent sets}

[Here the problem is whether every $\kappa$-c.c.\ Boolean algebra has free
caliber $\lambda$; the case of being ``near a strong limit singular $\mu$ of
cofinality $<\kappa$'' was considered in \cite{Sh:575}, we deal with the case
$\mu =\mu^{<\kappa}<\lambda <2^\mu$.  Here we make the set theoretic
preparation for a proof of the consistency of a negative answer with strong
violation of GCH.  We use Boolean algebras generated by $x_\alpha$'s freely
except $x_\alpha\cap x_\beta\cap x_\gamma =0$ for $\{\alpha ,\beta
,\gamma\}\in W$ for some set $W$ of triples with intersection having at most
one element.  The point is that the properties of ``$\bar\eta$ is a
$\lambda$-sequence for $\bar I$'' with such ideals $I$ (unlike the ones
associated with the Erd\"os-Rado theorem) are preserved by adding many Cohens
to $\mu$, (where $\mu\ll |{\rm Dom }(I_i) |$ etc.)]
\medskip

\n {\bf \S12. Constructing ccc Boolean algebras with no large independent sets}

[We complete the consistency results for which the ground was prepared in
\S11.  We construct the relevant Boolean algebra using a
$(\lambda,J)$-sequence for $\bar I$, $\bar I$ as there, using as building
blocks Boolean algebras generated e.g.\ from the triple system.  So we will
give sufficient conditions for the $\kappa$-c.c. and other properties of the
Boolean algebra.]  
\medskip

\n {\bf \S13. The singular case}

[We continue \S11, \S12 by dealing here with the case $\lambda$ is singular
but $(\forall\alpha <\lambda )$ $(|\alpha |^{<\kappa}<\lambda )$, note that
the forcing from \S12 essentially creates only such cases.]
\medskip

\n {\bf \S14. Getting free caliber for regular cardinals}

[We continue dealing with $\kappa$-c.c.\ Boolean algebras, giving a sufficient
condition for $\lambda$ being a free caliber, hence a consistency follows
(complementing \S11 and \S12; together this solves \cite[Problems 32, 33]{M2}
in the case we are not near a strong limit singular cardinal; thus
together with \S6 this gives a solution).]

\medskip

\n {\bf \S15. On irr: The invariant of the ultraproduct, greater than the
ultraproduct of invariants}

[We prove the consistency of ${\rm irr}\big(\prod_{n<\omega}\B_n/\DD\big)>
\prod_{n<\omega} {\rm irr}(\B_n)/\DD$ where $\DD$ is a nonprincipal
ultrafilter on $\omega$ and ${\rm irr}(\B )={\rm irr}_\omega (\B )$ and ${\rm
irr}_n(\B)=\sup\left\{ X:X\subset\B\right.$ and if $x_0, x_1,\ldots ,x_m$ are
distinct members of $X$, $m<n$ then $\left. x_0\notin \lan x_1,\ldots
,x_n\ran_\B\right\}$.  The way is to build $\B_n$ with ${\rm irr}_n(\B_n)=
\lambda^+$, ${\rm irr}_i(\B_n)=\lambda$, $\lambda=\lambda^{\aleph_0}$.  Our
earlier tries fail as the approximation to $\B_n$ did not work. So the point
is a version of $n$-graded independence phrased as $\lan F_\ell :\ell
<n\ran$, then solve [M2, Problem 26].
We then deal with $s(\; )$, hL$(-)$, hd$(-)$ and Length$(-)$, using the
construction of \S12 in ZFC, and solving \cite[Problems 22,  46,
51, 55]{M2}.] 
\medskip

\section{Introduction}

Our original aim was to construct special subsets of 
$\displaystyle{\prod_{i<\delta} \lambda_i}$, concentrating particularly on the
case when $\lambda_i$ converge to a strong limit singular.

This continues \cite{Sh:575} (so \cite{Sh:g}, \cite{Sh:462}, Ros{\l}anowski and
Shelah \cite{RoSh:534}), but as these are essentially notes from the author's
lectures in Madison, they are self contained. (\S1, \S4 just represent old
material, adding an illustration for Maharam algebras).

Some sections improve the general existence theorems. The main new point is
the case when we use
$$
I_i=\prod_{\ell <n_i}J^{{\rm bd}}_{\lambda_{\ell,i}}\quad \mbox{with the 
$\lambda_{\ell,i}$'s are regular decreasing}
$$
(as well as the case of the nonstationary ideal). We shall discuss this below
and give the definition after we first fix some notation. \medskip

\n {\bf Notation}:
\smallskip

1. $I$ denotes an ideal on a set $\mbox{Dom}(I)$, which means that $I$ is a
subset of $\PP(\mbox{Dom}(I))$ closed under (finite) unions and subsets,
$\hbox{Dom}(I)\not\in I$, and usually for simplicity, all singletons are
assumed to belong to $I$.
\smallskip

$I$ is $\kappa$-complete if it is closed under unions of $<\kappa$ elements.
\smallskip

2. $I,J$ denote ideals.
\smallskip

3. $I^+=\{ A\subseteq\mbox{ Dom}(I): A\not\in I\}$.
\smallskip

4. If $A$ is a set of ordinals with no last member we let
$$
J_A^{{\rm bd}}=\{B\subseteq A: B\;\mbox{ a bounded subset of $A$}\}.
$$

5.  The completeness of the ideal $I$, $\mbox{comp}(I)$ is the maximal
$\kappa$ such that $I$ is $\kappa$-complete (it is necessarily a well-defined
regular cardinal).
\smallskip

6. $[A]^\kappa =\{a\subseteq A:|a|=\kappa\}$, $[A]^{<\kappa}=\{ a\subseteq
A:|a|<\kappa\}$ etc.
\smallskip

7. ${\rm cov}(\lambda ,\mu ,\theta, \sigma )={\rm Min}\big\{|\PP|:\PP\subseteq
[\lambda ]^{<\mu}$, and for every $a\in [\lambda ]^{<\theta}$ there are
$\alpha <\sigma$ and $a_i\in\PP$ for $i<\alpha$ such that
$a\subseteq\bigcup_{i<\alpha}a_i\big\}$.  \medskip
\medskip

\centerline{* * * * * * * *}
\medskip

\n {\bf Definition}:
\smallskip
{\it We say $\bar\eta =\langle\eta_\alpha :\alpha <\lambda\ran$ is a $(\lambda
, I, J)$-sequence for $\bar I=\lan I_i:i<\delta\ran$ if
\begin{enumerate}
\item[a)] $I$ is an ideal on $\lambda$ (if not mentioned, we assume
$I=J^{{\rm bd}}_\lambda$).

\item[{}] $I_i$ is an ideal on $\mbox{Dom}(I_i)$.

\item[b)] $J$ is an ideal on $\delta$ (if not mentioned, we assume $J=J^{{\rm
bd}}_\delta$).

\item[c)] $\eta_\alpha\in \Pi_{i<\delta}{\mbox{Dom}}(I_i)$

\item[d)] If $X\in I^+$ $\underline{\rm then}$
$$
\left\{ i<\delta :\{\eta_\alpha (i):\alpha\in X\}\in I_i\right\}\in J.
$$
\end{enumerate}}

By \cite{Sh:575}, if $I_i$ is $\kappa_i$-complete, $\kappa_i>\sum_{j<i}
\kappa_j$, $\mu =\sum_{i<\delta}\kappa_i$ strong limit,
$|\mbox{Dom}(I_i)|<\mu$ and $2^\mu =\mu^+ =\lambda$, then there is such a 
sequence. We recall this in \S1.

As an example of the application of such $\bar\eta$, we presented the
following (presented in 1.13): Suppose that $\cBB$ is a Maharam measure
algebra of dimension $\ge\mu$, ${\rm cf}(\mu )=\aleph_0$.  Then we can find
$a_\alpha\in\cBB$ for $\alpha <\lambda$ such that Leb$(a_\alpha )>0$ and
$$
(\forall X\in [\lambda ]^\lambda )(\exists n)(\forall \alpha_0<\cdots
<\alpha_n\in X) \bigcap_{i\le n} a_{\alpha_i}=0.
$$
A ``neighborhood'' of $\mu$ being strong limit of cofinality $\aleph_0$ is
necessary.

\n Our usual case, which we call normal is: $\kappa_i>\prod_{j<i}
|\mbox{Dom}(I_j)|$ (this was not used in the measure algebra application, but
it is still good to have).
\smallskip

\n {\bf Main point}: The main new point of this paper is to build a $(\lambda,
I, J)$-sequence $\bar\eta$ for certain $\bar I$ {\it without using $2^\mu
=\mu^+$}. We describe the cases of $\bar{I}$ which we can handle.
\medskip

\n {\bf Case 1}: The easiest case of $I_i:I_i=J^{{\rm bd}}_{\lambda_i}$,
$\lambda ={\rm cf}\Big(\prod_{i<\delta}\lambda_i/J\Big)$.  We only need to
translate from the known pcf results.
\medskip

\n {\bf Case 2}:
$$
I_i=\prod_{\ell <n_i} J^{{\rm bd}}_{\lambda_{\ell ,i}},
$$
where $\lambda_{\ell ,i}$ are regular {\it increasing\/} with $\ell$ and $i$,
and $J$ is an ideal on $\{ (i,\ell ):i<\delta,\ell<n_i\}$ such that
$$
(\forall X\in J)(\exists^{ (J^{{\rm bd}}_\delta )^+}i)(\bigwedge_{\ell <n_i}
(i,\ell)\not\in X),
$$
and where for ideals $J_m$ $(m<n)$
$$ 
\prod_{m<n}J_m{\stackrel{\rm def}{=}}\left\{
\begin{array}{r}
X\subseteq\times_{m<n}\mbox{Dom}(J^m):\neg\exists^{J^+_0} x_0\exists^{J^+_1}
x_1\cdots \exists^{J^+_{n-1}}x_{n-1}\\
(\lan x_0,\ldots ,x_{n-1}\ran\in X)
\end{array}\right\}
$$
Starting from reasonable pcf assumptions and working a little, we can handle
this case as well.
\medskip

\n {\bf Main Case 3}:
$$
I_i=\prod_{\ell <n_i} J^{{\rm bd}}_{\lambda_{\ell ,i}},
$$
$\lambda_{\ell ,i}$ regular {\em decreasing} with $\ell$.
\medskip

We prove: $\underline{\rm If}$ $\bigwedge_in_i=n$, and $\lambda_\ell
=\mbox{tcf}\big(\prod_{i<\delta}\lambda_{\ell ,i}/J'\big)$ for $\ell <n$,
$\underline{\rm then}$ we can find $\lan\eta_{\bar\alpha}:\bar\alpha
<\prod_{i<\delta , \ell <n}\lambda_{\ell, i}\ran$ such that
$$
\mbox{if }\; X\in\Big(\prod_{\ell <n}J^{{\rm bd}}_{\lambda_\ell}\Big)^+
{\stackrel{\rm def}{=}}\left( J^{{\rm bd}}_{\lan\lambda_\ell:\ell<n\ran}
\right)^+
$$
$$
\mbox{then }\;\left\{ i<\delta:\{\eta_{\bar\alpha}(i):\bar\alpha\in X\}\in
\left(J^{{\rm bd}}_{\lan\lambda_{\ell,i}:\ell <n\ran}\right)\right\}\in J'.
$$
\smallskip

\n {\bf Interesting instances}:\quad $\lambda_{\ell ,i}$ decreasing with
$\ell$ and $i<j\Rightarrow\lambda_{i,\ell}<\lambda_{j,n}$.
\medskip

\n {\bf Case 4}: Like Case 3, but using the nonstationary ideal, or
nonstationary ideal restricted so some ``large subset'' of $\lambda_{\ell,i}$
instead of $J^{\rm bd}_{\lambda_{\ell,i}}$.
\medskip

\n {\bf Case 5}: Like Case 3 but using a suitable club guessing ideal $({\rm
id}^a(\overline C^{\ell ,i}))$.  \bigskip
 
On history, background etc.\ and on Boolean algebras, see Monk \cite{M1},
\cite{M2}.  This works continues \cite{Sh:575} and it evolved as
follows. Getting the thesis of Carrieres, which was based on \cite{Sh:92}, we
started thinking again on ``free calibers", this time on measure algebras. We
noted that \cite{Sh:575} gives the answer if e.g.\ $\lambda=(\beth_\omega)^+=
\beth_{\omega+1}$, and started to think of what is called here ``there is a
$(\lambda,J)$-sequence for $\bar{I}$''. We started to lecture on it (\S1, \S4,
then \S5,\S9,\S10).  Meanwhile Mirna D\v zamonja asked me doesn't this solve a
problem from her thesis. This was not actually the case, but it became so in
\S3. Then she similarly brought me p. 256 of Monk \cite{M2} and this
influenced the most of the rest of the paper, while later I also looked at
pages 255, 257 of \cite{M2}, but not so carefully. Lastly, \S15 is looking
back at the problems from \cite{RoSh:534}. Some of the sections are (revisions
of) notes from my lectures. So I would like to thank Christian Carrieres,
Donald Monk and the participants of the seminar in Madison for their
influence, and mainly Mirna D\v zamonja for god-mothering this paper in many
ways, and last but not the least Diane Reppert for typing the paper, and even
more for correcting and correcting and to David Fremlin who lately informs me
that 1.13 was well known and 3.7, 3.12 have already appeared in Plebanek
\cite{Pl1}, \cite{Pl2}.   

\section{The framework and an illustration}

We are considering a sequence $\lan I_i:i<\delta\ran$ of ideals, and we would
like to find a sequence $\bar\eta =\lan\eta_\alpha :\alpha <\lambda\ran$ of
members of $\displaystyle{\prod_{i<\delta}}\mbox{Dom}(I_i)$ which ``runs
away'' from $\bar A=\lan A_i:i<\delta\ran$ when $A_i\in I_i$ (see definition
1.1 below).

When $I_i$ is $\kappa_i$-complete, $\kappa_i>\displaystyle{\Pi_{j<i}}|
\mbox{Dom}(I_j)|$, $\mu=\displaystyle{\sum_{i<\delta}}\kappa_i$ strong limit
singular, $\lambda=\mu^+=2^\mu$, this is easy. We present this (all from
\cite{Sh:575}) and, for illustration, an example.
\medskip

\n {\bf 1.1 Definition}. {\it 1. We say that $\bar\eta$ is a $(\lambda,I,
J)$-sequence for $\bar I$ if: 
\smallskip
\begin{enumerate}
\item[a)]  $J$ is an ideal on $\delta$ and $I$ is an ideal on $\lambda$.
\smallskip

\item[b)] $\bar I=\lan I_i:i<\delta\ran$ where
\item[{}] $I_i$ is an ideal on Dom$(I_i)$.
\smallskip

\item[c)] $\bar\eta =\lan\eta_\alpha :\alpha<\lambda\ran$ where
$\eta_\alpha\in\displaystyle{\prod_{i<\delta}}\mbox{Dom}(I_i)$.
\smallskip
\item[d)] If $X\in I^+$ then
$$
\left\{ i<\delta :\{\eta_\alpha (i):\alpha \in X\}\in I_i\right\}\in J.
$$
\end{enumerate}

2. We say $\bar\eta$ is a weakly $(\lambda, I, J)$-sequence for $\bar I$ if we
weaken clause (d) to
\smallskip
\begin{enumerate} 
\item[d$^-$)] if $X\in I^+$ then
$$
\left\{ i<\delta :\{\eta_\alpha (i):\alpha\in X\}\in I^+_i\right\}\in
J^+.
$$
\end{enumerate}

3.  We may omit $J$ if $J=J^{{\rm bd}}_\delta$, we may omit $I$ if $I=J^{{\rm
bd}}_\lambda$, and then we may say ``$\bar\eta$ is a $\lambda$-sequence for
$\bar I$.''}

We can replace $\lambda$ by another index set.
\medskip

\n {\bf 1.2 Definition}: {\it 1. We say $\bar\eta$ is normally a $(\lambda
,I,J)$-sequence for $\bar I$ (or in short, ``$\bar\eta$ is normal", when $\bar
I, I, J$ are clear) if:
\begin{enumerate}
\item[($*$)] for every $i<\delta$,
$$
\mbox{comp}(I_i)>|\{ \eta_\alpha\restriction i:\alpha<\lambda\}|.
$$
\end{enumerate}

2.  We say $\bar I=\lan I_i:i<\lambda\ran$ is normal if}
$$
\mbox{comp}(I_i)>\prod_{j<i} |\mbox{Dom}(I_i)|.
$$

\medskip

\n {\bf 1.3 Claim}: If $\bar I=\lan I_i:i<\delta\ran$ is normal and $\bar\eta$
is a $(\lambda , I, J)$-sequence for $\bar I$ {\it then\/} $\bar\eta$ is a
normal (i.e.\ normally a $(\lambda ,I,J)$-sequence for $\bar I$).
\medskip

\n {\bf Proof}:  As for each $i<\delta$
$$
|\{\eta_\alpha\restriction i:\alpha<\lambda\}|\le |\prod_{j<i}\mbox{Dom}(I_j)|
=\prod_{j<i}|\mbox{Dom}(I_j)|<\mbox{comp}(I_i).
$$
\medskip

\n {\bf 1.4 Discussion}: Why is normality (and $\eta$-sequences in general) of
interest?  Think for example, of having for each $i<\delta$, a colouring
$\bfcc_i$, say a function with domain $[\mbox{Dom}(I_i)]^2$ (or even
$[\mbox{Dom}(I_i)]^{<\aleph_0}$), call its range the set of colours.  These
colourings are assumed to satisfy ``for every $X\in I^+_i$ we can find some
$Y\subseteq X$ with $Y\in I^+$, such that $\bfcc_i\restriction [Y]^2$ (or
$[Y]^{<\aleph_0}$) is of some constant pattern''.  Now using $\bar\eta$ we can
define a colouring $\bfcc$ on $[\lambda ]^2$ (or $[\lambda ]^{<\aleph_0}$)
``induced by the $\lan\bfcc_i:i<\delta\ran$'', e.g.
$$\begin{array}{rl}
\bfcc(\{\alpha ,\beta\})&=\bfcc_{i(\alpha,\beta)}(\{\eta_\alpha(i(\alpha,
\beta)),\eta_\beta(i(\alpha ,\beta ))\})\cr
&\cr
&\mbox{where $i(\alpha ,\beta )=\mbox{Min}\{ i:\eta_\alpha(i)\not= \eta_\beta
(i)\}$}.\end{array} 
$$

Now, normality (or weak normality) is a natural assumption, because of the
following:
\smallskip

\n {\bf 1.5 Claim}: If $\bar\eta$ is a normally $(\lambda, I,J)$-sequence for
$\bar I$, (or weakly so) and $X\in I^+$, $\underline{\rm then}$ the following
set is $=\delta\mbox{ mod }J$ (or $\not=\emptyset\mbox{ mod }J$):
$$\begin{array}{rl}
Y=\big\{ i<\delta :&\mbox{for some
$\nu\in\displaystyle{\prod_{j<i}}\mbox{Dom}(I_j)$ and $X_i\in I^+_i$}\cr
&\mbox{we have: $(\forall x\in X_i)(\exists\alpha\in X)[\nu =\eta_\alpha
\restriction i$ \& $x=\eta_\alpha(i)]$}\big\}.\end{array}
$$
\smallskip

\n {\bf Proof}.  Let $X_i=\{ \eta_\alpha (i):\alpha\in X\}$, by the
definitions it is enough to prove
\smallskip

\begin{enumerate}
\item[($*$)] if $X_i\in I^+_i$ then $i\in Y$.
\end{enumerate}

Let $Z_i=\{\eta_\alpha\restriction i:\alpha <\lambda\}$, so
$Z_i\subseteq\displaystyle{\prod_{j<i}}\mbox{Dom}(I_j)$ and
$|Z_i|<\hbox{comp}(I_i)$ by the normality of $\bar{I}$ and Claim 1.3. Now for
each $\nu\in Z$ let us define
$$
X^i_\nu =\{\eta_\alpha (i):\alpha\in X\;\;\mbox{ and }\;\;\eta_\alpha
\restriction i=\nu\}.
$$

Clearly $X_i=\bigcup\{ X^i_\nu:\nu\in Z_i\}$, and $I_i$ is $|Z_i|^+$-complete
(as $|Z_i|<\mbox{comp}(I_i)$).  As $X_i\in I^+_i$, necessarily for some
$\nu\in Z$ we have $X^i_\nu\in I^+_i$.  This exemplifies that $i\in Y$, as
required.{\hfill$\square_{\mbox{1.5}}$}
\medskip

\n {\bf 1.6 Conclusion}: Assume
\begin{enumerate}
\item[a)] $\bar\eta$ is a normal weak $(\lambda ,I,J)$-sequence for $\bar I$.
\item[b)] $\bfcc_i$ is a function from ${}^{\omega_>}(\mbox{Dom}(I_i))$ to a
set $\C$ of colours (or from $[\mbox{Dom}(I_i)]^{<\aleph_0}$).
\item[c)] $\bfdd$ is a function from ${}^{\omega >}\varepsilon (*)$ (or from
$[\varepsilon (*)]^{<\aleph_0}$) to $\C$. 
\item[d)] $\bfcc_i$ exemplifies $I_i\not\to (\bfdd )$ which means
\item[{}] ($*$) ~for every $X\in I^+_i$ we can find distinct $x_\zeta\in X$
for $\zeta <\varepsilon (*)$ such that: 

\hskip1.25truein if $n <\omega $ and  $\zeta_0<\cdots<\zeta_{n-1}<\varepsilon
(*)$ then 
$$
\bfcc_i (\lan x_{\zeta_0},\ldots ,x_{\zeta_{n-1}}\ran )=\bfdd
(\lan \zeta_0,\ldots ,\zeta_{n-1}\ran )
$$

\hskip1.25truein (or
$$
\bfcc_i(\{ x_{\zeta_0},\ldots , x_{\zeta_{n-1}}\} )=\bfdd
(\{\zeta_0,\ldots ,\zeta_{n-1}\} )).
$$
\item[e)] We define the colouring $\bfcc$ such that for all $n<\omega$
$$
\bfcc (\lan \alpha_0,\ldots ,\alpha_{n-1}\ran )=\bfcc_i(\lan\eta_{\alpha_0}(i),
\ldots,\eta_{\alpha_{n-1}}(i)\rangle)
$$
(or $\bfcc(\{\alpha_0,\ldots ,\alpha_{n-1}\})=\bfcc_i(\{\eta_{\alpha_0}(i),
\ldots,\eta_{\alpha_{n-1}}(i)\})$
\end{enumerate}

\hskip1.25truein if $\ell <m<n\Rightarrow i=\mbox{Min}\{j<\delta :\eta_{
\alpha_\ell}(j)\not=\eta_{\alpha_m}(j)\}$.
\smallskip

\n $\underline{\rm Then}$ $\bfcc$ exemplifies $I\not\to (\bfdd)$.

[Why? If $X\in I^+$, let $Y$ be the set as in Claim 1.5, hence $Y\in
J^+$. Pick an $i\in Y$, so there is $X_i\in I_i^+$ and $\nu$ exemplifying that
$i\in Y$. Let $\{x_\zeta:\,\zeta<\varepsilon(\ast)\}$ exemplify that
$I_i\not\to (\bfdd)$. For $\zeta<\varepsilon(\ast)$, let $\alpha_\zeta\in X$ be
such that $\eta_{\alpha_\zeta}\restriction i=\nu$ and $\eta_{\alpha_\zeta}(i)=
x_\zeta$. Hence for all $n<\omega$ and
$\zeta_0<\ldots\zeta_{n-1}<\varepsilon(\ast)$ we have 
$$
\bfcc(\langle \alpha_{\zeta_0},\ldots \alpha_{\zeta_{n-1}}\rangle)=
\bfcc_i(\langle x_{\zeta_0},\ldots x_{\zeta_{n-1}}\rangle)=
\bfdd(\langle \zeta_0,\ldots \zeta_{n-1}\rangle).
$$
\medskip

\n {\bf 1.6A Comments}. (1) Of course in 1.6 we can restrict ourselves to
colouring of pairs.  Note that the conclusion works for all $\bfdd$'s
simultaneously. Also, additional properties of the $\bfcc_i$'s are
automatically inherited by $\bfcc$, see 1.7 below.

(2) We can also be interested in colours of $n$-tuples, $n>3$, where
$i<\delta$ as in clause (e) of 1.6 does not exist.

(3)  What is the gain in the conclusion?

A reasonable gain is ``catching'' more cardinals, i.e.\ if $I_i=J^{{\rm
bd}}_{\lambda_i}$, $I=J^{{\rm bd}}_\lambda$, then in addition to having an
example for $\lambda_i$ we have one for $\lambda$.  A better gain is when $I$
is simpler than the $I_i$'s.  The best situation is when we essentially can
get $I=J^{{\rm bd}}_\lambda$, $J=J^{{\rm bd}}_\delta$ for all normal $\bar I$
with $\lan |\mbox{Dom}(I_i)|:i<\delta\ran$ increasing with limit $\mu$.
Assuming a case of G.C.H. this is trivially true.

\medskip
\centerline{* * *}
\medskip

Normally we can find many tuples for which there is $i<\delta$ as in clause
(e) of 1.6.

\medskip

\n {\bf 1.7 Fact}: In 1.6 if $\theta =(2^{|\delta |})^+$, or at least $\theta
={\rm cf}(\theta$) \& $(\forall \alpha <\theta )(|\alpha |^{|\delta|}<\theta
)$ then:

\begin{enumerate}
\item[($*$)] for every $X\in [\lambda ]^\theta$, we can find $Y\in [X]^\theta$
and $i<\delta$ and a $1$-to-$1$ function $h$ from $Y$ into $\mbox{Dom}(I_i)$
such that
$$
\bfcc (\lan \alpha_0,\ldots ,\alpha_{n-1}\ran)=\bfcc_i(\lan h(\alpha_0),\ldots
,h(\alpha_{n-1})\ran ) 
$$
for $\alpha_0,\ldots ,\alpha_{n-1}\in Y$ (actually $h(\alpha )=\eta_\alpha
(i)$, where for all $\alpha$ we have $\eta_\alpha\restriction i=\nu$ for some
$\nu\in\displaystyle{\prod_{j<i}}\mbox{Dom}(I_j)$).
\end{enumerate}

\n {\bf Proof}: By the $\Delta$-system lemma applied to
$\left\{\{\eta_\alpha\restriction i:i<\delta\} :\alpha\in X\right\}$.  More
elaborately, let $M\prec (\HH(\chi), \in <^*_\chi)$, where $\chi$ is large
enough, $\{\theta , X,I,J,\bar I, \bar\eta\}\subseteq M$ and
$M^\delta\subseteq M$, while $||M||=\theta$ and $M\cap\theta$ is an ordinal
$<\theta$.  If we choose $\alpha\in X\backslash M$, then we can choose
$i<\delta$ such that $\eta_\alpha\restriction i\in M$,
$\eta_\alpha\restriction (i+1)\not\in M$ (exists as $M^\delta\subseteq
M$). Now notice that for some such $\alpha$ and $i$ the set $Z\deq\{\eta_\beta
(i):\beta \in X, \eta_\beta\restriction i=\eta_\alpha\restriction i\}$ has
cardinality $\theta$. Let $h:Z\to X$ be such that $\gamma\in
Z\Rightarrow\eta_{h(\gamma )}\restriction i=\eta_\alpha\restriction i$ and
$\eta_{h(\gamma )}(i)=\gamma$.  Lastly let $Y={\rm Rang}(h)$.
\medskip

\n {\bf 1.8 Lemma}: Assume
\begin{enumerate}
\item[(a)] $I_i$ is a $\kappa_i$-complete ideal on $\lambda_i$ for $i<\delta$,
and $\delta$ is a limit ordinal
\item[(b)] $\kappa_i={\rm cf}(\kappa_i)>\sum_{j<i} \kappa_j$
\item[(c)] $\mu=\sup_{i<\delta}\kappa_i=\sup_{i<\delta}\lambda_i$
\item[(d)] ${\rm cf}(I_i,\subseteq )\le \mu^+$ (usually in applications it is
$<\mu$ as usually $2^{\lambda_i}<\mu$; the cofinality is that of a partially
ordered set)
\item[(e)] $\lambda =\mu^+=\mu^{|\delta|}$ (so $\lambda =\lambda^{|\delta|}$;
note that $\mu^{|\delta|}\ge \mu^{{\rm cf}(\mu )}\ge\mu^+$ always)
\end{enumerate}
Then some $\bar \eta$ is a $\mu^+$-sequence for $\lan I_i:i<\delta\ran$. 
\medskip

\n {\bf 1.9 Remark}. (1) We shall focus on the case $\mu$ as strong limit
singular, $\delta ={\rm cf}(\mu )$.

(2) We can weaken the requirement $\lambda =\mu^+$, but not now and here.
\medskip

\n {\bf Proof of 1.8}. Let $\YY_i\subseteq I_i$ be cofinal,
$$
|\YY_i|\le\lambda .
$$
So $\big|\prod_{i<\delta}\YY_i\big|\le\lambda^{|\delta|}=\lambda$, and we can
list $\prod_{i<\delta}\YY_i$ as $\left\langle\langle A^\zeta_i:i<\delta\rangle
:\zeta <\lambda\right\rangle$, where $A^\zeta_i\in \YY_i$.

For $\zeta <\lambda$, let $\lan \beta (\zeta ,\varepsilon ):\varepsilon <\mu
)$ list $\{ \beta :\beta <\max\{\mu ,\zeta\}\}$ (or $\{\beta
:\beta\le\zeta\}$).

Now we choose for $\zeta <\lambda$, a function $\eta_\zeta\in\prod_{i<\delta}
\lambda_i$.  Let $\eta_\zeta (i)$ be any member
of
$$
\lambda_i\backslash\bigcup\{ A^{\beta(\zeta,\varepsilon)}_i:\varepsilon
<\sum_{j<i} \kappa_j\}.
$$

\n [Why can we choose such $\eta_\zeta (i)$? Because $A^{\beta (\zeta ,\e
)}_i\in I_i$ and $I_i$ is $\kappa_i$-complete and $\kappa_i
>\sum_{j<i}\kappa_j$].

We claim that $\bar\eta\deq\langle\eta_\zeta:\,\zeta<\lambda \rangle$ is as
required.  Let $X$ be unbounded $\subseteq\lambda$, we need to show $Y$ is
co-bounded in $\delta$, where
$$
Y{\stackrel{\rm def}{=}}\{ i:\{\eta_\alpha (i):\alpha\in X\}\in I_i^+\}.
$$
Let $A^*_i=\{\eta_\alpha (i):\alpha\in X\}$ for every $i\not\in Y$.  Let
$A^*_i{\stackrel{\rm def}{=}}\emptyset$ for $i\in Y$.  Let $A_i\in\YY_i$,
$A_i\supseteq A_i^*$.  Let $\zeta <\lambda$ be such that $\lan A_i:i<
\delta\ran=\lan A^\zeta_i:i<\delta\ran$.  So for every $\alpha\in X\backslash
(\zeta +1)$, for every $i<\delta$ large enough $\eta_\alpha(i)\not\in A_i$.

[Large enough means: Just that letting $\e=\e_{\alpha,\zeta}<\mu$ be such that
$\zeta =\beta (\alpha ,\e )$ and letting $i^*=i^*_{\alpha ,\zeta}$ be such
that $\sum_{j<i^*}\kappa_j>\e$, then $i\in [i^*,\delta )\Rightarrow\eta_\alpha
(i)\notin A_i$].{\hfill $\square_{\rm 1.8}$}
\medskip

\n {\bf 1.10 Example}. $\lambda =\mu^+ =2^\mu$, $\mu$ strong limit of
cofinality $\aleph_0$.  Let $\mu =\sum_{n<\omega}\mu_n$.  Without loss of
generality $\mu_{n+1}>\beth_{n+7}(\mu_n)$.  Let $D_n=[\beth_{n+3}(\mu_n)^+]^n$.
$$
I_n\deq\left\{X\subseteq D_n:\,
\begin{array}{l}
\mbox{there is } h:X\to 2^{\mu_n}\\
\mbox{$\qquad$such that for no infinite }A \mbox{ is }
h\restriction [A]^n \mbox{ constant}
\end{array}
\right\}
$$
\smallskip

\n {\bf 1.10A Fact}.  $I_n$ is an ideal.
\smallskip

\n {\bf 1.10A Fact}. The ideal $I_n$ is not trivial (so $D_n\not\in
I_n$). [Why? By the Erd\"os-Rado Theorem, see 1.14-1.15 for a detailed
explanation.] 
\medskip

\n {\bf 1.12 Fact}. $I_n$ is $\mu^+_n$-complete.  [Why? If $h_i:D_n\to
2^{\mu_n}$ $(i<\mu_n)$, then there is $h:D_n\to 2^{\mu_n}$ such that
$h(x)=h(y)\Rightarrow\bigwedge_i h_i(x)=h_i(y)$].
\medskip

\n {\bf 1.12A Conclusion}. So, By Lemma 1.8, there is $\bar\eta
=\lan\eta_i:i<\mu^+=\lambda\ran$ which is a $\lambda$-sequence for $\langle
I_n:\,n<\omega\rangle$.
\medskip

We apply Conclusion 1.12A to measure algebras getting a well known result:
\medskip

\n {\bf 1.13 Application}. Assume $\lambda =\mu^+$ and $\mu$ is a strong limit
singular of cofinality $\aleph_0$ (i.e.\ as in 1.10). If $\cBB$ is a measure
algebra (Maharam) of dimension $\ge\mu$, we can find $a_\alpha\in \cBB$ for
$\alpha<\lambda$ with Leb$(a_\alpha )>0$ for each $\alpha$, and such that for
every $X\in [\lambda ]^\lambda$ we can find $n^*<\omega$, $\alpha_1,\ldots
,\alpha_{n^*}\in X$ such that
$$
\cBB\vDash \bigcap^{n^*}_{\ell =1} a_{\alpha_i}=0.
$$

\n {\bf Proof}.  Let $\bar\eta$ and $I_n$ be as in conclusion 1.12A (all in
the content of Example 1.10).  Let $\lan x_{n, \alpha}:n<\omega , \alpha
<\beth_{n+3} (\mu_n)^+\ran$ be independent in the sense of measure, all
elements of $\cBB$ and of measure $1/2$.

For any $\eta\in\prod_{n<\omega} D_n$, let
$$
y_{\eta ,n}=y_{\eta (n)}=1-\bigcap_{\beta\in\eta(n)}x_{n,\beta} -
\bigcap_{\beta\notin\eta (n)}(1-x_{n,\beta} ).
$$
Note that $\bigcap_{\beta\in\eta (n)}x_{n,\beta}$ has measure $2^{-n}$.  So
Leb$(y_{\eta ,n})=1-2\cdot 2^{-n}$ (hence Leb$(y_{\eta (n)})>0$ if $n\ge 2$).
Let $$ y_\eta =\bigcap_{n\ge 5}y_{\eta ,n}\in\cBB.
$$
So Leb$(y_\eta )\ge 1-2\cdot\sum_{n\ge 5}2^{-n}=1-2\cdot 2^{-4} =1-2^{-3}>
1/2$.  We let $a_\alpha =y_{\eta_\alpha}$ for $\alpha <\lambda$.  We check
that $\lan a_\alpha :\alpha <\lambda\ran$ is as required. Suppose $X\in
[\lambda ]^\lambda$.  So, as $\bar I$ is normal, for some $n>5$ and
$\nu\in\prod_{\ell <n}D_\ell$ we have
$$
Y_\nu=:\{\eta_\alpha (n):\alpha\in X,\eta_\alpha\restriction n=\nu\}\in
I^+_n. 
$$
(Note that $\nu$ is not really needed for the rest of the proof.)

So there is $\{\gamma_\ell :\ell <\omega\}\subseteq\beth_{n+3}(\mu_n)^+$
increasing such that
$$
[\{\gamma_\ell :\ell <\omega\} ]^n\subseteq Y_\nu .
$$
We use just $\lan\gamma_\ell :\ell <2n-1\ran$.

For $u\in [\{\gamma_\ell :\ell <2n-1\}]^n$ let $\alpha(u)\in X$ be such that 
$$
\eta_{\alpha (u)}(n)=u.
$$
It is enough to show that in $\cBB$
$$
\bigcap_u y_{\eta_{\alpha (u)}}=\bigcap_u a_{\alpha
(u)}=0 .
$$
So suppose that there is $z\in \cBB$ with $\mbox{Leb}(z)>0$ and such that
$z\le\bigcap_u y_{\eta_{\alpha (u)}}$. Then without loss of generality
$$
\ell <2n-1\Rightarrow z\le x_{n ,\gamma_\ell}\vee z\le 1-x_{n,\gamma_\ell}.
$$
\smallskip

\n {\bf Case 1}. $|\{\ell :z\le x_{n,\gamma_\ell}\}|\ge n$.  Let $u\in
[\{\gamma_\ell :\ell <2n-1\} ]^n$ be such that
$$
\bigwedge_{\gamma_\ell\in u}(z \le x_{n ,\gamma_\ell}).
$$
So $z\le \bigcap_{\gamma_{\ell\in u}} x_{n,\gamma_\ell}$.  But $z\le
y_{\eta_{\alpha (u)}}\le 1-\bigcap_{\gamma\in u} x_{n,\gamma}$ contradiction.
\medskip

\n {\bf Case 2}. Not 1. So $|\{\ell :z\le 1-x_{n ,\gamma_\ell}\}|\ge n$ and
continue as above using $1-x_{n,\gamma_\ell}$.

~{\hfill$\square_{\mbox{1.13}}$}
\medskip

Let us elaborate on the ideals used above.
\medskip

\n {\bf 1.14 Definition}: {\it For $n,\lambda ,\e$ let
$$
\mbox{ERJ}_\lambda^{n,\e}= J^{n,\e}_\lambda=
\left\{\begin{array}{l}
 A\subseteq [\lambda ]^n: \mbox{ there is no
$w\subseteq\lambda$,
$\mbox{opt}(w)=\e$}\\
\mbox{\hspace*{.75truein} and }\; [w]^n\subseteq
A\end{array}\right\}
$$
$$
\mbox{ERI}^{n,\e}_{\lambda,\mu}= I^{n,\e}_{\lambda ,\mu}
=\left\{\begin{array}{l}
A\subseteq [\lambda ]^n:\mbox{ there are $A_i\in J^{n
,\e}_\lambda$ for
$i<i(*)<\mu$}\\
\mbox{\hspace*{.75truein}such that
$\displaystyle{A=\bigcup_{i<i(*)}A_i}$.}\end{array}\right\}
$$
}
\smallskip

\n {\bf 1.15 Fact}: 1) $I^{n,\e}_{\lambda ,\mu}$ is a ${\rm cf}(\mu
)$-complete ideal on $[\lambda ]^n$, not necessarily proper (see (2)) .
$J^{n,\e}_\lambda$ is not necessarily an ideal.
\smallskip

2) $I^{n,\e}_{\lambda ,\mu}$ is a proper ideal i.e.\ $[\lambda ]^n\notin
I^{n,\e}_{\lambda ,\mu}$ $\underline{\rm iff}$
$$
\chi <\mu \Rightarrow\lambda\rightarrow (\e )^n_\chi.
$$

3) $I_n=I^{n,\omega }_{\beth_{n+7}(\mu_n)^+, (2^{\mu_n})^+}$ (where $I_n$ and
$\langle \mu_n:\n<\omega\rangle$ are from 1.10).

4) In the proof of 1.13 we could have used less, for example
$$
I_n=I^{n,2n+1}_{\beth_{n+1} (\mu_n)^+, \mu^+_n}
$$
as $\beth_{n+1}(\mu_n)^+\rightarrow (\mu_n^+)^n_{\mu_n}$ for $n\ge 1$. 
\medskip

\n {\bf Proof}. (3) $\underline{\rm First ~direction}$. Let $A\in I_n$, so
there is $h:A\to 2^{\mu_n}$ witnessing it.  Let $A_i=h^{-1}(i)$ for
$i<2^{\mu_n}$ Now $X\subseteq\lambda$, $|X|\ge \aleph_0\Rightarrow
[X]^n\not\subseteq A_i$, by the choice of $A$.

Hence
$$
A_i\in J^{n,\omega}_{\beth _{n+7}(\mu_n)^+}.
$$
Hence
$$
A\in I^{n,\omega}_{\beth_{n+7}(\mu_n)^+, (2^{\mu_n})^+}.
$$

\n $\underline{\rm Second ~direction}$: Let $A\in I^{n,\omega}_{
\beth_{n+7}(\mu_n)^+,(2^{\mu_n})^+}$, so there are $A_i(i<i(*)<
(2^{\mu_n})^+)$ such that $A_i\in J^{n,\omega}_{\beth_{n+7}(\mu_n)^+}$ and
$A=\bigcup_{i<i(*)} A_i$. 

Renaming, without loss of generality $i(*)\le 2^{\mu_n}$, and let
$$
A'_i=\left\{\begin{array}{ll}
\displaystyle{A_i\backslash\bigcup_{j<i} A_i}&\quad\mbox{if
$i<i(*)$}\\
&\\
\emptyset&\quad\mbox{otherwise, i.e.\ if $i\in
[i(*),2^{\lambda_n})$.}\end{array}\right.
$$
So $\lan A'_i:i<i(*)\ran$ is a partition of $A$. As $A_i\in
J^{n,\omega}_{\beth_{n+7}(\mu_n)^+}$, we know that $\neg (\exists
X\subseteq\beth_{n+7}(\mu_n)^+$ infinite) $([X]^n\subseteq A_i)$.  Hence,
letting $\kappa =\beth_{n+7}(\mu_n)^+$

$$
\neg (\exists X\subseteq\kappa\;\mbox{ infinite})\;\; ([X]^n\subseteq A_i]). 
$$
Define
$
h:A\to 2^{\mu_n}$ by
$$
h(\bar\alpha )=i\;\;\underline{\rm iff}\;\;\bar\alpha\in A'_i,
$$
so $h$ witnesses $A\in I_n$.{\hfill$\square_{\rm 1.15}$}
\medskip

\n {\bf 1.16 Definition}: {\it 1) A set $W\subseteq [\lambda ]^{<\aleph_0}$ is
called a ccc base if
$$
\mbox{($*$)}\;\;\mbox{for $u\not= v$ in $W$, $|u\cap v|<|u|/2$.}
$$

2) For $W\subseteq [\lambda ]^{<\aleph_0}$ let
$$\begin{array}{l}
I_{\lambda}[W]=\{ A\subseteq\lambda:W\cap [A]^{<\aleph_0}=\emptyset\}\\
I_{\lambda ,\kappa}[W]=\big\{ A\subseteq\lambda :\mbox{ $A$ is
the union of $<\kappa$ members of $I_\lambda [W]\big\}$.}
\end{array}
$$

3) For a Boolean algebra $\B$ we define $I_{\B ,\kappa}$ by letting $:X\in
I_{\B,\kappa}$ iff $X\subseteq\B\backslash\{ 1\}$ is the union of $<\kappa$
ideals of $\B$.} 
\medskip

\n {\bf 1.17 Claim}: 1) Assume

(a) $\bar\eta$ is a $(\lambda ,J)$-sequence for $\bar I=\lan I_i:i<\delta
\ran$, and ${\rm cf}(\lambda)>\delta$.

(b) for $i<\delta$, the function $h_i:{\rm Dom}(I_i)\to\lambda_i$ satisfies
$$
\alpha <\lambda_i\Rightarrow\{x\in{\rm Dom}(I_i):h_i(x)<\alpha\}\in I_i.
$$
Let $\bar h=\lan h_i:i<\delta\ran$ and let $f_\alpha =\bar h\circ\eta_\alpha
\deq\lan h_i(\eta_\alpha (i)):i<\delta\ran$ $\big(\in\prod_{i<\delta}
\lambda_i\big)$.

\underline{Then}

(c) $(\forall f\in\prod \lambda_i)$ $(\forall^{J_\lambda^{\rm bd}}\gamma
<\lambda)$ $(f<_J f_\gamma )$.

(d) for some club $E$ of $\lambda$, we have

$\quad$(d)$_E$ if $\alpha<\e\le\beta <\lambda$ and $\e\in E$ then $f_\alpha 
<_J f_\beta$ and

(So if $X\in [\lambda]^\lambda$, $(\forall\delta\in E)$ $|X\cap (\delta ,\min
(E\backslash (\delta +1)]|\le 1$ $\underline{\rm then}$ $\lan f_\alpha:\alpha
\in X\ran$ is $<_J$-increasing cofinal in $\prod_{i<\delta}\lambda_i$.
\medskip

2) If $\bar f=\lan f_\alpha :\alpha<\lambda\ran$, $E$ satisfies (d)$_E$ (and
of course $\sup_{i<\delta}\lambda_i<\lambda$) and $\mu <\lambda$
$\underline{\rm then}$ without loss of generality for $X$ as in (d)$_E$ the
sequence $\bar f\restriction X $ is $\mu$-free (see below).  Moreover $\bar f$
is $(\mu ,E)$-free (see below clause (1) of 1.18), if ($*$) or just the weaker
($*$)$'$ or just ($*$)$''$ below holds where:
\medskip

\n {\bf 1.18 Definition}: 
\medskip

(1) $\bar f$ is $\mu$-free if for $X\in [\lambda]^{<\mu}$ we can find $\bar
s=\lan s_\alpha :\alpha\in X\ran$, $s_\alpha\in J$ such that
$$[\alpha<\beta\,\&\,\,\alpha\in X\,\&\,\,\beta\in X\,\&\,\,i\in J\backslash 
s_\alpha\backslash s_\beta]\Rightarrow f_\alpha (i)<f_\beta (i).$$
\smallskip

(2) $\bar f$ is $(\mu ,E)$-free if for $X\in [\lambda]^{<\mu}$ we can find
$\bar s=\lan s_\alpha :\alpha\in X\ran$, $s_\alpha\in J$ such that
$$[\alpha\le\delta<\beta\,\,\&\,\,\alpha\in X\,\&\,\,\delta\in E
\,\&\,\,\beta\in X\,\&\,\,i\in \delta\backslash s_\alpha\backslash s_\beta]
\Rightarrow f_\alpha (i)<f_\beta (i).$$

(3) For $\lambda >\mu$ and $\bar\lambda=\lan\lambda_i:i<\delta\ran$ we
consider the conditions 
\medskip

$\quad$($*$) ~$\lambda=\chi^+$, $\chi={\rm cf}(\chi)\ge \mu$ for some
$\chi$
\medskip

$\quad$($*$)$'$ $\mu=\lim_J\lambda_i$ and $\{\delta <\lambda :{\rm
cf}(f(\delta ))<\mu\}\in I[\lambda ]$
\medskip

$\quad$($*$)$''$ ~there is $\bar f'=\lan f'_\alpha :\alpha <\lambda\ran$
which is $<_J$-increasing cofinal in $(\prod_{i<\delta}\lambda_i, <_J)$ and
is $\mu$-free.  
\medskip

\n {\bf Remark}. This applies to the construction in \S4, \S5, etc.\ (e.g.\
construction from $\lambda =\prod_{i<\delta}\lambda_i/J^{\rm bd}_\delta$).
\bigskip

\section{There are large free subsets}
The reader may wonder if really something like $\lambda ={\rm cf}(\lambda )\in
(\mu ,2^\mu ]$ for $\mu$ strong limit singular, is necessary for 1.13.  As in
\cite{Sh:575}, the answer is yes, though not for the same reason.

Of course, in what follows, Maharam measure algebra can be replaced by any
measure algebra.  The interesting case is $(\exists \chi )(\chi <\lambda
\le\chi^{\aleph_0})$.
\medskip

\n {\bf 2.1 Fact}: Let $\cBB$ be a Maharam measure algebra.  If $\beth_2\le\mu
=\mu^{\aleph_0}<\mbox{ cf}(\lambda )\le \lambda\le 2^\mu$ and $a_j\in\cBB^+$,
(so $\mbox{Leb}(a_j)>0$) for $\alpha<\lambda$ are pairwise distinct,
$\underline{\rm then}$ for some $X\in [\lambda ]^\lambda$ we have:
\begin{enumerate}
\item[($*$)] any nontrivial Boolean combination of finitely many members of
\item[{}] $\{ a_\alpha :\alpha \in X\}$ has positive measure.
\end{enumerate}
\medskip

\n {\bf Proof}: Let $\{ x_i:i<i(*)\}$ be a basis of the Maharam measure
algebra (so each $x_i$ has measure $1/2$ and $x_i$'s are measure-theoretically
independent). So for each $\alpha <\lambda$ we can find ordinals $i(\alpha
,n)<i(*)$ for $n<\omega$, and a Boolean term $\tau_\alpha$ such that $a_\alpha
=\tau_\alpha (x_{i(\alpha ,0)}, x_{i(\alpha ,1)},\ldots )$. Note that this
equality is only modulo the ideal of null sets.
\medskip

Without loss of generality, each $\tau_\alpha$ is a countable intersection of
a countable union of finite Boolean combinations of the $x_\alpha$'s.  Again
without loss of generality, $\lan i(\alpha ,n):n<\omega\ran$ is with no
repetition. Note that without loss of generality
$$
i(*)=\{ i(\alpha , n):\alpha <\lambda\;\;{\rm and}\;\; n<\omega\}.
$$
Hence without loss of generality $i(*)\le\lambda$, hence without loss of
generality $i(*)=\lambda$.  By Engelking Karlowicz Theorem \cite{EK}, we can
divide $\lambda$ to $\mu$ sets $\lan X_\zeta :\zeta <\mu\ran$ such that
\begin{enumerate}
\item[($*$)$_1$] the sets $A_{\zeta ,n}{\stackrel{\rm def}{=}}\{ i(\alpha
,n):\alpha\in X_\zeta\}$ for each $\zeta$ satisfy: $\lan A_{\zeta,n}:n<\omega
\ran$ are pairwise disjoint. 
\end{enumerate}
\smallskip

As the number of possible terms $\tau_\alpha$ is $\le 2^{\aleph_0}\le\mu$,
without loss of generality

\begin{enumerate}
\item[($*$)$_2$] if $\alpha ,\beta\in X_\zeta$ then $\tau_\alpha =\tau_\beta$,
call it $\tau^\zeta$.

Note also

\item[($*$)$_3$] if $Y\subseteq X_\zeta$ then
$$\begin{array}{l}
\mbox{ind}(Y)=:\Big\{\alpha\in Y:\mbox{ for no $m<\omega$ and
$\beta_0,\ldots ,\beta_{m-1}\in Y\cap\alpha$}\\
\mbox{do we have: $a_\alpha\in$ the complete subalgebra generated by}\\
\mbox{$\{x_{i(\beta_\ell ,n)}:\ell <m, n<\omega\}\Big\}$,}
\end{array}
$$
satisfies $|\mbox{ind}(Y)|+2^{\aleph_0}\ge |Y|$.
\end{enumerate}

\n [Why?  We can prove by induction on $\alpha\not\in \mbox{ind}(Y)$ that for
some $m<\omega$ and $\beta_0,\beta_1,\ldots,\beta_{m-1}\not\in\mbox{ind}(Y)
\cap\alpha$ we have $a_\alpha \in$ the complete subalgebra of $\cBB$ generated
by $\{x_{i(\beta_\ell ,n)}:\ell<m,n<\omega\}$, using the transitive character
of this property.  Now for each $m<\omega$ and $\beta_0,\ldots,\beta_{m-1}\in
\mbox{ind}(Y)$, the number of $a_\alpha$ such that $a_\alpha\in$ the
subalgebra generated by $\{ x_{i(\beta_\ell,n)}:\ell<m,n<\omega\}$ is at most
continuum.] 

As ${\rm cf}(\lambda )>\mu$, for at least one $\zeta <\mu$, $|X_\zeta
|=\lambda$, hence by ($*$)$_3$ we have $|\mbox{ind}(X_\zeta )|=\lambda$. So,
without loss of generality

\begin{enumerate}
\item[($*$)$_4$] (a) the sets $A_n=\{i(\alpha ,n):\alpha <\lambda\}$ are
pairwise disjoint, 
\smallskip

\item[{}] (b) $\tau_\alpha =\tau$ for $\alpha <\lambda$,
\smallskip

\item[{}] (c) for no $m<\omega$ and $\beta_0<\cdots<\beta_m<\alpha$ do we
have $a_{\beta_m}\in$ the complete subalgebra generated by $\{ x_{i(\beta_\ell
,n)}:\ell <m, n<\omega\}$.
\end{enumerate}
Now for each $\alpha <\lambda$ we define an ideal $I_\alpha$ on $\omega$: it
is the ideal generated by the sets
$$
Z_{\alpha ,\beta} =:\{ n<\omega :i(\beta ,n)=i(\alpha,n)\}\quad\mbox{for
$\beta <\alpha$.} 
$$
and (where $\mbox{ch}_A(n)$ is $1$ if $n\in A$ and $0$ if $n\notin A$)
\begin{eqnarray*}
J &=&\big\{ A\subseteq\omega :\tau (z_0, z_2, \ldots,z_{2n},\ldots )\\
&=&\tau(z_{0+\mbox{ch}_A(0)},
z_{2+\mbox{ch}_A(1)},\ldots ,z_{2n+\mbox{ch}_A(n)},\ldots)\}.
\end{eqnarray*}
As $\{ x_i:i<i(*)\}$ is free (in the measure theoretic sense), for $A\in J$,
and $\{\alpha_m:\,n<\omega\}$ and $\{\beta_n:\,n<\omega\}$ such that
$\alpha_n<i(*)$ are with no repetition and $\beta_n<i(*)$ with no repetition,
we have the following:

If $(\forall m, n<\omega )[\alpha_n=\beta_m\Rightleftarrow n=m\hbox{ \& }
n\not\in A]$ then $\tau (x_{\alpha_0},\ldots)=\tau (x_{\beta_0},\ldots )$.
(just apply the definition of $J$ to $\langle x_{\alpha_0},x_{\beta_0},
x_{\alpha_1},\ldots\rangle$).  By transitivity of equality we get $(\forall
n<\omega )$ $[n\not\in A\Rightarrow \alpha_n=\beta_n]\Rightarrow \tau
(x_{\alpha_0},\ldots )=\tau (x_{\beta_0},\ldots )$.  Hence $J$ is closed 
under subsets and (finite) unions.  By clause (c) of ($*$)$_4$ we know that
$\omega\notin I_\alpha$; so $I_\alpha$ is an ideal on $\omega$ though it is
possible that singletons are not in $I_\alpha$ (a violation of a convention in
\S0).  [In fact we could have eliminated this violation, but there is no
reason to put extra work for it.] Also $I_\alpha\subseteq J$.

Now, the number of possible ideals on $\omega$ is at most $\beth_2\le \mu
<{\rm cf}(\lambda )$ so it suffices to prove
\begin{enumerate}
\item[($*$)$_5$] if $Y\subseteq\lambda$, $\alpha\in Y\Rightarrow
I_\alpha\subseteq I$, where $I$ is an ideal on $\omega$ (so $\omega\notin I$
but singletons may or may not belong to $I$) extending $J$, $\underline{\rm
then}$ any finite Boolean combination of $\{a_\alpha :\alpha\in Y\}$ has
positive measure.
\end{enumerate}

\n {\bf Proof of ($*$)$_5$}:  Let $\beta_0<\cdots<\beta_{m-1}$ be from $Y$.
Let 
$$
A=\{ n<\omega:\mbox{ for some $\ell < k<m$ we have $i(\beta_\ell, n)=i(
\beta_k, n)$}\}.
$$
By the definition of $Z_{\alpha,\beta}$, clearly $A\in I$.  For $Z\subseteq
i(*)$ let $\cBB^*[Z]$ be the complete subalgebra of $\cBB$ generated by $\{
x_\beta :\beta\in Z\}$. We let $\cBB^*\deq\cBB^\ast [Z]$ if $Z=\{ i(\beta_\ell
,n):\,\ell<m,n\in A\}$.  Let $\cBB^\ast_\ell\deq\cBB^\ast
[\{i(\beta_\ell,n):\,n\in A\}]$.

As $\cBB^*_\ell$ is complete, for each $\ell <m$ we can find $b^-_\ell,
b^+_\ell\in\cBB^*_\ell$ such that
\begin{enumerate}
\item[(i)] $b^-_\ell\le a_{\beta_\ell}\le b^+_\ell$ \smallskip
\item[(ii)] if $c\in\cBB^*_\ell$ then $c\le a_{\beta_\ell}\Rightarrow c\le
b^-_\ell$ and $c\ge a_{\beta_\ell}\Rightarrow c\ge b_\ell^+$.
\end{enumerate}
By the  definition of $\cBB^*$ and  assumptions on $\lan x_i:i<i(*)\ran$ and  
$\lan a_\alpha :\alpha <\lambda\ran$ clearly
\begin{enumerate}
\item[($*$)$_6$] if $\{ i(\beta_\ell ,n):n\in A\}\subseteq Z$ and $\{i(
\beta_\ell,n):n\in\omega\backslash A\}\cap Z=\emptyset$ and $Z\subseteq i(*)$
then   

(ii)$_Z$ if $c\in\cBB^\ast [Z]$, then $c\le a_{\beta_\ell}\Rightarrow c\le
b^-_\ell$ and $c\ge a_{\beta_\ell}\Rightarrow c\ge b_\ell^+$.
\end{enumerate}
Obviously, for some Boolean terms $\tau^-_\ell$, $\tau^+_\ell$ we have
$$\begin{array}{l}
b^-_\ell =\tau^-_\ell (\ldots , x_{i(\beta_\ell ,n)},\ldots)_{n\in A}\\
b^+_\ell =\tau^+_\ell (\ldots , x_{i(\beta_\ell ,n)},\ldots)_{n\in A}.
\end{array}
$$
Now, clearly $\tau^-_\ell =\tau^-$ and $\tau^+_\ell =\tau^+$ for some fixed
$\tau^-$ and $\tau^+$.  Also $b^-_\ell <b^+_\ell$ as otherwise
$\omega\backslash A\in J$.  Let $b_\ell =b^+_\ell -b^-_\ell$ so
$\mbox{Leb}(b_\ell )>0$, and for some term $\tau^*$, $b_\ell =\tau^*(\ldots
,x_{i(\beta_l ,n)},\ldots )_{n\in A}$, and let $b=\bigcap_{\ell <m} b_\ell
\in\cBB^*$.

Clearly
\begin{enumerate}
\item[($*$)$_7$] Leb$(b)>0\Rightarrow$ any Boolean combination of the
$a_{\beta_\ell}(\ell <m)$ have positive measure.
\end{enumerate}
[Why? prove it on $\{ a_{\beta_\ell} :\ell <m'\}$ by induction on $m'\le m$
using ($*$)$_6$.]

For proving $\mbox{Leb}\Big(\bigcap_{\ell <m}b_\ell\Big)>0$, we define an
equivalence relation $E$ on $\omega$:
$$\begin{array}{l}
\mbox{$n_1E\; n_2$ iff for every $\ell <k<m$ we have}\\
\qquad\mbox{$i (\beta_\ell ,n_1)=i (\beta_k,n_1)\Leftrightarrow i(\beta_\ell
,n_2)=i (\beta_k,n_2)$.}\end{array}
$$
Clearly $E$ has finitely many equivalence classes, say $A_0, A_1,\ldots
,A_{k(*)-1}$.  For $k_1\le k(*)$ and $\bar\ell =\lan \ell_k:k_1\le
k<k(\ast)\ran$ let
$$\begin{array}{l}
Z_{k_1,\bar\ell}=\left\{ \tau^*(\ldots ,x_{i(\gamma_n,n)},\ldots
):\mbox{ for every $k<k(\ast)$, for some $\ell <m$}\right.\\
\mbox{we have $\lan \gamma_n:n\in A_k\ran =\lan
i(\beta_n,\ell ):n\in A_k\ran$ but if $k\ge k_1$ then $\ell
=\ell_k$}\big\}.\end{array}
$$
We prove by induction on $k<k(*)$ that for any $\bar\ell$
$$
c_{\bar{\ell}}\deq\mbox{Leb}(\bigcap\{ b:b\in Z_{k,\bar\ell}\})>0.
$$
(In fact the measure does not depend on $\bar\ell$.)

For $k=k(*)$ we have $\{ b_\ell :\ell <m\}\subseteq Z_{k,\lan\;\ran}$ so we
are done.
\smallskip

\n {\it The case $k=0$}:
\smallskip

It is trivial:  $Z_{0,\bar\ell}$ is a singleton $\{\tau^*(\ldots,x_{i(
\gamma_n,n)},\ldots )\}$, where $\gamma_n\in A_n$ so obviously is not zero.
\medskip

\n {\it The case $k_1+1$}:
\smallskip

So let $\bar\ell =\lan\ell_k:k_1+1\le k<k(*)\ran$, and we know that for each
$n<\omega$ the element $d_\ell=c_{(n)\whp\bar\ell}$ is $>0$.  For $\zeta <m$
let $f_\zeta$ be a function from $Y=\big\{ i(\beta_\eta,\ell):\ell<\omega$
such that if $\ell\in A_k$ then $k\in [k_1+1, k(*)]\Rightarrow n=\ell_k$ and
$k=k_1\Rightarrow n=0$ and $k<k_1\Rightarrow n<m\big\}$ into $\lambda$,
$f_\zeta$ is one to one, $f_\zeta$ is the identity on $Y^*=\{ i(\beta_n,\ell
)\in Y:\ell\not\in A_{k_1}\}$ and $\lan{\rm Rang}(f_\zeta\restriction
(Y:Y^*_0)):\zeta <m\ran$ are pairwise disjoint and
$$
\ell \in A_{k_1}\Rightarrow f_\zeta (i(\beta_0,\ell ))=i(\beta_\zeta ,\ell ).
$$
Now we can imitate the proof of $(*)_T$ and get $\bigcap_{n<m} d_n>0$.  Let 
$Y_\zeta ={\rm Rang}(f_\zeta )$, and note that $f_0$ is the identity and
$Y_0=Y$. Clearly $f_\zeta$ induces an isomorphism from $\cBB [Y_0]$ onto $\cBB
[Y_\zeta]$. Call it $\hat f_\zeta$ and easily $d_\zeta\deq\hat f_\zeta (d_0)$.
So we can imitate the proof of $(*)_7$ and get $\bigcap_{n<m}d_n>0$. But
$$
c_i=\bigcap_{n<m} c_{\lan n\ran\whp\bar\ell}=\bigcap_{n<m} d_n>0
$$
as required.{\hfill$\square_{\rm 2.1}$}
\medskip

\n {\bf 2.2 Discussion}.
\smallskip

(1) The proof of 2.1 gives more, almost a division to $\le\mu$ subfamilies of
independent elements (in the Boolean algebra sense), see 2.16 below.

(2) We may wonder if ``$\mu\ge \beth_2$'' is necessary.  Actually it almost is
not (see 2.5 below) but ${\rm cf}\lambda >2^{\aleph_0}$ is essential (see 3.10
below).

We shall see below (in 2.5) what we can get from the proof of 2.1.
\medskip

\n {\bf 2.3 Definition}: {\it For a Boolean algebra $\B$ we say $\lan\lan
a_\alpha,b_\alpha \ran :\alpha <\alpha^*\ran$ is an explicitly independent
sequence of intervals in in $\B$ if
\begin{enumerate}
\item[(a)] $\B\vDash a_\alpha <b_\alpha$
\item[(b)] if $u_0, u_1\subseteq\alpha^*$ are finite disjoint then
$$
\B\vDash\bigcap_{\alpha\in u_0}b_\alpha\cap\bigcap_{\alpha\in u_1}(-a_\alpha
)>0. 
$$ 
\end{enumerate}}
\medskip

\n {\bf 2.4 Claim}:  Assume
\smallskip

\n ($*$)$_{Y}[X]$
\begin{enumerate}
\item[(a)] $|X|=\chi$ and $\cBB (X)$ is a Maharam measure algebra with free
basis $\{ x_i:i\in X\}$. For $Z\subseteq X$ we let $\cBB (Z)$ be the complete
subalgebra of $\cBB (X)$ generated by $\{x_i:i\in Z\}$.

\item[(b)$_Y$] $a_\alpha\in\cBB^+$ (i.e.\ $\mbox{Leb}(a_\alpha )>0$) for
$\alpha\in Y$ and $\beta <\alpha\Rightarrow a_\beta\not= a_\alpha$, while
$|Y|=\lambda$.
\end{enumerate}
\smallskip

1) If $\lambda={\rm cf}(\lambda ) >\aleph_1$ \underline{then} for some $Y'\in
[Y]^\lambda$, $Z\in [X]^{<\lambda}$ and $a_\alpha^-\le a_\alpha^+$ from $\cBB
(Z)$ we have:
\begin{enumerate}
\item[(i)] for $c\in\cBB (Z)$ we have $c\le a_\alpha\Rightarrow c\le 
a^-_\alpha$ and $a_\alpha\le c\Rightarrow a^+_\alpha\le c$
\item[(ii)] if $u\in [Y']^{<\aleph_0}$, $\eta\in {}^u2$ and 
$$
\cap\{ a^+_\alpha :\alpha\in u,\eta (\alpha )=1\}\cap\{1-a^-_\alpha :\alpha\in
u, \eta (\alpha )=0\}\not= 0, 
$$
then $\bigcap_{\alpha\in u} a_\alpha^{[\eta (\alpha )]}\neq 0$, where $c^{[0]}
=-c$, $c^{[1]}=c$.
\end{enumerate}
\smallskip

2) Assume $\inf\{\mbox{Leb} (a_\alpha \triangle b):b\in\lan a_\beta :\beta
<\alpha\ran\}_{\cBB} >0$ for $\alpha\in Y$.  $\underline{\rm Then}$ in part
(1) we can demand $a^-_\alpha < a^+_\alpha$.  Hence
\smallskip
\begin{enumerate}
\item[($*$)] there is $Y''\in [Y']^\lambda$ such that $\lan a_\alpha
:\alpha\in Y'')$ is independent iff there is $Y''\in [Y']^\lambda$ such that
$\lan (a^-_\alpha , a^+_\alpha ):\alpha \in Y''\ran$ is explicitly
independent. (See Definition 2.3 above.)
\end{enumerate}

3) If $|Y|=\lambda >|X|=\chi$ and $\chi_1<\chi$, $\sigma = {\rm cov}$ $(\chi ,
\chi_1^+, \aleph_1, 2)<\lambda$ $\underline{\rm then}$ $Y$ can be represented
as the union of $\le\sigma$ subsets $Y'$ such that for each there is $Z\in
[\chi]^{\le \chi_1}$ satisfying $\{ a_\alpha :\alpha\in Y'\}\subseteq\cBB
(Z)$.
\medskip

4) If the clause ($\alpha$) below holds then we can represent $Y$ as the union
of $\le \mu$ subsets $Y'$ each satisfying (c) below (and (b)$_{Y'}$),
\begin{enumerate}
\item[(c)$_{Y'}$] $a_\alpha=\tau (\ldots, x_{i(\alpha ,n)},\ldots )_{n
<\omega}$, $n\not= m\Rightarrow i(\alpha, n)\not= i(\alpha ,m)$ and the sets
$A_n(Y')=\{ i(\alpha ,n):n<\omega\}$ are pairwise disjoint where

\item[($\alpha$)] (i) $2^{\aleph_0}\le\mu=\mu^{\aleph_0}$ and
$2^\mu\ge\lambda$ or at least

\item[{}] (ii) $2^{\aleph_0}\le \mu$ and the density of the ($<\aleph_1$)-base
product ${}^\omega\chi$ is $\le \mu$
\end{enumerate}

5) If $Y'$ is as in (4) i.e.\ satisfies clause (c), $\underline{\rm then}$ any
finite intersection of $a_\alpha$'s (for $\alpha\in Y'$) is not zero.
\smallskip

6) If $Y'$ is as in (4) i.e.\ satisfies clause (c) $\underline{\rm then}$ $Y'$
is the union of $\le \beth_2$ subsets $Y''$, such that
\begin{enumerate}
\item[($*$)$_{Y''}$] there is an algebra $M$ with universe $Y''$ and $\le
\beth_1$ functions (with finite arity, of course) such that:
\begin{center}
if [$u\subseteq Y''$, $\alpha\in u\Rightarrow\alpha\notin
\mbox{cl}_M\{u\cap\alpha\}$], then $\lan a_\alpha :\alpha\in u\rangle$ is
independent.
\end{center}
\end{enumerate}
\medskip

\n {\bf Proof}. Straight and/or included in the proof of 2.1.
{\hfill$\square_{\rm 2.4}$}
\medskip

\n {\bf 2.5 Claim}.  In 2.1 we can weaken ``$\mu\ge \beth_2$'' to ``$\mu\ge
2^{\aleph_0}$'' or even ``cf $(\lambda) >2^{\aleph_0}$'' except possibly when
$\lambda$ is singular but $\boxtimes$ below fails:
\begin{enumerate}
\item[{$\boxtimes$}] for any countable set $\aaa$ of regulars,
$|\mbox{pcf}(\aaa )|\le\aleph_0$ or ($*$) from 2.6.
\end{enumerate}
\smallskip

\n {\bf Proof}: Without loss of generality we assume ($*$)$_4$ from the proof
of 2.1 (as the proof of 2.1 up to that point works here too).  Let $J$ be as
there, so $J$ is an ideal on $\omega$, so
\begin{enumerate}
\item[($+$)] $J$ is an ideal on $\omega$ and $\lan i(\alpha,n):n<\omega\ran
/J$ for $\alpha <\lambda$ are pairwise distinct;
\end{enumerate}
by the following observation 2.6 for some ideal $I$ on $\omega$ extending $J$
and $X\in [\lambda]^\lambda$, we have
$$
\alpha\in X\mbox{ \& } \beta\in X,\alpha\not=\beta\Rightarrow \{ n:i(\alpha
,n)=i(\beta ,n)\}\in I.
$$
This is enough for continuing with the old proof of 2.1.{\hfill$\square_{\rm 
2.5}$}
\medskip

\n {\bf 2.6 Fact}. 1) ~If $J$ is an ideal on $\kappa$, $\lan f_\alpha/J:\alpha
<\lambda\ran$ are pairwise distinct functions in ${}^\kappa{\rm Ord}$ and
$\theta ={\rm cf}(\lambda )>2^\kappa$ $\underline{\rm then}$ for some ideal
$I$ on $\kappa$ extending $J$ and $X\in [\lambda]^\lambda$ we have: $\alpha\in
X$ \& $\beta\in X$ \& $\alpha\not=\beta\Rightarrow f_\alpha\not=_I f_\beta$
except possibly when
\begin{enumerate}
\item[($*$)] $\lambda$ is singular and $\neg\boxtimes_\kappa$ where 
\item[$\boxtimes_\kappa$] for any set $\aaa$ of regular cardinals $>\kappa$ we 
have $|\aaa |\le \kappa\Rightarrow |{\rm pcf}(\aaa )|\le \kappa$, 
\end{enumerate}

2) We can replace ($*$) by 
\begin{enumerate}
\item[$(*)'$] $\lambda$ is singular and $\neg\boxtimes^+_{\kappa ,\lambda}$ or 
$\neg\boxtimes^-_{\kappa ,\lambda}$, where

\item[$\boxtimes^+_{\kappa,\lambda}$] for no set $\aaa$ of regular cardinals 
$>\kappa$, do we have $|\aaa |\le \kappa$ and $\lambda=\sup (\lambda\cap {\rm 
pcf}(\aaa ))$

\item[$\boxtimes^-_{\kappa,\lambda}$] there are no $\chi$, ${\rm
cf}(\lambda)=
\theta<\chi<\lambda$ and increasing sequences $\bar\lambda^\zeta =\lan
\lambda^\zeta_i:i<\kappa \ran$ of regular cardinals $\in (2^\kappa ,\chi)$
such that $\lan\max {\rm pcf}\{ \lambda^\zeta_i:i<\kappa\}:\zeta <\theta \ran$
is increasing with limit $\lambda$ but for every ultrafilter $\DD$ on $\kappa$
we have
$$
\sup\Big\{{\rm tcf}\Big(\prod_{i<\kappa} \lambda^\zeta_i/\DD \Big):\zeta
<\theta\Big\}
 <\lambda.
$$
\end{enumerate} 
\medskip

\n {\bf Proof}.
\smallskip

1) ~Follows by (2).
\smallskip

2) The proof is split to cases.

{\bf Case 1}: $\lambda$ is regular.  We apply 7.3 here which is
\cite[6.6D]{Sh:430} (or in more details \cite[6.2]{Sh:513}). 
\smallskip

{\bf Case 2}: $\lambda$ singular.  First note
\smallskip

\n {\it 2.6A Fact}: $\neg\boxtimes^+_{\lambda ,\kappa}\Rightarrow\boxtimes^-_{
\lambda,\kappa}$.
\smallskip

\n [Why?  Let $\aaa$ exemplify $\neg\boxtimes^+_{\lambda,\kappa}$, let
$\theta_\e \in {\rm pcf}(\aaa )\setminus\{\lambda\}$ be increasing for $\e<
\theta$ with limit $\lambda$.  Let $\bbb_\e\subseteq \aaa$ be such that
$\theta_\e ={\rm max~pcf}(\bbb_\e)$ and let
$\lan\lambda_\zeta:\zeta<\kappa 
\ran$ list $\aaa$ and let $\lambda^\e_\zeta$ be: $\lambda_\zeta$ if
$\lambda_\zeta\in\bbb_\e$ and $(2^\kappa )^+$ if $\lambda_\zeta\not\in
\bbb_\e$. Now $\bar\lambda^\e=\lan\lambda^\e_\zeta :\zeta <\kappa\ran$
exemplifies $\neg\boxtimes_{\lambda ,\kappa}$.  First ${\rm max ~pcf}\{
\lambda^\e_\zeta:\zeta <\kappa\}=\theta_\e <\lambda$ and $\theta_\e$
is increasing with limit $\sup(\lambda\cap{\rm pcf}(\aaa))$.
Secondly, for every ultrafilter $\DD$ on $\kappa$ for each $\e$ we have
${\rm tcf}\Big(\prod_{\zeta <\kappa}\lambda^\e_\zeta/\DD\Big)$ is $(2^\kappa
)^+$ or is ${\rm tcf}\Big(\prod_{\zeta<\kappa}\lambda_\zeta /\DD \Big)$.
(Simplify the first case if $\{\zeta<\kappa:\lambda_\zeta\not\in\bbb_\e\}\in
\DD$ and the second case if $\{\zeta<\kappa:\lambda_\zeta\in\bbb_\e\}\in\DD$.)
So now if ${\rm tcf}\Big(\prod_{\zeta <\kappa}\lambda_\zeta/\DD \Big)\ge
\lambda$ implies ${\rm tcf}\Big(\prod_{\zeta<\kappa}\lambda^\e_\zeta/\DD\Big)=
(2^\kappa )^+$ as the later is $\le\theta_\e<\lambda$, so really there is no
ultrafilter $\DD$ on $\kappa$ for which ${\rm sup}\Big\{{\rm tcf}\Big(\prod_{
\zeta<\kappa}\lambda^\e_\zeta /\DD\Big):\e<\theta\Big\}<\lambda$, so the
second demand in $\boxtimes^-_{\lambda ,\kappa}$ holds. {\hfill$\square_{\rm 
2.6A}$}

Now we assume $\boxtimes^-_{\lambda ,\kappa}$. For every regular $\sigma\in
(2^\kappa ,\lambda )$ we apply 7.3 to $\lan f_\alpha :\alpha <\sigma\ran$, so
we can find $A_\sigma\subseteq\kappa$ and $\lan \gamma_{\sigma,i}:i<\kappa
\ran$ such that
\begin{enumerate}
\item[($*$)$_0$] for every sequence $\lan\beta_i:i\in A_\sigma\ran$ satisfying 
$\beta_i<\gamma_{\sigma ,i}$ there are $\sigma$ ordinals $\alpha <\sigma$ for 
which
$$\begin{array}{l}
i\in A_\sigma\Rightarrow\beta_i<f_\alpha (i)<\gamma_{\sigma ,i}\\
i\in\kappa\backslash A_\sigma\Rightarrow f_\alpha (i)=\gamma_{\sigma ,i}
\end{array}
$$
\item[($*$)$_1$] $B\in J\Rightarrow\sigma\in {\rm pcf}\{{\rm cf}(\gamma_{
\sigma ,i}):i<\kappa , i\in A_\sigma, i\notin B\}$.
\end{enumerate}
Let $J_\sigma =\{ B\subseteq\kappa :{\rm max ~pcf}\{{\rm cf}(\gamma_{\sigma
,i}):i\in \kappa \backslash A_\sigma\; {\rm and}\; i\in B\}<\sigma\}$, so
clearly $\sigma ={\rm tcf}\Big(\prod_{i<\kappa}{\rm cf}(\gamma_{\sigma,i})/
J_\sigma\Big)$ and $J\subseteq J_\sigma$. Let $A'_\sigma$ be such that
$A'_\sigma\subseteq A_\sigma$, and $\sigma ={\rm max ~pcf}\{{\rm cf}(\gamma_{
\sigma ,i}):i\in A'_\sigma\}$.  Also, as $\theta ={\rm cf}(\lambda )>
2^\kappa$, for some $A'\subseteq\kappa$ (infinite) the set $\Theta=\{\sigma:
\,2^\kappa <\sigma ={\rm cf}(\sigma )<\theta\; {\rm and}\; A'_\sigma =A'\}$ is
unbounded in $\lambda$.  Let $\lan\sigma_\e :\e<\theta\ran$ be an increasing
unbounded sequence of members of $\Theta$, such that its limit is
$\lambda$. Apply 7.3 (see Case 1) to $\lan g_\e \restriction A':\e<\theta\ran$
where $g_\e (i)=\gamma_{\sigma_\e ,i}$, and get $\lan\beta^*_i:i\in A'\ran$
and $B'\subseteq A'$ such that 
\begin{enumerate}
\item[($*$)] if $\lan\beta_i:i\in A'\ran$ satisfies 
$i\in A'\Rightarrow\beta_i<\beta^*_i$ $\underline{\rm then}$ for 
unboundely many ordinals $\e<\theta$
$$\begin{array}{l}
i\in B'\Rightarrow\beta_i<\gamma_{\sigma_\e,i}<\beta^*_i\\
i\in A'\backslash B'\Rightarrow\gamma_{\sigma_\e,i}=\beta^*_i.
\end{array}
$$
\end{enumerate}
Can $B'=\emptyset$? This would mean that for some unbounded $X\subseteq\theta$
we have $\e\in X\Rightarrow (\forall i\in A')[\gamma_{\sigma_\e,i}=\beta^*_i]$
hence $\{\sigma_\e :\e\in X\}\subseteq {\rm pcf}\{{\rm cf}(\beta^*_i):i\in
A'\}$, so $\{{\rm cf}(\beta^*_i):i\in A'\}$ has pcf of cardinality $\ge\theta
>2^\kappa$ whereas $|A'|\le \kappa$, contradiction, so really
$B'\not=\emptyset$. 

As we are assuming $\neg\boxtimes^-_{\kappa ,\lambda}$, there is an ultrafilter
$\DD$ on $A'$ such that
$$
\lambda\le\sup\Big\{ {\rm tcf}\Big(\prod_{i\in A'}\gamma_{\sigma_\e,i}/
\DD\Big):\e<\theta\Big\}.
$$ 
Clearly
$$
{\rm tcf}\Big(\prod_{i\in A'}\gamma_{\sigma_\e,i}/\DD\Big) \le\sigma_\e
<\lambda 
$$
(by the choice of $A'_{\sigma_\e}=A'$).  Without loss of generality $\sigma_\e
>\theta$ for each $\e <\theta$.  So we can choose, for each $\e$, a function
$h_\e\in\prod_{i\in A'}\gamma_{\sigma_\e,i}$ such that
\begin{enumerate}
\item[($*$)] if $\zeta <\theta$ and $\zeta\not=\e$, while $\lan\gamma_{
\sigma_\zeta,i}:i\in A'\ran\le_\DD\lan\gamma_{\sigma_\e ,i}:i\in A'\ran$
$\underline{\rm then}$
$$
\lan\gamma_{\sigma_\zeta ,i}:i\in A'\ran<_\DD h_\e.
$$
(Note that $\lan\gamma_{\sigma_\zeta ,i}:i\in A'\ran\neq_\DD\lan\gamma_{
\sigma_\e,i}:i\in A'\ran$ because of the cofinalities of the respective
ultraproducts.) So, considering $\DD$ as an ultrafilter on $\kappa$:
$$\begin{array}{rl}
X_\e=&\left\{\alpha <\sigma_\e:h_\e <_\DD f_\alpha 
<_\DD\lan\gamma_{\sigma_\e,i}:i<\kappa\ran\;\;\mbox{but}\right.\\
&\left.\beta <\alpha\Rightarrow\neg (f_\alpha\le_\DD f_\beta 
<_\DD\lan\gamma_{\sigma_\e,i}:i<\kappa\ran )\right\}
\end{array}
$$
has cardinality $\sigma_\e$.  So $X=\bigcup_{\e <\sigma} X_\e$ is as required.
\end{enumerate}
{\hfill$\square_{\rm
2.6}$}
\medskip

We may wonder whether we can remove or at least weaken the assumption ($*$): 
the answer is:
\medskip

\n {\bf 2.7 Claim}: 1) If $\kappa\le\lambda$ and $\theta ={\rm cf}(\lambda
)<\lambda$, and $\boxtimes^-_{\kappa,\lambda}$ (from 2.6) $\underline{\rm
then}$ for some $f_\alpha\in {}^\kappa\lambda$ (for $\alpha <\lambda$) the
conclusion of 2.6(1) fails.
\medskip

\n {\bf Proof}: 1) ~Let $\chi$, $\lambda^\zeta_i$ $(i<\kappa, \zeta<\theta)$
be as in $\boxtimes^-_{\kappa,\lambda}$.  
\medskip

Let $\aaa_\zeta =:\{\lambda^\zeta_i:i<\kappa\}$, and $\sigma_\zeta ={\rm max
~pcf}(\aaa_\zeta )$. Without loss of generality $\lan\sigma_\zeta :\zeta
<\theta\ran$ is increasing with limit $\lambda$.  By \cite[II, \S3]{Sh:g} for
each $\zeta <\theta$ we can find $\lan f^\zeta_\alpha:\alpha<\sigma_\zeta\ran$
be such that:
$$
\bbb\subseteq\aaa_\zeta\Rightarrow |\{ f^\zeta_\alpha\restriction\bbb :\alpha 
<{\rm max ~pcf}(\aaa_\zeta )\}|={\rm max ~pcf}(\bbb ).
$$
Define $\lan f_\alpha :\alpha <\lambda\ran$ by: $f_\alpha (i)=f_\alpha^{\zeta 
(\alpha )}(\lambda^\zeta_i)$ where $\zeta (\alpha )=\min\{\zeta:\sigma_\zeta>
\alpha\}$.  Now check. {\hfill$\square_{\rm 2.7}$}
\medskip

\n {\bf 2.8 Discussion}.  1) ~So if $2^\kappa <\lambda$, $\theta ={\rm cf}(
\lambda )$ then 2.7 shows that 2.6 is the best possible.  (Of course we still
do not know if $\boxtimes^-_{\lambda ,\kappa}$ is possible). See more in 3.11.

2) ~Note: If ${\rm cf}(\lambda )>2^\kappa$, and $(\forall\aaa )$ $(\aaa
\subseteq {\rm Reg ~\& ~}|\aaa |\le\kappa <\min(\aaa)\to |{\rm pcf}(\aaa )|\le
|\aaa|$), $\underline{\rm then}$ $\square^-_{\lambda ,\kappa}$ cannot occur as
without loss of generality 
$$
J_\zeta =\{ A\subseteq\kappa :{\rm max ~pcf}\{\lambda^\zeta_i:i\in A\}< {\rm 
max ~pcf}\{\lambda^\zeta_i:i<\kappa\}\}
$$
does not depend on $\zeta$.

\section{Strong independence in Maharam measure algebras}

\n {\bf 3.1 Claim}:  Assume
\begin{enumerate}
\item[(a)] $I_i$ is a $\kappa_i$-complete ideal on $\lambda_i$ for $i<\delta$.
\item[(b)] $\kappa_i>\sum_{j<i}\kappa_j$.
\item[(c)] $\mu=\sup_{i<\delta}\kappa_i$ is strong limit singular.
\item[(d)] $\lambda_i<\mu$.
\item[(e)] $\lambda =\mu^+=2^\mu$.
\end{enumerate}
$\underline{\rm Then}$ there is $\bar\eta$ a super $\lambda$-sequence for
$\lan I_i:i<\delta\ran$, where
\medskip

\n {\bf 3.2 Definition}. {\it We say $\bar\eta$ is a super $(\lambda,
J)$-sequence for $\lan I_i:i<\delta\ran$ if, in addition (to the demands in
1.8)
\begin{enumerate}
\item[($*$)] for every $n <\omega$ and $\beta_{\alpha,\ell}<\lambda$ (for
$\alpha <\lambda , \ell <n$) increasing with $\ell$, pairwise distinct we have
$$
\Big\{ i<\delta :\{\lan \eta_{\beta_{\alpha ,\ell}}(i):\ell<n\ran :\alpha
<\lambda\}\in \prod_{\ell<n}I_i\Big\}\in J.
$$
\end{enumerate}
Moreover
\begin{enumerate}
\item[($*$)] if $n<\omega$, $\beta_{\alpha ,\ell}<\lambda$ (for $\alpha
<\lambda$, $\ell <n$), $\beta_{\alpha , \ell}<\beta_{\alpha , \ell +1}$, and
the $\beta_{\alpha ,\ell}$ are pairwise distinct $\underline{\rm then}$ for
some $A\in J$ we have:
\end{enumerate}
{\it if $m<\omega$, $i_0<i_1<\cdots i_{m-1}$ belong to $\delta\backslash A$,
then}
$$
\left\{\bigr\langle\lan\eta_{\beta_{\alpha ,\ell}}(i_t):\ell<n\ran:t<m\bigr
\rangle:\alpha <\lambda\right\}\in\Big(\prod_{t<m}\Big(\prod_{\ell<n}I_{i_t}
\Big)\Big)^+.
$$
}
\smallskip

\n {\bf Proof}.  Like the proof of 1.8.  \hfill{$\square_{3.1}$}
\medskip

\n {\bf 3.3 Example}.  $\lambda =\mu^+=2^\mu$, $\mu=\sum_{i<\kappa}\lambda_i$,
$i<j\Rightarrow\delta=\kappa <\lambda_i <\lambda_j <\mu$ and each $\lambda_i$
is measurable with a $\big(\aleph_0+\sum_{j<i}\lambda_j\big)^+$-complete
normal (or just Ramsey for $n_i$) ultrafilter $\DD_i$ on $\lambda_i$.

Let $\bar n =\lan n_i:i<\kappa\ran$, $i<n_i<\omega$, (if $\kappa =\aleph_0$,
$n_i=i$ we may omit it)
$$
I_i=\{ A\subseteq [\lambda_i]^{n_i}:\mbox{ for some $B\in\DD_i$ we have
$[B]^{n_i}\cap A=\emptyset$}\}.
$$

Then
\begin{enumerate}
\item[($*$)$_1$] claim 3.1 applies.
\item[($*$)$_2$] for every $m<\omega$ and $X\in\prod_{\ell <m}I_i$ we can find
$A\in\DD_i$ such that:
$$
\{\bar s:\bar s=\lan s_\ell :\ell< m\ran, s_\ell\in [A]^{n_i},s_\ell<s_{\ell
+1}\}\cap X=\emptyset. 
$$
\end{enumerate}
\smallskip

\n {\bf 3.4 Definition}: {\it 1) For a Boolean algebra $\B$ we say $\lan
(a_\alpha , b_\alpha ):\alpha <\alpha^*\ran$ is a strongly independent
sequence of intervals if
\begin{enumerate}
\item[(a)] $\B\vDash a_\alpha < b_\alpha$
\item[(b)] if $\B'$ is a Boolean algebra extending $\B$ and $n <\omega$,
$\alpha_0 <\alpha_1<\cdots <\alpha_{n-1}<\alpha^*$ and $\B'\vDash$
``$a_{\alpha_\ell}\le x_\ell \le b_{\alpha_\ell}$'' for $\ell <n$, then any
non-trivial Boolean combination of $\lan x_\ell :\ell <n \ran$ is non-zero (in
$\B'$).
\end{enumerate}

2) We say, for a Boolean algebra $\B$ that $\lan (a_\alpha , b_\alpha ):\alpha
<\alpha^*\ran$ is a $\lambda$-anti independent sequence of intervals if
\begin{enumerate}
\item[(a)] $\B\vDash a_\alpha\le b_\alpha$
\item[(b)] if $\B'$ is a Boolean algebra extending $\B$ and $X\in
[\alpha^*]^\lambda$ and $\B'\vDash \mbox{``}a_\alpha\le x_\alpha \le
b_\alpha$'' for $\alpha\in X$, then there are $n <\omega$ and
$\alpha_0<\alpha_1<\cdots <\alpha_{n-1}$ from $X$ such that some non-trivial
Boolean combination of $\lan x_{\alpha_\ell}:\ell <n\ran$ is zero.
\end{enumerate}

3) We say $\lan (a_\alpha , b_\alpha ):\alpha <\alpha^*\ran$ is an independent
sequence of intervals in the Boolean algebra $\B$ if letting $\B'$, $x_\alpha$
be as in 3.5 below, we have $\lan x_\alpha :\alpha <\alpha^*\ran$ is
independent (in $\B'$).

4) We say $\lan (a_\alpha , b_\alpha ):\alpha <\alpha^*\ran$ is a strongly
$\lambda$-anti-independent sequence of intervals for the Boolean algebra $\B$
if:
\begin{enumerate}
\item[(a)] $\B\vDash a_\alpha\le b_\alpha$.
\item[(b)] if $\B'$, $X, x_\alpha (\alpha\in X)$ are as in 3.4(2)(b) above,
then the Boolean subalgebra of $\B'$ generated by $\{ x_\alpha :\alpha\in X\}$
contains no free subset of cardinality $\lambda$.
\end{enumerate}

5) We say $\lan (a_\alpha ,b_\alpha ):\alpha <\alpha^*\ran$ is mediumly
$\lambda$-anti independent (sequence of intervals of the Boolean algebra $\B$)
if
\begin{enumerate}
\item[(a)] $\B\vDash a_\alpha\le b_\alpha$
\item[(b)] if $\B'$ is the free extension of $\B$ for $\lan(a_\alpha,b_\alpha)
:\alpha <\alpha^*\ran$ (see 3.5), then the Boolean subalgebra of $\B'$
generated by $\{ x_\alpha :\alpha <\alpha^*\}$ contains no free subalgebra of
cardinality $\lambda$.
\end{enumerate}
}

\n {\bf 3.5 Definition}. {\it We say that $\B'=\B'(\B ,\lan (a_\alpha b_\alpha
):\alpha <\alpha^*\ran)$, or $\B'$ is the free extension of $\B$ for $\lan
(a_\alpha ,b_\alpha ):\alpha <\alpha^*\ran$, if
\begin{enumerate}
\item[($*$)] $\B'$ is the algebra freely generated by $\B\cup\{ x_\alpha
:\alpha <\alpha^*\}$ except for the equations:
\begin{enumerate}
\item[(a)] the equations which $\B$ satisfies
\item[(b)] $a_\alpha\le x_\alpha\le b_\alpha$, for $\alpha <\alpha^*$.
\end{enumerate}
\end{enumerate}}

\n {\bf 3.5A Observation}: 1) ~In 3.4(3), if $\B\subseteq\cBB (\alpha_0)$,
$\alpha_0+\omega +\alpha^*\le\alpha_1$ $\underline{\rm then}$ we can embed
$\B'$ into $\cBB (\alpha_1)$ over $\B$.

2) ~There are obvious implications among the notion from Definition 3.4 and
some equivalences: independent (3.4(3)) with explicitly independent; and
stronger independent with ``(a) of 3.4(1) and if $\alpha_1,\ldots , \alpha_2$,
$\beta_1,\ldots ,\beta_m<\alpha^*$ with no repetition,
$$
\B\vDash\bigcap^n_{\ell=1}a_{\alpha_\ell}\cap\bigcap^m_{\ell=1}(
-b_{\beta_\ell})>0."
$$
\smallskip

\n {\bf 3.6 Lemma}.  Assume $\mu$ is strong limit singular of countable
cofinality and $\lambda =\mu^+=2^\mu$.

$\underline{\rm Then}$ in $\cBB (\mu )$, (the Maharam measure algebra of
dimension $\mu$) we can find a sequence $\lan (a_\alpha , b_\alpha ):\alpha
<\mu\ran$ such that:
\begin{enumerate}
\item[(a)] $\cBB (\mu )\vDash a_\alpha <b_\alpha$
\item[(b)] $\lan (a_\alpha , b_\alpha ):\alpha <\lambda\ran$ is strongly
$\lambda$-anti independent.
\end{enumerate}
\medskip

\n {\bf Remark}.  What is the difference with 1.13?  Note that 3.4(ii)(b)
speaks on ``no free subset of the Boolean algebra'', not just of the set.
\medskip

\n {\bf Proof}. 

1) ~ Let $\mu =\sum_{n<\omega}\lambda^0_n$, (we may demand $\beth_{n+8}(
\lambda^0_n)<\lambda^0_{n+1}<\mu)$ and let $I_n$ be $\mbox{ERI}^{n,h}_{
\beth_{n-1}(\lambda^0_n)^+,(\lambda^0_n)^+}$ (see Definition 1.14, they were
used in the proof of 1.13). Let $\bar\eta=\lan\eta_\alpha:\alpha<\lambda\ran$
be as guaranteed by 3.1 (so $\mbox{lg}(\eta_\alpha )=\omega$, $\eta_\alpha(n)
\in [\lambda_n]^n$, where $\lambda_n=\beth_{n-1}(\lambda^0_n)^+$. So $I_{n+1}$
is $|\mbox{Dom}(I_n)|^+$-complete, (we could also have $\lan I_n:n<\omega\ran$
is normal). Renaming, let $x^n_\alpha$ (for $n<\omega$, $\alpha <\lambda_n$) be
the free generators of the Maharam algebra.

Define for $\alpha <\lambda$ and $n<\omega$
$$\begin{array}{l}
\mbox{$a^*_{\alpha ,n}=\cap\{ x^m_\beta :\beta\mbox{ appears in }\eta_\alpha
(m)\}$}\\ 
\mbox{$b^*_{\alpha ,n}=\cup\{(1- x^m_\beta :\beta\mbox{ appears in }
\eta_\alpha (m)\}$.}
\end{array}
$$
We define by induction on $n$, the elements $a_{\alpha ,n}$, $b_{\alpha ,n}$
as follows: for $n<5$ let $a_{\alpha ,n}=0$, $b_{\alpha ,n}=1$.  For $n\ge 5$
we let $a_{\alpha ,n}=a_{\alpha ,n-1}\cup (a^*_{\alpha ,n}\cap b_{\alpha ,n})$
and $b_{\alpha ,n}=b_{\alpha ,n+1}\cap (b^*_{\alpha ,n}\cup a_{\alpha ,n})$.
We can prove by induction on $n<\omega$ that $a_{\alpha ,n-1}\le a_{\alpha,n}
\le b_{\alpha,n}\le b_{\alpha,n-1}$.  We can compute the measure, e.g.~let
$(b_{\alpha ,n}-a_{\alpha ,n})=\prod\{ 1-2^{-(\ell -1)}:5\le\ell \le n\}$.

Let $a_\alpha =\bigcup_{n<\omega} a_{\alpha ,n}\in\cBB (\mu)$, $b_\alpha
=\bigcap_{n<\omega}b_{\alpha ,n}\in\cBB (\mu )$.

So clearly $\cBB (\mu )\vDash a_\alpha \le b_\alpha$, and by the measure
computations above, $\cBB (\mu )\vDash a_\alpha <b_\alpha$. So $\lan
(a_\alpha ,b_\alpha ):\alpha <\lambda\ran$ is a sequence of intervals.
Suppose $\B$, $c_\alpha$ (for $\alpha <\lambda$), is a counterexample to the
conclusion so there is an independent subset $\{ d_\alpha :\alpha <\lambda\}$
of $\lan c_\alpha :\alpha <\lambda\ran_{\B}\subseteq\B $.

So, for each $\alpha <\lambda$ for some $k_\alpha <\omega$ and a Boolean term
$\tau =\tau_\alpha (x_0,\ldots ,x_{k_\alpha -1})$ and some $\beta_{\alpha
,0}<\beta_{\alpha ,1}<\cdots <\beta_{\alpha,k_\alpha-1}$ we have $d_\alpha
=\tau_\alpha (c_{\beta_{\alpha ,0}}, c_{\beta_{\alpha,1}},\ldots,c_{\beta_{
\alpha,k_\alpha-1}})$.

As we can replace $\{ d_\alpha :\alpha <\lambda \}$ by any subset of the same
cardinality without loss of generality $\tau_\alpha =\tau$, so let $k_\alpha
=k(*)$.

Similarly, by the $\Delta$-system argument without loss of generality for some
$k<k(*)$ we have $\ell <k\Rightarrow\beta_{\alpha ,\ell}=\beta_\ell$ and
$\alpha (1)<\alpha (2)\Rightarrow \beta_{\alpha (1), k(*)-1}<\beta_{\alpha
(2),k}$.

Let $X_n=\{ \lan \eta_{\beta_{\alpha ,\ell}}(n):k\le\ell<k(*)\ran :\alpha
<\lambda\}\subseteq {}^{(k(*)-k)}([\lambda_n]^n)$.  So we know that $B=\{
n<\omega : n\ge k(*)-k\;\mbox{and}\; X_n\in(\prod^{k(*)-1}_{\ell =k} I_n)^+\}
\in J^+$. Let $n\in B$.  We can find a function $h:X_n\to\lambda$ such that
$\bar t\in X_n\mbox{ \& } h(\bar t)=\alpha\Rightarrow\bar t=\lan\eta_{\beta_{
\alpha ,\ell}}(n):k\le \ell <k(*)\ran$.  Let $m(*)<\omega$ be large enough, a
power of 2 for simplicity. 

As $X_n\in (\prod^{k(*)-1}_{\ell =k}I_n)^+$, we can find $\lan S_\ell :\ell\in
[k,k(*)]\ran$ and $\lan u_{\bar s}:\bar s\in S_\ell\ran$ for $\ell\in [k,
k(*))$ such that
\begin{enumerate}
\item[(a)] $S_k=\{\lan\quad\ran\}$.
\item[(b)] $u_{\bar s}\in [\lambda_n]^{m(*)}$.
\item[(c)] the $u_{\bar s}$'s are pairwise disjoint.
\item[(d)] $S_{\ell +1}=\{\bar s\;\hat{\phantom{A}}\; \lan w\ran :\bar s\in
S_\ell , w\in [u_{\bar s}]^n\}$.
\item[(e)] $S_{k(*)}\subseteq X_n$.
\end{enumerate}
(We just do it by induction on $\ell$ using the definition of $\prod^{k(*)-
1}_{\ell =k} I_n$ and the definition of $I_\ell$.) So it suffices to show that
$\lan d_{h(\bar t)}:\bar t\in S_{k(*)}\ran$ is not independent. For this just
note: 

\begin{enumerate}
\item[($\otimes$)] for every $\epsilon\in\bbr^{>0}$ if $n$ is large enough
compared to $k(*), 1/\epsilon$, and $m(*)$ is large enough compared to $n$
$\underline{\rm then}$ for every ultrafilter $\DD$ on $\cBB (\mu )$ we can by
downward induction on $\ell =k,\ldots ,k(*)-1$ find $u^-_{\bar s}\in [u_{\bar
s}]^{m(*)/2^{k(*)-\ell}}$ and $\eta_{\bar s}\in {}^{\{k,\ldots ,k(*)-1\}}2$
for $\bar s\in S_\ell$ such that: $\bar s\unlhd\bar t\in S_{\ell_1}$ and
$\ell\le\ell_1<k(*)$ and $\alpha\in u^-_{\bar t}\Rightarrow
[x^n_\alpha\in\DD\equiv\eta_{\bar s}(\ell_1)=1]$.
\end{enumerate}
Now let $\eta^*=\eta_{\lan\;\ran}$ (i.e. $\eta_{\bar s}$ for the unique $\bar
s\in S_0$) and for $m<k(*)$ letting $S'_m=\Big\{\bar s\in S_m:$ if $\ell <m$
then $\bar s(\ell)\in [u^-_{\bar s\restriction \ell}]^n\Big\}$, we have $\bar
s\in S'_{k(*)}\Rightarrow d_{h(\bar s)}\in\DD$ or $\bar s\in
S'_{k(*)}\Rightarrow d_{h(\bar s)}\notin \DD$.

So to prove that $\lan d_\alpha :\alpha <\lambda\ran$ is not independent it
suffices to find $S\subseteq S_{k(*)}$ such that
\begin{enumerate}
\item[$\otimes_S$] $\displaystyle{\bigcap_{\alpha\in S} 
d_\alpha\cap\bigcap_{\alpha\in S_{k(*)}\backslash S} d_\alpha =0}$
\end{enumerate}
or equivalently
\begin{enumerate}
\item[$\otimes'_S$] for no ultrafilter $\DD$ on $\cBB (\mu )$ do we have
$$
\alpha\in S_{k(*)}\Rightarrow [d_\alpha\in\DD\equiv\alpha\in S].
$$
\end{enumerate}
By the argument above it will suffice to have
\begin{enumerate}
\item[$\otimes''_S$] $\underline{\rm if}$ $\lan u^-_{\bar s}:\bar s\in\cup\{
S'_\ell :\ell <k(*)\ran$ satisfies: $S'_0=S_{0}$, $S'_\ell\subseteq S_{\ell}$,
$\bar s\in S'_\ell\Rightarrow u^-_{\bar s}\in [u_{\bar s}]^{m(*)/2^{2k(*)-
\ell}}$ and $S'_{\ell +1}=\{\bar s^\frown\lan w\ran :\bar s\in S'_\ell\;
\mbox{and}\; w\in [u^-_{\bar s}]^n\}$ $\underline{\rm then}$ $S\cap S'_{k(*)}
\not\in \{\emptyset ,S\}$.
\end{enumerate}
Now, not only that this is trivial by the probabilistic existence proof \'a la
Erd\"os but the proof gives much more than enough.{\hfill$\square_{\rm 3.6}$}
\medskip

\n {\bf 3.7 Claim}: Assume
\begin{enumerate}
\item[($*$)] $\lambda$ is regular $>\aleph_0$ and $\lan (a_\alpha, b_\alpha):
\alpha<\lambda\ran$ is a strongly (or just mediumly)
$\lambda$-anti-independent sequence of pairs from $\cBB (\lambda)$ satisfying 
$a_\alpha <b_\alpha$.
\end{enumerate}
Then:
\begin{enumerate}
\item[(a)] There is $\B'$, such that:
\smallskip

($\alpha$) $\B'$ is a subalgebra of $\cBB (\lambda )$,
\smallskip

($\beta$) $\B'$ has cardinality $\lambda$ and even dimension $\lambda$,
\smallskip

($\gamma$) there is no subset of $\B'$ of cardinality $\lambda$ which is
independent
\smallskip

\item[(b)] Let $\B'$, $x_\alpha (\alpha <\lambda )$ be as in 3.5, then the
Boolean algebra in clause (a) can be chosen isomorphic to $\lan x_\alpha
:\alpha <\lambda\ran_{\B'}$.
\end{enumerate}
\smallskip

\n {\bf Proof}: Straight. Clause (a) follows from clause (b).  For clause (b)
apply Definition 3.4(5) and 3.5A. (Note: we can use $\cBB'\subseteq\cBB
(\lambda +\lambda)$). It has already been done by Plebanek \cite{Pl1}. 
{\hfill$\square_{\rm 3.7}$} 
\medskip

\n {\bf 3.8 Conclusion}: For $\lambda$ as in 3.6 (i.e.\ $\lambda=\mu^+=2^\mu$,
$\mu$ strong limit $>{\rm cf}(\mu )=\aleph_0$) or just as in ($*$) of 3.7, we
have 
\begin{enumerate}
\item[($*$)] there is a topological space $X$ which is Hausdorff, compact zero
dimensional, with a measure $\mbox{Leb}$ on the Borel sets such that it has
dimension $\lambda$, so as a measure space is isomorphic to $\cBB (\lambda )$
but there is no homomorphism from $X$ onto ${}^\lambda 2$.
\end{enumerate}
\smallskip

\n {\bf Proof}: By 3.6(1) ($*$) of 3.7 holds so we can restrict ourselves to
this case.  So by 3.7 we know that clause (a) of 3.7 holds.  Now it follows
that ($*$) holds, more specifically, that the Check-Stone compactification of
$\B'$ (i.e.\ the set of ultrafilters of $\B'$ with the natural topology) and
the measure of $\B'$ (which is just the restriction of the one on $\B(\lambda
)$) satisfies ($*$) of 3.8.{\hfill$\square_{\rm 3.8}$
\medskip

\n {\bf 3.9 Example}: Assume $\cBB$ is a Maharam measure algebra of dimension
$\mu$ and free basis $\lan x_\alpha :\alpha <\mu\ran$, $\mu\ge\lambda >{\rm
cf}(\lambda )=\aleph_0$. Then $(*)_{2,\lambda}$ below holds, where

$(*)_{2,\lambda}$\quad there are positive pairwise distinct members $a_\alpha$
of $\cBB (\mu )$ for $\alpha <\mu$, such that for every $X\in [\lambda
]^\lambda$ for some $\alpha\not=\beta$ from $X$, $a_\alpha \cap a_\beta =0$.
\medskip

\n {\bf Proof}: Trivial: let $\lambda =\sum_{n<\omega}\lambda_n$, $\lambda_n
<\lambda_{n+1}$ and for $\alpha\in \big(\bigcup_{\ell <n}\lambda_\ell
,\lambda_n\big)$ we let $a_\alpha =x_{\omega +\alpha}\cap
(x_n-\bigcup_{m<n}x_m)$. {\hfill$\square_{\rm 3.9}$}
\medskip

\n {\bf 3.10 Fact}. Suppose $\aleph_0<{\rm cf}(\lambda )<\lambda$ and there
are positive $b_\alpha\in\cBB({\rm cf}(\lambda ))$ for $\alpha <{\rm
cf}(\lambda )$ such that for every $X\in [{\rm cf}(\lambda )]^{{\rm
cf}(\lambda)}$ for some $m<\omega$ and $\beta_0,\ldots ,\beta_m\in X$ we have
$\mbox{Leb}\Big(\bigcap_{\ell\le m} b_{\beta_\ell}\Big)=0$ and
$\mu\ge\lambda$. $\underline{\rm Then}$ we can find pairwise distinct
$a_\alpha \in\cBB(\lambda)$ for $\alpha <\lambda$ such that for every $X\in
[\lambda ]^\lambda$ for some $m<\omega$, $\beta_0,\ldots ,\beta_m\in X$ we
have $\mbox{Leb}\Big(\bigcap_{\ell\le m}a_{\beta_\ell}\Big)=0$ i.e.\ $\cBB
(\lambda )\vDash\bigcap_{\ell\le m} a_{\beta_\ell}=0$.  \medskip \n{\bf
Proof}: Like the proof of 3.9 replacing $x_n-\bigcup_{m<n} x_m$ (for
$n<\omega$) by $b_\alpha$ (for $\alpha< {\rm cf}(\lambda )$). (Just say that
if ${\rm cf}(\lambda )$ is a precaliber of $\cBB$ then so is $\lambda$.)
\medskip

\n {\bf 3.10A Remark}: 1) By 2.1 we have in 3.10 that necessarily ${\rm
cf}(\lambda )\le \beth_2$ is normally ${\rm cf}(\lambda)\le\beth_1$.

2) Note that 3.11 elaborates 2.7 above and 3.12 is complementary to \S2.
\medskip

\n {\bf 3.12 Example}: Assume $\aleph_0\le\sigma\le\theta={\rm cf}(\lambda)\le
2^\sigma\le\mu <\lambda$, 
$$
\lambda =\sup\{\max\mbox{ pcf}(\aaa ):\aaa\subseteq\mbox{Reg}\cap\mu\backslash
2^{\sigma},|\aaa |=\sigma , [\aaa ]^{<\sigma}\subseteq J_{<\max\mbox{
pcf}(\aaa )} [\aaa ], 
$$
\centerline{and ${\rm sup}({\rm pcf}(\aaa )\backslash\{{\rm max ~pcf}(\aaa 
)\})\le \mu\}$}
\smallskip

\noindent and there is $\AAA\subseteq [\sigma ]^\sigma$ such that $|\AAA
|\ge\theta$ and $[A\not= B\mbox{ \& } A\in \AAA\mbox{ \& } B\in\AAA\Rightarrow
|A\cap B|<\sigma]$. Or just for no uniform ultrafilter $\DD$ on $\sigma$ do we
have $|\DD\cap\AAA|\ge\sigma$.

$\underline{\rm Then}$ we can find ordinals $i(\alpha ,\e )$ for $\alpha
<\lambda ,\e<\sigma$ such that
\begin{enumerate}
\item[(a)] for $\alpha\not=\beta$, $\{\e:i(\alpha,\e)\neq i(\beta,\e)\}$ is
infinite.  Moreover 

\item[(a)$^+$] for any $\lambda'<\lambda$ for some ultrafilter $\DD$ on
$\sigma$, $\{\lan i(\alpha ,\e ):\e<\sigma\ran /\DD :\alpha <\lambda\}$ has
cardinality $\ge\lambda'$.

\item[(b)] for no ultrafilter $\DD$ on $\sigma$ do we have $\{\lan i(\alpha
,\e ):\e <\sigma\ran /\DD :\alpha <\lambda\}$ have cardinality $\lambda$.
\end{enumerate}
[Why?  Let
$$\begin{array}{l}
\displaystyle{\lambda=\sum_{\zeta <\theta}\lambda_\zeta,}\\
\displaystyle{\lambda_\zeta<\lambda,\lambda_\zeta=\max\mbox{ pcf}(\aaa_\zeta),
|\aaa_{i_\zeta}|=\sigma,}\\
\displaystyle{[\aaa_\zeta ]^{<\sigma}\subseteq J_{<\lambda_i}[\aaa_\zeta],}\\
\displaystyle{\mu\ge\sup(\mbox{pcf}(\aaa_\zeta)\backslash\{\lambda_\zeta\}).}
\end{array}
$$
Let $f^\zeta_\alpha\in \Pi\aaa_{\zeta}$ for $\zeta<\theta$, $\alpha<
\lambda_\zeta$ be such that $\lan f^\zeta_\alpha:\alpha<\lambda_\zeta\ran$ is
$< J_{<\lambda_\zeta }[\aaa_\zeta ]$-increasing cofinal and $\bbb\in J_{<
\lambda_\zeta}(\aaa_\zeta)\Rightarrow\mu\ge|\{f^\zeta_\alpha\restriction\bbb:
\alpha <\lambda_\zeta\}|$.  Let $\AAA =\{A_\zeta :\zeta <\theta\}$, let
$\aaa_\zeta =\{ \tau^\zeta_\e :\e\in A_\zeta\}$.  Lastly $i(\alpha ,\e )$ is
$$\begin{array}{ll}
\displaystyle{f^\zeta_\alpha (\e )}&\quad\mbox{if $\displaystyle{\bigcup_{\xi
\le\zeta}\lambda_\xi\le\alpha<\lambda_\zeta}$ \& $\e\in A_\zeta$}\\
\zeta &\quad\mbox{if $\displaystyle{\bigcup_{\xi<\zeta}\lambda_\xi\le\alpha
<\lambda_\zeta}$ \& $\e\notin A_\zeta$}\end{array}
$$
\smallskip

Now check.
\smallskip

\n {\bf 3.11A Remark}.  There are easy sufficient conditions: if $2^\sigma 
<\mu^1\le\mu$, ${\rm cf}(\mu^1)=\sigma$, ${\rm pp}(\mu^1)\ge\lambda$, 
($\forall\chi <\mu^1)({\rm cf}(\chi )\le\sigma\rightarrow {\rm pp}(\chi 
)<\mu^1$) and $\lambda <\mu^{+\omega}$ or at least $\lambda=\sup\big\{\chi:\mu
<\chi ={\rm cf}(\chi )<\lambda$ and $\neg (\exists\aaa )(\aaa\subset{\rm Reg}
\cap\chi\backslash\mu$ \& $|\aaa|\le\sigma$ \& $\chi\in {\rm
pcf}(\aaa))\big\}$. 
\medskip

\n {\bf 3.12 Example}.  Assume
\begin{enumerate}
\item[(a)] $\aleph_0<\theta = {\rm cf}(\lambda )\le 2^{\aleph_0}<\mu <\lambda$
\item[(b)]  there is a $\theta$-Luzin subset of ${}^\omega 2$.
\end{enumerate}
$\underline{\rm Then}$
\begin{enumerate}
\item[$\alpha$)] there are $a_\alpha\in\cBB (\mu )$ for $\alpha <\lambda$,
pairwise disjoint and for no $X\in [\lambda ]^\lambda$ is $\lan a_\alpha:
\alpha\in X\ran$ free.
\item[$\beta$)] Moreover for $X\in[\lambda ]^\lambda$ for some $n<\omega$ and
$\beta_0<\beta_1<\cdots <\beta_n$ from $X$ we have $\cBB (\lambda)\vDash
\bigcap_{\ell\le n} a_{\beta_\ell}=0$.  
\end{enumerate}
\smallskip

\n {\bf Proof}: (Has already appeared in Plebanek \cite{Pl1}.) By 3.10 it
suffices to prove its assumption. Let for $n<\omega$, $\lan c_{n,\ell}:\ell
<(n+1)^2\ran$ be a sequence of pairwise disjoint members of $\cBB (\omega)$
with union $1$, each with each with measure $1/n^2$.  For $\eta\in\prod_{n<
\omega}(n+1)^2$ let $b_\eta =\bigcap_{n<\omega}(1-c_{\eta ,\eta (\ell )})$.
Now suppose 
\begin{enumerate}
\item[($*$)] $X\subseteq {}^\omega 2$, $|X|=\theta$ and if $Y\in [X]^\theta$
then for some $n<\omega$ and $\nu\in\prod_{\ell <n}(\ell +1)^2$ we have
$$
\{\ell :\ell <(n+1)^2\}=\{\eta (n):\eta\restriction n=\mu,\eta\in Y\}.
$$
\end{enumerate}
So $\{ b_\eta :\eta\in X\}$ is as required. Lastly from clause (c) of the
assumption there is $X$ as required in ($*$) so, we are done.
{\hfill{$\square_{3.12}$}
\medskip

\n {\bf 3.13 Remark}. 1) So we can weaken clause (c) of the assumption to
($*$) from the proof, or variants of it.

2) Note that strong negation of (c) of 3.12 which is consistent, implies the
inverse situation.

\section{The interesting ideals and the direct pcf application}

Our problem, the existence of $(\lambda ,I,J)$-sequences for $\bar I$, depends
much on the ideals $I_i$ we use. Under strong set theoretic assumptions, there
are $\lambda$-sequences $\bar \eta$ by 1.8 (and 3.1); but we would like to
prove their existence (i.e.\ in ZFC). For some ideals, by \cite{Sh:g} we will
have many cases of existence, e.g.\ when $I_i$ is $J^{\rm bd}_{\lambda_i}$,
$\lambda_i$ regular.  But we are more interested in the existence for more
complicated ideals.  The first step up are $J^{{\rm bd}}_{\bar\lambda}$ with
$\bar\lambda$ a (finite) strictly increasing sequence of cardinals. The proof
for them is not much harder than with the $J^{{\rm bd}}_\lambda$'s. We then
consider the central ideal here: $J^{{\rm bd}}_{\bar\lambda}$ for
$\bar\lambda$ a (strictly) decreasing sequence of regular cardinals, and
explain why the existence of $\bar{\eta}$ for these ideals is more useful. We
also consider their strong relative which comes from the nonstationary
ideal. We would of course love to have even stronger ideals but there are
indications that for those which we considered and failed, the failure is not
completely due to incompetence, i.e. there are related independence results
(see later). We commence this section by reviewing some general definitions,
some of them used earlier in the paper.
\medskip

\n {\bf 4.1 Definition}. {\it 1) For a set $A$ of ordinals with no last
element (mainly $A=\lambda=\mbox{ cf}(\lambda )$)
$$
J^{{\rm bd}}_A=\{ B:B\subseteq A\;\mbox{ bounded}\}.
$$
2) If $A\subseteq\mbox{ Ord}$ is such that ${\rm cf}(\mbox{otp}(A))>\aleph_0$
and $A$ stationary in $\sup (A)$, we let
$$
J^{\mbox{nst}}_A=\{ B\subseteq A:\mbox{ $B$ is not a
stationary subset of $\sup (A)$}\}.
$$
3) If $A\subseteq \mbox{Ord}$, $\theta =\mbox{ cf}(\theta )<\mbox{cf}(
\mbox{otp}(A))$ and
$$
\{\delta <\sup (A):\delta\in A, \mbox{ cf}(\delta)=\theta\}
$$
is a stationary subset of $\sup (A)$, then let
$$
J^{\mbox{nst},\theta}_A=\{ B\subseteq A:\{\delta\in B:\mbox{ cf}(\delta
)=\theta\}\;\mbox{ is a nonstationary subset of $\sup (A)$}\}.
$$
\smallskip

\n {\bf 4.2 Definition}. {\it 1) For an ideal $J$ let $(\exists^{J^+}x)\varphi
(x)$ mean that
$$
\{ x\in\mbox{ Dom} (J):\varphi (x)\}\in J^+.
$$
2) For an ideal $J$ let $(\forall^Jx)\varphi (x)$ mean
$$
\{ x\in\mbox{ Dom}(J):\neg\varphi (x)\}\in J.
$$
}
\smallskip

\n {\bf 4.3 Definition}. {\it 1) $J=\prod_{\ell <n} J_\ell$ is the following
ideal on $\prod_{\ell <n}\mbox{ Dom}(J_\ell )$: for $X\subseteq\prod_{\ell
<n}\mbox{Dom}(J_\ell )$ we have
$$
X\in J^+\mbox{ iff } (\exists^{J^+_0}x_0)(\exists^{J^+_1}x_1)\cdots
(\exists^{J^+_{n-1}}x_{n-1})[\lan x_0,\ldots ,x_{n-1}\ran\in X].
$$
2) If $\bar\lambda=\lan\lambda_\ell :\ell <n\ran$ we let:
\begin{enumerate}
\item[a)] $J^{{\rm bd}}_{\bar\lambda}=\prod_{\ell <n}J^{{\rm
bd}}_{\lambda_\ell}$.
\item[b)] if ${\rm cf}(\lambda_\ell )>\aleph_0$ for $\ell<n$ then we let
$$
J^{\mbox{nst}}_{\bar\lambda}=\prod_{\ell<n}J^{\mbox{nst}}_{\lambda_\ell}.
$$
\item[c)] if ${\rm cf}(\lambda_\ell )>\theta={\rm cf}(\theta )$ for $\ell <n$
then we let 
$$
J^{\mbox{nst},\theta}_{\bar\lambda}=\prod_{\ell<n}J^{\mbox{nst},\theta}_{
\lambda_\ell}
$$
\item[d)] if $\bar\theta =\lan\theta_\ell :\ell <n\ran$ and ${\rm cf}(
\lambda_\ell )>\theta_\ell={\rm cf}(\theta_\ell )$ for $\ell <n$ then we let
$$
J^{\mbox{nst},\bar\theta}_{\bar\lambda} =\prod_{\ell<n}J^{\mbox{nst},
\theta_\ell}_{\lambda_\ell}.
$$
\end{enumerate}
}
\smallskip

\n {\bf 4.4 Claim}. If $\bar\lambda =\lan\lambda_\ell :\ell <n\ran$ is a
strictly increasing sequence of regular cardinals then the following
conditions (a)--(d) on $X\subseteq\prod_{\ell <n}\lambda_\ell =\mbox{
Dom}(J^{{\rm bd}}_{\bar\lambda})$ are equivalent:
\begin{enumerate}
\item[(a)] $X\in (J_{\bar\lambda}^{{\rm bd}})^+$.
\item[(b)] for no $\bar\alpha\in\prod_{\ell <n}\lambda_\ell$ do we have
$$
(\forall\bar\beta\in X)(\neg (\bar\alpha <\bar\beta))\;\mbox{ where $\bar\beta
<\bar\alpha=:\bigwedge_{\ell <n}\beta_\ell<\alpha_\ell$.}
$$
\item[(c)] We can find $\lan\alpha_\eta:\eta\in\bigcup_{m\le n}\prod_{\ell<m}
\lambda_\ell\ran$ such that:
\item[(i)] $\alpha_\eta <\lambda_{\mbox{lg}(\eta )}$
\item[(ii)] $\alpha_{\eta \hat{\mbox{~}}\lan i\ran}<\alpha_{\eta\hat{\mbox{~}}
\lan j\ran}$ for $i<j<\lambda_{\mbox{lg}(\eta )+1}$
\item[(iii)] $\eta\in\prod_{\ell <n}\lambda_\ell\Rightarrow\lan\alpha_{\eta
\restriction \ell}:\ell\le n\ran\in X$.
\item[d)] Like (c), adding
\item[(iv)] $\alpha_\eta =\alpha_\nu\Rightarrow\eta=\nu$.
\end{enumerate}
}
\medskip

\n {\bf Proof}. Straight.  For (b) $\Rightarrow$ (c) use induction on
$n=\lg(\bar{\lambda})$, see the proof at the end of the proof of 4.11, of
($*$) there.{\hfill$\square_{\rm 4.4}$}
\medskip

\n {\bf 4.4A Discussion}. From 4.4, we see that for for $X\in (J^{{\rm
bd}}_{\bar\lambda})^+$ there are patterns which necessarily occur as subsets
of $X$. These are essentially like the branches ($=$ maximal nodes) of a tree
with $n$ levels, with a linear order on each level and with no dependencies
between the different levels. These patterns were explored in \cite{Sh:462},
\cite{RoSh:534}, \cite{Sh:575}. The patterns considered there can be
represented as a set $\Delta\subseteq\prod_{\ell <n} B_\ell$, $B_\ell\subseteq
\mbox{ Ord}$ such that $\eta (i)=\nu (i)\Rightarrow\eta\restriction i=\nu
\restriction i$ (i.e.\ treeness). Now look at $J^{{\rm bd}}_{\bar\lambda}$,
where the gain is that $\Delta$ does not have a tree, that is, we have any
$\Delta\subseteq\prod_{\ell <n}B_\ell$, $B_\ell\subseteq\mbox{ Ord}$, so that
$\eta,\nu\in\Delta$ can have $\{\ell <n:\eta (\ell )=\nu (\ell )\}$ being
arbitrary (rather than being an initial segment), of course this depends on
the ideal.} 
\medskip

\n {\bf 4.5 Claim}: Assume $\bar J=\lan J_\ell :\ell <n\ran$ and $J_\ell$ is a
$\kappa_\ell$-complete ideal on $\lambda_\ell$. We also demand $\kappa_\ell
>\lambda_k$ when $\ell > k$. Let $J=\prod_{\ell <n}J_\ell$.

1) The following conditions on $X\subseteq\prod_{\ell <n}\lambda_\ell$ are
equivalent
\begin{enumerate}
\item[(a)] $X\in J^+$
\item[(b)] for no $\bar A=\lan A_\ell :\ell <n\ran$, $A_\ell\in J_\ell$ do we
have 
$$
\bar\beta\in X\Rightarrow \bigvee_\ell\beta_\ell\in A_\ell.
$$
\item[(c)] We can find $\lan\alpha_\eta :\eta\in \bigcup_{m\le n}\prod_{\ell
<m}\lambda_\ell\ran$ such that $\alpha_\eta <\lambda_{{\rm lg}(\eta )}$ and
\begin{enumerate}
\item[($*$)] for each $\nu\in\prod_{\ell <n}\lambda_\ell$ we have
$$
\lan\alpha_{\nu\restriction (\ell +1)}:\ell < n\ran\in X.
$$
\end{enumerate}
\end{enumerate}

2) If $[A\subseteq\lambda_\ell$ \& $|A|<\lambda_\ell]\Rightarrow A\in J_\ell$
then we can add
\begin{enumerate}
\item[(d)] like (c), but adding
\begin{enumerate}
\item[(iii)] $\alpha_{\nu\hat{\mbox{~}}\lan i\ran}<\alpha_{\nu\hat{\mbox{~}}
\lan j\ran}$ if $i<j<\lambda_{\mbox{lg}(\nu)+1}$.
\end{enumerate}
\end{enumerate}

\n {\bf Proof}. Similar to 4.4.
\medskip

\n {\bf 4.6 Claim}. Let $\bar\lambda =\lan\lambda_\ell :\ell <n\ran$ be a
decreasing sequence of regular cardinals
\begin{enumerate}
\item[(1)] If $\lambda_\ell >2^{\lambda_{\ell +1}}$ for $\ell <n$, then:
\begin{enumerate}
\item[($*$)] for every $A\in (J^{{\rm bd}}_{\bar\lambda})^+$, there are
$A_\ell\in (J^{{\rm bd}}_{\lambda_\ell} )^+$ such that $\prod_{\ell <n}
A_\ell\subseteq A$.
\end{enumerate}

\item[(2)] If $J=\prod_{\ell <n} J_\ell$ and $J_\ell$ is a $(2^{\lambda_{\ell
+1}})^+$-complete ideal on $\lambda_\ell$, $\underline{\rm then}$ $(\ast)$
holds, with $J$ in place of $J^{{\rm bd}}_{\bar\lambda}$ and $J_\ell$ in place
of $J^{\rm bd}_{\lambda_\ell}$.

\item[(3)] For every $A\in (J^{{\rm bd}}_{\bar\lambda})^+$ and $k<\omega$ we
can find $B_\ell\in [\lambda_\ell ]^k$ such that $\prod_{\ell <n}
B_\ell\subseteq A$.

\item[(4)] In (3), instead of $k$ and $J^{{\rm bd}}_{\lambda_\ell}$ (for $\ell
<n$) we can use any $\kappa$ and $((\lambda_{\ell +1})^\kappa )^+$-complete
ideal $J_\ell$ on $\lambda_\ell$ for $\ell <n$.
\end{enumerate}

\n {\bf Proof}. E.g. (3) We prove it by induction on $n$.

$\underline{n=1}$. Trivial, as singletons are in the ideal.

$\underline{n+1}$. Let $X_0=\left\{ \alpha<\lambda_0:\{\bar\alpha\in
\prod^{n-1}_{\ell =1}\lambda_\ell:\lan\alpha\ran\hat{\mbox{~}}\bar\alpha\in
A\}\in (\prod^{n-1}_{\ell =1}J^{\rm bd}_{\lambda_\ell})^+\right\}$.

Clearly, $X_0\in (J^{\rm bd}_{\lambda_0})^+$.

By the induction hypothesis, for each $\alpha\in X_0$, there is $\lan
B^\alpha_\ell : \ell =1,\ldots ,n-1\ran$, such that 
$$
B^\alpha_\ell\in [\lambda_\ell ]^k\;\;{\rm and}\;\; \prod^{n-1}_{\ell=1}
B^\alpha_\ell\subseteq\{\bar\alpha\in\prod^{n-1}_{\ell=1}\lambda_\ell:\lan
\alpha\ran^{\hat{~}}\bar\alpha\in A\}\deq\bar{B}^\alpha.
$$
So $X_0$ is the union of $\prod^{n-1}_{\ell =1}\lambda^k_\ell=\lambda_1$ sets
$X_0[\bar B]=\{\alpha\in X_0:\bar B^\alpha =\bar B\}$, so for some $\bar B$,
$|X_0[\bar B]|\ge k$ and let $B_0=$ first $k$ members of $X_{0,\bar
B}$.{\hfill$\square_{\rm 4.6}$}
\medskip

\n {\bf 4.7 Definition}. {\it For a partial order $P$ let
$\mbox{tcf}(P)=\lambda$ $\underline{\rm iff}$ there is an increasing cofinal
sequence of length $\lambda$ in $P$ ($\mbox{tcf}$ -- stands for true
cofinality);  so e.g.\ $(\omega , <)\times (\omega_1, <)$ has no true
cofinality, but $\mbox{tcf }\prod(\aleph_n, <)/\DD$ is well defined if $\DD$
is an ultrafiler on $\omega$.} 
\medskip

\n {\bf 4.9 Fact}. 1) ~If $J\supseteq J^{{\rm bd}}_\delta$ is an ideal,
$\lambda_i=\mbox{ cf}(\lambda_i)>\delta$, for $i<\delta$ and $\lambda
=\mbox{tcf}(\prod_{i<\delta}\lambda_i/J)$, $\underline{\rm then}$ there is a
$(\lambda ,J)$-sequence $\bar\eta =\lan\eta_\alpha :\alpha <\lambda\ran$ for
$\lan J^{{\rm bd}}_{\lambda_i}:i<\delta\ran$.

2) If $\lambda_i$ is increasing in $i$ then $\lan J^{\rm bd}_{\lambda_i}:i<
\delta\ran$ is normal (hence $\bar\eta$ is normal) provided that $\delta=
\omega$ or at least
\smallskip

($*$)$_1$ $\qquad\displaystyle{\lambda >\prod_{j<i}\lambda_j}$ for $i<\delta$
\smallskip

3) If we just ask $\bar\eta$ to be normal it suffices to demand
\smallskip

($*$)$_2$ $\qquad\displaystyle{\lambda_i>{\rm max ~pcf}\{\lambda_j :j<i\}}$ 
for $i<\delta$.
\medskip

\n {\bf Proof}. In $\prod_{i<\delta}\lambda_i /J$, there is a cofinal
increasing sequence $\lan f_\alpha :\alpha <\lambda\ran$. It is as required,
as we now show. Let $X\in [\lambda ]^\lambda$, let $X_i=\{ f_\alpha(i):\alpha
\in X\}$ for $i<\delta$. Define $f\in\prod_{i<\delta}\lambda_i$
$$
f(i)=\left\{\begin{array}{ll}
\sup (X_{i})+1&\quad\underline{\rm if}\,\, \sup (X_i)<\lambda_i\\
0&\quad\mbox{otherwise}.\end{array}\right.
$$
But $\lan f_\alpha :\alpha <\lambda\ran$ is cofinal, so for some
$\alpha_0<\lambda$, $f<_J f_{\alpha_0}$. Now $X\in [\lambda ]^\lambda$, so for
some $\alpha_1$, we have $\alpha_0<\alpha_1\in X$. As $\lan f_\alpha :\alpha
<\lambda\ran$ is increasing, $f_{\alpha_0} <_J f_{\alpha_1}$, hence $f<_J
f_{\alpha_1}$. So $A=\{ i:f(i)\ge f_{\alpha_1}(i)\}\in J$. But $f_{\alpha_1}
(i)\in X_i$, so $i\in\delta\backslash A\Rightarrow\lambda_i=\sup (X_i)$.

2) Easy

3) By \cite[II, 3.5]{Sh:g}.{\hfill$\square_{\rm 4.9}$}
\medskip

\n {\bf 4.10 Comment}. 1) ~This is good e.g.\ to lift a colouring of the
$\lambda_i$'s to one of $\lambda$. But we would like to have an upgrade as
well.

2) ~The kind of assumptions of 4.9 is the central interest in \cite{Sh:g}.
\medskip

\n {\bf 4.11 Claim}. Assume $\bar\lambda^i=\lan\lambda_{i,\ell}:\ell <n_i\ran$
is an increasing sequence of regulars $>\delta$ for $i<\delta$. Also assume
that $J$ is an ideal on $\{ (i,\ell ):i<\delta,\ell <n_i\}$ and
$$
\lambda=\mbox{tcf}(\prod_{i,\ell}\lambda_{i,\ell}/J),
$$
and for some ideal $J'$ on $\delta$, we have $J'\supseteq J^{{\rm bd}}_\delta$
and $J$ is generated by
$$
\big\{\{ (i,n):n<n_i, i\in A\} :A\in J'\big\}.
$$
$\underline{\rm Then}$ there is a $(\lambda, J')$-sequence $\bar\eta$ for
$\lan J^{{\rm bd}}_{\bar\lambda^i}:i<\delta\ran$.
\smallskip

2) $\lan J^{\rm bd}_{\bar\lambda^i}:i<\delta\ran$ is normal (hence $\bar\eta$
above is normal) if
\smallskip

($*$)$_1$ ~$\delta=\omega$ and $i<j<\delta\Rightarrow\lambda_{i,n_i-1}<
\lambda_{j,0}$.
\smallskip

$\underline{\rm or}$
\smallskip

($*$)$_2$ $\qquad\displaystyle{\prod\{\lambda_{i,\ell}:i<j,\ell<n_j\}<
\lambda_{j,0}}$
\smallskip

3) If we ask just $\bar\eta$ to be normal it suffices to demand
\smallskip

($*$)$_3$ $\qquad\displaystyle{{\rm max ~pcf}\{\lambda_{i,\ell}:i<j,\ell<n_j\}
<\lambda_{j,0}}$. 
 \medskip

\n {\bf Proof}. Again, let $\bar f=\lan f_\alpha :\alpha <\lambda\ran$ be
$<_J$-increasing cofinal. Let $\eta_\alpha (i)=\lan f_\alpha (i,\ell ):\ell
<n_i\ran\in\prod\bar\lambda^i$. Let $X\in [\lambda ]^\lambda$. Let
$X_i=\{\eta_\alpha (i):\alpha\in X\}$. If $X_i\in J^{{\rm bd}}_{\prod\bar
\lambda^i}$, then there is $\bar\alpha^i\in\prod\bar\lambda^i=\prod_{\ell<
n_i}\lambda_{i,\ell}$ such that
$$
\bar\beta\in X_i\Rightarrow \bigvee_{\ell <n_i}\beta_\ell <\alpha^i_\ell.
\leqno\mbox{($*$)}
$$
(We return to this at the end of the proof.)

So let $f\in\prod_{i,\ell}\lambda_{i,\ell}$ be given by $f(i,\ell)=
\alpha^i_\ell$. So, as before, for some $\alpha\in X$, $f<_J f_\alpha$. So
$$
A=\left\{ i:\bigwedge_{\ell <n_i} f((i,\ell ))\ge f_{\alpha}((i,\ell ))
\right\}\in J'.
$$
Now for $i\in\delta\backslash A$ we have $X_i\notin J^{{\rm bd}}_{\prod\bar
\lambda^i}$.
\medskip

\n [Why ($*$)? Prove the existence of $\bar\alpha^i$, for notational
convenience denoted here by $\bar\beta$, by induction on $n_i$. 
$\underline{\rm Here}$ we use ``increasing $\bar\lambda^i$''.
\medskip

$\underline{n_i=1}$. Clear
\medskip

$\underline{n_i=k+1}$. For $\alpha <\lambda_{i,0}$ define
$$
X_{i,\alpha}=\{\bar\beta\restriction [1,n_i):\bar\beta\in X_i\}.
$$
So we know that for some $\gamma_0<\lambda_{i,0}$
$$
\alpha \in [\gamma_0, \lambda_{i,0}]\Rightarrow X_{i,\alpha}\in
J^{{\rm bd}}_{\prod^{n-1}_{\ell =1}\lambda_{i,\ell}}.
$$
So for each such $\alpha$ we have $\bar\beta^\alpha\in\prod^n_{\ell=1}
\lambda_{i,\ell}$ as given by the induction hypothesis. Let
$$
\beta_\ell =\left\{\begin{array}{ll}
\gamma_{0}+1&\mbox{if $\ell =1$}\\
\cup\{\beta^\alpha_\ell :\alpha\in [\gamma_0,\lambda_{i,0})\}&
\mbox{otherwise.}\end{array}\right.
$$
Why is the latter $<\lambda_{i,\ell}$?  As $\lambda_{i,0}<\mbox{cf}(
\lambda_{i,\ell})$.]

~{\hfill$\square_{\rm 4.11}$}
\medskip

\n {\bf 4.12 Question}. Are there many cases fitting the framework of 4.11?
\medskip

\n {\bf 4.12A Answer}. Not so few. E.g.\ for any $\kappa$, for many $\lambda
={\rm cf}(\lambda )$ we have that $\lambda=\mbox{tcf}\Big(\prod_{i<\kappa}
\lambda_i/ J^{{\rm bd}}_\kappa\Big)$ for some sequence $\langle\lambda_i:\,i
<\kappa\rangle$. E.g.\ if $\aleph_0<\mbox{ cf}(\delta )=\kappa$ and $\kappa
<\mu =\beth_\delta <\lambda ={\rm cf}(\lambda)\le\beth_{\delta +1}$ or just
$\aleph_0<\kappa ={\rm cf }(\mu) <\lambda= {\rm cf}(\lambda)\le \mu^\kappa$
and $(\forall \chi <\mu )$ $[\chi^\kappa <\mu]$ $\underline{\rm then}$ there
is an increasing sequence of regulars $\lan\lambda_i:i<\kappa\ran$ with limit
$\beth_\delta$ or $\mu$ respectively as above. [Why? see \cite[VIII \S1,
2.6]{Sh:g}.] Even if $\kappa=\aleph_0$ this holds for many $\lambda$'s e.g.\
if $\mu <\lambda <\mu^{+\omega_1}$ or just $|\{\chi :\mu <\chi <\lambda\; {\rm
and}\; \chi=\aleph_\chi\}|<\mu$ see \cite[IX]{Sh:g} and use 4.13 below.

Note that by the $\mbox{pcf}$-theorem (see \cite[VIII, 2.6]{Sh:g})
\medskip

\n {\bf 4.13 Claim}. Assume $I$ to be an ideal on $\delta$, and $\lambda_{i,
\ell}=\mbox{ cf}(\lambda_{i,\ell})>|\delta |$ for $i<\delta$ and $\ell<n_i$
and $0< n_i<\omega$. $\underline{\rm Then}$ the following are equivalent
\begin{enumerate}
\item[(a)] for every $\lan k_i:i<\delta\ran\in\prod_{i<\delta}n_i$ we have
$$
\lambda =\mbox{tcf}\Big(\prod_{i<\delta}\lambda_{i,k_i}/I\Big).
$$

\item[(b)] letting
$$
I'=\Big\{ A\subseteq\bigcup_{i<\delta}\{ i\}\times n_i:\mbox{for some $B\in I$
we have $\displaystyle{A\subseteq\bigcup_{i\in B}\{ i\}\times n_i\Big\}}$},
$$
we have $\prod\lambda_{i,n}/I'$ has true cofinality $\lambda$.
\end{enumerate}
\medskip

{\bf Proof.} Let $A^*,B^*$ be a partition of $\bigcup\limits_{i<\theta}\{i\}
\times n$ such that $\lambda={\rm max~pcf}\{\lambda_{i,n}:(i,n)\in A^*\}$ and
$\lambda\notin {\rm pcf}\{\lambda_{i,n}: (i,n)\in B^*\}$ (exists by the pcf
theorem). Now:

\underline{(a)$\Rightarrow$(b)}

If $\prod_{i,n}\lambda_{i,n}/I'$ does not have true cofinality $\lambda$, then
for some $A\in (I')^+$ we have that $\prod_{(i,n)\in A}\lambda_{i,n}/I'$ has
true cofinality $\lambda'\neq\lambda$ (here we use the pcf theorem) and wlog
$A\subseteq A^*\ \vee\ A\subseteq B^*$, hence $\lambda\notin{\rm pcf}\{
\lambda_{i,n}:(i,n)\in A\}$. Let $B=\{i<\delta:\,(\exists n<n_i)[(i,n)\in A]
\}$, so by the definition of $I'$ we know $B\in I^+$. So, for $i\in B$ we can
choose $k_i\in\{ 0,\ldots ,n_i-1\}$ such that $(i,k_i)\in A$. So $\{ (i,k_i):
i\in B\}\subseteq A$ hence ${\rm pcf}\{\lambda_{i,k_i}:i\in B\}\subseteq {\rm
pcf}\{\lambda_{i,k}:(i,k)\in A\}$, but $\lambda$ does not belong to
the later, hence not to the former, contradicting (a).
\medskip

\n $\underline{\neg{\rm (a)}\Rightarrow\neg{\rm (b)}}$
\medskip

So there is $\lan k_i:i<\delta\ran\in\prod_{i<\delta}n_i$ such that $\neg[{\rm
tcf}(\prod\lambda_{i,k_i}/I)=\lambda ]$ hence by the pcf theorem, for some
$A\in (I)^+$, we have ${\rm max ~pcf}\{ \lambda_{i,k_i}:i\in A\}<\lambda$, let
$B=\{ (i,k_i):i\in A\}$, so clearly ${\rm max~pcf}\{\lambda_{i,k_i}:(i,k_i)\in
B\}<\lambda$. But by the definition of $I'$, we have $B\in (I')^+$ so we get
contradiction to (b).{\hfill$\square_{\rm 4.13}$}
\medskip

\n {\bf 4.14 Remark}. See more on related topics in \cite{Sh:589}. 

\section{$\lambda$-sequences for decreasing $\bar\lambda^i$ by pcf.}

\n {\bf 5.1 Discussion}. Our aim here is to get ``decreasing $\bar\lambda$''
from ``increasing $\bar\lambda$'' (for $J^{{\rm bd}}_{\bar\lambda}$), in some
sense, to ``make gold from lead''. We do this by using pcf assumptions, then
proving that these assumptions are very reasonable.

(Note: when we cannot materialize the pcf assumptions the situation is close
to SCH, and then we have other avenues for construction of $\lambda$-sequences
for some $I$, e.g. (1.8, 3.1).)
\smallskip

\centerline{ * * * }
\smallskip

In the following claim the interesting case is when $\lambda_\ell$ are
increasing, $\bar\lambda^i=\lan\lambda_{\ell ,i}:\ell <n\ran$ decreasing
sequence of regular cardinals, $\lambda_{\ell ,i}>\prod_{j<i\atop
m<n}\lambda_{m,j}$, or at least $\lambda_{\ell ,i}>\max\mbox{ pcf}\{
\lambda_{m,j}:m<n, j<i\}$.
\medskip

\n {\bf 5.2 Claim}. Assume
\begin{enumerate}
\item[(a)] $\bar \lambda=\lan\lambda_\ell :\ell < n\ran$,
$\bar\lambda^i=\lan\lambda_{\ell ,i}:\ell <n\ran$ for $i<\delta$.
\item[(b)] $I$ is an ideal on $\delta$.
\item[(c)] $\lambda_\ell =\mbox{tcf}\displaystyle{\Big(\prod_{i<\delta}
\lambda_{\ell ,i}/I\Big)}$ for $\ell <n$.
\item[(d)] $\bar f^\ell =\lan f_{\ell ,\alpha}:\alpha<\lambda_\ell\ran$ is
$<_I$-increasing and cofinal in $\displaystyle{\prod_{i<\delta}\lambda_{\ell
,i}}$. 
\item[(e)] $\delta <\lambda_{\ell ,i}=\mbox{cf}(\lambda_{\ell ,i})$.
\item[(f)] for $\displaystyle{\bar\alpha\in\prod_{\ell <n}\lambda_\ell}$ let
$f_{\bar\alpha}$ be defined by $f_{\bar\alpha}(i)=\lan f_{\ell,\alpha_\ell}
(i):\ell <n\ran\in\displaystyle{\prod_{\ell <n}\lambda_{\ell,i}}$.
\end{enumerate}
$\underline{\rm Then}$ for any $X\in (J^{{\rm bd}}_{\bar\lambda})^+$ we have
$$
\big\{ i:\{ f_{\bar\alpha}(i):\bar\alpha\in X\}\in J^{{\rm bd}}_{\bar
\lambda^i}\big\}\in I.
$$
\smallskip

\noindent {\bf Proof}. Let $X_i=\{ f_{\bar\alpha}(i):\bar\alpha\in X\}$ and
let $B=\{ i<\delta :X_i\in J^{{\rm bd}}_{\bar\lambda^i}\}$. 

Assume $B\in I^+$ and we shall get a contradiction. For each $i\in B$, $m<n$
and $\bar\alpha\in \displaystyle{\prod_{\ell <m}\lambda_{\ell ,i}}$, let
$$
X^i_{\bar\alpha}=\Big\{\bar\beta\in\prod^{n-1}_{\ell =m}\lambda_{\ell ,i}:
\bar\alpha\hat{\phantom{A}}\bar\beta\in X_i\Big\}
$$
and let $g_i(\bar\alpha )=\min\Big\{\gamma\le\lambda_{m,i}:$ if $\beta\in
[\gamma ,\lambda_{m,i})$ then $X^i_{\bar\alpha\hat{\mbox{~}}\lan\beta\ran}\in
J^{{\rm bd}}_{\prod_{\ell=m+1}^{n-1}\lambda_{\ell,i}}\Big\}$.

This definition just unravels the definition of $J^{{\rm
bd}}_{\bar\lambda^i}$; note
\begin{enumerate}
\item[($*$)] $X^i_{\lan\;\ran}=X_i\in J^{{\rm bd}}_{\bar\lambda^i}$
\item[($*$)$'$] if $X^i_{\bar\alpha}\in J^{{\rm bd}}_{\prod_{\ell \ge {\rm
lg}(\bar\alpha)}\lambda_{\ell ,i}}$ then $g_i(\bar\alpha )<\lambda_{{\rm
lg}(\bar\alpha )}$.
\end{enumerate}
Now we choose by induction on $m< n$ ordinals $\alpha_m<\lambda_m$ such that
for $m\le n$ we have
\begin{enumerate}
\item[($*$)$_m$] $B_m= :\big\{ i\in B:X^i_{\lan f_{\ell,\alpha_\ell}(i):\ell<
m\ran}\in J^{{\rm bd}}_{\prod_{\ell \ge m}\lambda_{\ell,i}}\big\} =B\mbox{ mod
}I$. 
\end{enumerate}
So, stipulating $J^{\rm bd}_{\prod_{\ell\ge n}\lambda_{\ell,i}}=\{\emptyset
\}$, the ideal on $\{\lan\;\ran\}$, we have that ($*$)$_0$ holds with $B=B_0$.

If ($*$)$_m$ is true, clearly
$$
\lan g_i(\lan f_{\ell ,\alpha_\ell}(i):\ell <m\ran ):i\in B_m\ran
$$
is in $\prod_{i<\delta}\lambda_{m,i}$. But $B_m\in I^+$ and $\lan f_{m,\alpha}
:\alpha <\lambda_m\ran$ is $<_I$-increasing cofinal in $\prod_{i<\delta}
\lambda_{m,i}$. So for some $\alpha_m$
$$
B'_m=\{ i\in B_m:g_i(\lan f_\ell,\alpha_\ell(i):\ell<m\ran)\ge\alpha_m\}\in I.
$$
Defining $B_{m+1}$ using this $\alpha_m$, we easily obtain
$$
B_{m+1}\supseteq B_m\backslash B'_m\;\mbox{ so we see that $(*)_{m+1}$ holds.}
$$
So
$$
\bar\alpha=\lan\alpha_\ell :\ell <n\ran\in\prod_{\ell <n}\lambda_\ell
$$
is well defined.

In the inductive definition of $\alpha_m$, any larger $\alpha'_m$ would serve
in place of $\alpha_m$ (of course it would influence the future choices). So,
in addition to ($*$)$_m$, we can demand
$$
\Big\{\bar\beta\in\prod^{n-1}_{\ell=m}\lambda_\ell:\lan\alpha_\ell:\ell<m\ran
\hat{\mbox{~}}\bar\beta\in X\Big\}\in (J^{\rm bd}_{\bar\lambda\restriction
[m,n)})^+.\leqno\mbox{($**$)$_m$} 
$$
So from $(\ast\ast)_n$ we get $\lan\alpha_\ell :\ell <n\ran\in X$ hence for
all $i$ we have $\lan\;\ran\in X^i_{\lan f_{\ell ,\alpha_\ell}(i):\ell<n
\ran}$, by the definition. But $B_n=\{ i\in B:X^i_{\lan f_{\ell,\alpha_\ell}
(i):\ell <n\ran}\in J^{\rm bd}_{\prod_{\ell\ge n}\lambda_{\ell,i}}\}=B \mbox{
mod }I$, so $B_n\not= \emptyset$, and if $i\in B_n$ this means $X^i_{\lan
f_{\ell ,\alpha_\ell}(i):\ell <n\ran}\in J^{\rm bd}_{\prod_{\ell\ge n}
\lambda_{\ell ,i}}=\{\emptyset\}$ so $X^i_{\lan f_{i,\alpha_\ell}(i):\ell<
\ran}=\emptyset$, contradicting the previous sentence.{\hfill $\square_{\rm
5.2}$}
\medskip

In fact, more generally,
\medskip

\n {\bf 5.3 Claim}. Assume
\begin{enumerate}
\item[a)] $\bar\eta^\ell =\lan \eta^\ell_\alpha :\alpha <\lambda\ran$ in an
$(I,J,\lambda)$-sequence for $\lan I_{i ,\ell}:i<\delta\ran$ for each $\ell
<n$.
\item[b)] $I_i=\prod_{\ell <n} I_{i,\ell}$
\item[c)] $\bar\eta =\lan\eta_\alpha :\alpha <\lambda\ran$ where
$\eta_\alpha\in\prod_{i<\delta}\mbox{ Dom}(I_i)$ and
$$
\eta_\alpha (i)=\lan\eta^\ell_\alpha (i):\ell <n\ran.
$$
\end{enumerate}
$\underline{\rm Then}$ $\bar\eta$ is an $(I,J,\lambda )$-sequence for $\lan
I_i:i<\delta\ran$.
\medskip

\n {\bf Proof}. Like the proof of 5.2.{\hfill$\square_{\rm 5.3}$}
\medskip

\n {\bf 5.4 Claim}. Assume
\begin{enumerate}
\item[a)] $\lambda =\mbox{ tcf}\Big( \prod_{i<\delta}\theta_{\ell ,i}/J\Big)$
for $\ell <n$ and $\theta_{\ell ,i}$ are increasing with $\ell$.
\item[b)] $\theta_{\ell ,i}=\mbox{ tcf}\Big(\prod_{\e <\e_i}\tau_{\ell
,i,\e}/J_i\Big)$and $\tau_{\ell ,i,\e}$ are regular decreasing with $\ell$
i.e.\ $\tau_{\ell ,i,\e}\ge \tau_{\ell+1,i,\e}$ (the interesting case is $>$).
\end{enumerate}
Let
$$\begin{array}{r}
J^*=\Big\{ A:A\subseteq\{ (\ell ,i,\e ):\ell <n, i<\delta
,\e <\e_i\}\\
\mbox{and }\;{\displaystyle{\bigwedge_\ell}}(\forall^Ji)(\forall^{J_i}\e)
[(\ell,i,\e )\notin A]\Big\}\end{array}
$$
Let
$$
I_{i,\e}=\prod_{\ell <n}J^{{\rm bd}}_{\tau_{\ell ,i,\e}}.
$$
$\underline{\rm Then}$
$$
\lambda =\mbox{ tcf}\Big(\prod_{i,\e}\tau_{\ell,i,\e}/J^*\Big)
$$
and we can find $\bar\eta_\alpha \in\prod_{i,\e} I_{i,\e}$ for $\alpha
<\lambda$ such $\lan\eta_\alpha :\alpha <\lambda\ran$ is a $(\lambda,
J^*)$-sequence for $\lan I_{i,\e}:i,\e\ran$.
\medskip

\n {\bf Proof}. Straight. (Using 5.3 and \cite[I,2.10]{Sh:g}.)
{\hfill$\square_{\rm 5.4}$
\medskip

\n {\bf 5.5 Example}. Assume
\begin{enumerate}
\item[($*$)] $\lan\lambda_i:i<\delta\ran$ is a strictly increasing sequence of
regulars, $\delta <\lambda_0$, $\lambda =\mbox{tcf}\big(\prod_{i<\delta}
\lambda_i/J^{{\rm bd}}_\delta\big)$.
\end{enumerate}

\n {\bf 5.5A Discussion}. This may seem a strong assumption, but getting such
representations is central in \cite{Sh:g}. If $\mu$ is strong limit singular
\begin{enumerate}
\item[$\otimes$] $\aleph_0<\kappa ={\rm cf}(\mu )<\mu <\lambda =\mbox{cf}(
\lambda )\le 2^\mu$,
\end{enumerate}
then there is such $\lan\lambda_i:i<\mbox{ cf}(\mu )\ran$, $\lambda_i<\mu
=\sup (\lambda_i)$. So without loss of generality $2^{\lambda_i}<\lambda_{i+
1}$ (see 4.12A).

Now fix $n$ for simplicity. Let
$$
\lambda_{\ell ,i}=\lambda_{n\times i+n-\ell}.
$$
So
$$
\bar\lambda^i=\lan\lambda_{\ell ,i}:\ell <n\ran\;\;\mbox{ is strictly
decreasing.} 
$$
In 4.13 an example is given for 5.2. For 5.4 we have e.g.
\smallskip

\n {\bf 5.6 Claim}. Assume
\begin{enumerate}
\item[(a)] $\mu$ is strong limit.
\item[(b)] $\aleph_0=\mbox{ cf}(\mu )<\mu$.
\item[(c)] $2^\mu\ge \mu^{+\omega +1}=\lambda$ [also $\mu^{+\omega_4+\omega
+1}=\lambda$ is OK, or just $\lambda =\mu^{+\delta+1}<{\rm pp}^+(\mu )$ and
${\rm cf}(\mu^{+\delta})<\mu$] 
\end{enumerate}
Then:
\smallskip

1) We can find $\lambda_{\ell ,i}$, $k_\ell$ such that $(\ell \le i<\omega )$
\begin{enumerate}
\item[(A)] $\lambda_{i,\ell}<\mu=\sum_{m,j}\lambda_{m,j}$.
\item[(B)] $2^{\lambda_{\ell +1 ,i}}<\lambda_{\ell ,i}$ and $2^{\lambda_{0
,i}}<2^{\lambda_{i+1, i+1}}$.
\item[(C)] $\mbox{tcf}\big(\prod_{i<\omega }\lambda_{\ell ,i}/J^{{\rm
bd}}_\omega\big) =\mu^{+k_\ell}$.
\item[(D)] $0<k_m<k_{m+1}<\omega$.
\item[(E)] $\lambda=\mbox{tcf}\Big(\prod_{m<\omega}\mu^{+k_m}/J^{{\rm
bd}}_\omega\Big)$
\end{enumerate}

2) For every $n<\omega$, we can find $J,\lambda'_{\ell ,i}$ $(\ell <n,
i<\omega )$ such that
\begin{enumerate}
\item[(i)] There is $\bar\eta$ a $\lambda$-sequence for
$$
\lan J^{{\rm bd}}_{\lan\lambda'_{\ell ,i}:\ell<n\ran}:i<\omega\ran.
$$
\item[(ii)] $2^{\lambda'_{\ell +1},i}<\lambda'_{\ell ,i}$
\item[(iii)] $2^{\lambda'_{0,i}}<\lambda'_{n-1,i+1}$
\item[(iv)] $(\forall A\in J)(\exists^\infty i)$ $[n\times\{i\}\cap A=
\emptyset]$.
\end{enumerate}
\smallskip

\n {\bf 5.7 Remark}. (1) This Claim can be used with no further reference to
$\mbox{pcf}$: just for any $\mu$ as in (a)--(c), we have $\bar\eta$ for which
we can construct colourings, objects, etc.

(2) There are theorems with $n$ increasing, they are somewhat cumbersome.

Of course, we can use
$$
I'_m=\prod^{n_m+1}_{i=n_m}J^{\rm bd}_{\lan\lambda_{\ell ,i}:\ell <n_m\ran}.
$$

(3) $\underline{\rm Note}$: $2^\mu\ge\mu^{\omega +1}$ is a strong negation of
$2^\mu =\mu^+$ which was very useful here. (Our general theme is: $\neg$SCH is
a good hypothesis) and we shall deal with closing the gap.

(4) $\underline{\rm Note}$: if $2^\mu =\mu^{+n(*)}$, we can prove nice things
with $I=J^{{\rm bd}}_{\lan\mu^{+n(*)-\ell}, \ell <n(*)\ran}$.

(5) If $\aleph_0<{\rm cf}(\mu )<\mu$ the parallel claim is even easier, and
$\mu$ being a strong limit is necessary only for (B).  \medskip

\n {\bf Proof of 5.6}. (1) We will just give a series of quotations.

First ${\rm cf}(\mu^{+\omega})=\aleph_0$, so by \cite[II 1.6]{Sh:g}, there is
an increasing sequence $\lan\theta_i:i<\omega\ran$ of regulars with limit
$\mu^{+\omega}$ such that
$$
\lambda =\mu^{+\omega +1}=\mbox{ tcf}\Big(\prod_{i<\omega}\theta_i/J^{{\rm
bd}}_\omega\Big),
$$
so for $i$ large enough $\theta_i>\mu$. So without loss of generality
$\bigwedge_i\theta_i>\mu$.

So let $\theta_i=\mu^{+k_i}$, $k_i\in (0,\omega )$ strictly increasing. By
\cite[IX, 5.9, p. 408]{Sh:g}, we have ${\rm pp}(\mu )>\mu^{+k_i}$. (We would
like to have ${\rm pp}(\mu )=2^\mu$, but only ``almost proved''.) This means
by the no hole theorem \cite[II, 2.3]{Sh:g} that for some countable set
$\aaa_\ell$ of regulars $<\mu$, $\mu =\sup (\aaa_\ell )$ and $\mu^{+k_\ell}
\in\mbox{pcf}(\aaa_\ell )$. So by the pcf theorem, without loss of generality
$\mu^{+k_\ell}=\max \mbox{ pcf}(\aaa_\ell )$ and $\mu^+,\ldots,\mu^{+(k_\ell
-1)}\notin \mbox{pcf}(\aaa_\ell )$ (alternatively use \cite[VIII, \S1]{Sh:g}). 

So necessarily
$$
\mu^{+k_\ell}=\mbox{ tcf}\big(\Pi \aaa_\ell /J^{{\rm bd}}_{\aaa_\ell}\big).
$$
Let $\mu =\sum_{n<\omega}\mu_n$, $\mu_n<\mu_{n+1}<\mu$. We start choosing
$\lambda_{\ell,i}$ by induction on $i$, for all $i$ by downward induction on
$\ell$, so that
$$
\lambda_{\ell ,i}>\mu_i,\;\;\lambda_{\ell ,i}\in \aaa_i,
$$
and (B) holds.  So, as $\lambda_{\ell,i}\in \aaa_i$ and $\lambda_{\ell ,i}$ is
increasing with $i$, with limit $\mu$, we have
$$
\mbox{tcf}\Big(\prod_i\lambda_{\ell ,i}/J^{{\rm bd}}_\omega\Big)=\mu^{+k_\ell}.
$$

(2) Let $h:\omega\to \omega$ be such that $(\forall m)(\exists^{\aleph_0}i)$
($h(i)=m$). Choose by induction on $i$, $\lambda'_{\ell,i}\in\{\lambda_{h(i),
m}:m<\omega\}$ such that (b) $+$ (c) of (2) hold.

For each $i$ we do this by downward induction on $\ell$.  Then apply the last
theorem.{\hfill$\square_{\rm 5.7}$}
\medskip

We may deal with all $n$'s at once, at some price. The simplest case is:
\medskip

\n {\bf 5.8 Claim}. Assume
\begin{enumerate}
\item[(a)] $\lan A_\ell :\ell <\omega\ran$ is a sequence of pairwise disjoint
sets.
\item[(b)] $\lambda=\mbox{tcf}\big(\prod_{n<\omega}\theta_n/J^{{\rm
bd}}_\omega\big)$
\item[(c)] $\theta_n=\mbox{tcf}\big(\prod_{\ell <\omega}\tau_{n,\ell}/J^{{\rm
bd}}_\omega\big)$, $\tau_{n,\ell}$ regular $>\aleph_0$.
\item[(d)] $h:\omega\to \omega$ is such that $|h^{-1}(\{ n\})|=\aleph_0$,
$J=\big\{ A\subseteq \omega\times \omega :(\forall^{J^{{\rm bd}}_\omega}n)$
$(\forall^{J^{\rm bd}_\omega}m)$ $(h(n)=\emptyset =A\cap\{ m\}\times
[h(n), 2h(n))\big\}$. 
\end{enumerate}
$\underline{\rm Then}$ there is a $(\lambda ,J)$-sequence for $\lan J^{{\rm
bd}}_{\tau_{n,\ell}}: (n,\ell )\in \omega\times \omega\ran$.
\medskip

\n {\bf Proof}. Straight.
\medskip

\n {\bf 5.9 Remark}. 1) ~We can replace $\lan \theta_n:n<\omega\ran$ by
$\lan\theta_i:i<\delta\ran$.
\smallskip

2) Another way to get an example for 5.4 is to have $\lan \mu_i:i<\kappa\ran$
increasing continuous, $\kappa={\rm cf}(\kappa )>\aleph_0$, $\kappa <\mu_0$,
$\mu=\mu_\kappa =\sum_{i<\kappa}\mu_i$, ${\rm cf}(\mu_i)\le |\delta|$, ${\rm
pp}_{|\delta|}(\mu_i)<\mu_{i+1}$, $\chi_i=|\mbox{Reg}\cap
[\mu_i,\mbox{pp}^+_{|\delta|}(\mu_i))|$, $S\subseteq \kappa$ stationary such
that for every $S'\subseteq S$ stationary we have $\prod_{i\in
S'}\chi_i>\chi_\kappa$.
\smallskip

3) In all the cases here we can get normality as in \S4.
\smallskip

4) See 1.16, 1.17.

\section{Products of Boolean Algebras}
Monk asks \cite[Problem 35, p.~15]{M2}.
\medskip

\n {\bf 6.1 Monk's Problem}. Does $\prod_{n<\omega} \mbox{FBA}(\beth_n)$ have
free caliber $\beth^+_\omega$? Here:
\medskip

\n {\bf 6.1A Notation}. $\mbox{FBA}(\beta )$ is the Boolean algebra freely
generated by $\lan x_\alpha :\alpha <\beta\ran$.
\medskip

\n {\bf 6.2 Definition}. {\it 1) We say that the cardinal $\lambda$ is a free
caliber of the Boolean algebra $\B$ $\underline{\rm if}$ for every $X\in
[\B]^\lambda$ there is $Y\in [\B ]^\lambda$ such that $Y$ is independent in
$\B$.

2) Free Cal$(\B )=\{\lambda\le |\B| :\lambda\;\mbox{ is a free caliber of
$\B$}\}$.}
\medskip

We show that e.g.\ if $\beth^+_\omega =2^{\beth_\omega}$ then the answer is
NO.
\medskip

\n {\bf 6.3 Claim}. Assume
\begin{enumerate}
\item[(a)] there is a normal\footnote{if $\bar I$ is normal,
i.e. $\kappa_{i+1}>\lambda_i$, the normality of $\bar\eta$ follows.}  super
$(\lambda,J)$-sequence $\bar\eta$ for $\bar I=\lan I_i:i<\delta\ran$,
\item[(b)] $I_i={\rm ERI}^2_{\lambda_i,\kappa_i}=:\big\{X\subseteq[
\lambda_i]^2:$ for some $h:\chi \to\kappa_i$, $|\mbox{Rang }h|<\kappa_i$, and
for no $u\in [\lambda_i]^{\aleph_0}$ do we have ($h\restriction [u]^2$
constant) \& $[u]^2\subseteq X\big\}$. 
\item[(c)] $\delta <\omega_1$.
\end{enumerate}
$\underline{\rm Then}$ $\lambda$ is not a free caliber of $\prod_{i<\delta}
\mbox{FBA}(\lambda_i)$.
\smallskip

\n {\bf 6.3A Remark}. By 3.1, if $\lambda =\mu^+=2^\mu$, $\mu$ strong limit
$>\aleph_0={\rm cf}(\mu)$, then we can find such $\kappa_i$, $\lambda_i<\mu$
and $\bar\eta$ for $\delta =\omega$.
\medskip

\n {\bf Proof}. By renaming without loss of generality
$$
\eta_\alpha
(i)\ge\sum_{j<i}\lambda_j.\leqno\mbox{($*$)$_1$}
$$
Let $f_\alpha (i)=\{ f^0_\alpha (i), f^1_\alpha (i)\}$, $f^0_\alpha(i)<
f^1_\alpha (i)$ ($<\lambda$). First we deal with the case $\delta =\omega$, as
its notation is simpler. Let $\B_n= {\rm FBA}(\lambda_n)$ be freely generated
by $\{ x^n_\alpha :\alpha <\lambda_n\}$.  We define $g^*_\alpha\in\prod_{n<
\omega}\B_n$ for $\alpha <\lambda$ by
$$
g^*_\alpha (\ell ) = \bigcap_{k<\ell}\Big(x^\ell_{f^0_\alpha (k)}-x^\ell_{
f^1_\alpha (k)}\Big).
$$
Note:
\begin{enumerate}
\item[$\otimes_1$] for $\alpha <\beta <\lambda$, we have $g^*_\alpha$,
$g^*_\beta$ are distinct elements of $\prod_{n<\omega}\B_n$.
\item[$\otimes_2$] if $f^0_n(\beta )=f^1_n(\alpha)$ and $m>n$ then $\B_m\vDash
g^*_\alpha (m)\cap g^\ast_\beta(m)=0.$
\end{enumerate}
[Why? as $x^m_{f^0_n(\alpha )}-x^m_{f^1_n(\alpha)}$ is disjoint to $x^m_{
f^0_n(\beta )}-x^m_{f^1_n(\beta )}$.]

\begin{enumerate}
\item[$\otimes_3$] if $n<\omega$ and for $i=1,2$ we have $\alpha_i,\beta_i<
\lambda$ and $f^0_n(\beta_i)=f^1_n(\alpha_i)$ and
$$
\bigwedge_{k<n} f^0_k(\alpha_1)=f^0_k(\alpha_2)\; {\rm and}\;
\bigwedge_{k<n}f^1_k(\beta_1)=f^1_k(\beta_2)
$$
then
$$
\prod_{n<\omega}\B_m\vDash g^*_{\alpha_1}\cap g^*_{\beta_1}=g^*_{\alpha_2}\cap
g^*_{\beta_2}.
$$
\end{enumerate}
[Why? Check each coordinate in the product, for $m>n$ use $\otimes_2$ to show
that both sides are zero, and if $m\le n$ use the last two assumptions.]

Now if $X\in[\lambda ]^\lambda$ then there are such $\alpha_1, \alpha_2,
\beta_1, \beta_2$ (using the choice of $\bar\eta$ and its
normality).{\hfill$\square_{\rm 6.3}$}
\medskip

\n {\bf 6.4 Claim}. Assume
\begin{enumerate}
\item[($*$)(a)] $\mu=\mu^\theta <\lambda=\mbox{ cf }(\lambda)\le 2^\mu$, and
$\lan \chi_i:i<\theta\ran$ a sequence of cardinals, or
\item[(b)] $2^\theta <\lambda =\mbox{ cf}(\lambda )$ and in the
($<\theta^+$)-base product topology on ${}^{\sup (\chi_i)}2$ the density is
$<\lambda$, or at least in the box product topology on $\prod_{i<\theta}(
{}^{\chi_i}{2})$ (where each ${}^{\chi_i}2$ has Tychonov topology) has density
$<\lambda$. 
\end{enumerate}
$\underline{\rm Then}$ $\prod_{i<\theta} \mbox{FBA}(\chi_i)$ has free caliber
$\lambda$.
\medskip

\n {\bf Proof}. As in \S2.{\hfill$\square_{\rm 6.3}$}
\medskip

Probably the choice of the product of $\lan{\rm FBA}(\beth_n):n<\omega\ran$ in
the original question was chosen just as the simplest case, as is often
done. But in this case the products of uncountably many free Boolean algebras
behave differently.
\medskip

\n {\bf 6.5 Claim}. Assume $\lambda=\mbox{ cf}(\lambda )>2^\theta$, ${\rm
cf}(\theta) >\aleph_0$ and $(\forall\alpha <\lambda )$ $(|\alpha|^{\aleph_0}<
\lambda)$. $\underline{\rm Then}$ $\prod_{i<\theta}\mbox{FBA}(\chi_i)$ has
free caliber $\lambda$. 
\medskip

\n {\bf Proof}. First assume a stronger assumption
\begin{enumerate}
\item[($*$)\phantom{$^-$}] $\lambda =\mu^+$, ${\rm cf}(\mu )=\theta >\aleph_0$
and $(\forall\alpha <\mu )$ $(|\alpha |^\theta <\mu )$,
\end{enumerate}
or alternatively
\begin{enumerate}
\item[($*$)$^-$] $\lambda ={\rm cf}(\lambda )$ and $\mu >2^\theta$ are as in
7.3 below and we assume $i<\theta\Rightarrow \chi_i\le\mu$.
\end{enumerate}
(This was our first proof. It possibly covers all cases under some reasonable
pcf hypothesis, and illuminates the method).

Let $g^*_\alpha\in\prod_{i<\theta}\mbox{FBA}(\chi_i)$ for $\alpha <\lambda$ be
pairwise distinct, and we should find $X\in [\lambda ]^\lambda$ such that
$\lan g^*_\alpha :\alpha\in X\ran$ is independent. Let
$$
g^*_\alpha (i)=\tau_{\alpha ,i}(x_{\beta_{\alpha ,i,0}}, x_{\beta_{\alpha ,
i,1}},\ldots , x_{\beta_{\alpha, i, m(\alpha ,i)-1}})
$$
where $\tau_{\alpha ,i}$ is a Boolean term. Without loss of generality no
$x_{\beta_{\alpha ,i,\ell}}$ is redundant, $\beta_{\alpha , i, m}$ increasing
with $m$. As $2^\theta <\lambda ={\rm cf}(\lambda )$ without loss of
generality $\tau_{\alpha ,i}=\tau_i$ and so $m(\alpha ,i)=m(i)$. Let
$f_\alpha$ be the function with domain $\theta$, $f_\alpha (i)=\lan
\beta_{\alpha ,i,\ell}:\ell <m(i)\ran$. Let $f_\alpha^{[\ell]}(i)=
\beta_{\alpha ,i,\ell}$, so ${\rm Dom}(f_\alpha^{[\ell]})=\{ i<\theta :\ell
<m(i)\}$.

If ($*$) holds then by 7.1 and 7.2 (see below) we have
\begin{enumerate}
\item[$\circledcirc$] there are $u^*$, $m^*$, $v$, $\bar\beta^*$, $X$ such
that
\end{enumerate}

\begin{enumerate}
\item[(a)] $u^*\in [\theta ]^\theta$ and $X\in [\lambda ]^\lambda$
\item[(b)] $i\in u^*\Rightarrow m(i)=m^*$
\item[(c)] $v\subseteq m^*$ but $v\not= m^*$
\item[(d)] $\bar\beta^*=\lan \beta^*_{\ell ,i}:\ell <m^*, i\in u^*\ran$
\item[(e)] $\ell\in v\Rightarrow \lan f_\alpha^{[\ell ]}\restriction
u^*:\alpha\in X\ran$ is $<_{J^{\rm bd}_{u^*}}$-increasing and cofinal in
$\prod_{i\in u^*}\beta^*_{\ell ,i}$.
\item[(f)] $\ell\in m^*\backslash v\Rightarrow f_\alpha^{[\ell ]}\restriction
u^*=\lan \beta^*_{\ell ,i}:i\in u^*\ran$
\item[(g)] for every $\bar\gamma\in\prod_{\ell\in v\atop i\in
u^*}\beta^*_{\ell ,i}$ for $\lambda$ ordinals $\alpha\in X$ we have,
$$
i\in u^*\;\mbox{ \& }\; \ell\in v\Rightarrow\gamma_{\ell ,i}<f_\alpha^{[\ell
]}(i)<\beta^*_{\ell ,i}.
$$
\item[(h)] if $\ell\in v$, $\alpha\in X$, $i\in u^*$ then $f_\alpha^{[\ell
]}(i)>\sup\left\{\beta^*_{\ell_1, i_1}:\beta^*_{\ell_1,i_1}<\beta^*_{\ell,
i}\right.$ where $\ell_1<m^*$ and $i_1<\theta\big\}$ and $\alpha <\beta\in X$
implies: for every $i\in u^*$ large enough we have $f_\beta^{[\ell ]}(i)>\max
\left\{ f^{[\ell_1]}_\alpha (i_1):\beta^*_{\ell_1,i_1}=\beta^*_{\ell
,i}\right.$ and $\ell_1<m^*$ and $\left. i_1<\theta\right\}$ (the interesting
case is $i_1=i$).
\end{enumerate}
Now for any $n<\omega$, and $\alpha_0<\cdots <\alpha_{n-1}$ from $X$, we have
\begin{enumerate}
\item[\phantom{$\otimes$}$\otimes$\phantom{$\otimes$}] for every $i\in u^*$
large enough $$\begin{array}{l} \lan f_{\alpha_0}(i), f_{\alpha_1}(i),\ldots,
f_{\alpha_{n-1}}(i)\ran =\\ \\ \qquad \big\langle\lan\beta_{\alpha_0,i,\ell}:
\ell <m^*\ran ,\lan\beta_{\alpha_1,i,\ell}:\ell<m^*\ran,\ldots ,\lan
\beta_{\alpha_{n-1},i,\ell}:\ell<m^*\ran\big\rangle\end{array}
$$
is as in a $\Delta$-system, in fact
$$
\beta_{\alpha_{k(1)}, i,\ell (1)}=\beta_{\alpha_{k(2)},i,\ell (2)}\Rightarrow 
(k(1), \ell (1))=(k(2), \ell (2))\vee (\ell (1)=\ell (2)\in v).
$$
As $v\not= \{0,1,\ldots ,m^\ast-1\}$ and in $\tau$ no variable is redundant
clearly

\item[$\otimes$] for every $i\in u^*$ large enough, $\lan\tau
(x_{\beta_{\alpha_0,i,0}},\ldots ), \tau (x_{\beta_{\alpha_1,i,0}},\ldots
),\ldots\ran$ is independent.
\end{enumerate}
This implies that $\lan g^*_{\alpha_\ell}:\ell <n\ran$ is independent (in
$\prod_{i<\theta}\mbox{FBA}(\chi_i)$) as required.

If we do not have ($*$) or ($*$)$^-$, by $(\forall\alpha <\lambda )$ $(|\alpha
|^{\aleph_0}<\lambda )$ and $2^\theta <\lambda ={\rm cf}(\lambda )$ without
loss of generality for some $\tau =\tau (x_{1},\ldots , x_{n-1})$ and infinite
$u\subseteq \theta$, and some $X\in [\lambda ]^\lambda$ we have: $\lan
f_\alpha\restriction u:\alpha\in X\ran$ is with no repetition, $\tau_{\alpha,
i}=\tau$ for $\alpha\in X$, $i\in u$. So without loss of generality
$u=\theta$. Then we can find an ultrafilter $\DD$ on $\theta$ as in 7.4 below
and then the proof above works.{\hfill$\square_{\rm 6.5}$}
\medskip

\n {\bf 6.6 Comment}. Before we use 7.4, we wonder if ``$\chi_i\le\mu$'' is
necessary in ($*$)$^-$ of 6.5. This is quite straight. We can omit it if
$$
\aaa\subseteq\mbox{ Reg }\cap\lambda\backslash\mu,\;\; |\aaa|\le\theta
\Rightarrow \max\mbox{ pcf}(\aaa )<\lambda.
$$
\smallskip

\n {\bf 6.7 Problem}. 1) Which of the following statements is consistent with
ZFC:
\begin{enumerate}
\item[a)] $\mu$ is strong limit, ${\rm cf}(\mu )=\aleph_0$, and for every
$\lambda\in\mbox{ Reg}\cap (\mu ,2^\mu ]$ and cardinals $\chi_n$ such that
$\mu=\sum\limits_{n<\omega}\chi_n$, $\lambda$ is a free caliber of
$\prod\limits_{n<\omega}\mbox{FBA}(\chi_n)$ (What about ``some such
$\lambda$''? See 6.11 below.) 
\item[b)] the same for all such $\mu$.
\end{enumerate}

2) Can you prove in ZFC that for some strong limit $\mu$, $\theta =\mbox{
cf}(\mu )<\mu$ and for some set $\lan\aaa_i :i<\sigma\ran$ where $\sigma
=\theta^+$ or $\sigma =(2^\theta )^+$, pairwise disjoint there is $\lambda\in
(\mu ,2^\mu ]\cap\bigcap_{i<\sigma}\mbox{ pcf}(\aaa_i)$.

Now we turn to another of Monk's problems.
\medskip

\n {\bf 6.8 Claim}. Assume
\begin{enumerate}
\item[($*$)] $\kappa >\aleph_0$ is weakly inaccessible and $\lan 2^\mu :\mu
<\kappa\ran$ is not eventually constant.
\end{enumerate}
$\underline{\rm Then}$
\begin{enumerate}
\item[(a)] there is a $\kappa -c.c.$ Boolean algebra of cardinality
$2^{<\kappa}$, with no independent subset of cardinality $\kappa^+$.
\end{enumerate}
\medskip

\n {\bf Proof}. There are sequences $\lan (\II_i, \JJ_i):i<\kappa\ran$, $\lan
(\kappa_i, \lambda_i):i<\kappa\ran$ such that $\JJ_i$ is a dense linear order
of cardinality $\lambda_i$ and $\II_i\subseteq \JJ_i$ a dense subset of
$\JJ_i$ of cardinality $\kappa_i$, $\lan \kappa_i:i<\kappa\ran$ increasing
with limit $\kappa$, and $\lambda_j>\sum_{i<j}2^{\kappa_i}(\ge\sum_{i<j}
\lambda_i$), by \cite[3.4]{Sh:430}.

Let $\B_i$ be $\mbox{Intalg}(\JJ_i)$, the Boolean algebra of closed-open
intervals of $\JJ_i$. Let $\B$ be the free product of $\{ \B_i:i<\kappa\}$, so
$\B$ extends each $\B_i$ and each element of $\B$ is a Boolean combination of
finitely many elements of $\bigcup_{i<\kappa}\B_i$. It is straight to check
$\B$ is as required:
\begin{enumerate}
\item[($*$)$_1$] $|\B|=\sum_{i<\kappa}|\B_i|+\aleph_0=\sum_{i<\kappa}\lambda_i
=\sum_{1<\kappa}2^{\kappa_i}=2^{<\kappa}$.
\item[($*$)$_2$] $\B$ satisfies the $\kappa-c.c.$
\end{enumerate}
[Why? Let $a_i\in\B\backslash\{ 0\}$ for $i<\kappa$, so let $a_i=\tau_i(b_{i,
0},\ldots ,b_{i,n_i-1})$ for $i<\kappa$. $b_{i,\ell}\in\B_{\alpha_{i,\ell}}$.
As we can replace $a_i$ by any $a'_i$, $0<a'_i\le a_i$ without loss of
generality $a_i=\bigcap_{\ell <n_i} b_{i,\ell}$, $b_{i,\ell}\in\B_{\alpha_{i,
\ell}}\backslash\{ 0\}$. So without loss of generality $\alpha_{i,0}<
\alpha_{i,1}<\cdots <\alpha_{i,n_i}$. As $\kappa >\aleph_0$ is regular and as
we can replace $\lan a_i:i<\kappa\ran$ by $\lan a_i:i\in X\ran$ whenever $X\in
[\kappa ]^\kappa$, without loss of generality for some $m$, $\bigwedge_{\ell
<m}\alpha_{i,\ell}=\alpha_\ell$ and $i<j$ \& $\{\ell,k\}\subseteq [m,n]
\Rightarrow\alpha_{i,\ell}<\alpha_{j,k}$.  Let $a'_i=\bigcap_{\ell <m}b_{i,
\ell}$, so clearly
$$
a'_i\cap a'_j\not= 0\Leftrightarrow a_i\cap a_j\neq 0\Leftrightarrow
\bigwedge_{\ell <m}b_{i,\ell}\cap b_{j,\ell}\not= 0.
$$
But $\B_i$ satisfies the $\kappa$-Knaster condition (as $\kappa ={\rm
cf}(\kappa )>$ density $(\JJ_i)$), so can we finish.]
\begin{enumerate}
\item[($*$)$_3$] $\B$ has no independent subset of cardinality $\kappa^+$.
\end{enumerate}
[Why? Let $a_i\in\B$ for $i<\kappa^+$, let $a_i=\tau_i(b_{i,0},\ldots
,b_{i,n_i-1})$ and let $b_{i,\ell}\in\B_{\alpha_{i,\ell}}\backslash \{
0,1\}$. We can replace $\lan a_i:i<\kappa^+\ran$ by $\lan a_i:i\in X\ran$ for
$X\in [\kappa^+]^{\kappa^+}$, so without loss of generality $\tau_i=\tau$,
$n_i=n$ and $\alpha_{i,\ell}=\alpha_\ell$.  Let $b_{i,\ell}=\bigcup_{k\in
u_{i,\ell}}[x_{i,\ell,k}, x_{i,\ell, k+1})$ where $\bar x^{i,\ell}=\lan
x_{i,\ell , k}:k\le k_{i,\ell}\ran$ is an increasing sequence of elements of
$\{-\infty\}\cup\JJ_i\cup\{\infty\}$, $x_{i,\ell,0 }=-\infty$, $x_{i,\ell
,k_{i,\ell}}=\infty$, $u_{i,\ell}\subseteq k_{i,\ell}$.  We can find
$y_{i,\ell ,k}\in\II_i$ such that $x_{i,\ell ,k}<y_{i,\ell ,k}<x_{i,\ell
,k+1}$. Without loss of generality $k_{i,\ell}=k_\ell$, $y_{i,\ell ,k}=y_{\ell
,k}$, $u_{i,\ell}=u_\ell$.

Without loss of generality $y_{i,\ell ,k}=y_{\ell ,k}$. For a finite
$A\subseteq\B$ let at$(A)={\rm at}(A,\B )$ be the number of atoms in the
Boolean subalgebra of $\B$ which $A$ generates (all this was mainly for
clarity). Now for any finite $u\subseteq\kappa^+$
$$\begin{array}{l}
\displaystyle{\mbox{at}(\{ a_i:i\in u\}, \B)\le \mbox{at}(\{b_{i,\ell}:i\in u,
\ell <n\}, \B )\le }\\ 
\qquad\displaystyle{\prod_{\ell <n} \mbox{at}(\{
b_{i,\ell}:i\in u\}, \B_{\alpha_{i,\ell}})\le }\\
\quad\displaystyle{\prod_{\ell <n}\mbox{at}(\{ x_{i,\ell
,k}:i\in u, k<k_\ell\}, \B_{\alpha_{i,\ell}}\})\le }\\
\quad\displaystyle{\prod_{\ell <n}\Big(\sum_{k<k_\ell}
(|u|+1)\Big)\le k^*\times |u|^n }\end{array}
$$
for $k^{*}=\max\{ k_\ell +1:\ell <n\}$. So if $u$ is large enough this is
$<2^{|u|}$, showing non independence.  {\hfill$\square_{\rm 6.8}$}
\medskip

\n {\bf 6.9 Claim}. Let $\B$ be the completion of $\mbox{FBA}(\chi )$
\begin{enumerate}
\item[1)] $\lambda$ is not a free caliber of $\B$ if
\begin{enumerate}
\item[($*$)\phantom{$^-$}] $\lambda =\mu^+=2^\mu$, $\mu\le \chi$, $\mu$ strong
limit singular of cofinality $\aleph_0$
\end{enumerate}
\item[2)] $\lambda$ is a free caliber of $\B$ if
\begin{enumerate}
\item[($*$)] $\mu =\mu^{\aleph_0}<\lambda =\mbox{ cf }(\lambda)\le 2^\mu$,
$\chi\ge \lambda$, or at least
\item[($*$)$'$] $\chi\ge \mu$, $\mu <\lambda =\mbox{ cf }(\lambda)\le
2^\lambda$, $\mu$ strong limit singular of cofinality $\aleph_0$ and the
($<\aleph_1$)-box product topology on ${}^\chi\omega$ has density $<\lambda$.
\end{enumerate}
\end{enumerate}
\smallskip

\n {\bf Proof}. 1) By 6.3, 6.3A's proofs.

2) If ($*$) use 6.5, if ($*$)$'$ the proof is similar.  {\hfill$\square_{\rm
6.9}$}
\medskip

\n {\bf 6.9A Remark}. We can deal with singular cardinals similarly as in the
earlier proofs.
\medskip

\n {\bf 6.10 Claim}. In the earlier claims if
\begin{enumerate}
\item[$(*)_1$] $\lambda =\mu^{++}$,

\n or at least if

\item[$(*)_2$] $\mu <\lambda$, and $[\alpha <\lambda\Rightarrow\mbox{
cf}([\alpha ]^\theta ,\subseteq )<\lambda ]$, $\chi=\sup_{i<\theta}\chi_i$
\end{enumerate}
$\underline{\rm then}$ ``in the $(\le\theta^+)$-box product topology,
${}^\chi\theta$ has density $<\lambda$'' can be replaced by ``in the
$(<\theta^+)$-box product topology, ${}^\mu\theta$ has density $<\lambda$''.
\medskip

\n {\bf 6.11 Conclusion}. 1) ~Let $\ell \in\{1,2\}$ for simplicity. The
following questions cannot be answered in ZFC (assuming the consistency of
large cardinals).

Assume $\beth^{+\ell}_\omega\le \beth_{\omega +1}$
\begin{enumerate}
\item[a)$_\ell$] Does $\prod_{n<\omega}\mbox{FBA}(\beth_m)$ have free caliber
$\beth^{+\ell}_\omega$?
\item[b)$_\ell$] Does the completion of $\mbox{FBA}(\beth_\omega)$ have free
caliber $\beth^{+\ell}_\omega$?
\item[c)$_\ell$] Does the completion of $\mbox{FBA}(\beth_\omega^{+\ell})$
have free caliber $\beth^{+\ell}_\omega$?
\end{enumerate}

2) Moreover we can add
$$
\mbox{for $x\in \{ a,b,c\}$ even $(*)_1+(*)_2$, and $\neg (*)_1+\neg (*)_2$.}
$$
\smallskip

\n {\bf Proof}. 1) Let $\ell =2$. By Gitik, Shelah \cite{GiSh:597} it is
consistent with $ZFC$ that with the $(<\aleph_1)$-box product topology,
${}^{(\beth_\omega)}\omega$ has density $\le \beth^+_\omega$, so we can use
6.3, 6.4(i) (using 6.11 of course). For the other direction by Gitik, Shelah
\cite{GiSh:597} the necessary assumptions for 6.2, 6.9(i) are consistent.

For $\ell =1$, if $\beth^+_\omega =2^{\beth_\omega}$ then the answer is NO by
6.2, 6.9.

To get consistency for $\lambda =\beth^+_\omega$ we need dual: in ${}^\mu
\omega$, for every $\mu^+$ open sets there is a point belonging to $\mu^+$ of
them (this is phrased in 6.12 below). This too is proved consistent in
\cite{GiSh:597}.

2) Similarly.{\hfill$\square_{\rm 6.11}$}
\medskip

\n {\bf 6.12 Definition}. {\it $\mbox{Pr}_{\theta, \sigma}(\lambda ,\mu )$
means:
\begin{quotation}
\noindent if $f_\alpha$ is a partial function from $\mu$ to $\theta$ such
that $|\mbox{Dom}(f_\alpha )|<\sigma$ for $\alpha <\lambda$,

\noindent $\underline{\rm then}$ some $f\in {}^\mu\theta$ extends $\lambda$ of
the functions $f_\alpha$.
\end{quotation}
If $\sigma =\theta$ we may omit it.}
\medskip

\n {\bf 6.13 Claim}. In Claim 6.11 the assumption on the density of box
products can be replaced by cases of Definition 6.12:
\begin{enumerate}
\item[(a)] [2.1] Assume $\B=\cBB (\chi )$ is a Maharam measure algebra of
dimension $\chi$, ${\rm cf}(\lambda )>2^{\aleph_0}$ and ${\rm cf}(\lambda
)>\beth_2\vee \lambda ={\rm cf}(\lambda)\vee\boxtimes^-_{\lambda
,\aleph_0}$. If ${\rm Pr}_{\aleph_0}(\lambda ,\chi)$ $\underline{\rm then}$
$\B$ has $\lambda$ as a free caliber.

\item[(b)] [6.4] Assume $2^\theta <\lambda ={\rm cf}(\lambda )$, $\chi
=\sup_{i<\theta}\chi_i$. If ${\rm Pr}_\theta (\lambda ,\chi )$ then
$\prod_{i<\theta}{\rm FBA}(\chi_i)$ has free caliber $\lambda$.
\end{enumerate}

\n {\bf Proof}. Straight. {\hfill{$\square_{6.13}$}}
\medskip

In fact cases of Pr are essentially necessary and sufficient conditions.
\medskip

\n {\bf 6.14 Claim}. 1) Assume $\lambda ={\rm cf}(\lambda )>2^{\aleph_0}$, and
$\chi_n$ are cardinal. The following conditions are equivalent
\begin{enumerate}
\item[(a)] $\displaystyle{\prod_{n<\omega}{\rm FBA}(\chi_n)}$ has free caliber
$\lambda$.
\item[(b)] If for $\alpha <\lambda$, $i<\omega$, $(u^\alpha_i, v^\alpha_i)$ is
a pair of disjoint finite subsets of $\chi_i$ $\underline{\rm then}$ for some
$X\in [\lambda ]^\lambda$ we have
$$
i<\omega\Rightarrow \bigcup_{\alpha\in X}u^\alpha_i\cap \bigcup_{\alpha\in X}
v^\alpha_i=\emptyset,
$$
(i.e.\ if $f^\alpha_i$ is a finite function from $\chi_i$ to $\{ 0,1\}$ for
$i<\omega$, $\alpha <\lambda$, then for some $\lan f_i:i<\omega\ran$
$$
(\exists^\lambda\alpha <\lambda )(\forall i<\omega)f^\alpha_i\subseteq f_i.
$$
\end{enumerate}

\n {\bf Proof}. Straight. {\hfill{$\square_{6.14}$}}
\medskip

\n {\bf 6.15 Discussion}. For measure, the parallel seems cumbersome. We still
may like to be more concrete on the dependencies appearing. Note
\begin{enumerate}
\item[$\otimes_1$] in 3.6, we can have $\bar x=\lan x_\alpha:\alpha
<\lambda\ran$ satisfies
 
($*$)$_{\B,\bar x}$ for every $X\in [\lambda ]^\lambda$, $m<\omega$, and
$\beta (\alpha ,k)<\lambda$ for $\alpha <\lambda$, $k<2m$ pairwise distinct,
for every $n$ large enough there are pairwise distinct $\alpha_0,\ldots,
\alpha_{2n-1}\in X$ such that
$$
0=\bigcap_{\ell <n}\Big(\bigcup_{k<m}\Big(x_{\beta (\alpha_{2\ell},k)}\Delta
x_{\beta (\alpha_{2\ell +1},k)}\Big)\Big).
$$
\item[$\otimes_2$] if ($*$)$_{\B ,\bar x}$ holds then the Boolean algebra
$\B'=\lan x_\alpha :\alpha <\lambda\ran_\B$ has no independent subset of
cardinality $\lambda$. Moreover, if $x'_\alpha\in\B'$ for $\alpha <\lambda$
are distinct, then ($*$)$_{\B', \lan x'_\alpha :\alpha <\lambda\ran}$.
\end{enumerate}  

\section{A nice subfamily of functions exists}

We expand and continue on \cite[6.6D]{Sh:430}, \cite[6.1]{Sh:513}. 

\n {\bf 7.0 Claim}. Assume
\begin{enumerate}
\item[(A)] $\lambda =\mbox{ cf}(\lambda )\ge\mu >2^\kappa$
\item[(B)] $\DD$ is a $\mu$-complete\footnote{in parts (0), (1),
$\mu=(2^\kappa)^+$ is O.K.} filter on $\lambda$ containing the
co-bounded subsets of $\lambda$
\item[(C)] $f_\alpha :\kappa\to\mbox{ Ord}$ for $\alpha <\lambda$
\item[(D)] $\alpha<\mu\Rightarrow |\alpha|^{\le\kappa}<\lambda$.
\end{enumerate}
$\underline{\rm Then}$

0) We can find $w\subseteq \kappa$ and $\bar\beta^*=\lan\beta^*_i:i<
\kappa\ran$ such that: $i\in \kappa\setminus w\Rightarrow {\rm cf}(\beta^*_i)
>2^\kappa$ and for every $\bar\beta\in\prod_{i\in \kappa\setminus w}\beta^*_i$
for $\lambda$ ordinals $\alpha<\lambda$ (even a set in $\DD^+$) we have
$\bar\beta<f_\alpha\restriction (\kappa\backslash w)<\bar\beta^*\restriction
(\kappa\backslash w)$, $f_\alpha \restriction w=\bar\beta^*\restriction w$,
and $\sup\{\beta^*_j:\beta^*_j<\beta^*_i\}<f_\alpha (i)<\beta^*_i$.

1) We can find a partition $\lan w^*_\ell :\ell <2\ran$ of $\kappa$, $X\in
\DD^+$ and $\lan A_i:i<\kappa\ran$, $\lan\bar\lambda_i:i<\kappa\ran$, $\lan
h_i:i<\kappa\ran$, $\lan n_i:i<\kappa\ran$ such that:
\begin{enumerate}
\item[(a)]   $A_i\subseteq\mbox{ Ord}$,
\item[(b)] $\bar\lambda_i=\lan\lambda_{i,\ell}:\ell <n_i\ran$ and $2^\kappa
<\lambda_{i,\ell}\le \lambda_{i,\ell +1}\le\lambda$
\item[(c)] $h_i$ is an order preserving function from $\prod_{\ell
<n_i}\lambda_{i,\ell}$ onto $A_i$ so $n_i=0\Leftrightarrow |A_i|=1$. (The
order on $\prod_{\ell <n_i}\lambda_{\ell ,i}$ being lexicographic $<_{\ell
x}$.) 
\item[(d)] $i<\kappa$ \& $\alpha\in X\Rightarrow f_\alpha (i)\in A_i$, and we
let $f^*_\alpha (i,\ell )=[h^{-1}_i(f_\alpha (i))](\ell )$, so
$f^*_\alpha\in\prod\limits_{\scriptstyle i<\kappa\atop\scriptstyle\ell <n_i}
\lambda_{i,\ell}$,   
\item[(e)] $i\in w^*_0\Rightarrow n_i=0$ (so $|A_i|=1$)
\item[(f)] if $i\in w^\ast_1$ then $|A_i|\le\lambda$, hence $|\bigcup_{i\in
w^\ast_1 }A_i|\le\lambda$
\item[(g)] if $\displaystyle{g\in\prod_{i<\kappa\atop\ell<n_i}\lambda_{i,
\ell}}$ then
$$
\{\alpha\in X:g<f^*_\alpha\}\in\DD^+
$$
\item[(h)] $\mu \leq {\rm max ~pcf}\{\lambda_{i,\ell}:i\in w_1^\ast\; {\rm and }\; \ell
<n_i\}\le\lambda$ when $w_1^*\neq\emptyset$ (so e.g.\ under GCH
$${\rm max ~pcf}\left\{ {\rm cf}(\lambda_{i,\ell}):i\in w^*_1\right.\mbox{
and }\left.\ell <n_i\right\}=\lambda).$$ 
\end{enumerate} 

2) In part (1) we can add $(*)_1$ to the conclusion if (E) below holds,
\begin{enumerate}
\item[($*$)$_1$] if $\lambda_{i,\ell}\in [\mu ,\lambda )$ then $\lambda_{i,n}$
is regular.
\end{enumerate}

\begin{enumerate}
\item[(E)] For any set $\aaa$ of $\le\kappa$ singular cardinals from the
interval $(\mu ,\lambda )$, we have $\max\mbox{ pcf} \{ {\rm cf}(\chi):\chi
\in\aaa\}<\lambda$
\end{enumerate}

3) Assume in part (1) that (F) below holds. Then we can demand ($*$)$_2$.

\begin{enumerate}
\item[($*$)$_2$] $\lambda^i_\ell \ge\mu$ for $i\in w_2$, $\ell< n_i$.
\end{enumerate}
\begin{enumerate}
\item[(F)] ${\rm cf}(\mu )>\kappa$ and $\alpha <\mu\Rightarrow\DD$ is $[|
\alpha |^{\le\kappa}]^+$-complete.
\end{enumerate}
4) If in part (1) in addition (G) below holds, $\underline{\rm then}$ we can
add
\begin{enumerate}
\item[$(*)_3$]\quad $\lambda\in{\rm pcf}_{\sigma{\rm -complete}}\{
\lambda^i_\ell:i\in w^*_1\; {\rm and }\; \ell <n_i\}$ if $w_1^*\neq\emptyset$, 
\end{enumerate}
moreover 
\begin{enumerate}
\item[$(*)_4$] if $\ell_i< n_i$ for $i\in w^*_1$ then $\lambda\in
{\rm pcf}_{\sigma{\rm -complete}}\{{\rm cf}(\lambda^i_{\ell_i}): i\in w^*_1\}$ 
\end{enumerate}

\begin{enumerate}
\item[(G)]
\begin{enumerate}
\item[(i)] $(\forall\alpha <\lambda )(|\alpha |^{<\sigma }<\lambda )$ and
$\sigma ={\rm cf}(\sigma )>\aleph_0$ 
\item[(ii)] $\DD$ is $\lambda$-complete.
\end{enumerate}
\end{enumerate}

5) If in part (1) in addition (H) below holds then we can add  
\begin{enumerate}
\item[$(*)_5$]\quad if $m<m^*$, $A\in J_m$ and $\ell_i< n_i$ for $i\in
\kappa\setminus A$ (so $w^*_0\subseteq A$) \underline{then}
$\lambda\in {\rm pcf}\{\lambda^i_{\ell_i}: i\in \kappa\setminus A\}$
\end{enumerate}
\begin{enumerate}
\item[(H)]
\begin{enumerate}
\item[(i)] $m^*<\omega$ and $J_m$ an $\aleph_1$-complete ideal on
$\kappa$ for $m<m^*$
\item[(ii)] ${\cal D}$ is $\lambda$-complete.
\end{enumerate}
\end{enumerate}
\smallskip

\n {\bf Remark}. 1) If $\lambda_{i,\ell}$ is singular we can replace it with a
sequence $\lan\gamma_{i,\ell_1}:\zeta<{\rm cf}(\lambda_{i,\ell})\ran$, and the
index set $\lan\lan\alpha\ran :\alpha <\lambda_{i,\ell}\ran$ by $\lan (\zeta
,\gamma ):\zeta <{\rm cf}(\lambda_{i,\ell})$ and $\gamma<\gamma_{i,\ell_2}
\ran$, and $\gamma_{i,\ell_1}$, are replaced by sequences of regular
cardinals. Not clear if all this helps.

2) The reader may concentrate on the case (F) + (G)(ii) holds. 
\medskip

\n {\bf Proof}. 0) By part (1).

1) Let $\chi$ be regular large enough. Choose $N$ such that \begin{enumerate}
\item[(i)] $N\prec({\cal H}(\chi),\in )$
\item[(ii)] $2^\kappa +1\subseteq N$ and $\| N\|=2^\kappa$
\item[(iii)] $\kappa , \mu ,\lambda ,\DD$ and $\lan f_\alpha :\alpha
<\lambda\ran$ belong to $N$ \item[(iv)] $N^\kappa\subseteq N$.
\end{enumerate}
Next choose $\delta(*)<\lambda$ which belongs to $B^*=\cap\{ B\in\DD :B\in
N\}$, which is the intersection of $\le 2^\kappa <\mu$ members of $\DD$.
Necessarily $B^\ast\in \DD$ so $\delta (*)$ exists. For each $i<\kappa$ let
$Y_i=:\{ A\in N: A$ a set of ordinals and $f_{\delta (*)}(i)\in A\}$,
clearly $Y_i\not=\emptyset$ as $\bigcup_{\gamma <\lambda}( f_\gamma (i)+1)\in N$,
hence there is a set $A_i\in Y_i$ of minimal order type. As $N^\kappa
\subseteq N$ clearly $\bar A=:\lan A_i:i\in\kappa\ran$ belongs to $N$.

Let us define:
$$\begin{array}{l}
w^*_0=:\{ i<\kappa :|A_i|=1\}\\
w^*_1=:\{ i<\kappa :2^\kappa <|A_i|\le\lambda\}.
\end{array}
$$
Now note
\begin{enumerate}
\item[($*$)$_1$] $A_i\not= \emptyset$. 
\item[{}] [Why? As $A_i\in Y_i$ hence $f_{\delta (\ast)}(i)\in A_i$.]
\item[($*$)$_2$] $|A_i|=1$ iff $A_i=\{ f_{\delta (*)}(i)\}$ iff $f_{\delta 
(\ast)}(i)\in N$ (iff $i\in w^*_0$).
\item[{}] [Why? Think.]
\item[($*$)$_3$] W.l.o.g. $A_i\subseteq\{ f_\alpha (i):\alpha <\lambda\}$
\item[{}] [Why? As $\{ f_\alpha (i):\alpha <\lambda\}\in Y_i$ and $A_i\cap\{ 
f_\alpha (i):\alpha <\lambda\}\in Y_i$.]

Hence
\item[($*$)$_4$] If $i\in \kappa\backslash w^\ast_0$ then $|A_i|\le\lambda$.
 \end{enumerate}
Let 
$$\begin{array}{rl}
K_i=\big\{ (\bar\lambda ,\bar \beta )\in N:&\mbox{for some $n$, $\bar\lambda 
=\lan\lambda_\ell :\ell <n\ran\in N$,}\\
&\mbox{and $\bar\beta =\lan\beta_\eta :\eta\in\prod_{\ell<n}\lambda_\ell\ran
\in N$ and $\beta_\eta\in\{ f_\alpha (i):\alpha<\lambda\}$}\\
&\mbox{and $f_{\delta (\ast)}(i)\in\{\beta_\eta :\eta\in\prod_{\ell 
<n}\lambda_\ell\}$ and:}\\
&\mbox{if $\eta <_{\ell x}\nu$ are from $\prod_{\ell <n}\lambda_\ell$ then 
$\beta_\eta\le \beta_\nu\big\}$.}
\end{array}
$$
We define a partial order $<^*$ on $K_i$
$$
(\bar\lambda^1, \bar\beta^1)<^* (\bar\lambda^2, \bar\beta^2)\;\;{\rm iff:}
$$
$$
\{\beta^1_\eta:\eta\in\prod_\ell\lambda^1_\ell\}\subseteq\{\beta^2_\eta:\eta
\in\prod_\ell\lambda^2_\ell\}
$$
and:
$$
{\rm otp} (\prod\limits_\ell \lambda^1_\ell, \leq_{\ell x}) < {\rm otp}
(\prod\limits_\ell \lambda^2_\ell, \leq_{\ell x}) \quad$$
 $$\underline
{\rm or} \quad
{\rm otp}(\prod\limits_\ell \lambda^1_\ell , \leq_{\ell x}) =
{\rm otp}(\prod\limits_\ell \lambda^2_\ell, \leq_{\ell x})\quad{\rm and}\quad 
{\rm lg}(\bar\lambda^1)<
{\rm lg}(\bar\lambda^2)$$
$$\underline {\rm or}\ {\rm otp}(\prod_\ell
\lambda^1_\ell, \leq_{\ell x}) = {\rm otp}(\prod_\ell
\lambda^2_\ell, \leq_{\ell x}), {\rm lg}(\bar\lambda^1)=
{\rm lg}(\bar\lambda^2) \ {\rm and}
$$
$$
\bigvee_{k<{\rm lg }(\bar\lambda^1)}\Big[ \lambda^1_{{\rm
lg}(\bar\lambda^1)-1-k}<\lambda^2_{{\rm lg}(\bar\lambda^2)-1-k}\; {\rm and}\;
\bigwedge_{\ell <k}\lambda^1_{{\rm lg}(\bar\lambda^1)-1-\ell}=\lambda^2_{{\rm
lg}(\bar\lambda^2)-1-\ell}\Big].
$$
\begin{enumerate}
\item[($*$)$_5$] $(K_i, \le^*)\subseteq  N$ is a partial order which is a well quasi
order (i.e.\ no strictly decreasing $\omega$-chains)

\item[{}] [Why? Reflect.]
\item[($*$)$_6$] there is $(\bar\lambda ,\bar\beta )\in K_i$ such that
$\bigwedge_{\ell <{\rm lg }(\bar\lambda)}\lambda_\ell\le |A_i|$,
\item[{}] [Why? By ($*$)$_7$ below.]
\item[($*$)$_7$] ${\rm otp}(A_{i})\le {|A_i|^n}$ for some $n<\omega$
\item[{}] [Why? By Dushnik-Milner \cite{DM}, we can find $A_{i,n}\subseteq
A_i$ for $n<\omega$ such that $A_i=\bigcup_{n<\omega}A_{i,n}$ and ${\rm
otp}(A_{i,n})\le |A_i|^n$. So as $A_i\in N$ there is such sequence $\lan A_{i,
n}:n<\omega\ran$ in $N$ so $A_{i,n}\in N$ hence for some $n$ we have $f_{
\delta (*)}(i)\in A_{i,n}\in N$, so by the choice of $A_i$ clearly ${\rm otp}(
A_i)\le |A_i|^n$.]
\end{enumerate}
So we can find a $<^*$-minimal $(\bar\lambda^i, \bar\beta^i)\in K_i$ and let
$n_i={\rm lg}(\bar\lambda^i)$. Note: 
\begin{enumerate}
\item[($*$)$_{\rm 8}$] we can above in the choice of $A_i$ demand
$A_i=\{\beta^i_\eta :\eta\in\prod\limits_{\ell <n_i}\lambda^i_\ell\}$ where
$(\bar\lambda^i, \bar\beta^i)$ is a $<^*$--minimal in $K_i$.
\item[($*$)$_9$] $\lambda^i_\ell\le \lambda^i_{\ell +1}\le\lambda$ for $\ell 
<n_i$.
\item[{}] [Why? The second inequality by ($*$)$_4$ and ($*$)$_6$, the first
inequality as otherwise by renaming we can omit $\lambda^i_{\ell+1}$ and contradict the $<^*$--minimality of $(\bar{\lambda}^i,\bar{\beta}^i)$.]
\end{enumerate}
Let $\lan\eta^*_i:i<\kappa\ran$ be such that $\beta^i_{\eta_i^*}=f_{\delta
(*)}(i)$ and $\eta^*_i\in\prod_{\ell <\eta_i}\lambda^i_\ell$.
\begin{enumerate}
\item[($*$)$_{\rm 10}$] $\lambda_{i,\ell}>2^\kappa$ [trivial]
\end{enumerate}
Let $Y=\big\{\alpha <\lambda :$ for every $i<\kappa$ we have $f_\alpha (i)\in
A_i\big\}$, as $\bar f\in N$ and $\lan A_i:i<\kappa\ran\in N$ necessarily
$Y\in N$. Also $Y\in\DD^+$ because $\delta (*)\in Y$ and the choice of $Y$.
So for $\alpha\in Y$ we can choose $\lan\eta^\alpha_i:i<\kappa\ran$ such that
$\eta^\alpha_i\in\prod_{\ell < n_i}\lambda^i_\ell$ and $f_\alpha
(i)=\beta^i_{\eta^\alpha_i}$. We now define $f^*_\alpha\in\prod_{i<\kappa\atop
\ell <n_i}\lambda^i_\ell$ for $\alpha <\lambda$: $f^*_\alpha (i,\ell
)=\eta^\alpha_i(\ell )$.

Note: 
\begin{enumerate}
\item[($*$)$_{\rm 11}$] $\lan\bar\lambda^i:i<\kappa\ran$,
$\lan\bar\beta^i:i<\kappa\ran$ and $\bar f$, hence $\bar f^*=\lan f^*_\alpha
:\alpha \in Y\ran$ belong to $N$.

\item[($*$)$_{\rm 12}$] $\eta^\alpha_i(\ell )=f^*_\alpha(i,\ell )\in [\sup (N\cap
\lambda^i_\ell ),\lambda^i_\ell]$ and $\alpha\in Y\Rightarrow f^*_\alpha
(i,\ell )<\lambda^i_\ell$.

\item[{}] [Why? Reflect.]

\item[($*$)$_{\rm 13}$] for every $g\in \prod_{i<\kappa\atop \ell
<n_i}\lambda^i_\ell$ and $X\in [Y]^\lambda\cap N$ such that $\delta (*)\in X$
$\underline{\rm there ~is}$ $\alpha \in X$ such that
$$
g<f^*_\alpha\;\;{\rm i.e.}\;\; i<\kappa\; {\rm and}\; \ell <n_i\Rightarrow
g(i,\ell )<f^*_\alpha (i,\ell ).
$$
[Why? If not, there is such $g$, so as $\lan (\bar\lambda^i,\bar\beta^i):i<
\kappa \ran$, $\bar f=\lan f_\alpha :\alpha <\lambda\ran$ and $X$, $Y$ belong to
$N$ also $\bar f^*=\lan f^*_\alpha :\alpha \in X\ran$ belongs to $N$, so all
the requirements on $g$ are first order with parameters from $N$, so without
loss of generality $g\in N$. Now $\delta (*)\in X$ cannot satisfy the
requirement hence there are $i<\kappa$, $\ell <n_i$ such that $g(i,\ell
)>f^*_{\delta (*)}(i,\ell )$ contradicting ($*$)$_{\rm 12}$.]
\end{enumerate}

Let $Z_i=\{\eta\in\prod_{i<n_i}\lambda^i_\ell:\;{\rm if}\;\nu\in\prod_{\ell<
n_i}\lambda^i_\ell\;{\rm and}\;\nu<_{\ell x}\eta\;{\rm then}\;\beta^i_\nu<
\beta^i_\eta\}$
$$
Z^+_i=\{\eta\restriction k:\eta\in Z_i\; {\rm and }\; k\le n_i\}.
$$
As $(\bar\lambda^i, \bar\beta^i)\in N$ clearly also $Z_i$, $Z^+_i\in N$.

\begin{enumerate}
\item[($*$)$_{\rm 14}$] if $i<\kappa$, $k<n_i$ then
$$
\lambda^i_k={\rm otp}\{\eta (k):(\eta^*_i\restriction k)\triangleleft \eta\in
Z_i\}.
$$
[Why? Let $Z'_i=\{\eta\in Z_i:\lambda^i_k>{\rm{otp}}\{\nu\in Z_i:\eta\restriction k\vartriangleleft\nu\in Z_i\}\}$. So $\eta^*_i\in Z'_i\in N$, by renaming 
$$\eta\in Z'_i\Longrightarrow\lambda^i_k>{\rm{sup}}\{\nu\in Z_i:\eta\restriction k\vartriangleleft\nu\in Z_i\},$$
and if $\lambda^i_k$ is regular we get a contradiction to $(*)_8$ as in the proof of $(*)_9$. If $\lambda^i_k$ is singular, we by renaming get the desired equality.]

Hence
\item[($*$)$_{\rm 15}$] without loss of generality $\lan\beta^i_\eta:\eta\in
\prod_{\ell<n_i}\lambda^i_\ell\ran$ is increasing (with $<_{\ell x}$)
\item[($*$)$_{\rm 16}$] $\mu\le\sup\{ \lambda^i_\ell :i\in w^*_1\; {\rm and}\; 
\ell <n_i\}$.

\item[{}] [Why? Otherwise let $\mu >\mu_0=\sup\{\lambda^i_\ell :i\in w^*_1\;
{\rm and}\; \ell <n_i\}$, and so $B^*\deq\{\beta^i_\eta :i<\kappa,\eta\in
\prod_{\ell <\eta_i} \lambda^i_\ell\}$ has cardinality $\mu_0$ so there
is $\PP\in N$, $|\PP |<\lambda$, $\PP\subseteq [\mu_0]^{\le\kappa}$ and $\PP$
is cofinal in $([\mu_0]^{\le\kappa}, \subseteq)$.

\item[{}] (Why? By assumption ($D$)). Note that if for some $X\in (\DD
+Y)^+$, $\bar f\restriction X$ is constant we are done. Otherwise
$a\in\PP\Rightarrow\{\alpha <\lambda: {\rm Rang}(f_\alpha )\subseteq a\}=
\emptyset \; {\rm mod }\;\DD$ but $\DD$ is $\mu$-complete hence $X^*=:\{\alpha
\in Y :(\exists a\in\PP )[{\rm Rang}(f_\alpha )\subseteq a]\}=\emptyset\; {\rm
mod }\DD$ and $X^*\in N$ and $\delta (*)\in X^*$, contradicting the choice of
$X^*$.] 

\item[($*$)$_{17}$] ${\rm max ~pcf}\big\{\lambda^i_\ell :i\in w^*_i$ and $\ell
<n_i\big\}\le\lambda$.
\item[{}] [Why? By ($*$)$_{13}$]
 
\item[$(*)_{18}$] $\lambda^i_\ell$ has cofinality $> 2^\kappa$.
\item[{}] [Why? Otherwise we can decrease it, getting a contradiction to the $<^*$--minimality of $(\bar{\lambda}^i,\bar{\beta}^i)$.]
\end{enumerate}
The conclusion can be checked easily

2) Let $\aaa =\{{\rm cf}(\lambda^i_\ell ):\lambda^i_\ell\;{\rm
is~singular~and}\; \mu\leq \lambda^i_\ell< \lambda\}$
and use (E). 

3) Easy.

4) Assume that the desired conclusion fails. For this we choose not just one
model $N$ but an $(\omega+1)$-tree of models. More precisely, we choose by
induction on $i\le\omega$ a sequence $\lan N_\eta:\,\eta\in T_i\ran$ such that
\begin{description}
\item{(a)} $T_i\subseteq {}^i\lambda$,
\item{(b)} $j<i\,\,\&\,\,\eta\in T_i\Rightarrow \eta\restriction j\in T_j$,
\item{(c)} $|T_i|<\lambda$,
\item{(d)} $N_\eta\prec (\HH(\chi),\in)$ satisfies (i)--(iv) from the
proof of part (1)
\item{(e)} For $\eta\in T_i$, $\eta\in N_\eta$ and $\lan N_\nu:\,\nu\in
\bigcup_{j<i}T_j\ran\in N_\eta$ and
$$
\nu\lhd\eta\Rightarrow N_\nu\prec N_\eta\,\,\&\,\,N_\nu\in N_\eta.
$$
\item{(f)} if $i=0$, then $T_i=\{\lan\ran\}$,
\item{(g)} If $i$ is $\omega$, then
$T_i=\{\eta\in {}^i\lambda:\,(\forall j<i)(\eta\restriction j\in T_i)$,
\item{(h)} If $i=j+1$, $\eta\in T_j$ and $\lan a_{\eta,\e}:\,
\e<\e_\eta<\lambda\ran$ list $[\sup(N_\eta\cap\lambda)]^{<\sigma}$,
\underline{then}
\[\{\nu\in T_i:\,\eta\lhd\nu\}=\{\eta\whp\lan\alpha\ran:\,\alpha<\e_\eta\},\]
and $a_{\eta,\e}\in N_{\eta\whp\lan\e\ran}$,
\item{(i)} $T=
\bigcup_{i\le\omega} T_i$.
\end{description}
There is no problem to carry the definition (note that
$\varepsilon_\eta< \lambda$ by assumption (G)(i) and
$|T_{m+1}|<\lambda$ as in addition $\lambda$ is regular, and
$|T_\omega|<\lambda$ by assumption (G)(i) as $\sigma > \aleph_0$). Now
\[B^\ast=\bigcap\{B\in \DD:\,\mbox{for some }\eta\in T\mbox{ we have }
B\in N_\eta\}\] 
being the intersection of $\leq |T|+2^\kappa<\lambda$ sets in $\DD$,
belongs to $\DD$ (using assumption (G)(ii)), so 
choose $\delta(\ast)\in B^\ast$. Now we choose by induction on $k<\omega$,
$\eta_k\in T_k$ and $w_0^k$, $w_1^k$, $\langle (\bar{\lambda}^{i,k},
\bar{\beta}^{i,k}):\,i<\kappa\ran\in N_{\eta_k}$ as in the proof of (1) for
$N_{\eta_k}$, such that $w^k_0\subseteq w^{k+1}_0$, $\eta_k\lhd\eta_{k+1}$ and
$(\forall i\in w_1^k)[(\bar{\lambda}^{i,k+1},\bar{\beta}^{i,k+1})<^*(\bar{
\lambda}^{i,k}, \bar{\beta}^{i,k})]$. The last assertion holds by the
assumption toward contradiction and basic pcf. 

If $\cup_{k<\omega} w_0^k=\kappa$, then $f_{\delta(\ast)}\in N_{\cup_k
\eta_k}$, hence $\delta(\ast)\in N_{\cup_k\eta_k}$, contradiction. If $i\in
\kappa\setminus \cup_{k<\omega} w^k_0$, then $\lan(\bar{\lambda}^{i,k},\bar{
\beta}^{i,k}):\,k<\omega\ran$ is strictly decreasing in $K_i$ (more exactly in
$\bigcup_{k<\omega}K_i[N_{\eta_k}]$), contradicting a parallel of
$(\ast)_{11}$. 

(5) We choose by induction on $t\in \omega$ the objects $N_t$, $\delta_t$, $\bar
    A^t = \langle A^t_i: i<\kappa\rangle$, $\langle (\bar \lambda^t_i,
    \bar \beta^t_i): i<\kappa\rangle$, $\langle h^t_i: i<\kappa\rangle$,
    $K^t_i$ such that 
\begin{enumerate}
\item[(a)] for each $t$, they are as required in the proof of part (1)
\item[(b)] $N_t\in N_{t+1}$, $K^t_i \subseteq K^{t+1}_i$ and
$(\bar \lambda^{t+1}_i, \bar \beta^{t+1}_i) \leq^* (\bar \lambda^t_i,
\bar \beta^t_i)$ in $K^{t+1}_i$
\item[(c)] for each $t$ for some $m_t< m^*$ we have $\{i<\kappa: (\bar
\lambda^{t+1}_i, \bar \beta^{t+1}_i) <^* (\bar \lambda^t_i, \bar
\beta^t_i)\}= \kappa\ {\rm mod}\ J_{m_t}$.
\end{enumerate} 
No problem to carry it by assumption toward contradiction. So for some
$m$, $\{t: m_t = m\}$ is infinite, contradicting ``$J_m$ is
$\aleph_1$-complete, and for each $i<\kappa$, $\bigcup_t K^t_i$ well ordered by $<^*$''.
{\hfill{$\square_{7.0}$}}
\medskip

We spell out a special case.
\medskip

\n {\bf 7.1 Fact}. Assume
$$
2^\theta <\mu,\;\mbox{ cf}(\mu )=\theta\quad {\rm and}\quad
(\forall\alpha <\mu )( |\alpha |^\theta <\mu ).\leqno(*)
$$
and $\lambda =\mu^+$.

$\underline{\rm Then}$:

1) For every sequence $\bar f=\lan f_\alpha :\alpha <\lambda\ran$ of functions
from $\theta$ to the ordinals, we can find $u^*\in [\theta ]^\theta$ and
$\bar{\beta}^*=\lan\beta^*_i:i\in u^*\ran$ such that one of the following
cases occurs:
\begin{enumerate}
\item[($*$)$_1$] for some $X\in [\lambda ]^\lambda$, $f_\alpha \restriction
u^*=\bar\beta^*$ for $\alpha\in X$,
\item[($*$)$_2$] $\beta^*_i$ is a limit ordinal (for every $i\in
u^*$), and $\lan{\rm cf}(\beta^*_i):i\in
u^*\ran$ is strictly increasing with limit $\mu$ and $\lambda=\mbox{tcf}\Big(
\prod_{i\in u^*}{\rm cf}(\beta^*_i)/J^{{\rm bd}}_{u^*}\Big)$ and for every
$\bar\gamma\in\prod_{i\in u^*}\beta^*_i$ for $\lambda$ ordinals $\alpha
<\lambda$ we have $(\forall i\in u^*)(\gamma_i<f_\alpha(i)<\beta^*_i)$,
\item[($*$)$_3$] $\beta^*_i$ is a limit ordinal of cofinality $\lambda$ for
$i\in u^*$ and for some $X\in [\lambda ]^\lambda$ we have: $i\in u^*
\Rightarrow\lan f_\alpha (i):\alpha\in X\ran$ is strictly increasing with
limit $\beta^*_i$ and for $\alpha\in X$, the interval $[f_\alpha (i),
\beta^*_i)$ is disjoint to 
$$
\{ f_\beta (j):\beta\in X\;{\rm and}\; j\in u^*\backslash\{ i\}\; \&
\; \beta_j \neq \beta_i \; 
{\rm or}\; \beta <\alpha \; {\rm and}\; j\in u^*\}.
$$
\end{enumerate}

2) Assume $\theta>\aleph_0$. For every sequence $\bar f=\lan f_\alpha :\alpha
<\lambda\ran$ of pairwise distinct functions from $\theta$ to ${}^{\omega
>}\mbox{Ord}$ such that $|\{f_\alpha (i):\alpha <\lambda\}|<\lambda$ for
$i<\theta$, we can find $u^*\in [\theta ]^\theta$ and $n(*)\in [1,\omega )$
and $v\subseteq n^*$ nonempty and $\bar\beta^*=\lan\beta^*_{\ell ,i}:\ell
<n^*, i\in u^*\ran$ such that for each $i$
\begin{enumerate}
\item[(a)] for $\ell\in v$ we have that $\beta^*_{\ell ,i}$ is a limit
ordinal, $\lan{\rm cf}(\beta^*_{\ell ,i}):i\in u^*\ran$ is strictly increasing
with limit $\mu$ and $\lambda =\mbox{tcf}(\Pi_{i\in u^*} \mbox{ cf}(\beta^*_{\ell
,i})/J^{{\rm bd}}_{u^*})$, and also for $i<j$ in $u^*$, and $\ell,
k\in v$ we have
${\rm cf}(\beta^*_{\ell ,i})<{\rm cf}(\beta^*_{k,j})$
\item[(b)] for every $\bar\gamma\in\prod_{\ell ,i}\beta^*_{\ell ,i}$ for
$\lambda$ ordinals $\alpha <\lambda$ we have
$$
\begin{array}{l}
\displaystyle{(\forall i<u^*)(\forall\ell \in v)[\gamma_{\ell,i}<(f_\alpha
(i))(\ell )<\beta^*_{\ell ,i}]\;\; {\rm and}}\\ 
\displaystyle{(\forall i\in u^*)(\forall\ell\in n^*\backslash v) [f_\alpha
(i))(\ell )=\beta^*_{\ell ,i}].} 
\end{array}
$$
\end{enumerate}

3) In part (2), we can replace $u^*\in [\theta ]^\theta$ by $u\in J^+$ for any normal ideal
$J$ on $\theta$. Moreover if $\{\delta <\theta :(\forall\alpha < {\rm cf }(\delta)
)(|\alpha |^{<\sigma}<{\rm cf}(|\delta|)\}$ is stationary then ${\rm Rang}(f_\alpha
)\subseteq{}^{\sigma >} {\rm Ord}$ is fine.  If we omit the assumption $|\{f_\alpha(i):\alpha<\lambda\}|<\lambda$, instead of $v$ we have a partition $(v_1, v_2, v_3)$ of $\{\ell:\ell<n^*\}$ such that clause (a) holds for $\ell\in 2$, clause (b) holds with $\ell\in v_2\cup v_3$, $\ell\in v_1$ instead of $\ell\in v$, $\ell\in n^*\setminus v$, and the parallel of $(*)_3$ holds for $\ell\in v_3$.
\medskip

\n {\bf Proof}. 1) By 7.0(0), we know that
\begin{enumerate}
\item[$\otimes$] there is $\lan \beta^*_i:i<\theta\ran$ and $w^*\subseteq\theta$ such
that letting $u^*=\theta\backslash w^*$ we have:
\begin{enumerate}
\item[a)] for every $\bar \gamma\in \prod_{i\in u^*}\beta^*_i$ for $\lambda$
ordinals $\alpha <\lambda$ we have
$$\begin{array}{l}
\displaystyle{i\in w^*\Rightarrow f_\alpha (i)=\beta^*_i}\\
\displaystyle{i\in u^*\Rightarrow\gamma_i<f_\alpha(i)<\beta^*_i}\end{array}
$$
\end{enumerate}
\end{enumerate}
\medskip

\n {\bf Case 1}. $|u^*|<\theta$.

So for some $X\in [\lambda]^\lambda$ we have $\lan f_\alpha\restriction
w^*:\alpha\in X\ran$ is constant. Easily ($*$)$_1$ holds.
\medskip

\n {\bf Case 2}. $u'=\{ i\in u^*:\mbox{cf}(\beta^*_i)=\lambda\}$ is unbounded
in $\theta$. 

Clearly ($*$)$_3$ holds and we get $X$ by ``thinning'': choose by induction on
$\gamma <\lambda$ the $\gamma$-th member $\alpha_\gamma <\lambda$ of $X$.
\medskip

\n {\bf Case 3}. Neither case 1 nor case 2.

Let $\mu =\sum_{i<\theta}\mu_i$, $\mu_i<\mu$ increasing with $i$. Choose
$j_i\in u^*$ such that
$$
j_i>\bigcup_{\zeta <i}j_\zeta\;\;\mbox{ \& }\;\;{\rm
cf}(\beta^*_{j_i})>\mu_i
+ \sum\limits_{\zeta< i} \mbox{ cf}(\beta^*_{j_\zeta}).
$$
This is possible as $\mu >\mbox{pcf}(\{{\rm cf}(\beta^*_\zeta):\zeta\in u^*,\;
\mbox{ cf}(\beta^*_i)<\mu_i\}$ because $(\mu_i)^\theta < \mu$ by
assumption.

Let us choose $u^*=\{j_i: i< \theta\}$, and $X=\{\alpha: i\in w^*
\Rightarrow f_\alpha (i) =\beta^*, \mbox{ and } i\in u^* \Rightarrow
f_\alpha(i) < \beta^*_i\}$. Now $\langle \mbox{ cf}(\beta^*_j): j\in
u^*\rangle$ is strictly increasing with limit $\mu$ hence $\Pi_{j\in
u^*}\beta^*_j / J^{\rm bd}_{u^*}$ is $\mu$-directed and even
$\mu^+$-directed. By clause (a), $\{f_\alpha \restriction u^*:
\alpha\in X\}$ is cofinal in $(\Pi_{j\in u^*} \beta^*_j,
<_{J^{\rm bd}_{u^*}})$ and $|X|=\lambda$. Together $(\Pi_{j\in u^*}
\beta^*_j, <_{J^{\rm bd}_{u^*}})$ has true cofinality $\lambda =
\mu^+$. Now all the demands in $(*)_2$ hold, so we are done.

2) First without loss of generality ${\rm lg}(f_\alpha (i))=n^*$, i.e.\ does not
depend on $\alpha$, secondly e.g.\ by successive applications of part (1).

3) Similar.{\hfill$\square_{\rm 7.1}$}
\medskip

\n {\bf 7.2 Conclusion}. For

1) In 7.1(1), ($*$)$_2$ and ($*$)$_3$ implies
\begin{enumerate}
\item[($*$)$'_2$] there are $u^*$, $\beta^*=\lan \beta^*_i:i\in u^*\ran$ and
$X$ such that
\begin{enumerate}
\item[(a)] $u^*\in [\theta ]^\theta$,
\item[(b)] $X\in [\lambda ]^\lambda$
\item[(c)] $\lan f_\alpha\restriction u^*:\alpha\in X\ran$ is $<_{J^{\rm
bd}_{u^*}}$-increasing
\item[(d)] for every $\bar\gamma\in\prod_{i\in u^*}\beta^*_i$ there are
$\lambda$ ordinals $\alpha\in X$ such that
$$
i\in u^*\Rightarrow \gamma_i<f_\alpha (i)<\beta^*_i.
$$
\item[(e)] if ($*$)$_3$ $\underline{\rm then}$:
\begin{enumerate}
\item[(i)] $\alpha <\beta$ from $X\Rightarrow f_\alpha\restriction
u^*<f_\beta\restriction u^*$ 
\item[(ii)] if $i\not= j$ are in $u^*$ and $\beta^*_i<\beta^*_j$ then $\alpha
\in X\Rightarrow f_\alpha (j)>\beta^*_i$
\item[(iii)] if $i,j\in u^*$, $\beta^*_i=\beta^*_j$ and $\alpha <\beta$ are
from $X$ then $f_\alpha (i)<f_\beta (j)$.
\end{enumerate}
\end{enumerate}
\end{enumerate}

2) Similarly for 7.1(2), getting $\circledcirc$ from the proof of 6.5.
\medskip

\n {\bf Proof}. Straight. Choose the $\gamma$-th member of $X$ for
$\gamma< \lambda$, by induction on $\gamma$.
{\hfill$\square_{\rm 7.2}$
\medskip

Similarly we can prove
\medskip

\n {\bf 7.3 Claim}. Assume
\begin{enumerate}
\item[(A)] $\lambda =\mbox{ cf}(\lambda )>2^\theta$,
\item[(B)] $\mu =\min\{\mu :\mu^\theta\ge\lambda\}$, ${\rm cf}(\mu )=\theta
>\aleph_0$. 
\item[(C)] if $\aaa\subseteq\mbox{ Reg }\cap\mu\backslash 2^\theta$,
$|\aaa|\le\theta$, $\lambda\in\mbox{ pcf}_{\theta{\rm -complete}}(\aaa )$,
then for some $\bb\subseteq\aaa$, $\lambda=\mbox{ tcf} (\prod\bb /[\bb
]^{<\theta})$. (Note: this holds if ${\frak d}\subseteq{\rm Reg}\setminus
2^\theta\ \&\ |{\frak d}|\leq\theta\ \ \Rightarrow\ \ |\mbox{pcf }({\frak
d})|\le\theta$. Why? Now $\langle{\frak b}_\theta[\aaa]:\theta\in{\rm
pcf}(\aaa)\rangle$ is well defined and $\lambda\in{\rm pcf}_{\theta{\rm -complete}}(
\aaa)$ so letting ${\rm pcf}(\aaa)\cap\lambda$ be $\langle\theta_\zeta:\zeta<
\theta\rangle$, choose $\mu_\zeta\in{\frak b}_\lambda[\aaa]\setminus
\bigcup\limits_{\xi<\zeta}{\frak b}_{\theta_\xi}[\aaa]$, and let ${\frak b}=
\{\mu_\zeta:\zeta<\theta\}$). 
\end{enumerate}
\smallskip

\n $\underline{\rm Then}$ the conclusions of 7.1, 7.2 hold.
\medskip

\n {\bf Proof}. Similar. {\hfill{$\square_{7.3}$}}
\medskip

\n {\bf 7.4 Fact}. 1) Assume
\begin{enumerate}
\item[(A)] $\lambda =\mbox{ cf}(\lambda )>2^\theta$ and $n<\omega$
\item[(B)] $f_\alpha^\ell\in {}^\theta\mbox{Ord}$ for $\ell <n$, $\alpha
<\lambda$.
\item[(C)] $\alpha\not=\beta\Rightarrow\lan f^\ell_\alpha :\ell <n\ran\not=
\lan f^\ell_\beta :\ell <n\ran$.
\item[(D)] $(\forall\alpha <\lambda )(|\alpha |^{\aleph_0}<\lambda )$.
\end{enumerate}
$\underline{\rm Then}$ we can find an ultrafilter $\DD$ on $\theta$ (possibly
a principal one) and $X\in [\lambda ]^\lambda$, $v\subseteq n$ and
$f_\ell\in{}^\theta\mbox{Ord}$ for $\ell <n$ such that
\begin{enumerate}
\item[(a)] for $\ell\in n\backslash v$ and $\alpha\in X$ we have
$f_\alpha^\ell/\DD= f_\ell /\DD$
\item[(b)] for $\alpha <\beta$ from $X$ and $\ell$, $m\in v$ such that $f_\ell /\DD
=f_m/\DD$ (e.g.\ $\ell =m$) we have $f^\ell_\alpha /\DD< f^m_\beta/\DD$.
\item[(c)] if $\ell$, $m<n$ and $f_\ell / \DD < f_m / \DD$ and
$\alpha$, $\beta$ are from $X$ then $f^\ell_\alpha / \DD <
f^\ell_\beta / \DD$.
\end{enumerate}

2) Assume 
\begin{enumerate}
\item[(a)] $\lambda = \mbox{ cf}(\lambda) > 2^\theta$ and $(\forall \alpha<\lambda)(|\alpha|^{<\sigma}<\lambda)$ and $\aleph_1 +
|\varepsilon(*)|^+ \leq \sigma = \mbox{ cf}(\sigma)$ and 
\item[(a)] $f^\varepsilon_\alpha \in {}^\theta \Ord$ for $\varepsilon
< \varepsilon(*)$ and $\alpha< \lambda$
\item[(c)] $I$ is a $\sigma$-complete filter on $\theta$ 
\item[(d)] $\DD$ is a $\lambda$-complete filter on $\lambda$ to which
all cobounded subsets of $\lambda$ belong.
\end{enumerate}
\underline{Then} we can find $X$, $v$, $f_\varepsilon$ (for
$\varepsilon < \varepsilon(*)$) and $\bar w$, $J$ such that 
\begin{enumerate}
\item[$(\alpha)$] $X\in [\lambda]^\lambda$
\item[$(\beta)$] $f_\varepsilon \in {}^\theta \Ord$ for $\varepsilon <
\varepsilon(*)$
\item[$(\gamma)$] $J$ is a $\sigma$-complete ideal on $\theta$
extending $I$ 
\item[$(\delta)$] $\bar w = \langle w_\varepsilon: \varepsilon <
\varepsilon(*)\rangle$, $w_\varepsilon \subseteq \theta$
\item[$(\varepsilon)$] if $\alpha \in X$ and $\varepsilon <
\varepsilon(*)$ then $f^\varepsilon_\alpha \restriction w_\varepsilon
= f_\varepsilon \restriction w_\varepsilon$
\item[$(\zeta)$] if $\alpha < \beta$ are from $X$ then $\varepsilon <
\varepsilon(*) \Rightarrow f^\varepsilon_\alpha < f^\varepsilon_\beta
\ {\rm mod}\ (J+ w_\varepsilon)$ moreover 
$$
\begin{array}{ll}
\{ i< \theta: & \mbox{ for some } \zeta, \xi< \varepsilon(*)\mbox{ we
have } i\notin w_\zeta, i\notin w_\xi \mbox{ and }\\
\ & f_\zeta(i) \leq f_\xi(i) \mbox{ but }f^\zeta_\alpha(i) \geq
f^\xi_\beta(i)\}\in J
\end{array}
$$
\item[$(\eta)$] if $\alpha\in X$ and $i< \theta$, $\zeta$, $\xi<
\varepsilon(*)$, $f_\zeta(i) < f_\xi(i)$ then $f_\zeta(i) <
f^\xi_\alpha(i)$
\item[$(\theta)$] if $2^{|\varepsilon(*)|} < \sigma$ then $\varepsilon
< \varepsilon(*) \Rightarrow w_\varepsilon \in J \ \vee\ \theta
\setminus w_\varepsilon \in J$.
\end{enumerate}

3) We can combine 7.0(1) with part (2) (having $\langle
   \lambda^\varepsilon_{i, \ell} : \ell< \ell^\varepsilon_i\rangle$).
\medskip

\n \underline{Remark.} can prove also the parallel of 7.0(5).
\medskip

\n {\bf Proof}. 1) Like the proof of 7.0(4) or by part (2) for $\sigma=\aleph_1$.

2) We repeat the proof of 7.0(4) except that $T\subseteq
   \bigcup\limits_{i<\sigma} {}^i\lambda$. After defining $B^*\in \DD$
   and choosing $\delta(*)$, for $\eta\in T$, $\varepsilon<
   \varepsilon(*)$ and $i< \theta$ we let $\beta_{\varepsilon, i,
   \eta} = \min[N_\eta \cap \Ord\setminus
   f^\varepsilon_{\delta(*)}(i)]$ and $w_{\varepsilon, \eta} =
   \{i<\theta: f^\varepsilon_{\delta(*)}(i)\in N_\eta\}$. 

So clearly
\begin{enumerate}
\item[$(*)_1$] $\eta\vartriangleleft \nu\in T \Rightarrow (\forall
   \varepsilon < \varepsilon(*)) (\forall i<\theta)
   (\beta_{\varepsilon, i, \eta} \geq \beta_{\varepsilon, i, \nu}\ \&\
   w_{\varepsilon, \eta} \subseteq w_{\varepsilon, \nu})$
\end{enumerate}
and
\begin{enumerate}
\item[$(*)_2$] $i\notin w_{\varepsilon, \eta} \Rightarrow \mbox{ cf}
(\beta_{\varepsilon,i,\eta}) > 2^\theta$.
\end{enumerate} 

Let $J_\eta$ is the $\sigma$-ideal on $\theta$ generated by 
$$
\begin{array}{ll}
I\cup
\{w\subseteq \theta: & \mbox{ for some }\varepsilon <
\varepsilon(*)\mbox{ we have } w\subseteq \theta\setminus
w_{\varepsilon, \eta}\mbox{ and }\\
\ & \lambda>\mbox{max pcf } \{\mbox{ cf}(\beta_{\varepsilon, i, \eta}): i<w\}\}
\end{array}
$$

If for some $\eta$, $\theta\notin J_\eta$ then we are done
(choosing the $\alpha$-th member of $X$  by induction on $\alpha$). So
assume that $\eta\in T \Rightarrow \theta \in J_\eta$. We now choose
by induction on $\zeta< \sigma$, a sequence $\eta_\zeta \in T_i$ such
that $\xi< \zeta \Rightarrow \eta_\xi = \eta_\xi \restriction j$ and
$$
\zeta = \xi+1 \Rightarrow \{i<\theta: (\exists \varepsilon <
\varepsilon(*))(\beta_{\varepsilon, i, \eta_\zeta} >
\beta_{\varepsilon, i, \eta_\xi})\}= \theta\; {\rm mod}\;  I.
$$
For some $\varepsilon<\varepsilon(*)$ and infinite $Y\subseteq\theta$ we have:
$$\xi\in Y\Rightarrow Z_\xi=\{i<\theta:\beta_{\varepsilon,i,\xi}>\beta_{\varepsilon,i,\xi+1}\}=\theta\;{\rm{mod}}\;I.$$
But for $\xi<\zeta$ we have $\beta_{\varepsilon,i,\xi+1}\geq\beta_{\epsilon,i,\zeta}$ by $(*)_1$.  Without loss of generality ${\rm{otp}}(Y)=\omega$. As $I$ is $\sigma$--complete and $\sigma\geq\aleph_0$, there is an $i\in\cap\{Z_\xi:\xi\in Y\}$, and $\langle\beta_{\varepsilon,i,\xi}:\xi\in Z\rangle$ is strictly decreasing, contradiction.

Now for $\zeta =0$, $\zeta$ limit there are no ``serious'' demands and
for $\zeta$ successor ordinal we use $\theta\in J_\eta$.

3) Left to the reader (and not used).
{\hfill$\square_{\rm 7.4}$}    
\medskip

\n {\bf 7.4A Fact} Assume
\begin{enumerate}
\item[(A)] $\lambda = \mu^+$, $\mu> 2^\theta$, $\theta = \mbox{ cf}(\mu)
>\aleph_0$
\item[(B)] $|\epsilon(*)|^+ +\aleph_0
<\theta$
\item[(C)] $f^\varepsilon_\alpha\in {}^\theta \Ord$ for $\varepsilon <
\varepsilon(*)$, $\alpha< \lambda$
\item[(D)] $(\forall \alpha< \mu)(|\alpha|^\theta< \mu)$ 
\end{enumerate}
\underline{Then} we can find a stationary $S\subseteq \{\delta<
\theta: \mbox{ cf}(\delta)\geq |\varepsilon(*)|^+ +\aleph_0\}$ and unbounded
subset $X^\prime$ of $\lambda$ and $S_\varepsilon \subseteq S$ and
$f_\varepsilon \in {}^S\Ord$ for $\varepsilon < \varepsilon(*)$
\begin{enumerate}
\item[(a)] for $\varepsilon< \varepsilon(*)$ we have $\alpha \in X
\Rightarrow f^\varepsilon_\alpha \restriction S_\varepsilon =
f_\varepsilon \restriction S_\varepsilon$
\item[(b)] for $\varepsilon_1 < \varepsilon(*)$ and $\alpha< \beta$
from $X$ if $S_{\varepsilon, \zeta} = \{i\in S: f_\varepsilon(i) \leq
f_\zeta(i)\} \setminus S_\varepsilon \setminus S_\zeta$ is unbounded
in $\theta$ then $f^\varepsilon_\alpha \restriction S_{\varepsilon,
\zeta} < f^\zeta_\beta \ {\rm mod}\ J^{\rm bd}_{S_{\varepsilon, \zeta}}$
\item[(c)] if $\zeta$, $\varepsilon < \varepsilon(*)$, $f_\zeta (i) <
f_\varepsilon (i)$, and $\alpha\in X$ then $f_\zeta(i)<
f^\varepsilon_\alpha(i)$
\item[(d)] if $2^{|\varepsilon(*)|} < \theta$ then $\varepsilon <
\varepsilon(*) \Rightarrow S_\varepsilon \in \{\emptyset, S\}$.
\end{enumerate}
\medskip

\n {\bf Proof}. 
Let $\bar f^\varepsilon =\langle f^\varepsilon_\alpha: \alpha<
\lambda\rangle$, let $\chi$ be large enough and $\langle
\lambda_\varepsilon: \varepsilon < \theta\rangle$ be increasing
continuous with limit $\mu$, and choose by induction on $\zeta<
\theta$, an elementary submodel $N_\zeta$ of $({\cal H}(\chi), \in,
<^*_\zeta)$ of cardinality $(\lambda_\zeta)^\theta$ such that
$(\lambda_\zeta)^\theta \subseteq N_\zeta$, ${}^\theta (N_\zeta)
\subseteq N_\zeta$, $\{\bar f^\varepsilon: \varepsilon<
\varepsilon(*)\}\in N_\zeta$, and $\langle N_\xi: \xi< \zeta\rangle\in
N\zeta$.

Choose $\delta(*)\in \lambda\setminus \bigcup\limits_{\zeta< \theta}
N_\zeta$, possible as $|\bigcup\limits_{\zeta< \theta} N_\zeta| =
|\sum_{\varepsilon<\theta}(\lambda_\varepsilon)^\theta| = \mu < \lambda$. For each
$\zeta< \theta$, $\varepsilon < \varepsilon(*)$ and $i< \theta$ let
$\beta^*_{\varepsilon, i, \zeta} = \min(N_\zeta \cap \Ord \setminus
f^\varepsilon_{\delta(*)}(i))$.

For each limit $i< \theta$ of cofinality $> |\varepsilon(*)|$ look at
$\langle \beta^*_{\varepsilon, i, \zeta}: \zeta< i\rangle$, it is a
non increasing sequence of ordinals hence is constant on some end
segment, i.e. for some $j_{\varepsilon, i} < i$ we have 
$$
j_{\varepsilon, i} \leq \zeta < i \Rightarrow \beta^*_{\varepsilon, i,
\zeta} = \beta^*_{\varepsilon, i, j_{\varepsilon, i}}.
$$
As $\mbox{ cf}(i) > |\varepsilon(*)|$, necessarily $j_i =
\sup\{j_{\varepsilon, i}: \varepsilon < \varepsilon(*)\}$ is $< i$,
hence for some $j(*) < \theta$ we have 
$$
S=\{i<\theta: \mbox{ cf}(i) > |\varepsilon(*)|, i\mbox{ a limit
ordinal}\}\mbox{ is stationary.}
$$
The rest should be clear.
{\hfill$\square_{\rm 7.4A}$}    

\medskip

\n {\underline{Remark}} We can demand $S\subseteq S^*$ in 7.4 if $S^*\subseteq\{\delta<\theta:{\rm{cf}}(\delta)\geq |\varepsilon(*)|^++\aleph_0\}$ is stationary.

\centerline{* * * * * * *}
\medskip

\n {\bf 7.5 Discussion} We may wonder what occurs for ultraproducts of free
Boolean algebras $\prod_{i<\theta}\mbox{FBA}(\chi_i)/\DD$ (or even reduced
products, recall $\mbox{FBA}(\chi_i)$ is the free Boolean algebra
generated, say, by $\{x_\alpha: \alpha< \chi_i\}$ freely). Now
\begin{enumerate}
\item[($*$)$_1$] if $\DD$ is $\aleph_1$-complete, the situation is as in the
$\theta >\aleph_0$ case for products.
\item[($*$)$_2$] if
$$
(\exists A_0 ,A_1,\ldots )\Big(\bigcap_{n<\omega} A_n\in \DD \;\mbox{ \&
}\;\bigcap_{n<\omega} A_n=\emptyset\Big)
$$
the situation is as in the $\theta =\aleph_0$ case.
\end{enumerate}
\medskip

\n {\bf 7.6 Claim}. Assume
\begin{enumerate}
\item[(a)] $\lambda =\mu^{++}$, $\mu >2^\theta$
\item[(b)] $f_\alpha :\theta\to\mbox{ Ord}$ for $\alpha <\lambda$.
\end{enumerate}
$\underline{\rm Then}$ we can find $\bar u^*=\lan u^*_0, u^*_1, u^*_2\ran$,
$\bar\beta^*$, $X$ such that
\begin{enumerate}
\item[(a)] $\lan u^*_0, u^*_1, u^*_2\ran$ is a partition of $\theta$
\item[(b)] $\bar\beta =\lan\beta^*_i:i<\theta\ran$
\item[(c)] $X\in [\lambda ]^\lambda$ (we can use an appropriate ideal $J$ on
$\lambda$ and demand $X\in J^+$)
\item[(d)] $\alpha\in X\Rightarrow f_\alpha\restriction u^*_0=\lan\beta^*_i:
i\in u^*_0\ran$
\item[(e)] if $i\in u^*_1$ $\underline{\rm then}$ $\lan f_\alpha (i):\alpha\in
X\ran$ is strictly increasing with limit $\beta^*_i$ (so ${\rm cf}(\beta^*_i)
=\lambda$)
\item[(f)] $i\in u^*_2\Rightarrow 2^\theta <\mbox{ cf}(\beta^*_i)<\mu$
\item[(g)] for every $\bar\gamma\in\prod_{i\in u^*_2}\beta^*_i$ for $\lambda$
ordinals $\alpha\in X$ we have $i\in u^*_2\Rightarrow\gamma_i<f_\alpha(i)<
\beta^*_i$
\item[(h)] $\lambda=\max\mbox{ pcf}\{{\rm cf}(\beta^*_i):i\in u^*_2\}$ if
$u^*_2\not=\emptyset$.
\end{enumerate}
\smallskip

\n {\bf Proof}. Let $\bar C=\lan C_\alpha :\alpha <\mu^+\ran$ be such that
$\mbox{otp}(C_\alpha )\le\theta^+$, $[\beta\in C_\alpha\Rightarrow C_\beta =
C_\alpha\cap\beta ]$, $C_\alpha$ a set of successor ordinals and $S^*=\{
\delta <\lambda^+:\mbox{ cf} (\delta )=\theta^+$ and $\alpha =\sup (C_\alpha
)\}$ is stationary (exists by \cite[\S1]{Sh:420}.)

Let $\bar f=\lan f_\alpha :\alpha <\lambda\ran$ be given.  Let $\chi$ be
strong limit such that $\bar f\in \HH(\chi )$.  We choose $M_\alpha$ by
induction on $\alpha <\mu^+$ such that

\begin{enumerate}
\item[($\alpha$)] $M_\alpha \prec (\HH (\chi ),\in ,<^*_\chi)$
\item[($\beta$)] $\| M_\alpha\|=2^\theta$ and $2^\theta +1\subseteq M_\alpha$
and ${}^{\theta\ge}(M_\alpha )\subseteq M_\alpha$
\item[($\gamma$)] $\lambda$, $\bar f$, $\bar C$ and $\alpha$ belong to
$M_\alpha$
\item[($\delta$)] $\lan M_\beta :\beta <\alpha\ran$ belongs to $M_{\alpha}$
and $\beta\in C_\alpha\Rightarrow M_\beta \prec M_\alpha$
\end{enumerate}
Now for every $\beta\in\lambda\backslash\bigcup_{\alpha <\mu^+}M_\alpha$ we
define a $g^\beta_\alpha\in {}^\theta(M_\alpha\cap\mbox{ Ord})$ and a function
$F_\beta$ from $\mu^+$ to $\mu^+$, as follows 
$$
g^\beta_\alpha (i)=\min (M_\alpha\cap\chi\backslash f_\beta
(i))\leqno\mbox{($*$)$_2$}
$$
[Why is it well defined? As $\bar f\in M_\alpha$ also $\cup\{ f_\gamma
(i)+1:\gamma <\lambda\}\in M_\alpha$ and $f_\beta (i)$ is smaller than that
ordinal.]
\smallskip

We let
$$
\begin{array}{l}
u^\beta_{\alpha ,0}=\{ i<\theta :f_\beta (i)\in M_\alpha\}\\
u^\beta_{\alpha ,1}=\{ i<\theta :f_\beta (i)\not\in M_\alpha\; {\rm and}\; 
{\rm cf}(g^\beta_\alpha (i))=\lambda\}\\
u^\beta_{\alpha,2}=\{ i<\theta: \mbox{ cf}(g^\beta_\alpha (i))\le\mu^+\;\;
\mbox{and}\;\;f_\beta (i)\notin M_\alpha\}.
\end{array}
$$

\n {\bf Note}: $f_\beta(i)\notin M_\alpha\Rightarrow \lambda\ge {\rm cf }
(g_\beta(i))>2^\theta$.

[Why? If $i\in\theta\backslash u^\beta_{\alpha ,0}$ and $\lambda <{\rm
cf}(g^\beta_\alpha (i))$, then
$$
\cup\{ f_\gamma (i):\gamma <\lambda\; {\rm and}\; f_\gamma (i)<g^\beta_\alpha
(i)\}
$$
belongs to $M_\alpha$ and contradicts the choice of $g^\beta_\alpha (i)$.  If
$i\in\theta\backslash u^\beta_{\alpha ,0}$ and ${\rm cf}(g^\beta_\alpha(i))
\le 2^\theta$ then $g^\beta_\alpha(i)=\sup(M_\alpha\cap g^\beta_\alpha(i))$.] 

Let
$$
J= J_{\beta ,\alpha}=\{ u\subseteq u^\beta_{\alpha ,2}:\mu^+>\max\mbox{
pcf}\{{\rm cf}(g^\beta_\alpha (i)):i\in u\}\}.
$$
By the pcf theorem (\cite[VIII 2.6]{Sh:g}) there is $u^\beta_\alpha\subseteq
u^\beta_{\alpha,2}$ such that: 
$$
\begin{array}{l}
\mu^+\notin\mbox{ pcf}\{{\rm cf}(g^\beta_\alpha (i)):i\in u^\beta_{\alpha,2}
\backslash u^\beta_\alpha\}\\ \\ \mu^+\ge {\rm max~pcf}\{{\rm cf}(
g^\beta_\alpha (i)):i\in u^\beta_\alpha\}.
\end{array}
$$
If $u^\beta_\alpha\notin J$ let $\bar h_{\beta ,\alpha}=\lan h_{\alpha,\beta,
\gamma}:\gamma <\mu^+\ran \in M_\alpha$ be $<_{J\restriction u^\beta_\alpha
}$-increasing and cofinal in $\prod_{i\in u^\beta_\alpha}g^\beta_\alpha (i)$. 
[Why it exists? as ${\rm Dom}(g^\beta_\alpha )\subseteq M_\alpha$ and ${\rm
Rang} (g^\beta_\alpha )\subseteq M_\alpha$]. Hence for some $\gamma =\gamma
(\alpha ,\beta )<\mu^+$, 
$$ 
f_\beta\restriction u^\beta_\alpha<h_{\alpha,\beta,\gamma(\alpha,\beta)}\ 
\mbox{ mod}(J\restriction u^\beta_\alpha).
$$
In fact any $\gamma'\in [\gamma (\alpha ,\beta ),\mu^+)$ will do, and now we
let $F_\beta (\alpha )=\gamma (\alpha ,\beta )$. If $u^\beta_\alpha\in
J_{\beta ,\alpha}$ we let $F_\beta (\alpha )=\alpha +1$.

So the set $E_\beta =\{\delta <\mu^+:$ $\delta$ a limit ordinal such that
$(\forall\alpha <\delta )F_\beta (\alpha )<\delta\}$ is a club of
$\mu^+$. Hence there is $\delta =\delta_\beta\in S^*\cap {\rm acc}(E_\beta )$
(i.e.\ $\delta =\sup (E_\beta\cap\delta )$ and $\delta\in S^*$). Now for each
$i<\theta$ the sequence $\lan g^\beta_\alpha (i):\alpha\in
C_{\delta_\beta}\ran$ is nonincreasing as $\lan M_\alpha :\alpha \in
C_{\delta_\beta}\ran$ is increasing, hence is eventually constant, say for
$\alpha\in C_{\delta_\beta}\backslash\alpha^*(\beta ,i)$. But ${\rm
otp}(C_{\delta_\beta})=\theta^+$ so $\alpha^*(\beta )=\sup\{\alpha^*(\beta
,i):i<\theta\}$ is $<\delta_\beta$, and reflection shows that
$$
\alpha\in C_{\delta_\beta}\backslash\alpha^*(\beta)\Rightarrow u^\beta_\alpha
\in J_{\beta ,\alpha}.
$$
Choose such $\alpha^\otimes_\beta$. So for some $\alpha^\otimes$,
$\delta^\otimes$ we have
$$
X=\Big\{\beta <\lambda :\beta\notin\bigcup_{\alpha<\mu^+}M_\alpha\;\;
\mbox{and}\;\; \alpha_\beta^\otimes =\alpha^\otimes,\delta_\beta =
\delta^\otimes\Big\}
$$ 
belongs to $[\lambda ]^\lambda$. Now we continue as in 7.0.{\hfill$\square_{\rm 7.4}$}
\medskip

\n {\bf 7.7 Claim}. 1) In 7.6 we can replace $\lambda =\mu^{++}$, by $\lambda
=\tau^+$, $\tau =\mbox{ cf}([\tau ]^{\le\mu},\subseteq )$ using
\cite[\S2]{Sh:420}. 

2) Also if $\lambda$ is weakly inaccessible $>\beth_\omega$, $(\forall\alpha
<\lambda )$ $[\lambda >\mbox{ cf}([\alpha ]^{\le\mu}, \subseteq )]$ we can get
7.6.

\section{Consistency of ``$\PP (\omega_1)$ has a free caliber" and discussion 
of pcf}
\smallskip

This solves [M2, Problem 37].

\n {\bf 8.1 Claim}: Assume for simplicity GCH and $P$ is adding
$\aleph_{\omega_1}$ Cohen reals. In $V^P$ we have $2^{\aleph_0}=\aleph_{
\omega_1}$, $2^{\aleph_1}=\aleph_{\omega_1+1}$ and
\begin{enumerate}
\item[($*$)] there is no complete Boolean algebra $\B$ of cardinality
$2^{\aleph_1}$ with FreeCal$(\B )= \emptyset$.  In fact for any complete
Boolean algebra $\B$ of cardinality $2^{\aleph_1}$ we have
$\aleph_{\omega_1+1}\in\mbox{FreeCal}(\B )$.
\end{enumerate}
\smallskip

\n {\bf Proof}. Clearly (as if the Boolean algebra $\bf{B}$ has cardinality $2^{\aleph_1}=\aleph_{\omega_1+1}$ and satisfies the ccc then $(*)$ holds, i.e., $\aleph_{\omega_1+1}\in{\rm{FreeCal}}({\bf{B}})$, because
$V^P\vDash (\aleph_{\omega_1})^{\aleph_0}=\aleph_{\omega_1}$, see and otherwise we can reduce to the case ${\bf{B}}={\mathcal{P}}({\omega_1})$) it is enough to
show
\begin{enumerate}
\item[($*$)$_1$] $V^P\vDash\aleph_{\omega_1+1}\in\mbox{FreeCal} (\PP
(\omega_1))$.
\end{enumerate}

\n So let $p^*\in P$
$$
p^*\Vdash_P\mbox{`` }\langle\name{a}_\alpha:\alpha<\aleph_{\omega_1+1}\ran
\mbox{ is a sequence of distinct elements of $\PP(\omega_1)$ ''}.
$$

{\bf Note}: $P=\{ f:f\mbox{ is finite function from $\aleph_{\omega_1}$ to
$\{0,1\}$\}}$ so without loss of generality $p^*=\emptyset$.

For each $\alpha <\aleph_{\omega_1+1}$ and $i<\omega_1$ there is a maximal
antichain $\lan f_{\alpha ,i,n}:n<\omega\ran$ of $P$ and sequence of truth
values $\lan \bt_{\alpha ,i,n}:n<\omega\ran$ such that
$$
f_{\alpha ,i,n}\Vdash_P\mbox{ ``}i\in\name{a}_\alpha\; \mbox{ $\underline{\rm
iff}$ }\; \bt_{\alpha ,i,n}\mbox{''}.
$$
Let $A_\alpha=\bigcup_{i<\omega_1,n<\omega}\mbox{Dom}(f_{\alpha ,i,n})\cup
{\rm Dom}(p^*)$, so $A_\alpha\in [\aleph_{\omega_1}]^{\le\aleph_1}$. Let
$A_\alpha=\{\gamma_{\alpha ,j}:j<j_\alpha\}$, $\gamma_{\alpha ,j}$ strictly
increasing with $j$.

As $V\vDash 2^{\aleph_1}<\aleph_{\omega_1+1}$, without loss of generality
\begin{enumerate}
\item[($*$)$_2$] (a) ~~$j_\alpha =j^*$
\item[{}] (b) ~~the truth value of ``$\gamma_{\alpha,j}\in\mbox{Dom}(f_{
\alpha ,i,n})$" does not depend on $\alpha$.
\end{enumerate}

Let $\name{a}$ be the Mostowski collapse of the name, i.e.\ $\name{a}
=OP_{j^*,A_\alpha}(\name{a}_\alpha )$ for each $\alpha$ (without loss of
generality it does not depend on $\alpha$). [Remember $OP_{A,B}(\beta
)=\alpha$ iff $\alpha\in A$, $\beta\in B$, ${\rm opt}(\beta\cap B)={\rm otp}
(\alpha\cap A)$.] 

We call $\lan (g^0_i, g^1_i,\xi_i) :i<\omega_1\ran$ a witness above $f^*$ if:
\begin{enumerate}
\item[(i)] $f^*, g^0_i, g^1_i\in P$ and $p^*\leq f^*$
\item[(ii)] $f^*\le g^0_i$, 
\item[(iii)] $f^*\le g^1_i$
\item[(iv)] $\mbox{Dom}(g^\ell_i)\subseteq j^*$
\item[(v)] $\lan\mbox{Dom}(g^0_i)\cup\mbox{Dom}(g^1_i)\backslash
\mbox{Dom}(f^*):i<\omega_1\ran$ are pairwise disjoint
\item[(vi)] $g^0_i\Vdash \mbox{``}\xi_i\in\name{a}$''
\item[(vii)] $g^1_i\Vdash \mbox{``}\xi_i\notin \name{a}$''
\item[(viii)] $\xi_i<\omega_1$ and $\xi_i\not=\xi_j$ for $i\not= j$.
\end{enumerate}

\n {\bf Fact}: There are $f^*$ and a witness $\lan
(g^0_i,g^1_i,\xi_i):i<\omega_1\ran$ 
above $f^*$ and $X\subseteq\aleph_{\omega_1+1}$ unbounded and an ideal
$J\supseteq J^{\rm bd}_{\omega_1}$ on $\omega_1$ such that: letting
$$
u_{\alpha ,i}=\mbox{ OP}_{A_\alpha ,j^*}\big(\mbox{Dom}(g_i^0)\cup
\mbox{Dom}(g^1_i)\backslash\mbox{Dom}(f^*)\big)
$$
\begin{enumerate}
\item[$\oplus$] if $\alpha <\beta$ are in $X$ then
$$
\{ i: u_{\alpha ,i}\cap u_{\beta ,i}\not=\emptyset\}\in J.
$$
\end{enumerate}
\smallskip

We show how to finish the proof assuming the fact, and then prove the fact.
For some unbounded $X\subseteq \aleph_{\omega_1+1}$ we have $\alpha\in
X \Rightarrow f^{**}=\mbox{OP}_{A_\alpha
,j^*}(f^*)$ i.e.\ does not depend on $\alpha\in X$. (As there are $\le
|P|=\aleph_{\omega_1}< \aleph_{\omega_1+1}$ possibilities.)

We shall prove
$$
\begin{array}{ll}
f^{**}\Vdash_P & \mbox{ ``$\lan\name{a}_\alpha :\alpha\in X\ran$ is
independent (as a family of subsets of $\omega_1$),}\\
\ & \mbox{ even modulo $J^{\rm bd}_{\omega_1}$''}.
\end{array}
$$
This is more than enough.

If not then for some $n<\omega$ and pairwise distinct $\alpha_1,\ldots
,\alpha_{2n}\in X$, we have:
$$
\neg(f^{**}\Vdash\mbox{``$\displaystyle{\bigcap^n_{\ell=1}\name{a}_{
\alpha_\ell}\cap\bigcap^{2n}_{\ell =n+1}(\omega_1\backslash\name{a}_{
\alpha_\ell})}$ is unbounded in $\omega_1$''})
$$
So for some $f^1$, $f^{**}\leq f^1\in P$, and $\zeta <\omega_1$ we
have 
$$
f^1\Vdash_P\mbox{``$\displaystyle{\bigcap^n_{\ell =1}\name{a}_{\alpha_i}\cap
\bigcap^n_{\ell =n+1}(\omega_1\backslash\name{a}_{\alpha_i})\subseteq
\zeta}$''}.\leqno\mbox{$\boxtimes$}
$$
Now letting
\begin{eqnarray*}
g^0_{\alpha ,i}&=&g^0_i\circ {\rm OP}_{j^*,A_\alpha}\;\;{\rm and}\\
g^1_{\alpha ,i}&=&g^1_i\circ {\rm OP}_{j^*,A_\alpha}
\end{eqnarray*}
we have
$$
\mbox{Dom}(g^0_{\alpha, i}) \cup \mbox{Dom}(g^1_{\alpha, i}) \subseteq
\{\gamma_{\alpha, j}: j< j^*\}.
$$
Let
\begin{eqnarray*}
B&=:&\{ i<\omega_1:\xi_i <\zeta\}\in J,\\
B_{\ell ,m}&=:& \{ i:u_{\alpha_\ell ,i}\cap
u_{\alpha_m,i}\not=\emptyset\}\in J\;\;\mbox{ for $\ell\not=m$,}\\
B_\ell &=:&\{ i:\mbox{Dom}(f^1)\cap(\mbox{Dom}(g^0_{\alpha,i})\cup\mbox{Dom}
(g^1_{\alpha ,i}))\not=\mbox{Dom}(f^{**})\}\in J\\
&&\mbox{ (in fact is finite).}
\end{eqnarray*}
So we can find $i\in \omega_1\backslash \bigcup_{\ell\not= m} B_{\ell,m}
\backslash \bigcup_\ell B_\ell\backslash B$ (because the set of inappropriate $i$'s
is in $J$).

So $\displaystyle{f^2=f^1\cup\bigcup^n_{\ell =1} g^0_{\alpha_\ell ,i}
\cup\bigcup^{2n}_{\ell =n+1}g^1_{\alpha_\ell,i} \in P}$ forces that the
intersection from $\boxtimes$ is not $\subseteq\zeta$, contradicting the
choice of $f^1$.{\hfill$\square$}
\medskip

\n {\bf Proof of the Fact}: We divide the proof into two cases, depending on
the answer to:
\medskip

\n {\bf Question}: Are there $\zeta <\omega_1$ and $j^{**}<j^*$ such that: for
no $g^0$, $g^1\in P_{j^*}$ above $p^*$ and $\xi \in [\zeta ,\omega_1]$ do we have
$$\begin{array}{l}
g^0\restriction j^{**}=g^1\restriction j^{**}\\
g^0\Vdash\mbox{``}\xi\in\name{a}\mbox{''}\\
g^1\Vdash\mbox{``}\xi\notin\name{a}\mbox{''?}
\end{array}
$$
\smallskip

\n {\bf Case A}: The answer is YES.

For some unbounded $X\subseteq \aleph_{\omega_1+1}$ and $\gamma^{**}$ we
have $j\leq j^{**} \ \&\ \alpha\in X \Rightarrow \gamma_{\alpha, j} =
\gamma^{**}_j$. 

So $\name{a}$ is actually a $P_{j^{**}}$-name.

So for $\alpha\in X$, $\name{a}_\alpha$ depends only on
$\{f\in\name{G}_P:\mbox{Dom}(f)\subseteq 
\{ \gamma_{\alpha ,i}:i<j^{**}\}\}$.

Hence there are $\le |\gamma^{**}|^{\aleph_1}<\aleph_{\omega_1}$ such names,
contradiction.
\medskip

\n {\bf Case B}: The answer is NO.

Let $\lan j(\zeta ):\zeta <\omega_1\ran$ be increasing continuous with limit
$j^*$ (because as stated above, $\bigcup_{j<j^*}\gamma^*_j=\aleph_{\omega_1}$
so necessarily ${\rm cf}(j^*)=\aleph_1$.

So for every $\zeta <\omega_1$, we have $\lan\xi^*_\zeta ,g^0_\zeta,g^1_\zeta
\ran$ giving the counterexample for $\zeta, j^{**}=j(\zeta )$,
w.l.o.g. $\mbox{Dom}(g^0_\zeta) = \mbox{Dom}(g^1_\zeta)$. As $\lan
j(\zeta ):\zeta <\omega_1\ran$ is increasing continuous, by Fodor's lemma we can
find $S\subseteq \omega_1$ stationary and $j^\otimes <j^*$ and $n^*$ such that
$$
\mbox{Dom}(g^0_\zeta\restriction j(\zeta))\cup\mbox{Dom}(g^1_\zeta
\restriction j(\zeta ))\subseteq j^\otimes,
$$
$\mbox{Dom}(g^\ell_\zeta\setminus j^\otimes)$ has $n^*$ elements and
$\zeta_1< \zeta\in S \Rightarrow \xi^*_{\zeta_1} < \zeta$.
 
We can find $X'\subseteq \aleph_{\omega_1 +1}$ unbounded such that for
all $\alpha\in X' $, $\lan 
\gamma_{\alpha ,j}:j<j^\otimes\ran
=\lan\gamma^\otimes_j:j<j^\otimes\ran$. 
Let $\langle \varepsilon(i): i<\omega_1\rangle$ be a (strictly)
increasing sequence listing $S$, and $\xi_i =
\xi^*_{\varepsilon(i)}$. For $\ell< n^*$, $\alpha\in X^\prime$ and
$i<\omega_1$ we let $f^\ell_\alpha(i)$ be the $\ell$-th member of
$\{\gamma_{\alpha, j}: j\in \mbox{Dom}(g^0_{\varepsilon(i)})\setminus
j^\otimes\}$. Now we apply Fact 7.4(1) and get an ideal $J$ (the dual
to $\DD$) and $X\subseteq X^\prime$ unbounded in $\aleph_{\omega_1+1}$
which are as required.
{\hfill$\square_{\rm 8.1}$}
\medskip

\n {\bf 8.2 Discussion}:
\smallskip
1) Clearly we can replace $\aleph_1$, $\aleph_{\omega_1+1}$ by any $\theta$,
$\lambda$ as in 7.4.

2) Normally if $\mu$ is strong limit singular of cofinality $\theta$, (at least
large enough), we can find long intervals $\aaa_i$ of the $\mbox{Reg}\cap\mu$
for $i<\theta$, $i<j\Rightarrow\sup (\aaa_i)<\min (\aaa_j)$ such that
$(\forall\bar\lambda\in\prod_i\aaa_i)$ [max pcf $(\mbox{Rang
}(\bar\lambda))=\lambda^*$] for some $\lambda^*\in [\mu ,2^\mu ]$, usually
$\mbox{cf}(2^\mu )$. This is a strong indication that $\lan I^n_{\sup
(\aaa_i), \min (\aaa_i)}:i<\theta\ran$ will have a $\lambda$-sequence, so for
example there is a $(2^\theta )^+$-c.c. Boolean algebra of cardinality
$\lambda$ having no independent subset of cardinality $\lambda$, for which
even $\lambda$-Knaster property fails.

To make this happen for no $\mu$, we need a very special pcf structure in the
universe. But we do not know even if the following simple case is consistent.
\medskip

\n {\bf 8.3 Question}: Is it consistent that
\begin{enumerate}
\item[($*$)] for every set $\aaa$ of odd (or even) regular cardinals with
$|\aaa |<\mbox{Min}(\aaa )$ we have $\mbox{max pcf}(\aaa )$ is odd (or even
respectively) (we may moreover ask $(\forall\alpha )$
$2^{\aleph_\alpha}=\aleph_{\alpha +2}$)?
\end{enumerate}

Essentially by \cite[\S5]{Sh:430}:
\medskip

\n {\bf 8.4 Lemma}: Assume $\mu >\theta=\mbox{cf}(\mu )$, $\mu$ strong limit,
$\mu =\sum_{i<\theta}\mu_i$, $\mu_i<\mu$ strong limit, $\mbox{cf}(\mu_i)=
\sigma_i$ and $2^{\mu_i}=\mu^+_i$, $\mu_i=\sum_{\zeta <\sigma_i}\mu_{i,\zeta}$,
$n_{i,\zeta}<\omega$, $\lambda =\mbox{tcf}\big(\prod_{i<\theta}\mu^+_i/J^*)$,
$J^{\rm bd}_\theta\subseteq J^*$.

Let $I_{i,\zeta}=\mbox{ERI}^{n_{i,\zeta}}_{\beth_{n_{i,\zeta}}(\mu_{i,
\zeta})^+, \mu^+_{i,\zeta}}$ and
$$
J=\sum_{J^*}J^{\rm bd}_{\sigma_i}.
$$
$\underline{\rm Then}$
\begin{enumerate}
\item[a)] there is a $(J,\lambda)$-sequence $\bar\eta$ for
$$
\lan I_{i,\zeta}:i<\theta ,\zeta <\sigma_i\ran.
$$
\item[b)] if $i<\theta\Rightarrow\sigma_i=\theta$ then we can find $\zeta
(i)<\theta$ for $i<\theta$ such that there is a $\lambda$-sequence $\bar\eta$
for $\lan I_{i,\zeta (i)}:i<\theta\ran$.
\end{enumerate}
\smallskip

\n {\bf 8.5 Remark}. So if $S=\{\mu :\mu\mbox{ strong limit }, \mbox{cf}(\mu
)=\aleph_0, 2^\mu =\mu^+\}$ is unbounded, then for a class of cardinals $\mu$
which is closed unbounded
\begin{enumerate}
\item[($*$)] (a) ~ $\mu$ strong limit and $\mu=\sup (S\cap\mu )$
\item[{}] (b) if $\mbox{cf}(\mu )=\aleph_0$ then we can find $\lambda\in (\mu
,2^\mu ]\cap\mbox{Reg}$ and $\mu_n<\mu =\sum_n\mu_n$, $\mu_n<\mu_{n+1}$ and
there is a $\lambda$-sequence $\bar\eta$ for $\lan
I^n_{\beth_n(\mu_n)^+,\mu^+_n)}: n<\omega\ran$.
\end{enumerate}
\bigskip

\section{Having a $\lambda$-sequence for a sequence of non-stationary ideals}
\smallskip

\n {\bf 9.1 Lemma}. Assume
\begin{enumerate}
\item[(a)] $\mu$ is a strong limit singular of cofinality $\theta$
\item[(b)] $\lambda =2^\mu =\mbox{cf}(\lambda )$
\item[(c)] $\lambda_i$ regular increasing for $i<\delta$ with limit $\mu$,
$\delta <\mu$ (usually $\delta =\theta$),
\item[(d)] $J$ is an ideal on $\delta$ extending $J^{\rm bd}_\delta$
\item[(e)] $\lambda=\mbox{tcf}(\prod_{i<\delta}\lambda_i/J)$,
\item[(f)] $\lan A_\zeta :\zeta <\zeta (*)\ran$ is a partition of $\delta$ (so
$A_\delta$ pairwise disjoint) each $A_\zeta$ in $J^+$ (otherwise not
interesting), 
\item[(g)] $|\delta |<\sigma =\mbox{cf}(\sigma )<\lambda_0$.
\end{enumerate}
$\underline{\rm Then}$ there is a sequence $\bar\eta =\lan \eta_\alpha :\alpha
<\lambda\ran$, $\eta_\alpha\in\prod_{i<\delta}\lambda_i$, ${\rm cf}(
\eta_\alpha (i))=\sigma$, satisfying
\begin{enumerate}
\item[($*$)] For any sequence $\lan F_{\zeta ,i}:\zeta <\zeta (*),
i<\delta\ran$ of functions, for every large enough $\alpha <\lambda$ we have
\item[($**$)] if $\zeta<\zeta (*)$, $F_{\zeta,i}(\eta_\alpha\restriction
\bigcup_{\xi <\zeta}A_\xi )$ is a club of $\lambda_i$ for $i<\delta$ (really
$i\in A_\zeta$), 
\end{enumerate}
$\underline{\rm then}$
$$
\{ i\in A_\zeta :\eta_\alpha (i)\notin F_{\zeta,i}(\eta_\alpha\restriction
\bigcup_{\xi <\zeta}A_\zeta )\}\in J.
$$
Moreover
\begin{enumerate}
\item[($**$)$^+$] if $\zeta <\zeta (*)$, $n <\omega$ and $\beta_0 ,\ldots
,\beta_{n-1}<\alpha$, and for each $i\in A_\zeta$ we have: $F_{\zeta,i}(
\eta_{\beta_0},\beta_0\ldots,\eta_{\beta_{n-1}},\beta_{n_1},\eta_\alpha
\restriction \bigcup_{\xi<\zeta}A_\xi )$ is a club of $\lambda_i$, then
$$
\{ i\in A_\zeta :\eta_\alpha (i)\not\in F_{\zeta ,i}(\eta_{\beta_0}, \beta_0
\ldots ,\eta_{\beta_{n-1}},\beta_{n-1}, \eta_\alpha\restriction\bigcup_{\xi
<\zeta}A_\xi )\}\in J.
$$
\end{enumerate}
\smallskip

\n {\bf 9.2 Discussion}. For a given $\mu$ as in ($a$), clause ($b$) may fail,
but then we will have another lemma. What about ($e$)?

If $\theta >\aleph_0$ there are such $\lan \lambda_i:i<\theta\ran$ even for
$J=J^{\rm bd}_\theta$ (see \cite[VIII, \S1]{Sh:g}. If $\theta=\aleph_0$ we do
not know, but we know that the failures are ``rare''. E.g.\
$$
\{\delta <\omega_1:\beth_\delta\;\mbox{fails ($e$) i.e.}\;
\neg[\beth_{\delta+1} =^+{\rm pp}(\beth_\delta)]\}
$$
is not stationary. About pp$_{J^{\rm bd}_\omega}$ e.g.\ if $|\aaa|\le
\aleph_0\Rightarrow |{\rm pcf}(\aaa )|\le\aleph_0$ we then can get it, see 
\cite{Sh:E12}, Part C, on \cite[IX]{Sh:g}.
\medskip

\n {\bf 9.3 Remark}. 1) This can be rephrased as having a $(\lambda
,J)$-sequence for $\lan\Pi J^{{\rm nst},\sigma}_{\lambda_{i,n}}:i<\delta\ran$
with $\lambda_{i,n}$ decreasing.

So compared to earlier theorems, the $\lambda, \lambda_i$ for which the Lemma
applies are fewer, but the result is stronger: nonstationary ideal and we get
also the ``super'' version see ($**$).

2) Of course another variant is to start with $I_i=J_{\lambda_i}^{{\rm
nst},\sigma}$ and get $J=J_\lambda^{{\rm nst},\sigma}$.

3) Considering functions with finitely many $\eta_\beta$'s, $\beta <\alpha$ as
parameters (i.e. $(\ast\ast)^+$); thinning $\lan f_\alpha :\alpha
<\lambda\ran$ the conclusion follows. 

4) In ($**$)$^+$ instead $n<\omega$ we can ask $n<\sigma$ if $(\forall\alpha
<\lambda)$ ($|\alpha |^{<\sigma}<\lambda$) \medskip

\n {\bf Proof of 9.1}. For simplicity we concentrate on ($**$) (in 10.1 we
concentrate on the parallel of ($**$)$^+$). List the possible $\lan F_{\zeta
,i}:i<\delta,\zeta<\zeta (*)\ran$, i.e.\ sequence with each $F_{\zeta,i}$
having the ``right" domain and range, which are clear from the statement, as
$\big\langle \lan F^\beta_{\zeta ,i}:i<\delta ,\zeta <\zeta (*)\ran :\beta
<\lambda\big\rangle$. Let us define $\eta_\alpha\in\prod_i\lambda_i$ by
induction on $\alpha$.

For a given $\alpha$ we choose $\eta_\alpha\restriction A_\zeta$ by induction
on $\zeta <\zeta (*)$.

Define for $i\in A_\zeta$, $\beta <\alpha$
$$
C^\beta_i=\left\{\begin{array}{ll}
\displaystyle{F^\beta_{\zeta,i}\big(\eta_\alpha\restriction\bigcup_{\xi
<\zeta}A_\zeta\big)}&\quad\mbox{if this set is a club of $\lambda_i$}\\
\lambda_i &\quad\mbox{otherwise}.\end{array}\right.
$$
So we need
\medskip

\n {\bf Fact}. There is $\eta\in\prod_{i\in A_\zeta}\lambda_i$ such that
$$
\begin{array}{l} 
\displaystyle{\bigwedge_{\beta <\alpha}\{ i\in A_\zeta:\eta (i)\notin
C^\beta_i \}\in J}\\ 
\displaystyle{i\in A_\zeta\Rightarrow {\rm cf}(\eta(i))=\sigma}.
\end{array}
$$
\smallskip

\n {\bf Proof of the Fact}. We shall choose by induction on $\e <\sigma$ a
function $g_\e\in\prod_{i\in A_\zeta}\lambda_i$ such that $\e_1<\e\Rightarrow
g_{\e_1}<g_\e$ (in all coordinates) and
$$
(\forall \beta <\alpha )(\forall^Ji\in A_\zeta )[(g_\e (i),
g_{\e +1}(i))\cap C^\beta_i\not= \emptyset].
$$
Why is this enough?

Let $\nu=\eta\restriction A_\zeta$ be defined by
$$
\nu(i)=\bigcup_{\e<\sigma} g_\e (i).
$$
Now $\nu (i)<\lambda_i$ as $g_\e(i)<\lambda_i$ and $\sigma <\lambda_i={\rm
cf}(\lambda_i)$. (We can also say something for $\sigma\ge\mu$, but not now.)
Also $\lan g_\e (i):\e <\sigma\ran$ is strictly increasing, so $\mbox{cf}(\nu
(i))=\sigma$.

Now let $\beta <\alpha$ and define
$$
B^*_\beta =\{ i\in A_\zeta :\nu (i)\notin C^\beta_i\}.
$$
We would like to have $B^*_\beta\in J$. For each $i\in B^*_\beta$, the
sequence $\lan g_\e (i):i<\sigma\ran$ is a strictly increasing sequence of
ordinals with limit not in $C^\beta_i$.

So for some $\e_{\beta ,i}<\sigma$
$$
C^\beta_i\cap (g_{\e_{\beta ,i}}(i), \nu (i)) =\emptyset.
$$
So
$$
\bigwedge_{\e\ge\e_{\beta ,i}}(g_{\e (i)}, g_{\e +1}(i))\cap
C_{\beta ,i}=\emptyset.
$$
Let $\e_\beta =\sup_{i<\delta}\e_{\beta ,i}$.

Now $\e_{\beta ,i}<\sigma$ \& $\sigma =\mbox{cf}(\sigma )>|\delta |\ge
|A_\zeta |$, so $\e_\beta <\sigma$. So
$$
\bigwedge_{i\in B^*_\beta}(g_{\e_\beta}(i), g_{\e_{\beta +1}}(i))\cap
C^\beta_i=\emptyset
$$
hence $B^*_\beta\in J$ as required, i.e.~$\nu$ is the required $\eta$.

Why is the choice of the $g_\e$ possible?
\medskip

\n {\bf Construction}.

$\underline{\e=0}$. Trivial.
\smallskip

$\underline{\e {\rm ~limit}}$. $g_\e (i)=\bigcup_{\e_1<\e}
g_{\e_1}(i)<\lambda_i$ (as $\e<\sigma <\lambda_i=\mbox{cf}(\lambda_i)$).
\smallskip

$\underline{\e+1}$. For $\beta <\alpha$ define $h_{\beta ,\e}\in\prod_{i\in
A_\zeta}\lambda_i$ by
$$
h_{\beta ,\e}(i){\stackrel{\rm def}{=}}\min\{\gamma <\lambda_i:(g_\e (i),
\gamma )\cap C^\beta_i\not= \emptyset\}.
$$
So $\{ h_{\beta ,\e}:\beta <\alpha\}$ is a subset of $\prod_{i\in
A_\zeta}\lambda_i$ of cardinality $<\lambda$, but $\prod_{i<\delta}\lambda_i/
J$ hence $\prod_{i\in A_\zeta}\lambda_i/(J\restriction A_\zeta )$ has true
cofinality $\lambda$ (as if $A_\zeta\in J$ there is nothing to prove). So
there is $g'_\e\in\prod_{i\in A_\zeta}\lambda_i$ which is $a<_{J\restriction
A_\zeta}$-upper bound of $\{ h_{\beta ,\e}:\beta<\alpha\}$.

Let $g_{\e +1}(i)=\max\{ g'_\e (i), g_\e (i)+1\}$, clearly it is as
required.{\hfill$\square_{\rm 9.1}$}
\medskip

\n {\bf 9.4 Claim}. 1) Assume
\begin{enumerate}
\item[(a)] $\bar\eta=\lan\eta_\alpha :\alpha <\lambda\ran$, where
\item[{}] $\displaystyle{\eta_\alpha\in\prod_{i\in{\rm Dom}(J)}\mbox{Dom}
(I_i)}$ and $J$ is an ideal on $\delta$ extending $J_\delta^{\rm bd}$, each
$I_i$ and ideal and $I$ an ideal on $\lambda$ extending $J^{\rm bd}_\lambda$
\item[(b)] $\lan A_\zeta:\zeta <\zeta (*)\ran$ is a partition of ${\rm Dom}
(J)$, $A_\zeta\notin J$
\item[(c)] for every $\bar F=\lan F_i:i\in{\rm Dom}(J)\ran$, for the
$I$-majority of $\alpha <\lambda$, for every $\zeta <\zeta(*)$ $\underline{\rm
if}$ $F_i(\eta_\alpha\restriction\bigcup_{\xi <\zeta} A_\xi )\in I_i$ for
$i\in A_\zeta$, $\underline{\rm then}$ 
$$
(\forall^Ji\in A_\zeta )[(\eta_\alpha (i)\notin
F_i\big(\eta_\alpha\restriction\bigcup_{\xi <\zeta}A_\xi\big)].
$$
\item[(d)] $\displaystyle{I^*_j=\prod_{\ell< n_j}I_{i(j,\ell)}}$ for
$j<\delta^*$ where $i(j,\ell )<\delta$ 
\item[(e)] $J^*=\big\{A\subseteq\delta^*:\mbox{ for some }B\subseteq\delta, \
\ \displaystyle{\bigwedge_\zeta (B\cap A_\zeta )\in J}\ \mbox{ and }\
\;\displaystyle{\bigwedge_{i\in A}\bigvee_{\ell<n_j}i(j,\ell)\in B}\}$\quad is
an ideal on $\delta^*$. 
\item[(f)] $\eta^*_\alpha$ is defined by
$$
\eta^*_\alpha(j)=\lan\eta_\alpha (i(j,\ell )):\ell <n_j\ran
$$
\end{enumerate}
$\underline{\rm Then}$\footnote{so we have dealt here with the case of $J^{\rm
bd}_{\bar\lambda}$, $\bar\lambda$ decreasing} $\lan\eta^*_\alpha :\alpha
<\lambda\ran$ is a $(\lambda ,J^*, I)$-sequence for $\lan I^*_j:j<\delta\ran$.
\smallskip

2) If we strengthen clause (c) to the parallel of $(\ast\ast)^+$ in 9.1,
\underline{then} $\langle \eta_\alpha^\ast:\,\alpha<\lambda\rangle$ is a super
$(\lambda, J^\ast, I)$-sequence for $\langle I^\ast_i:\,i< \delta\ran$.

\smallskip
\n {\bf Proof}: Straightforward.{\hfill$\square_{\rm 9.4}$}
\medskip

\n {\bf 9.4A Conclusion}. Assume (a)--(g) of 9.1 (see 9.2) and (a), (e) of
9.4. Then there is a super $(\lambda , J^*)$-sequence for $\lan
I^*_j:j<\delta\ran$.
\medskip

\n {\bf 9.5 Conclusion}. Assume $\mu >\mbox{ cf}(\mu )=\aleph_0$ is a strong
limit, and
$$
\lambda=2^\mu=\mbox{cf}(2^\mu)=\mbox{tcf}\Big(\prod_{n<\omega}\lambda_n/J^{\rm
bd}_\omega\Big),
$$
$\lambda_n$ regular $<\mu$.

Let $\lan k_n:n<\omega\ran$ be such that
$$
(\forall k)(\exists^\infty n)\; (k_n=k),
$$
and e.g.\ $\theta =(2^{\aleph_0})^+$.

For $n<\omega$ and $k<k_n$ let $\ell(n,k)=\sum\{k_m:m<n\}+k$ and let
\begin{eqnarray*}
I_n&=& \prod_{k<k_n} J^{{\rm nst},\theta}_{\lambda_{\ell(n,k)}}\\ J&=&\{
A\subseteq \omega :\sup_{n\in A}k_n<\omega\}. 
\end{eqnarray*}

\n $\underline{\rm Then}$ there is a $(\lambda ,J)$-sequence for $\lan
I_n:n<\omega \ran$ (even a super one)..
\medskip

\n {\bf Proof}. By Lemma 9.1 and Claim 9.4, we choose in 9.4 the parameters
$\delta =\omega$, $\zeta (*)=\omega$ and let
$$
A_\zeta =\big\{\sum_{m\le n}k_m-\zeta:k_n>\zeta\big\}.
$$
 ~{\hfill$\square_{\rm 9.5}$}
\medskip

We may wonder on the ``tcf'' assumption; at the expense of using ``some $J$''
this can be overcome:
\medskip

\n {\bf 9.6 Claim}. Assume $\mu >\mbox{ cf}(\mu )=\aleph_0$ strong limit
singular,
\begin{eqnarray*}
\lambda &=& 2^\mu =\mbox{ cf}(2^\mu )\in\mbox{ pcf}\{ \lambda_n:n<\omega\}\\
\lambda_n &=& \mbox{cf}(\lambda_n)<\mu
\end{eqnarray*}
and $\lan k_n:n<\omega\ran$ is as in 9.5.  $\underline{\rm Then}$ we can find
$i(n,\ell )<\omega$, $\ell <k_n$ with no repetitions,
$$
i(n,0)>i(n-1,k_{n-1}-1)>\cdots >i(n-1,0),
$$
and letting
$$
I_n=\prod_{\ell <k_n} J^{{\rm nst},\theta}_{\lambda_{i(n,\ell)}},
$$
we have: for some ideal $J\supseteq J^{\rm bd}_\omega$ on $\omega$, there is a
$(\lambda ,J)$-sequence for $\lan I_n:n<\omega\ran$.
\medskip

\n {\bf Proof}. Let
$$
\mbox{pcf}_{J^{\rm bd}_\omega}(\{\lambda_n:n<\omega\}) =\left\{
\begin{array}{l}
\displaystyle{\chi :\mbox{cf}(\chi )=\chi =\mbox{tcf}\big(\prod_{n\in A}
\lambda_n/J^{\rm bd}_A\big)}\\ 
\qquad\qquad\mbox{for some infinite $A\subseteq \omega$}\end{array}\right\}.
$$
\smallskip

By a pcf claim:
\medskip

\n {\bf 9.6A Fact}. We can find increasing $\lan\chi_\e :\e <\e (*)\ran$,
$\e(*)<\omega_1$, a limit ordinal, $J^*$ an ideal $\supseteq J^{\rm bd}_{\e
(*)}$, such that
$$
\chi_\e \in \mbox{pcf}_{J^{\rm bd}_\omega}(\{\lambda_n:n<\omega\}),
$$
say
$$
\chi_\e =\mbox{tcf}\Big(\prod_{n\in B_\e}\lambda_n/ J^{\rm bd}_{B_\e}\Big),
$$
$\lan B_\e :\e <\e (*)\ran$ is a partition of $\omega$, and
$$
\lambda= {\rm tcf}\Big(\prod_{\e<\e (*)}\lambda_\e /J^*\Big).
$$
\smallskip

\n {\bf Continuation of The proof of 9.6}.  Let $\lan k_n:n<\omega\ran$ be as
before. Choose
$$
\lan i(n,\ell ):\ell <k_n\ran\;\; \mbox{for each $n$}
$$
such that
\begin{enumerate}
\item[(a)] $i(n,\ell )>i(n,\ell +1)$
\item[{}] $i(n,\ell_1)<i(n+1,\ell_2)$, and
\item[(b)] for every $k$ and $\e_0,\ldots ,\e_{k-1}$, for infinitely many $n$
we have
$$
k_n=k,\qquad\quad i(n,\ell )\in B_{\e_\ell}.
$$
\end{enumerate}

\n Let
$$
A_\ell =\{ i(n,\ell ):n<\omega, k_n>\ell\}.
$$
So
\begin{enumerate}
\item[{}] $\lan A_\ell :\ell <\omega\ran$ is a sequence of pairwise disjoint
subsets of $\omega$ such that $|A_\ell\cap B_{\e_\ell}|=\aleph_0$.
\end{enumerate}
We apply 9.1 for
$$
\lan A_n:n<\omega\ran,\;\; \lan\lambda_n:n<\omega\ran, \lambda ,\mu.
$$
~{\hfill$\square_{\rm 9.6}$}
\medskip

\n {\bf 9.6B Remark}. If $\mu >{\rm cf}(\mu )>\aleph_0$, $2^\mu$ regular, the
parallel to 9.5 always occurs.{\hfill$\square_{\rm 9.6}$}
\medskip

If we use $\bar A=\lan A_0\ran$, $A_0=\delta$ in 9.1:
\medskip

\n {\bf 9.6C Conclusion}
\smallskip

In 9.1 we get:
\begin{enumerate}
\item[{}] there is a $(\lambda ,J)$-sequence for $\lan I_i:i<\delta\ran$, even
a super one.
\end{enumerate}
\smallskip

\n {\bf Remark 9.7}. By the proofs in \cite[\S1]{Sh:420} we can replace $\lan
S^{\lambda_i}_\theta :i<\delta\ran$, $S^\lambda_\theta =\{\delta <\lambda_i:
{\rm cf}(\delta )=\theta\}$ by some large enough $\bar S =\lan
S_i:i<\delta\ran$, where $S_i\in I[\lambda_i]$, see below.

Also if $\lan f_\alpha :\alpha <\lambda\ran$ is $<_J$-increasing cofinal in
$\prod_{i<\delta}\lambda_i/J$, continuous when it can be, then for some club
$E$ of $\lambda$ we have $\lan f_\delta :\delta \in E$, $\mbox{cf}(\delta
)=\theta$, $\bar f\restriction\delta$ has an exact least upper bound lub$\ran$
is OK.  Probably more interesting is to strengthen $I^{{\rm nst},\theta}_{
\lambda_i}$ by club guessing, as follows.
\medskip

\n {\bf 9.8 Definition}. {\it For $\bar C=\lan C_\delta :\delta\in S\ran$,
$S\subseteq\lambda$, stationary
$$
{\rm id}^a(\bar C)=\left\{\begin{array}{l}
A\subseteq\lambda:\mbox{ for some club
$E$ of $\lambda$ the set}\\
\mbox{\hspace*{.75truein} $\{\delta\in S: C_\delta\subseteq
E\}$ is not stationary}\\
\mbox{\hspace*{.75truein} (so as we can shrink $E$,
equivalently, empty)}\end{array}\right\}.
$$}
\smallskip

\n {\bf 9.9 Lemma}. Assume
\begin{enumerate}
\item[(a)] $\mu$ is a strong limit singular $\theta$
\item[(b)] $\lambda =2^\mu =\mbox{cf}(\lambda )$
\item[(c)] $\lambda_i$ regular increasing for $i<\delta$ with limit $\mu$, 
$\delta <\mu$ (usually $\delta ={\rm cf}(\mu )$),
\item[(d)] $J$ is an ideal on $\delta$ extending $J^{\rm bd}_\delta$
\item[(e)]  $\lambda=\mbox{tcf}(\prod_{i<\delta}\lambda_i/J)$,
\item[(f)] $\lan A_\zeta :\zeta <\zeta (*)\ran$ is a partition of $\delta$ (so
pairwise disjoint)
\item[(g)] $\sigma ={\rm cf}(\sigma )<\mu$, moreover $\sigma <\lambda_0$ and
satisfies

$\otimes^{\sigma ,\delta}_J$ we have $\sigma >\delta$ (or at least if $A_\e
\in J$ for $\e<\sigma$ then
$$
\{ i<\delta:i\in A_\e\;\mbox{for every large enough $\e<\sigma$}\} \in J)
$$
\end{enumerate}
\underline{Then}

1) ~For $\theta\in {\rm Reg}\cap (\sigma ,\lambda_0)$ we can find $\lan
S_i:i<\delta\ran$, $\lan \bar C^i:i<\delta\ran$, $\bar I=\lan
I_i:i<\delta\ran$, $\bar\eta =\lan\eta_\alpha :\alpha <\lambda\ran$ such that
\begin{enumerate}
\item[($\alpha$)] $S_i\in I[\lambda_i]$ is stationary, and $\delta\in 
S_i\Rightarrow {\rm cf}(\delta )=\sigma$
\item[$(\beta)$] $\bar C^i=\lan C^i_\delta :\delta\in S_i\ran$, $C^i_\delta$ a
club of $\delta$,
\item[$(\gamma)$] $I_i={\rm id}^a(\bar C^i)=\big\{A\subseteq\lambda_i:\;$ for
some club $E$ of $\lambda_i$\ \ we have:\ \ $\delta\in S\cap A_i$ implies
$\sup (C^i_\delta\backslash E)<\delta\big\}$
\item[$(\delta)$] 
\begin{enumerate}
\item[($*$)] For any sequence $\lan F_{\zeta ,i}:\zeta <\zeta (*),
i<\delta\ran$ of functions, for every large enough $\alpha <\lambda$ we have
\item[($**$)] if $\zeta<\zeta (*)$, $F_{\zeta,i}(\eta_\alpha\restriction
\bigcup_{\xi <\zeta}A_\xi )$ a member of ${\rm id}^a(\bar C^i)$ for $i<\delta$
(really $i\in A_\zeta$), 
\end{enumerate}
\end{enumerate}
$\underline{\rm then}$
$$
\{ i\in A_\zeta :\eta_\alpha (i)\in F_{\zeta,i}(\eta_\alpha\restriction
\bigcup_{\xi <\zeta}A_\zeta )\}\in J.
$$
Moreover
\begin{enumerate}
\item[($**$)$^+$] if $\zeta <\zeta (*)$, $n <\omega$ and $\beta_0 ,\ldots
,\beta_{n-1}<\alpha$ and for each $i\in A_\zeta$ we have: $F_{\zeta,i}(
\eta_{\beta_0},\beta_0\ldots,\eta_{\beta_{n-1}},\beta_{n-1},\eta_\alpha
\restriction \bigcup_{\xi<\zeta}A_\xi )$ in a member of ${\rm id}^a(\bar C^i)$
then 
$$
\big\{i\in A_\zeta :\eta_\alpha (i)\in F_{\zeta ,i}(\eta_{\beta_0},\beta_0,
\ldots,\eta_{\beta_{n-1}},\beta_{n-1},\eta_\alpha\restriction\bigcup_{\xi
<\zeta}A_\xi)\big\}\in J.
$$
\end{enumerate}
\smallskip

\n {\bf Remark}: 1) Included in the proof are imitations of proofs from 
\cite[\S1]{Sh:420} and of 9.1.

2) We have a bit of flexibility in the proof.

3) In ($**$)$^+$, we can replace $n<\omega$ by $n<\tau$ when $(\forall\alpha
<\lambda)$ ($|\alpha |^{<\tau}<\lambda$)
\medskip

\n {\bf Proof}: Let $\theta =2^\sigma$. By \cite[\S1]{Sh:420} we can find
$\bar e^i$ such that:
\begin{description}
\item[(i)] for $i<\delta$, $\bar e^i=\lan e^i_\alpha:\alpha\in S_i\ran$,
$S_i\in I[\lambda]$, 
\item[(ii)] $e^i_\alpha$ a club of $\alpha$ of order type $\sigma$ such that
\quad $\alpha \in S_i\Rightarrow {\rm cf}(\alpha )=\sigma$, 
\item[(iii)] for $\chi$ large enough, $x\in\HH (\chi )$, we can find $\lan
N_i:i\le \sigma\ran$ such that $x\in N_\e\prec (\H (\chi ),\in,<^*_\chi )$,
$\lan N_\zeta :\zeta\le\e\ran\in N_{\e+1}$, $N_\e$ increasing continuous, $\|
N_\e\|=\theta$, $\theta +1\subseteq N_\e$ and
$$
i<\delta\quad\Rightarrow\quad\sup e^i_{\sup (N_{\sigma}\cap\lambda_i)}\in S_i.
$$
\end{description}
For $\bar d\in\bigcup\big\{\prod\limits_{i<\delta}e_i:e_i$ a club of
$\sigma\big\}$ let $\bar e^{i,\bar d}=\lan e^{i,\bar d}_\alpha:\alpha\in S_i
\ran$, $e^{i,\bar{d}}_\alpha =\lan\beta\in e^i_\alpha:{\rm otp}(e^i_\alpha
\cap\beta )\in d_i\ran$. For each such $\bar d$ we repeat the proof of 9.1, so
we choose $\eta_\alpha=\eta^{\bar d}_\alpha$ by induction on $\alpha<\lambda$,
and for each $\alpha$, choose $\eta_\alpha\restriction\big(\bigcup\limits_{\e
<\zeta}A_\e\big)$ by induction on $\zeta\le\zeta(*)$. If we succeed fine, so
assume we fail. So for some $\alpha=\alpha[\bar d]$, $\zeta=\zeta[\bar d]$ the
situation is: $\lan\eta^{\bar d}_\beta:\beta<\alpha\ran$ and $\eta^{\bar
d}_\alpha\restriction \big(\bigcup\limits_{\e<\zeta} A_\e \big)$ are defined,
but we cannot define $\eta^{\bar d}_\alpha\restriction A_\zeta$ and as there
we can compute a family $\E =\E^i_{\bar d}$ of cardinality $<\lambda$ whose
members has the form $\bar B=\lan B_i:i<\delta\ran$, $B_i\in {\rm id}^a (\bar
e^{i,d})$ and let $E^i_{B_i}$ be a club of $\lambda_i$ exemplifying $B_i\in
{\rm id}^a(\bar e^{i,d})$; let $\E^i_{\bar d}=\{\lan E^i_{B_i}:i<\delta\ran:
\bar B=\lan B_i:i<\delta\ran\in \E\}$. Let $\lan N_i:i\leq\sigma\ran$ be as in
$\otimes$(iii) for $x=\big\{\lan\lan\E^1_{\bar d},\bar d\ran :d_i\subseteq
\sigma\;\mbox{ a club for $i<\delta$}\ran$, $\bar\lambda$, $\lan\bar
e^i:i<\delta\ran\big\}$.

As in the proof of 9.1 quite easily:
$$
\e\le\sigma\;\mbox{\&}\; \bar B=\lan B_i:i<\delta\ran\in\bigcup_{\bar d}\E_{
\bar d}\quad \Rightarrow\quad \{ i<\delta:\sup (N_\e\cap \lambda_i)\not\in
E^i_{B_i}\}\in J. 
$$
Let $d_i=\big\{{\rm otp}(e^i_{\sup (N_\sigma\cap\lambda_i)}\cap\sup (N_\e\cap
\lambda_\ell )): \e<\sigma$ and $\sup (N_\e\cap\lambda_i)\in e^i_{\sup
(N_\sigma\cap\lambda_i)}\big\}$. Clearly $d_i$ is a club of $\sigma$ and let
$\bar d=\lan d_i:i<\delta\ran$. Now $\lan\sup (N_\sigma\cap \lambda_i):i\in
A_{\zeta[\bar d]}$ is as required.
\medskip

\centerline{* * *}
\medskip

\n {\bf 9.9A Conclusion}. 1)\ \ In 9.9 we get:
\begin{quotation}
\noindent for some function $c:[\lambda]^2\longrightarrow\sigma$, for every
$X,Y\in [\lambda]^\lambda$ and $\zeta<\sigma$, for some $\alpha\in X$,
$\beta\in Y$ we have $\alpha>\beta$ and $c(\{\alpha,\beta\})=\zeta$.
\end{quotation}
2)\ \ In 9.9 we can add:
\begin{quotation}
\noindent if e.g.~$\chi=(2^\lambda)^+$, for every $X\subseteq 2^\mu$ for every
$\alpha<\lambda$ large enough, for $\zeta<\zeta(*)$, there is a sequence
$\langle N_{\varepsilon}:\varepsilon <\sigma\rangle$ as in the proof of
9.9 such that
\begin{description}
\item[$(\boxtimes)$] \quad $\{i\in A_\zeta: \eta_\alpha(i)\neq\sup( N_{
\delta}\cap\lambda_1)\}\in J$
\end{description}
\end{quotation}

\n {\bf 9.9B Remark}: In 9.9A(1) we get even $Pr_1(\lambda, \lambda,
\sigma, \sigma )$.
\medskip

\n {\bf Proof}: 1)\ \ \ We relay on part 2).

\n 2)\ \ \ For $\alpha<\beta$ let $c(\{\alpha,\beta\})=\zeta$ if 
$$
\{i\in A_0: f_\beta(i)\geq f_\alpha(i)\ \mbox{ or }\ f_\beta(i)<f_\alpha(i)\
\&\ \zeta\neq{\rm otp}(e^i_{f_\alpha(i)}\cap\beta)\}\in J,
$$
and zero if there is no such $\zeta$.

\n Let $X,Y\in [\lambda]^\lambda$. take $\alpha\in X$ large enough, so that we
can find $\langle N_\varepsilon:\varepsilon\leq\sigma\rangle$ as there, with
$(\boxtimes)$ for part (2). We can find $\beta\in N_{\zeta+1}\cap Y$ such
that $\langle
\sup(N_\zeta\cap\lambda_i): i<\delta\rangle <_J \eta_\beta$ (as $Y\cap
N_{\zeta+1}$ is unbounded in $\lambda\cap N_{\zeta+1}$). Now $\alpha>\beta$
are as required.{\hfill$\square_{\rm 9.9A}$}
\medskip

\n {\bf 9.10 Claim}. In 9.1 
\begin{enumerate} \item[1)] Instead of ``$\mu>\theta ={\rm cf}(\theta )>|\delta
|$'' we can assume only
\begin{enumerate}
\item[$\otimes_1$] $\mu >\theta=\mbox{cf}(\theta )$ and if $\lan u_\zeta
:\zeta <\theta\ran$ is a sequence of members of $J$ then
$$
\big\{ i<\delta :\theta =\sup\{ \zeta: i\notin u_\zeta\}\big\} =\delta \;{\rm
mod }J.
$$
\end{enumerate}
\item[2)] Weakening the conclusion of 9.1 to ``weak $(J,\lambda )$-sequence'',
we can replace ``$\theta=\mbox{cf}(\theta )>|\delta |$'' by
\begin{enumerate}
\item[$\otimes_2$] $\theta=\mbox{cf}(\theta )$ and if $\lan u_\zeta :\zeta
<\theta\ran$ is a sequence of members of $J$ then
$$
\big\{ i<\delta :\theta =\sup\{\zeta :i\notin u_\zeta\}\big\}\in J^+.
$$
\end{enumerate}
\item[3)] In part (1) and (2), if $\theta >\aleph_0$, then we can find $\bar
C^i=\lan C^i_\delta :\delta\in S_\theta^{\lambda_i}\ran$ with $C^i_\delta$ a
club of $\delta$ such that: we can replace $I^{{\rm nst}, \theta}_{\lambda_i}$
by ${\rm id}^a_{\lambda_i}(\bar C^i)$, see 9.9. above.
\end{enumerate}

\section{The power of a strong limit singular is itself singular: existence}

\n {\bf 10.1 Lemma}. Assume
\begin{enumerate}
\item[(a)] $\mu$ strong limit singular.
\item[(b)]$2^\mu$ is singular, $\lambda ={\rm cf}(2^\mu )$
(so $2^\mu >\lambda >\mu$)
\item[(c)] $\mu >\sigma ={\rm cf}(\sigma )>{\rm cf}(\mu )$.
\item[(d)] $2^\mu ={\rm pp}(\mu )$ (see discussion in \S9)
\end{enumerate}

\n $\underline{\rm Then}$
\begin{enumerate}
\item[($\alpha$)] we can find $J$, $J^*$, $\bar\theta^i=\lan \theta^i_\zeta:
\zeta<\mbox{ cf}(\mu )\ran$ for $i<\lambda$ and $\bar\lambda$ such that
\begin{enumerate}
\item[(i)] $\bar\theta^i$ is an increasing sequence of regular cardinals
$<\mu$ with limit $\mu$ for $i<\lambda$.
\item[(ii)] $\bar\lambda =\lan\lambda_\alpha :\alpha <\lambda\ran$ is an
increasing sequence of regulars $\in (\mu+\lambda,2^\mu )$ with limit $2^\mu$.
\item[(iii)] $J\subseteq J^*$ are ideals on $\mbox{cf}(\mu )$, $\mbox{cf}(\mu
)$-complete.
\item[(iv)] $\lambda_\alpha =\mbox{tcf}\big(\prod_{\zeta }\theta^\alpha_\zeta
/J\big)$
\item[(v)] $\lan\bar\theta^\alpha :\alpha <\lambda\ran$ is
$<_{J^*}$-increasing, i.e.\ $\alpha <\beta\to \{\zeta <\mbox{cf}(\mu 
):\theta^\alpha_\zeta \ge\theta^\beta_\zeta\}\in J^*$, with $<_{J^*}$-exact
upper bound $\lan\theta^*_\zeta :\zeta <\mbox{cf}(\mu )\ran$ and
($\theta^*_\zeta$ is a cardinal $<\mu$, normally singular) $\mu =\lim$
$\lan\theta^*_\zeta :\zeta<{\rm cf}(\mu )\ran$ and
$$
\bigwedge_{\alpha <\lambda\atop\zeta <{\rm cf}(\mu )}\theta^\alpha_\zeta
<\theta^*_\zeta.
$$
\item[(vi)] If $J\not= J^{\rm bd}_{{\rm cf}(\mu )}$, then ${\rm cf}(\mu
)=\aleph_0$ and ${\rm pp}_{J^{\rm bd}_{{\rm cf }\mu}}(\mu )<2^\mu$ and $J$
as in 9.6 so for most such $\mu$ we have the conclusion of (1), see
\cite{Sh:E12} and \S4.
\end{enumerate}
\item[($\beta$)] If $J$, $\bar\theta^\alpha (\alpha <\lambda )$, $\bar\lambda$
are as in clause ($\alpha$) $\underline{\rm then}$ we can find $\bar\eta
=\lan\eta_\alpha :\alpha <\lambda\ran$ such that
\begin{enumerate}
\item[(i)] $\bar\eta =\lan\eta_\alpha:\alpha <\lambda\ran$, $\eta_\alpha\in
\prod_{\zeta <{\rm cf}(\mu)}\theta^*_\zeta \subseteq {}^{{\rm cf}(\mu )}\mu$.
Moreover $\eta_\alpha\in\prod_{\zeta <{\rm cf}(\mu)}\theta^\alpha_\zeta$ and
$\sigma ={\rm cf}(\eta_\alpha (i))$ for $\alpha<\lambda$, $i<{\rm cf}(\mu)$.
\item[(ii)] if $\bar C=\lan C_\zeta:\zeta<{\rm cf}(\mu )\ran$,
$\theta^\alpha_\zeta\cap C_\zeta$ a club of $\theta^\alpha_\zeta$ for $\alpha
<\lambda$, $\zeta <{\rm cf}(\mu )$, $\underline{\rm then}$ for some
$\alpha^*=\alpha^*_{\bar C}$ we have 
$$
\alpha \in [\alpha^*,\lambda )\Rightarrow (\forall^J\zeta<{\rm cf}(\mu))[
\eta_\alpha (\zeta )\in C_\zeta].
$$
\end{enumerate}
\item[($\gamma$)] Assume $\lan A_\e :\e <\e^*\ran$ is a partition of
$\mbox{cf}(\mu )$ to sets not in $J$.  $\underline{\rm Then}$ we can add
\end{enumerate}
\begin{enumerate}
\item[(ii)$^+$] For any sequence of functions
$$
F=\lan F_\zeta :\zeta <{\rm cf}(\mu )\ran,
$$
for some $\alpha^*=\alpha^*_{\bar F}$, for every $\alpha\in
[\alpha^*,\lambda )$ we have
\end{enumerate}
\begin{enumerate}
\item[($*$)] if $\e <\e^*$, $n<\omega$, $\beta_\ell <\alpha$ for $\ell <n$
$\underline{\rm then}$
$$
\left\{\begin{array}{r}
\zeta <{\rm cf}(\mu ):F_\zeta (\ldots ,\beta_\ell,\eta_{\beta_\ell},\ldots
,\eta_\alpha\restriction \bigcup\limits_{\xi<\e}A_\zeta )\cap
\theta^\alpha_\zeta\; \mbox{ is a club of $\theta^\alpha_\zeta$}\\
\mbox{but $\eta_\alpha (\zeta )\notin F_\zeta(\beta_\ell,\eta_{\beta_\ell},
\ldots,\eta_\alpha\restriction \bigcup\limits_{\xi <\e} A_\xi)\cap
\theta^\alpha_\zeta$}\end{array}\right\}
$$
belongs to $J$. (If we use constant $F$ this reduces to (ii).)
\end{enumerate}
\smallskip

\n {\bf Proof}. Of clause ($\alpha$):

First choose $\lan\lambda_\alpha^0 :\alpha <\lambda\ran$ as demanded in clause
(ii) (but we will manipulate it later, possible by clause (e)). Now as in
9.6, for each $\alpha$ there are
$$
J_\alpha ,\bar\theta^\alpha =\lan\theta^\alpha_\zeta :\zeta <{\rm cf}(\mu
)\ran
$$
as there, so satisfying (i), (iii), (iv), (vi).

As $\lambda=\mbox{ cf}(\lambda )>\mu >2^{{\rm cf}(\mu)}$, we can replace
$\bar\lambda$ by a subsequence, so without loss of generality $J\subseteq
J^*$, so $J^*$ is ${\rm cf}(\mu )$-complete and $\bar\theta^\alpha$ is
$<_J$-increasing, see 7.0. So $\lan\bar\theta^\alpha :\alpha <\lambda\ran$ has
$<_{J^*}$-exact upper bound $\bar\theta^*$, without loss of generality
$$
\bigwedge_{\alpha ,\zeta}\theta^\alpha_\zeta <\theta^*_\zeta.
$$
So clause (v) holds.
\medskip

\n $\underline{\rm Note}$: If ${\rm cf}(\mu )>\aleph_0$ we
have $J=J^{\rm bd}_\mu$.

\n {\bf Proof of ($\beta$) + ($\gamma$)}. (Here $\mbox{cf}(\mu )$ can be
replaced by any $\delta\le\mu$ such that $\mbox{cf}(\delta )=\mbox{cf}(\mu
)$.)

List all relevant $\bar F=\lan F_\zeta :\zeta <\delta\ran$ with values subsets
of $\mu$. So there are $\le 2^\mu$ of them, list them as $\lan\bar
F^i:i<2^\mu\ran$ with
$$
\bar F^i=\lan F^i_\zeta :\zeta <\delta\ran.
$$
We choose $\eta_\alpha\in\prod_{\zeta <{\rm cf}(\mu )}\theta^*_\zeta$ by
induction on $\alpha$.

For a given $\alpha <\lambda$ we choose $\eta_\alpha\restriction A_\e$ by
induction on $\e<\e^*$.  We will choose $\eta_\alpha\restriction A_\e$ such
that
\begin{enumerate}
\item[($*$)] if $n<\omega$, $\beta_0$, $\beta_1,\ldots ,\beta_{n-1}<\alpha$
and $i<\sup\{\lambda_\beta :\beta <\alpha\}$ (necessarily $<\lambda_\alpha$),
\end{enumerate}
$$
\left\{\begin{array}{rl} 
\zeta\in A_\e :&F^i_\zeta (\ldots \beta_\ell,\eta_{\beta_\ell},\ldots
,\eta_\alpha\restriction\bigcup_{\xi <\e} A_\xi )\cap\theta^\alpha_\zeta\\ 
&\mbox{is a club of $\theta^\alpha_\zeta$ but $\eta_\alpha (\zeta )$ does not
belong to it}
\end{array}\right\}\in J. 
$$
But in 9.1's proof we have shown that this is possible.{\hfill$\square_{\rm
10.1}$}
\smallskip

\centerline{* * * * * * * * *}
\smallskip

We have conclusions variants similar to the case $2^\mu$ is regular.
\smallskip

\section{Preliminaries to the construction of ccc Boolean algebras with no
large independent sets} 

Monk \cite{M2} asks:
\medskip

\n {\bf Problem 33}. Assume $\mbox{cf}(\mu )\le \kappa <\mu <\lambda\le
\mu^{{\rm cf}(\mu )}$. Is it possible in ZFC that there is a Boolean algebra
of cardinality $\lambda$, satisfying the $\kappa$-cc with no independent
subset of cardinality $\lambda$? 
\medskip

\n This is closely related to the problem of ``is $\lambda$ a free caliber of
such Boolean algebra'' (see also in Monk \cite{M2}).
\smallskip

Why in ZFC? Because of earlier results under ``$\mu$ strong limit, $2^\mu
=\mu^+$'', I think.
\smallskip

The real problem seems to me is for $\lambda$ regular, and we shall prove that
``almost always'' there is such a Boolean algebra, so we prove the consistency
of failure.

We shall use $\lan J^{\rm bd}_{\lan\lambda_{i,0},\lambda_{i,1}\ran}:i<\delta
\ran$ with regular $\lambda_{i,0}>\lambda_{i,1}$, but we use Boolean algebras
whose existence is only consistent.

So we shall use $\bar\eta$ a $(\lambda,J)$-sequence for $\lan J^{\rm bd}_{\lan
\lambda_{i,0},\lambda_{i,1}\ran}:i<\delta\ran$, if $\delta =\omega$ the
Boolean algebra $\B$ will have a dense subalgebra $\B^*$ which will be the
free product of $\{ \B_n:n<\omega\}$, $x^-_t, x^+_t\in \B_n$ for $t\in {\rm
Dom}(I_n)$ and $\B =\lan \B^*, y_\alpha :\alpha <\lambda\ran$ where
$y_\alpha\in$ completion of $\B^*$ is defined from $\lan x^-_{\eta_\alpha
(n)}, x^+_{\eta_\alpha (n)}:n<\omega\ran$. We need special properties of
$\B_n$, $x^-_t$, $x^+_t$ $(t\in {\rm Dom}(I_n))$. The construction continues
\cite[\S3]{RoSh:534}. Concerning the parallel to 6.13 see later.
\smallskip

For the case $\mu$ strong limit we can use instead subalgebras of the measure
algebra. See \S2. Now we have consistency (and independence) for $\lambda$,
$\mu <\lambda\le 2^\mu$, $\mu$ strong limit singular, hence we concentrate on
the other case where the behavior is different i.e.\ when for some $\chi$ we
have ${\rm cf}(\mu )\le\kappa <\chi =\chi^{<\kappa}<\mu <\lambda <\mu^{{\rm
cf}(\mu )}\le 2^\chi$. The proof here uses ideals which are ``easier'' and can
be generalized to get ``non-$n$-independent subset of $\B$ of cardinality
$\lambda$ for some specific $n$". For this we need to start with ``there
is a
$\lambda_n$-complete uniform filter $\DD_n$ on $\lambda^{+n}_n$''. 
\medskip

\n {\bf 11.1 Definition}. {\it We say $(\B_1, \bar x^+, \bar
x^-)$ witness $(I,\TT )$ if
\begin{enumerate}
\item[(a)] $\TT$ is a set of Boolean terms written as $\tau =\tau (x_1,\ldots
,x_{n_\tau})$
\item[(b)] $I$ is an ideal
\item[(c)] $\B$ is a Boolean algebra
\item[(d)] $\bar x^+=\lan x^+_t:t\in{\rm Dom}(I)\ran$, $x^+_t\in\B$
\item[(e)] $\bar x^-=\lan x^-_t:t\in {\rm Dom}(I)\ran$, $x^-_t\in\B$
\item[(f)] $x^-_t<x^+_t$
\item[(g)] If $X\in I^+$ and $\B\subseteq \B'$ and
$$
\B'\vDash x^-_t\le y_t\le x^+_t,\;\;\mbox{ for $t\in X$}
$$
$\underline{\rm then}$ for some $\tau (x_1,\ldots ,x_n)\in\TT$ and pairwise
distinct $t_1,\ldots ,t_n\in X$ we have
$$
\B'\vDash \tau (y_{t_1}, y_{t_2},\ldots ,y_{t_n})=0.
$$
\end{enumerate}}
\smallskip

\n {\bf 11.2 Explanation}. We think of having $\bar\eta$ a $(\lambda
,J)$-sequence for $\lan I_i:i<\delta\ran$, and having $(\B_i,\bar x^+_i,\bar
x^-_i)$ witnessing $(I_i,\TT)$ for $i<\delta$ and using the sequence of
intervals $\lan(x^-_{i,\eta_\alpha (i)}, x^+_{i,\eta_{\alpha (i)}}):i<\delta
\ran$ as a sequence of approximations for an element $x_\alpha$ of the desired
Boolean algebra $\B$ of cardinality $\lambda$.

But we may think not only of ``$\{ x_\alpha:\alpha<\lambda\}$ has no
independent subset of cardinality $\lambda$'' but of other subsets of $\B$. So
sometimes we use
\smallskip

\n {\bf 11.3 Definition} {\it 1) We say that $(\B ,\bar x^-, \bar x^+)$
strongly witnesses $(I,\TT )$ if:
\begin{enumerate}
\item[(a)-(f)] as before, and
\item[(g)$^+$] If $\B\subseteq\B'$,
$$
\B'\vDash x^-_t\le y_t\le x^+_t\;\;\mbox{for $t\in {\rm Dom}(I)$,}
$$
$\lan b_\ell :\ell\le m\ran$ is a sequence of pairwise disjoint non-zero
members of $\B'$, $m<\omega$ and
$$
X\in \Big(\prod^m_{\ell =1}I\Big)^+,
$$
and $u\subseteq [1,m]$, $\underline{\rm then}$ we can find $n$, $\tau (x_1,
\ldots ,x_n)\in\TT$ and distinct $\bar t^1,\ldots ,\bar t^n\in X$, so $\bar
t^r=\lan t^r_\ell :\ell =1,\ldots ,m\ran$, such that $\tau (c_{\bar
t^1},\ldots ,c_{\bar t^n})=0$ where
$$
c_{\bar t}=b_0\cup\bigcup_{\ell\in [1,m]\atop\ell \in u}(b_\ell\cap
y_{t_\ell})\cup\bigcup_{\ell\in [1,m]\atop\ell\notin u} (b_\ell -y_{t_\ell})
$$
\end{enumerate}

2) We say that $(\B ,\bar x^+, \bar x^-)$ witness $(I,\TT )$ {\it
$m$-strongly\/} if we restrict ourselves to this $m$. Similarly
$[m_1,m_2]$-strongly.}
\medskip

Next we need our specific $(\B,\bar x^-,\bar x^+,I)$. The following is
essentially from \cite[p.244--246]{Sh:126}.
\smallskip

\n {\bf 11.4 Claim}.
\begin{enumerate}
\item[1)] If $\mu =2^\lambda =\lambda^+$, (or just $\mu\nrightarrow
[\mu]^2_\mu$) and $2^\mu =\mu^+$, $\underline{\rm then}$ we can find $\bar
F=\lan F_\alpha:\alpha <\mu^{+}\ran$ such that:
\end{enumerate}
$(*)^\mu_{\bar F}$
\begin{enumerate}
\item[(a)] $F_\alpha :[\mu]^2\rightarrow \alpha\times \mu$ is one to one.
\item[(b)] If $A\in (J^{\rm bd}_{\lan\mu^+,\mu\ran})^+$, then for some
$$\begin{array}{c}
(\alpha ,i_0),\; (\alpha ,i_1),\; (\beta , i_2)\in A\\
\mbox{we have }\;\; F_\alpha (\{ i_0, i_1\})=(\beta ,i_2).
\end{array}
$$
We write this also as
$$
F(\{\alpha , i_0\}, \{\alpha ,i_1\})=(\beta , i_2).
$$
We can add that for every $\beta $ we have $|{\rm Rang}(F_\alpha )\cap
(\{\beta\}\times\mu )|\le 1$ for $\alpha>\mu$. We do not strictly distinguish
$\bar F$ from $F$.
\end{enumerate}

2) The property $(*)^\mu_{\bar F}$ is preserved by forcing notions which have
the $(3,J^{\rm bd}_{\lan \mu^+,\mu\ran})^+$-c.c.\ (see 11.6 below).

3) Let $\B=\B_{\bar F}$ be the Boolean algebra freely generated by
$$
x^+_{\alpha ,i}=x^+_{(\alpha ,i)}; x^-_{\alpha ,i}=x^-_{(\alpha ,i)}(\mbox{for
}(\alpha ,i)\in\mu^+\times\mu )
$$
except that $x^-_{\alpha ,i}\le x^+_{\alpha ,i}$ and
$$
x^+_{(\alpha ,i)}\cap x^+_{(\alpha ,j)}\cap x^+_{F_\alpha (i,j)}=0.
$$
$\underline{\rm Then}$
\begin{enumerate}
\item[(i)] $(\B,\bar x^+, \bar x^-)$ witness $(J^{\rm bd}_{\lan\mu^+,\mu\ran},
\{ x_0\cap x_1\cap x_2=0\})$.
\item[(ii)] $\B$ satisfies the ccc.
\end{enumerate}
\smallskip

\n {\bf 11.5 Remark}. On more general Boolean algebras generated by such
equations see Hajnal, Juhasz, Szemintklossy \cite{HaJuSz}.
\smallskip

\n {\bf 11.6 Definition}. {\it For an ideal $J$ and forcing notion $P$, we say
that $P$ satisfies the $(n,J)$-c.c.\ if for $\lan p_t:t\in A\ran$, $A\in J^+$,
there is $B\subseteq A$, $B\in J^+$ such that any $n$ conditions in $\{
p_t:t\in B\}$ have a common upper bound.}
\medskip

\n {\bf 11.7 Fact}. If $P$ is the forcing notion $P_{\chi ,\theta}$ of adding
$\chi$ Cohens for $\theta$ and $\lambda^{<\theta}=\lambda$ then $P$ satisfies
$(n,J)$-c.c.\ for $n<\omega$, $J=J_{\lan \lambda^{++}, \lambda^+\ran}$.
\medskip

\n {\bf Proof of 11.4}. 1) Let $\{ A_\alpha :\alpha <\mu^+\}$ list all subsets
$A$ of $\mu^+\times \mu$ of cardinality $\mu$ such that for every $\beta
<\mu^+$ we have $|A\cap (\{\beta\}\times\mu )|\le 1$. For every $\alpha$ such
that $\mu<\alpha<\mu^+$ choose $H_\alpha:[\mu]^2\to\alpha$ such that $(\forall
X\in [\mu ]^\mu )$ $[H''_\alpha([X]^2)=\alpha]$.  For each $\alpha$, choose
$F_\alpha (i,j)\in \{\beta^\alpha_{\{i,j\}}\}\times\mu$ by induction on
$<^\otimes$, where $\{i,j\}<^\otimes\{i',j'\}$ iff $\max\{i,j\}<\max\{i',j'\}
\vee (\max\{i,j\}=\max\{i', j'\}$ \& $\min\{i,j\}<\min\{i',j'\})$, with
$\beta^\alpha_{i,j}$ with no repetition so that
$$
F_\alpha(i,j)\in\alpha\times\mu^+\backslash\cup\{\{\beta^\alpha_{i',j'}\}
\times\mu :\{i', j'\}<^\otimes\{i,j\}\},
$$
and if possible
$$
F_\alpha (i,j)\in A_{H_\alpha(\{i,j\}) },
$$
which occurs if $A_{H_\alpha (\{ i,j\})}\subseteq\alpha\times\mu$.
\smallskip

2) Trivial. Let $P$ be the forcing notion. Let $p^*\Vdash$ ``$\name{A}\in
(J^{\rm bd}_{\lan\mu^+,\mu\ran})^+$ and it exemplifies a contradiction to
$(*)^\mu_{\bar F}$''.  Let $A{\stackrel{\rm def}{=}}\{(\alpha,i):p^*\nVdash(
\alpha,i)\notin \name{A}\}$. So $A\subseteq\mu^+\times\mu$ and,
$$
p^*\Vdash ``A\supseteq\Name{A},\;\;\Name{A}\in (J^{\rm bd}_{\lan\mu^+,\mu
\ran})^+,\mbox{''}
$$
hence
$$
A\in (J^{\rm bd}_{\lan\mu^+,\mu\ran})^+.
$$
For $(\alpha ,i)\in A$ there is $p_{(\alpha ,i)}\ge p^*$ such that
$$
p_{(\alpha ,i)}\Vdash \mbox{``}(\alpha ,i)\in\Name{A}\mbox{''}.
$$
Apply (3, $J^{\rm bd}_{\lan \mu^+ ,\mu\ran})$-cc to $\lan p_{(\alpha
,i)}:(\alpha ,i)\in A\ran$, and obtain $B$ as in Definition 11.6. As $B\in
(J^{\rm bd}_{\lan \mu^+,\mu\ran})^+$, by $(*)^\mu_{\bar F}$ we can find
$(\alpha ,i_0)$, $(\alpha ,i_1)$, $(\beta ,i_2)\in B$ such that
$$
F_\alpha (\{i_0, i_1\})=(\beta ,i_2).
$$
But by the choice of $B$ there is $q\in P$ such that
$$
q\ge p_{(\alpha ,i_0)}, p_{(\alpha ,i_1)}, p_{(\beta ,i_2)}
$$
(hence $q\ge p^*$). So
$$
q\Vdash\mbox{``}(\alpha ,i_0), (\alpha,i_1), (\beta,i_2)\in\Name{A}\;\;
\mbox{and}\;\; F_\alpha (\{ i_0,i_1\})=(\beta,i_2)\mbox{''}. 
$$ 
But this contradicts the assumption on $p^*$, $\Name{A}$.

3) For clause (i), read the definition. For clause (ii):

Call $\ZZ \subseteq\mu^+\times \mu$ {\it closed} if $F(t_1, t_2)=t_3$ \& $|\{
t_1,t_2,t_2\}\cap \ZZ|>1\Rightarrow\{ t_1, t_2,t_3\}\subseteq\ZZ$. Now,
\begin{enumerate}
\item[($*$)] if $F(t_i, s_i)=r_i$ ~~ for $i=0,1$
\end{enumerate}
$\underline{\rm then}$
$$
\{ t_0, s_0, r_0\}\cap \{ t_1,s_1,r_1\}
$$
has $\le 1$ or 3 elements.
\smallskip

\n [Why? As each $F_\alpha$ is one to one and
$$
F=\bigcup_{\alpha <\mu^+} F_\alpha\restriction (\{\alpha\}\times [\mu ]^2)
$$
and
$$
\lan\{\alpha\}\times [\mu ]^2:\alpha <\mu^+\ran\; \mbox{ are pairwise
disjoint]}
$$
\begin{enumerate}
\item[($**$)] if $\ZZ\subseteq\mu^+\times\mu$, and $\B_\ZZ$ is defined
naturally: it is freely generated by $\{ x^+_t, x^-_t:t\in\ZZ\}$ except the
equations explicitly demanded on those variables, $\underline{\rm then}$
$\B_\ZZ\subseteq \B$ (even if $\ZZ$ is not closed).
\end{enumerate}
\smallskip

\n [Why? If $f:\{ x^-_t, x^+_t: t\in\ZZ\}\to\{ 0,1\}$ preserves the equations,
and we define
$$
f^*:\{ x^-_t, x^+_t: t\in\mu^+\times\mu\}\rightarrow\{ 0,1\}
$$
by
$$
f^*(y)\stackrel{\rm def}{=}\left\{
\begin{array}{ll}
f(y)&\quad\mbox{if $y=x^\pm_t$, $t\in\ZZ$}\\
0&\quad\mbox{if $y=x^\pm_t$, $t\not\in \ZZ$,}\end{array}\right.
$$
then $f^*$ preserves the equations.]
\begin{enumerate}
\item[($***$)] $\B\vDash {\rm ccc}$
\end{enumerate}
[Why? Let $\lan a_\zeta:\zeta<\omega_1\ran$ be a sequence of non-zero
elements. We can find finite $\ZZ_\zeta$ such that $a_\zeta\in\B_{
\ZZ_\zeta}$. Let $f_\zeta :\B_{\ZZ_\zeta}\rightarrow\{ 0,1\}$ be a
homomorphism such that $f_\zeta (a_\zeta )=1$. Let
$$
\ZZ_\zeta^+{\stackrel{\rm def}{=}}\ZZ_\zeta\cup\cup\{\{t_1,t_2,t_3\}:
F(t_1,t_2)=t_3, \; {\rm and }\; \{t_1,t_2,t_3\}\cap\ZZ_\zeta >1\}.
$$
Without loss of generality $\lan \ZZ^+_\zeta :\zeta <\omega_1\ran$ is a
$\Delta$-system with heart $\ZZ^+$.

Without loss of generality $f_\zeta\restriction \{ x^+_t:t\in\ZZ^+\}$ is
constant.

Without loss of generality $\ZZ_\zeta\cap\ZZ^+$ is constant.

So
\begin{enumerate}
\item[($*$)$_4$] If $\zeta\not= \xi <\omega_1$
$$
F(t_1,t_2)=t_3\;\mbox{ and }\;\{t_1,t_2,t_3\}\subseteq\ZZ_\xi\cup\ZZ_\zeta,
$$
$\underline{\rm then}$
$$
\{t_1,t_2,t_3\}\subseteq\ZZ_\zeta\;\mbox{ or }\;\{t_1,t_2,t_3\}\subseteq
\ZZ_\xi.
$$
\end{enumerate}
\smallskip
[Why? Without loss of generality
$$
|\{ t_1,t_2,t_3\}\cap\ZZ_\zeta|\ge 2
$$
So
$$
\{ t_1, t_2, t_3\}\subseteq\ZZ^+_\zeta.
$$
Now if $t_i\in \ZZ^+_\zeta\backslash \ZZ_\zeta$, then $t_i\not\in\ZZ_\xi$
(otherwise $t_i\in\ZZ^+_\zeta \cap\ZZ^+_\xi$, hence $t_i\in\ZZ^+$, but
$\ZZ_\zeta\cap\ZZ^+$ is constant). So $\{ t_1, t_2, t_3\}\subseteq\ZZ_\zeta$.]

Now $f_\zeta\cup f_\xi$ preserves the equations on $\ZZ_\zeta\cup\ZZ_\xi$ and
by the homomorphism it induces, $a_\zeta\cap a_\xi$ is mapped to $1$, so
$\B_{Z_\zeta\cup Z_\xi}\vDash {\rm ``}a_\zeta\cap a_\xi\not= 0$'' hence by
($**$) we have $\B\vDash {\rm ``}a_\zeta\cap a_\xi\not=
0$''.]{\hfill$\square_{\rm 11.4}$} \medskip

\n {\bf 11.7A Fact}. Assume
\begin{enumerate}
\item[(a)] $(\B ,\bar x^-, \bar x^+)$ is a witness for $(I,\TT)$.
\item[(b)] $y^-_t=-x^+_t$, $y^+_t=-x^-_t$ for $t\in {\rm Dom}(I_i)$,
\item[{}] $\bar y^-=\lan y^-_t:t\in {\rm Dom}(I)\ran$, $\bar y^+=\lan
y^+_t:t\in {\rm Dom}(I)\ran$
\item[(c)] $\TT'=\{ -\tau (-x_0, \ldots ,-x_{n-1}):\tau (x_0,\ldots
,x_{n-1})\in \TT\}$
\end{enumerate}
Then $(\B ,\bar y^-, \bar y^+)$ is a witness for $(I,\TT')$ (and is called the
dual of $(\B ,\bar x^-, \bar x^+)$).

We may consider
\medskip

\n {\bf 11.8 Definition}. {\it 1) Let $(*)^\mu_{\bar F, \bar H}$
mean
\begin{enumerate}
\item[(a)] $\bar F=\lan F_\alpha :\alpha <\mu^+\ran$, $F_\alpha$ is a partial
function from $[\mu ]^{2}$ into $\alpha\times\mu$
\item[(b)] $\bar H=\lan H_\alpha :\alpha <\mu^+\ran$, $H_\alpha$ is a partial
function from $[\mu]^{2}$ into $\{ 0,1\}$
\item[(c)] if $A\in (J^{\rm bd}_{\lan \mu^+,\mu\ran})^+$ and $\ell <2$ then
for some $(\alpha ,i_0)$, $(\alpha, i_1)\in A$ we have $F_\alpha (i_0, i_1)\in
A$ and $H_\alpha (i_0,i_1)=\ell$.
\item[(d)] the Boolean algebra $\B_{\bar F,\bar H}$ defined below satisfies
the c.c.c. We may write $F=:\bigcup_{\alpha <\mu^+}F_\alpha$,
$H=:\bigcup_{\alpha <\mu^+}H_\alpha$ instead of $\bar F$, $\bar H$
respectively.
\end{enumerate}

2) $\B_{\bar F, \bar H}$ is the Boolean algebra generated freely by $\{ x^-_t,
x^+_t:t\in\mu^+\times \mu\}$ except that $x^-_t\le x^+_t$ and $x^+_{t_0}\cap
x^+_{t_1}\cap x^+_{t_2}=0$ when $F(t_0, t_1)=t_2$, $H(t_0, t_1)=0$ and
$(-x^-_{t_0})\cap (-x^-_{t_1})\cap (-x^-_{t_2})=0$ when $F(t_0, t_1)=t_2$,
$H(t_0, t_1)=1$.}
\medskip

\n {\bf 11.9 Remark}. Of course $\B_{\bar F, \bar H}$ is defined from two sets
of triples, which are disjoint and no distinct two have $>1$ element in
common.
\medskip

\n {\bf 11.10 Claim}. Assume $(*)^\mu_{F_0}$ of 11.4(1) and e.g..\
$\mu=\lambda^+$, $\lambda^{<\theta}=\lambda$. 

1) For some $(\theta^{<\theta})^+$-c.c., $\theta$-complete, forcing notion $P$
of cardinality $\le\mu^+$ we have
$$
\Vdash_P \mbox{``}(*)^\mu_{F,H}\;\;\mbox{ for some $F,H$''}.
$$

2) If $(*)^\mu_{F,H}$ and $Q$ is a forcing notion satisfying the $(3, J^{\rm
bd}_{\lan \mu^+,\mu\ran})$-c.c.\ then in $V^Q$ we have $(*)^\mu_{F,H}$. If
$V=V^P_0$, $P$ as above it is enough that $P*\Name{Q}$ satisfies the $(3,
J^{\rm bd}_{\lan\mu^+,\mu\ran})$-c.c.
\medskip

\n {\bf Proof}. 1) Let
$$
P=\left\{ \begin{array}{rl}
(f,h):&\mbox{for some $u=u_{(f,h)}\subseteq \mu^+\times \mu$ of cardinality
$<\theta$ we have:}\\
&\mbox{$f,h$ are partial functions, $\mbox{Dom}(f)=\mbox{Dom}(h)\subseteq
(\mbox{Dom }F)\cap [u]^2$}\\
&\mbox{$f\subseteq F_0$ and Rang$(h)\subseteq\{0,1\}$ and $\B_{f,h}$ satisfies
the c.c.c.}\end{array}\right\} 
$$
where $\B_{f,h}$ is defined as in 11.8(2) (and see 11.9).

\n $\underline{\rm The ~order}$ $(f_1,h_1)\le (f_2,h_2)$ $\underline{\rm iff}$
\begin{enumerate}
\item[(i)] $u_{(f_1,h_1)}\subseteq u_{(f_2,h_2)}$,
\item[(ii)] $f_1=f_2\restriction [u_{(f_1,h_1)}]^2$
\item[(iii)] $h_1=h_2\restriction [u_{(f_1,h_1)}]^2$
\item[(iv)] $\B_{(f_1,h_1)}\subseteq \B_{(f_2,h_2)}$ moreover
\item[{}] $\B_{(f_1,h_1)}\lcirc\B_{(f_2,h_2)}$.
\end{enumerate}
The reader can check {\hfill$\square_{\rm 11.10}$}
\medskip

\n {\bf 11.11 Claim}. Assume $2^{\lambda^{+\ell}}=\lambda^{+\ell +1}$ for
$\ell <n$ and let $\lambda_\ell =\lambda^{n-\ell +1}$

1) We can find $W$ such that
\begin{enumerate}
\item[(a)] $W\subseteq\Big[\prod_{\ell <n}\lambda_\ell\Big]^n$
\item[(b)] if $u_1\not= u_2$ belongs to $W$ then $|u_1\cap u_2|\le 1$
\item[(c)] if $A\in (J^{\rm bd}_{\lan \lambda_\ell :\ell <n\ran})^+$ then
$[A]^n\cap W\not=\emptyset$
\item[(d)] $\lan\lambda_\ell :\ell <n\ran$ is a decreasing sequence of
regulars
\end{enumerate}

2) there is a forcing notion $Q$ of cardinality $\lambda^{+n}$,
$\lambda^+$-complete satisfying the $\lambda^+$-c.c.\ and even the $(n,J^{\rm
bd}_{\lan\lambda_\ell :\ell <m\ran})$-c.c. and adding $W$ satisfying (a), (b),
(c) of part 1 and
\begin{enumerate}
\item[(e)] $W$ is locally finite: if $A\subseteq \prod_{\ell <n}\lambda_\ell$
is finite, then for some finite $B$, $A\subseteq B\subseteq \prod_{\ell
<n}\lambda_\ell$ and $w\in W$ \& $|w\cap B|\ge 2\Rightarrow w\subseteq B$
\end{enumerate}

3) If $P$ is adding $\chi$ many $\theta$-Cohen reals, $\lambda=\lambda^\theta$
and in $V$, $\overline W$ satisfies (a), (b), (c), (d) and (e), then in $V^P$
still clause (c) holds (and trivially the other demands on $W$). (See
\cite{Sh:126}.) 
\medskip

\n {\bf Proof}. 1) We prove by induction on $n$ that for any such $\lambda$
satisfying $\ell <n\Rightarrow 2^{\lambda^{+\ell}}=\lambda^{+(\ell +1)}$ we
can find $(W,F)$ such that (a), (b), (c) of 11.11(1) hold for $W$,
$\lan\lambda^{+(\ell +1)}:\ell <n\ran$ and
\begin{enumerate}
\item[(f)] $F:W\to\lambda^+$ satisfies: if $A\in (J^{\rm bd}_{\lan\lambda^{+
(\ell+1)}:\ell<n\ran})^+$, then $\mbox{Rang}(F\restriction [A]^n)=\lambda^+$. 
\end{enumerate}
\smallskip

The induction step is as in the previous proof.

2) Similar to the proof of 11.10

3) Because $P$ satisfies the $(n,J^{\rm bd}_{\lan\lambda_\ell:\ell<n
\ran})$-c.c.{\hfill$\square_{\rm 11.11}$}
\medskip

\n {\bf 11.12 Claim}. Assume
\begin{enumerate}
\item[(A)] $W$, $\lan \lambda_\ell :\ell <n\ran$ satisfy (a), (b), (c), (d)
and (e) of Claim 11.11(1).
\item[(B)] $3\le m< n/2$, $n>6$.
\item[(C)] $\B$ is the Boolean algebra generated by $\{x^-_t,x^+_t:t\in
\prod_{\ell <n}\lambda_\ell\}$ freely except:
\item[{}] ($*$)$_1$ ~$x^-_t\le x^+_t$
\item[{}] ($*$)$_2$ ~if $w=\{t_0,\ldots,t_{n-1}\}\in W$, where $t_\ell$ is
increasing 
\item[{}] ~in the lexicographic order, and $u\subseteq n$, $|u|\ge m$ and
$n-|u|>m$, then
$$
\bigcap_{\ell\in u} x^+_{t_\ell}\cap\bigcap_{\ell <n,\ell\not\in u}
(-x^-_{t_\ell})=0.
$$
\item[(D)] $\TT=\TT_{n,m}=\{\bigcap_{\ell\in u} x_\ell\cap\bigcap_{\ell <n,
\ell\not\in u}(-x_\ell ):u\subseteq n, m\le |u|\le n-m\}$.

\n $\underline{\rm Then}$

\begin{enumerate}
\item[(i)] $\B\vDash {\rm ``}x^-_t<x^+_t$ \& $x^-_s\nleq x^+_t$'' for $t\not=
s$ in $\prod_{\ell <n}\lambda_\ell$
\item[(ii)] $(\B, \bar x^-, \bar x^+)$ is a witness for $(J^{\rm
bd}_{\prod_{\ell <n}\lambda_\ell},\TT )$
\item[(iii)] $\B$ satisfies the ccc.
\end{enumerate}
\end{enumerate}

\n {\bf Proof}. Clearly $\B\vDash x^-_t\leq x^+_t$ by the equation in
($*$)$_1$ and $\B\vDash \mbox{``}x^-_t\not= x^+_t$'' because the function
$f_0$ given by,
$$
f_0(x^-_s)=0,\qquad f_0(x^+_s)=\left\{\begin{array}{ll}
1&\quad s=t\\
0&\quad s\not= t\end{array}\right.
$$
preserves all the required equations (as $2\le m$). Taken together, $\B\vDash
x^-_t<x^+_t$. Also $\B\vDash x^-_t\nleq x^+_s$ when $t\not= s$ using $f_1$
defined by
$$
f_1(x^+_r)= f_1(x^-_r)=\left\{
\begin{array}{ll}
1&\quad\mbox{if $r=t$}\\
0&\quad\mbox{if $r\not= t$.}\end{array}\right.
$$
So clause (i) of the conclusion holds. Clause (ii) holds easily by the
equation in ($*$)$_2$ and assumption (A) i.e.\ (c) of 11.11(1).

We are left with verifying clause (iii) i.e.\ the c.c.c. So let
$a_\zeta\in\B\backslash\{ 0\}$ for $\zeta <\omega_1$. For every $\zeta$ we can
find a finite set $Z_\zeta\subseteq\prod_{\ell <n}\lambda_\ell$ such that
$a_\zeta\in \lan x^-_t, x^+_t:t\in Z_\zeta\ran$. By 11.11, i.e. by clause (A),
without loss of generality
\begin{enumerate}
\item[($*$)] if $w\in W$ \& $|w\cap Z_\zeta |\ge 2\Rightarrow w\subseteq
Z_\zeta$. 
\end{enumerate}
Let $f^*_\zeta :\{ x^-_t, x^+_t:t\in Z_\zeta\}\rightarrow \{ 0,1\}$ be such
that it preserves all the equations (from ($*$)$_1$ $+$ ($*$)$_2$) on these
variables and so the homomorphism it induces from $\B_{Z_\zeta}$ to $\{
0,1\}$, $\hat f^*_\zeta$ maps $a_\zeta$ to $1$. Without loss of generality
$\lan Z_\zeta:\zeta <\omega_1\ran$ is a $\Delta$-system with heart $Z$ and
$f^*_\zeta\restriction \{ x^-_t, x^+_t:t\in Z\}$ is constant.

Let $\zeta (1)<\zeta (2) <\omega_1$ and define $f_2$
\begin{eqnarray*}
f_2(x^-_t)&=&\left\{\begin{array}{ll}
f^*_{\zeta (1)}(x^-_t)&\quad\mbox{$\underline{\rm if}$ ~$t\in
Z_{\zeta (1)}$}\\
f^*_{\zeta (2)}(x^-_t)&\quad\mbox{$\underline{\rm if}$ ~$t\in
Z_{\zeta (2)}$}\\
0&\quad\mbox{$\underline{\rm if}$
~otherwise}.\end{array}\right.\\
&&\\
f_2(x^+_t)&=&\left\{\begin{array}{ll}
f^*_{\zeta (1)}(x^+_t)&\quad\mbox{$\underline{\rm if}$ ~$t\in
Z_{\zeta (1)}$}\\
f^*_{\zeta (2)}(x^+_t)&\quad\mbox{$\underline{\rm if}$ ~$t\in
Z_{\zeta (2)}$}\\
0&\quad\mbox{$\underline{\rm if}$ ~otherwise}.\end{array}\right.
\end{eqnarray*}
Clearly it is well defined and with the right domain. Does $f_2$ preserve all
the equations?
\medskip

\n {\bf Case 1}. $x^-_t\leq x^-_t$

if $t\not\in Z_{\zeta (1)}\cup Z_{\zeta (2)}$ trivial (both get value zero),

and if $t\in Z_{\zeta (\ell )}$ then trivial (as $f^*_{\zeta (\ell)}$
preserves this equation).
\medskip

\n {\bf Case 2}. $\bigcap_{\ell\in u} x^+_{t_\ell}\cap\bigcap_{\ell <n\atop
\ell\not\in u} (-x^-_{t_\ell})=0$.

If $\ell\in\{ 1,2\}$ and $\{ t_0,\ldots ,t_{n-1}\}\subseteq Z_{\zeta (1)}$
this holds as $f^*_{\zeta (\ell)}$ preserves this equation. So assume this
fails for $\ell =1,2$ so $|\{ t_0,\ldots ,t_{n-1}\}\cap Z_{\zeta (\ell )}|\le
1$ hence $2\ge |\{ t_0,\ldots ,t_{n-1}\}\cap (Z_{\zeta (1)}\cup Z_{\zeta
(2)})|$ so $\{ \ell :t_\ell\not\in Z_{\zeta (1)}\cup Z_{\zeta (2)}\}$
necessarily includes members of $u$, hence the equation
holds.{\hfill$\square_{\rm 11.12}$
\medskip

\n {\bf 11.13 Comment}. 1) ~If in addition we have $\kappa$-complete maximal
ideals $I_{n,\ell}$ on $\lambda_{n,\ell}$ extending $J^{\rm
bd}_{\lambda_{n,\ell}}$ and $\lan \lambda_{n,\ell}:\ell <n\ran$ as above for
$\bar\eta$ a $(\lambda ,J)$-sequence e.g.\ for $\lan I^*_n:n<\omega\ran$ where
$I^*_n=\prod J_{\lan \lambda_{n,\ell}:\ell <n\ran}$, we are in a powerful
situation as it can be applied to $n$-tuples rather than each one
separately. But above we prepare the proof for not using it by having strong
equations.

2) ~We can waive the ``locally finite'' demand proving as in the proof of
($***$) in the proof of 11.4.
\medskip

\section{Constructing ccc Boolean Algebras with no large independent sets} 

On such constructions see Ros{\l}anowski Shelah \cite[\S3]{RoSh:534}.
\smallskip

\n {\bf 12.1 Construction's Hypothesis}

We assume
\begin{enumerate}
\item[(a)] $\bar\eta$ is a normal $(\lambda ,J)$-sequence for $\lan
I_i:i<\delta\ran$
\item[(b)] $(\B_i, \bar x_i^-, \bar x^+_i)$ is a witness for $(I_i, \TT_i)$,
$\|\B_i\|=|\mbox{ Dom}(I_i)|$.
\item[(c)] $\lambda =\mbox{cf}(\lambda )$, $\sum\limits_{i<\delta}|
\mbox{Dom}(I_i)|<\lambda$.
\end{enumerate}
\medskip

\n {\bf 12.2 Remark}. Actually $\TT_i$ do not influence the construction, only
the properties of the Boolean algebra constructed. Similarly, the normality
and the fact that $\|\B_i\|=|\mbox{ Dom}(I_i)|$, as well as clause (c).

We define a Boolean algebra $\B$ and $y_\alpha\in\B$ $(\alpha <\lambda )$ as
follows:
\medskip

\n {\bf 12.3 The construction}.
\medskip

\n {\bf Case 1}. $\delta =\omega$.

Let $\B_*$ be the free product of $\{\B_i:i<\delta\}$ (so $\B_n =*_{i<n}\B_i$,
$\B^n\subseteq \B^{n+1}\subseteq\B_*$, so $\B_*=\lan\bigcup_{n<\omega}\B_n
\ran_{\B_*})$.

Let $\B^c_*$ be the completion of $\B_*$.

For each $i<\delta$ and $\eta\in\{\eta_\alpha\restriction i:\alpha
<\lambda\}\subseteq\prod_{j<i}\mbox{Dom}(I_j)$ we define $y_\eta^-<y^+_\eta$
in $\B^i$. This is done by induction on $i$.
\smallskip

\n $\underline{i=0}$. $y^-_\eta =0$, $y^+_\eta =1$.
\smallskip

\n $\underline{i=j+1}$. $y^-_\eta =y^-_{\eta \restriction j}\cup
(y^+_{\eta\restriction j}\cap x^-_{i,\eta (j)})$

$y^+_\eta =y^-_{\eta\restriction j}\cup (y^+_{\eta\restriction j}\cap
x^+_{i,\eta (j)})$.
\smallskip

\n So easily
$$
j<i\Rightarrow y^-_{\eta_\alpha\restriction j}\le y^-_{\eta_\alpha\restriction
i}<y^+_{\eta_\alpha\restriction i}\le y^+_{\eta_\alpha\restriction j}.
$$
Now let $y_\alpha$ be $\mbox{lub}\{ y^-_{\eta_\alpha\restriction
i}:i<\delta\}$. (Note: If $\B_i\vDash {\rm ``}0<x^-_{i,t}<x^+_{i,t}<1$'' for
$t\in\mbox{Dom }(I_i)$, then also $y_\alpha =$ maximal lower bound of $\{
y^+_{\eta_\alpha\restriction i}:i<\delta\}$. This will not be used.)
\smallskip

\n [Otherwise, the difference contains some member of $\B_*$, hence of some
$\B^i$ $(i<\delta )$, but there is none.]

Lastly $\B =\B_{\bar\eta, \bar I, \lan (\B_i, \bar x^-_i, \bar
x^+_i):i<\delta\ran}$ is the subalgebra of $\B^c_*$ generated by $\B_*\cup\{
y_\alpha :\alpha <\lambda\}$ (by the finitary operations, so it is not
complete).
\smallskip

\n {\bf Case 2}. $\delta >\omega$.

We find by induction on $i<\delta$, $\B^i$, $\{ (y_\eta^-, y^+_\eta
):\eta\in\{\eta_\alpha\restriction i:\alpha <\lambda\}\}$ such that
\begin{enumerate}
\item[(i)] $\B^i$ increasing (by $\subseteq$, even $\lcirc$)
\item[(ii)] $\B^i\vDash y^-_\eta <y^+_\eta$ (when $\bigvee_\alpha\eta
=\eta_\alpha\restriction i$)
$$
j<i\Rightarrow\B^i\vDash y^-_{\eta\restriction j}\le y^-_\eta\le y^+_\eta\le
y^+_{\eta\restriction j}.
$$
\item[(iii)] $\B^0$ is the trivial Boolean algebra.
\item[(iv)] if $i=j+1$ then $\B^i=\B^j*\B_j$ (free product) and for
$y\in\{\eta_\alpha\restriction i:\alpha <\lambda\}$
\begin{eqnarray*}
y^-_\eta &=&y^-_{\eta\restriction j}\cup
(y^+_{\eta\restriction j}\cap x^-_{j,\eta (j)})\\
y^+_\eta &=& y^-_{\eta\restriction j}\cup
(y^+_{\eta\restriction j}\cap x^+_{j,\eta (j)})
\end{eqnarray*}
\item[(v)] For $i$ limit, $\B^i$ is generated freely by
$$
\bigcup_{j<i}\B^j\cup \{ y^-_\eta, y^+_\eta :\eta\in\{\eta_\alpha\restriction
i:\alpha <\lambda\}
$$
except: the equations in $\B$ and
$$
y^-_{\eta\restriction j}\le y^-_\eta\le y^+_\eta\le y^+_{\eta\restriction
j}\;\mbox{ for }\; j<i, \eta\;\; \mbox{as above}.
$$
Lastly, $\B\subseteq$ completion $\Big(\bigcup_{i<\delta}\B^i\Big)$ is defined
as in case 1 using $y_\alpha\deq y^-_{\eta_\alpha}$.
\end{enumerate}
\smallskip

\n {\bf 12.3A The construction}: A variant
$$
\bar x^\pm_i=\lan x^\pm_{i,\eta}:\eta\in\{\eta_\alpha\restriction (i+1):\alpha
<\lambda\}\ran
$$
so we use $x_{i,\eta_\alpha\restriction (i+1)}$ instead of $x_{i,\eta_\alpha
(i)}$.
\medskip

\n {\bf 12.3B The construction}: A variant. It is like 12.3A but we are given
$(\B_i^\pm, \lan x_{i,\alpha}:\alpha <\lambda_i\ran )$ and we define by
induction on $i$, $\B^i$ (increasing with $i$), and follows:
\smallskip

\n {\bf Case 1}: $i=0$: $\B^i$ is the trivial Boolean algebra,
$y^-_{\eta_\alpha\restriction i}=0$, $y^+_{\eta_\alpha\restriction i}=1$.
\smallskip

\n {\bf Case 2}: $i=j+1$: $\B^i$ is generated by $\B^j\cup\{x^-_{\eta_\alpha
\restriction i}, x^+_{\eta_\alpha\restriction i}:\alpha<\lambda\}$ freely
except the equations in $\B^j$ and 
$$
\tau (\ldots x^-_{\eta_{\alpha_\ell}\restriction i},
x^+_{\eta_{\alpha_\ell}\restriction i},\ldots )_{\ell <n}=0
$$
whenever $\B_i\vDash\tau (\ldots ,x^-_{\eta_{\alpha_\ell}(j)},
x^+_{\eta_{\alpha_\ell}(j)},\ldots )_{\ell <n}=0 $; lastly defines
\begin{eqnarray*}
y^-_{\eta_\alpha\restriction i}&=& y^-_{\eta_\alpha\restriction j}\cup
\Big( y^+_{\eta_\alpha\restriction j}\cap x^-_{j,\eta_\alpha (j)}\Big)\\
y^+_{\eta_\alpha\restriction i}&=& y^-_{\eta_\alpha\restriction j}\cup
\Big( y^+_{\eta_\alpha\restriction j}\cap x^-_{j,\eta_\alpha (j)}\Big)
\end{eqnarray*}
\smallskip

\n {\bf Case 3}: $i$ limit

$\B^i$ is generated by $\bigcup_{j<i}\B^j\cup\big\{y^-_{\eta_\alpha
\restriction i}, y^+_{\eta_\alpha\restriction i}:\alpha <\lambda\}$ freely
except the equations in $\B^j$ for $j<i$ and $y^-_{\eta_\alpha\restriction
j}\le y^-_{\eta_\alpha\restriction i}\le y^+_{\eta_\alpha\restriction i}\le
y^+_{\eta_\alpha\restriction j}$ for $\alpha <\lambda$.
\medskip

\n {\bf 12.4 Comment}: Clearly 12.3A includes 12.3B as a special case, but
mostly there is no real difference in the uses. The reader may concentration
on 12.3B.
\medskip

\n {\bf 12.5 Discussion}. Usually the conclusions are of the form: among any
$\lambda$ elements of $\B$, something occurs. The first need is
$\|\B\|=\lambda$, a trivial thing.
\medskip

\n {\bf 12.6 Fact}. ($*$)$_3$ $\Rightarrow$ ($*$)$_2$ $\Rightarrow$
($*$)$_1$,\quad where
\begin{enumerate}
\item[($*$)$_1$] $\|\B\|=\lambda$
\item[($*$)$_2$] for every $\alpha <\beta <\lambda$
$$
\{i:\B_i\vDash\neg (\exists y)(x^-_{\eta_\alpha (i)}\le y\le x^+_{\eta_\alpha
(i)}\wedge x^-_{\eta_\beta (i)}\le y\le x^+_{\eta_\beta (i)})\}\not= \emptyset
$$
i.e.
$$
\{i:\B_i\vDash x^-_{\eta_\alpha (i)}\le x^+_{\eta_\beta (i)}\vee x^-_{
\eta_\beta (i)}\nleq x^+_{\eta_\alpha (i)}\}\not=\emptyset
$$
\item[($*$)$_3$] if $t\not= s$ are in $\mbox{Dom}(I_i)$ for some $i<\delta$,
$\underline{\rm then}$
$$
\B_i\vDash x^-_t\nleq x^+_s\vee x^-_s\nleq x^+_t.
$$
\end{enumerate}
\smallskip

\n {\bf Proof}. Easy.{\hfill$\square_{\rm 12.5}$}
\medskip

\n {\bf 12.7 Remark}. If not said otherwise, all examples satisfy ($*$)$_3$.

We will also be interested in stronger properties. In section 15 we will be
interested in the case $(\B ,\bar x^-,\bar x^+)$ the pairs $(x^-_\eta
,x^+_\eta )$, $(x^-_\nu ,x^+_\nu)$ were independent.
\medskip

\n {\bf 12.8 Claim}. Assume
\begin{enumerate}
\item[($*$)] $a_\alpha\in\B$ for $\alpha <\lambda$.
\end{enumerate}
$\underline{\rm Then}$ we can find in $\B$ a sequence $\lan b_\ell :\ell\le
m\ran$ a $\B$-partition of $1$ (i.e.\ a sequence of disjoint non-zero elements
with union 1), $m\ge 0$, and $X\in [\lambda ]^\lambda$ and $c\le b_0$ in $\B$
and $n$, and Boolean terms $\tau_\ell$ for $\ell =1,\ldots ,m$ with $n$
variables and ordinals $\gamma_{\alpha ,k}\in X$ for $\alpha\in X$, $k<n$ and
$\gamma_k$ for $k\in [n,n^*)$, where $n^*\ge n$ and $i^*<\delta$, $\nu_k$ for
$k<n^*$ such that 
\begin{enumerate}
\item[(i)] $n=0$ iff $m=0$ iff $\lan a_{\alpha}:\alpha\in X\ran$ constant
\item[(ii)] $\gamma_{\alpha,0}<\gamma_{\alpha,1}<\cdots<\gamma_{\alpha,n-1}$
and $\gamma_n<\gamma_{n+1}<\cdots <\gamma_{n^*-1}<\gamma_{\alpha,0}$
\item[(iii)] if $\alpha <\beta$ are in $X$ then $\gamma_{\alpha ,n-1}<
\gamma_{\beta ,0}$.
\item[(iv)] if $\alpha\in X$ then $a_\alpha\le\bigcup_{\ell\le m}b_\ell$,
$a_\alpha\cap b_0=c$ and $[\ell\in [1,m]\ \Rightarrow\ a_\alpha\cap b_\ell=
\tau_\ell(y_{\gamma_{\alpha ,0}}, y_{\gamma_{\alpha,1}},\ldots,y_{\gamma_{
\alpha ,n-1}})]$, and $[\ell\in [1,m]\ \Rightarrow\ 0< a_\alpha \cap
b_\ell <b_\ell]$ (so $\tau_\ell$ non-trivial).
\item[(v)] $\eta_{\gamma_{\alpha ,k}}\restriction i^*=\nu_k$ for $k<n$
\item[(vi)] $\{ b_\ell :\ell\le m\}\subseteq\lan\B_i\cup\{ y_{\gamma_k} :k\in
[n,n^*)\}\ran$ and $\eta_{\gamma_k}\restriction i^*=\nu_k$ for $k\in [n,n^*)$
\item[(vii)] $\lan\nu_k:k<n^*\ran$ is with no repetition.
\end{enumerate}
\smallskip

\n {\bf Proof}: By the $\Delta$-system lemma and Boolean algebra
manipulation.{\hfill$\square_{\rm 12.8}$} 
\medskip

\n {\bf 12.9 Claim}. A sufficient condition to
\begin{enumerate}
\item[$\otimes_0$] $\B$ has no independent subset of cardinality $\lambda$
\end{enumerate}
is
\begin{enumerate}
\item[$\otimes_1$] $\underline{\rm if}$ $a_\alpha$, $X$, $n$, $m$,
$\tau_\ell$, $\gamma_{\alpha ,k}$, $b_\ell$($\alpha\in X$, $k<n$, $\ell <n$)
are as above in 12.8, and $c_0=0$, $m=1$, $\underline{\rm then}$ $\{a_\alpha
:\alpha\in X\}$ is not independent,
\end{enumerate}
which follows from:
\begin{enumerate}
\item[$\otimes_2$] $\underline{\rm if}$ $a_\ell$, $X$, $n$, $m$, $\tau_\ell$, 
$\gamma_{\alpha ,k}$ $(\alpha\in X, k<n)$ are as above in 12.8, $c_0=0$,
$m=1$, then 
$$\left\{i:\,\begin{array}{rl}
&\mbox{for every $A,\B'$ and $y_t$, if $A\in
((I_i)^n)^+,\B_i\subseteq\B'$,}\\
&\mbox{$\B' \models x^-_t\le y_t\le x^+_t$ for $t\in\mbox{Dom}(I_i)$},\\
&\mbox{$\underline{\rm then}$ $\lan\tau_1(y_{t_0},\ldots,y_{t_{n-1}}):\lan
t_0,\ldots ,t_{n-1}\ran\in A\ran$}\\ 
&\mbox{is not strongly independent}\end{array}\right\}\in J^+.
$$
\end{enumerate}

\smallskip

\n {\bf 12.10 Remark}. If we ask more on $\bar\eta$, we can weaken
$\otimes_2$, like:
\smallskip

\n $\underline{\rm if}$ $n<\omega$, $\lan\gamma_{\alpha,k}:k<n\ran$ increasing
$\alpha<\beta\Rightarrow\gamma_{\alpha,n-1}<\gamma_{\alpha,0}$, then letting
$\eta'_\alpha =\lan\lan \eta_{\gamma_{\alpha,k}}(i):k<n\ran:i<\delta\ran$,
gives that $\bar\eta'=\lan\eta'_\alpha :\alpha <\lambda \ran$ is a $(\lambda
,J)$-sequence for $\lan (I_i)^n:i<\delta\ran$ as well as some weaker versions.
\medskip

\n {\bf Proof of 12.9}.  $\underline{\otimes_1\Rightarrow\otimes_0}$.

We choose by induction on $\ell\le m$ a sequence $\lan(\tau_\ell,\gamma^\ell_{
\alpha ,0},\ldots ,\gamma^\ell_{\alpha ,m(\ell )-1}):\alpha<\lambda\ran$ such
that 
\begin{enumerate}
\item[(i)] $\tau_\ell =\tau_\ell (x_ i,\ldots ,x_{m(\ell)-1})$ is a Boolean
term
\item[(ii)] $\gamma^\ell_{\alpha,0}<\gamma^\ell_{\alpha,1}<\cdots<
\gamma^\ell_{\alpha,m(\ell)-1}<\lambda$
\item[(iii)] $\alpha <\beta <\lambda\Rightarrow\gamma^\ell_{\alpha,m(\ell)-1}
<\gamma^\ell_{\beta ,0}$ when they are well defined. 
\item[(iv)] $\tau_\ell (a_{\gamma^\ell_{\alpha,0}},\ldots,\gamma^\ell_{\alpha 
,m(\ell)-1})\cap \bigcup_{\ell_1\le \ell}b_{\ell_1}=0$.
\end{enumerate}
$\underline{{\rm For} \;\ell =0}$ Let $\tau_\ell (x_0,x_1)=x_0-x_1$, so
$m(\ell)=2$
$$
\gamma^0_{\alpha ,0}=2\alpha ,\qquad \gamma^\ell_{\alpha ,1}=2\alpha +1.
$$
\smallskip

\n $\underline{{\rm For} \; \ell+1}$. For each $\alpha(*)<\lambda$, apply
$\otimes_1$ with $1-b_{\ell+1}, b_{\ell+1}$, $\lan a^\ell_{\alpha(*)+\alpha}:
\alpha<\lambda\rangle$, where $a^\ell_\alpha\deq\tau^\ell_\alpha
(a_{\gamma^\ell_{\alpha,0}},\ldots,a_{\gamma^\ell_{\alpha ,m(\ell)-1}})$
here standing for $b_0, b_1, \lan a_\alpha:\alpha<\lambda\ran$ there, and get
a Boolean term $\tau^{\ell +1}_{\alpha (*)}(x_0,\ldots,x_{m(\ell +1,\alpha
(*))-1})$, and ordinals $\beta^\ell_{\alpha^*(*),0}<\ldots<\beta^\ell_{\alpha
(*), m(\ell+1,\alpha(*))-1}$, all in the interval $[\alpha(*),\lambda)$, such
that 
$$
\tau^{\ell +1}_{\alpha(*)}(a^\ell_{\beta^\ell_{\alpha(*),0}},a^\ell_{
\beta^\ell_{\alpha(*),1}},\ldots,a^\ell_{\beta^\ell_{\alpha(*),m(\ell+1,
\alpha(*)-1)}})=0.
$$
Let $X\in [\lambda]^\lambda$ be such that 
\begin{enumerate}
\item[(a)] $\alpha\in X\Rightarrow\tau_\alpha^{\ell+1}=\tau^*_\ell$, $m(\ell,
\alpha )=m(\ell ,*)$
\item[(b)] $X$ is thin enough, i.e.~if $\alpha<\beta$ are in $X$ then
$\beta^\ell_{\alpha,0},\ldots,\beta^\ell_{\alpha,m(\ell,*)}<\beta$. 
\end{enumerate}
Now if $\e$ is the $\zeta$-th element of $X$ we let
$$
u^{\ell +1}_\zeta =\{ \gamma^\ell_{\beta ,m}:m<m(\ell )\; {\rm and}\;
\beta\in \{\beta^\ell_{\e ,0}, \ldots ,\beta^0_{\e ,m(\ell ,*)-1}\}\}.
$$
So $|u_\zeta^{\ell +1}|=m(\ell )\times m(\ell ,*)$ let $m(\ell +1)=m(\ell)
\times m(\ell ,*)$ let $\gamma^{\ell +1}_{\zeta ,0}<\gamma^{\ell+1}_{\zeta,1}
<\cdots <\zeta^{\ell +1}_{\zeta ,m(\ell +1)-1}$ list $u^{\ell +1}_\zeta$,
and it should be clear what is $\tau_{\ell +1}$. For $\ell=m$ we have finished.
\medskip

$\underline{\otimes_2\Rightarrow\otimes_1}$. Straight.{\hfill$\square_{\rm
12.9}$}
\medskip

\n {\bf 12.11 Fact}. 1) ~In 12.8 we can add (so in $\otimes_2$ of 12.9 we can
assume) that
\begin{enumerate}
\item[(viii)] $\tau_1 ( x_0,\ldots ,x_{n-1})\in\{ x_{k}, -x_k:k<n\}$
\end{enumerate}
$\underline{\rm if}$
\begin{enumerate}
\item[($*$)] for a set of $i<\delta$ from $J^+$ we have $\lan x^+_{i,t}-
x^-_{i,t}:t\in\mbox{ Dom }(I_i)\ran$ is a sequence of pairwise disjoint
(nonzero) elements of $\B_i$.
\end{enumerate}

2) ~Assume
\begin{enumerate}
\item[($*$)$^+$] for every $i<\delta$ we have $\lan x^+_{i,t}-x^-_{i,t}:t\in
{\rm Dom}(I_i)\ran$ is a sequence of pairwise disjoint (non zero) elements of
$\B_i$.
\end{enumerate}
\smallskip

\n $\underline{\rm Then}$
\begin{enumerate}
\item[(a)] in 12.8 above we can add:
$$
b_0,\ldots ,b_m=\bigcup_{i<\delta} B^i.
$$
\item[(b)] Under 12.3B we can add: for $k\in [1,m)$, if $i$ is large enough,
if $\alpha_0,\ldots ,\alpha_{n-1}\in X$ letting $b^\ell_k$ be the projection
of $b_\alpha$ in $\B^{i+1}$ (i.e.\ any element $b$ satisfying
$$
(\forall x\in\B^{i+1})(x\le b_k\to b\le b_k, x\ge b_k\to b\ge b_k)
$$
(there is a minimal and maximal such $b^i_k$ and they are in $\lan\B^i\cup
\big\{\rho: \rho=f_\alpha\restriction (x+1)$ for some $i,\neg
(\nu\triangleleft\rho )\big\}\ran$), $f_{\alpha_\ell}\restriction
i=f_{\alpha_0}\restriction i$, $\lan f_{\alpha_\ell}(i):\ell\le s\ran$ is with
no repetitions and $\tau (x_0,\ldots ,x_{s-1})$ is a Boolean term then
$$\begin{array}{l}
\B\restriction b_k\vDash \tau (b_k\cap y_{\alpha_0},\ldots ,b_k\cap 
y_{\alpha_{s-1}})=0\quad\Rightarrow\\
\B^{i+1}\vDash \tau (b^i_k\cap y_{\alpha_0},\ldots,b_k^i\cap
y_{\alpha_{s-1}})=0
\end{array}
$$
\end{enumerate}
(we can even be more explicit).
\medskip

\n {\bf Proof}. Straightforward.{\hfill$\square_{\rm 12.11}$}
\medskip

We can now phrase sufficient conditions for having free caliber $\lambda$ (for
$\TT$) and for having no $\TT$-free subset of $\B$ of cardinality $\lambda$.
\medskip

\n {\bf 12.12 Claim}. 1) Sufficient conditions for ``$\B$ satisfies the
$\kappa$-c.c.'' are ($\kappa$ is regular uncountable and):
\begin{enumerate}
\item[($*$)$_1$] $\delta =\omega$ and each $\B_i$ satisfies the
$\kappa$-Knaster condition.
\item[($*$)$_2$] each $\B_i$ satisfies the $\kappa$-Knaster condition and
$(\forall\alpha <\kappa )$ $(|\alpha |^{|\delta|}<\kappa)$
\item[($*$)$_3$] each $\B_i$ satisfies the $\kappa$-Knaster condition, $\kappa
>\delta$ and for every $A\in [\lambda ]^\kappa$, and limit ordinal
$\delta'\le\delta$ for some $B\in [A]^\kappa$ and $i<\delta$ we have 
$$
\alpha\in B,\; \beta\in B,\;\eta_\alpha\restriction\delta'\not=\eta_\beta
\restriction\delta'\Rightarrow\ell g(\eta_\alpha\cap\eta_\beta )=i
$$
\end{enumerate}
(follows from ``$\bar\eta$ is $\kappa^+$-free'', see 1.17 and Definition 1.18).
\smallskip

\n {\bf 12.13 Claim}. Assume
\begin{enumerate}
\item[(A)] $\bar\eta$ is a normal $(\lambda ,J)$-sequence for $\lan
I_i:i<\delta\ran$
\item[(B)] $(\B_i, \bar x^-_i, \bar x^+_i)$ is a witness for $(I_i, \{ x_0\cap
x_1\cap x_2=0\})$
\item[(C)] $\B$ is as constructed in 12.1, 12.3.

\n $\underline{\rm Then}$
\begin{enumerate}
\item[($\alpha$)] $\lambda$ is not a free caliber of $\B$.
\item[($\beta$)] $\B$ has cardinality $\lambda$ and satisfies the
$\kappa$-c.c. if $\kappa$ is as in 12.12.
\end{enumerate}
\end{enumerate}
\smallskip

\n {\bf Proof}. Straightforward.{\hfill$\square_{\rm 12.13}$}
\medskip

\n {\bf 12.14 Conclusion}. Assume for simplicity that $V\vDash \mbox{GCH}$,
$\theta =\theta^{<\theta}<\chi =\chi^{<\chi}$ and $P$ is the forcing notion of
adding $\chi$ $\theta$-Cohen subsets of $\theta$, i.e.\
$$
\mbox{$P=\{f:f$ is a partial function from $\chi$ to $\{0,1\}$ with domain of
cardinality $<\theta\}$}. 
$$
$\underline{\rm Then}$ (cardinal arithmetic on $V^P$ is well known) and
\begin{enumerate}
\item[($*$)] if $\mbox{cf}(\mu )<\theta <\mu <\chi$ $\underline{\rm then}$
there is a $(2^{{\rm cf}(\mu )})^+$-c.c.\ Boolean algebra $\B$ of cardinality
$\lambda =\mu^+$ such that $\lambda$ is not a free caliber of $\B$ (and even
satisfying the $\kappa$-c.c.\ if $\kappa$ is as in 12.12).
\end{enumerate}
\smallskip

\n {\bf Proof}. Use 12.13 and \S11.{\hfill$\square_{\rm 12.14}$}
\medskip

The problem of ``$\B$ with no independent subset of cardinality $\lambda$'' is
somewhat harder.
\medskip

\n {\bf 12.15 Claim}. Assume
\begin{enumerate}
\item[(A)] $\bar\eta$ is a normal $(\lambda ,J)$-sequence for $\lan
I_i:i<\delta\ran$.
\item[(B)] $(\B_i, \bar x^-_i, \bar x^+_i)$ is a witness for $(I_i,
\TT_{n_i,m_i})$ (on $\TT_{n_i,m_i}$ see 11.12 clause (D)).
\item[(C)] $3\le m_i<n_i/2$
\item[(D)] for every $k<\omega$, $\{ i:km_i<n_i\}\in J^+$
\item[(E)] $\B$ is as in construction 12.1, 12.3. 
\end{enumerate}
$\underline{\rm Then}$
\begin{enumerate}
\item[(i)] $\B$ does not have a free subset of cardinality $\lambda$.
\item[(ii)] $\B$ has cardinality $\lambda$ and satisfies the $\kappa$-c.c.\ if
$\kappa$ is as in 12.12.
\end{enumerate}
\smallskip

\n {\bf Proof}. Straightforward (using the criterion in 12.9).
{\hfill$\square_{\rm 12.15}$}
\medskip

\n {\bf 12.16 Conclusion}. Assume for simplicity $V\vDash \mbox{GCH}$, and
$\theta =\theta^{<\theta}<\chi =\chi^{<\chi}$ is the forcing notion of adding
$\chi$ $\theta$-Cohen reals. $\underline{\rm Then}$ cardinal arithmetic in
$V^P$ is well known and
\begin{enumerate}
\item[($*$)] if $\mbox{cf }(\mu) <\theta <\mu <\chi$ $\underline{\rm then}$
there is a $(2^{{\rm cf}(\mu )})^+$-c.c.\ Boolean algebra $\B$ of cardinality
$\lambda =\mu^+$ without an independent subset of cardinality $\lambda$.
\item[($**$)] we can demand that $\B$ satisfies the $({\rm cf}(\mu )^+)$-c.c.\ 
if 
$$
V\vDash\mbox{``}\{\delta <\mu^+:{\rm cf}(\delta )={\rm cf}(\mu )\}\in
I[\lambda]\mbox{''}
$$
\end{enumerate}
\smallskip

\n {\bf Proof}. By 12.15, where $(\B_i, \bar x^-_i, \bar x^+_i)$ is provided
by 11.12 (and $W$ for it by 11.11). For ($**$) see 1.17(2).
{\hfill$\square_{\rm 12.16}$}
\medskip

\centerline{* * *}
\medskip

We would also like sufficient condition for inequalities, for simplicity
$n=2$.
\medskip

\n {\bf 12.17 Claim}: 1) ~Assume 12.1, 12.3 and ($*$) of 12.11 and $n<\omega$
and $\tau^0=\tau^0(x_0,\ldots,x_{n-1})$ a Boolean term and $\tau^1=\tau^1(
-x_0,\ldots,-x_{n-1})$. $\underline{\rm Then}$ ($*$)$_1$ $\Rightarrow$
($*$)$_2$, where 
\begin{enumerate}
\item[$(*)_1$] if $\ell <2$ for a set of $i<\delta$ from $J^+$ we have: if
$X\in I^+_i$ $\underline{\rm then}$ for some $t_0,\ldots ,t_{n-1}\in X$,
pairwise distinct, we have
$$
\B_i\vDash \tau^\ell (x_{i,t_0},\ldots ,x_{i,t_{n-1}})=0,
$$
\item[$(*)_2$] if $a_\alpha\in\B$ for $\alpha<\lambda$ then for some
$k<\omega$ and $\alpha_{\ell,m}<\lambda$ for $\ell<k$, $m<n$ we have
$\alpha_{\ell,0}<\alpha_{\ell,1}<\ldots<\alpha_{\ell,m-1}<\alpha_{\ell+1,0}$
(for $\ell<k$) and for some $i(\ell)\in\{0,1\}$ for $\ell<k$ we have
$$
\B_i\vDash \bigcap_{\ell<k}\tau^{i(\ell)} (a_{\ell,0},\ldots,a_{\ell,m-1})=0.
$$
\end{enumerate}
2) Assume 12.1, 12.3 (using 12.3B) and ($*$) of 12.8 and for simplicity
$I_i=J^{\rm bd}_{\lambda_i}$ $\underline{\rm and}$ assume further $n<\omega$,
$\bf t$ a function from $\{ 0,\ldots ,n-1\}$ to $\{+1,-1\}$ and
$\tau^0=\tau^0(x_0,\ldots ,x_{n-1})$ a Boolean term, increasing in $x_\ell$ if
${\bf t}(\ell )=+1$, decreasing with $x_\ell$ if ${\bf t}(\ell )=-1$. Let
$\tau^1(x_0,\ldots ,x_{n-1})=\tau^0(-x_0,\ldots ,-x_{n-1})$. Assume also
$\tau_0(-x_0,\ldots,-x_{n-1})=0$ if $x_\ell\in \{ 0,1\}$ and $\bigwedge_\ell
(x_\ell =1\equiv {\bf t}(\ell)=1 )$ or $\bigwedge_\ell x_\ell =1\equiv {\bf
t}(\ell )=-1$. Then ($*$)$_3$ $\Rightarrow$ ($*$)$_4$, where
\begin{enumerate}
\item[($*$)$_3$] for a set of $i<\delta$ which belongs to $J^+$ the following
holds:\quad if $\gamma_{\alpha ,\ell}<\lambda_i$ and
$$
[\alpha<\beta <\lambda_i\Rightarrow \max_{\ell <n}\gamma_{\alpha ,\ell}
<\min_{\ell <n}\gamma_{\alpha ,\ell}
$$ 
then for some $\alpha(0)<\cdots <\alpha (n-1)$ we have, for every $\ell <n$:
$$\begin{array}{l}
\tau^0(x^{{\bf t}(\ell)}_{\gamma_{\alpha(0),\ell}},x^{{\bf t}(\ell)}_{
\gamma_{\alpha(1),\ell}},\ldots,x^{{\bf t}(\ell )}_{\gamma_{\alpha(n-1),
\ell}})=0\\
\tau^1(x^{-{\bf t}(\ell )}_{\gamma_{\alpha(0),\ell}}, x^{-{\bf t}(\ell)}_{
\gamma_{\alpha(1),\ell}},\ldots,x^{-{\bf t}(\ell)}_{\gamma_{\alpha(n-1),
\ell}})=0
\end{array}
$$
\item[($*$)$_4$] if $a_\alpha\in\B$ for $\alpha <\lambda$ then for some
$\alpha_0<\cdots <\alpha_{n-1}$ we have $\tau^0(a_{\alpha_0},\ldots,a_{
\alpha_{n-1}})=0$ 
\end{enumerate}
\medskip

\n {\bf Proof}. Easy. 
\medskip

\centerline{* * * * * * * * * *} 
\medskip

\n {\bf 12.18 Comments}.

1) This concludes the proof of the consistency of the existence, answering a
part of Monk's problem 33.

2) We can get ``$\B\vDash(\mbox{cf }(\mu) )^+$-c.c.'' when 12.12 provides one.

3) We may still like to get ``no $k$-independent set'' for some specific $k$
as done in 12.17. Probably also 11.13 will help but we have not really looked
into it. 

Clearly it is supposed to have, for a $J^+$-set of $i$'s:
\begin{enumerate}
\item[($*$)$_i$] for some function $F$, if $m<\omega$, and $X\subseteq ({\rm 
Dom }\; I_i)^m$ is $F$-large (i.e.\ if $k<\omega$, $\bar t^0,\ldots 
\bar t^{k-1}\in X$ and $F(\bar t^0,\ldots ,\bar t^{k-1})\in I$ then 
for some $\bar t\in X$, ${\rm Rang}\;\bar t\cap F(\bar t^0,\ldots 
,\bar t^{k-1})=\emptyset$).
\end{enumerate}
$\underline{\rm Then}$ for some distinct $\bar t^0,\ldots ,\bar t^{n-1}\in X$, 
we have
$$
\ell <m\Rightarrow \tau (t^0_\ell, t^1_\ell,\ldots ,t_\ell^{n-1})=0.
$$
See more in 15.11, 15.12.
\medskip

\section{The singular case}
We continue to deal with problem 33 of Monk \cite{M2}. This time we
concentrate on the case $\lambda$ is singular. Though a priori this looked to
be the side issue, we can get quite a coherent picture.

\n {\bf Note}: If $\kappa >\mbox{cf}(\lambda )$ there is such a Boolean
algebras (the disjoint sum of $\mbox{cf}(\lambda )$ Boolean algebras each of
cardinality $<\lambda$). Moreover
\medskip

\n {\bf 13.1 Claim}. Assume
\begin{enumerate}
\item[($*$)] $\lambda >\mbox{cf}(\lambda )=\theta$ and $(\forall\alpha
<\lambda )$ $(|\alpha |^{<\kappa}<\lambda )$ and $\lambda >\kappa
=\mbox{cf}(\kappa)>\aleph_0$.
\end{enumerate}

1) The following conditions are equivalent:
\begin{enumerate}
\item[(A)] there are $\B$ and $a_\zeta$ such that
\begin{enumerate}
\item[(a)] $\B$ is a $\kappa$-c.c.\ Boolean algebra
\item[(b)] $a_\zeta\in\B\backslash\{ 0\}$ for $\zeta <\theta$
\item[(c)] if $\lan w_\zeta :\zeta <\theta\ran$ is a sequence of pairwise
disjoint finite subsets of $\theta$ $\underline{\rm then}$ for some finite
$u\subseteq \theta$ we have
$$
\bigcap_{\zeta\in u}\bigcup_{\xi\in w_\zeta} a_\xi =0.
$$
\end{enumerate}
\item[(B)] there is a Boolean algebra $\B$ of cardinality $\lambda$ with no
independent subset of cardinality $\lambda$.
\end{enumerate}

2) The following conditions are equivalent
\begin{enumerate}
\item[(A)$'$] there are $\B$, $a_\zeta$ such that
\begin{enumerate}
\item[(a)] $\B$ is a $\kappa$-c.c. Boolean algebra
\item[(b)] $a_\zeta\in\B\backslash\{ 0\}$ for $\zeta <\theta$
\item[(c)] for any $X\in [\theta ]^\theta$ for some finite $w\subseteq X$ we
have $\bigcap_{\zeta\in w} a_\zeta =0$
\end{enumerate}
\item[(B)$'$] there is a Boolean algebra $\B$ of cardinality $\ge\lambda$
which does not have $\lambda$ as a free caliber.
\end{enumerate}
\smallskip

\n {\bf Proof}. 1) $\underline{\mbox{(A) $\Rightarrow$ (B)}}$.  The case
$\theta =\aleph_0$ is easier, so we leave it to the reader.

Without loss of generality $\B$ has cardinality $\theta$. Let $\lambda=\kappa
+\theta +\sum_{\zeta <\theta}\lambda_\zeta$ where $\lambda >\lambda_\zeta
>\kappa+\theta+\sum_{\xi <\zeta}\lambda_\xi$. Let $\B^*$ be the Boolean
algebra freely generated by $\B\cup\{x_{\zeta,\alpha}:\zeta<\theta,\alpha
<\lambda^+_\zeta\}$ except for the equations in $\B$ and
$$
x_{\zeta ,\alpha}\le a_\zeta\quad \mbox{(for $\zeta <\theta$, $\alpha
<\lambda^+_\zeta$).}
$$
Clearly $\B\subseteq\B^*$ and assume that $\{b_\gamma:\gamma<\lambda\}
\subseteq\B^*$ is independent. Then for each $\gamma$ there are $n(\gamma
)<\omega$ and Boolean terms $\tau_\gamma$ and $\zeta_{\gamma ,\ell}<\theta$,
$\alpha_{\gamma ,\ell}<\lambda_{\zeta_{\gamma ,\ell}}$ for $\ell <n_\zeta$ and
$c_{\gamma ,\ell}\in\B$ for $\ell <m(\gamma )$ such that $b_\gamma=\tau_\gamma
(x_{\zeta_{\gamma,0},\alpha_{\gamma,0}},\ldots,x_{\zeta_{\gamma,n(\gamma)-1},
\alpha_{\gamma,n(\gamma)-1}}$, $c_{\gamma,0},\ldots,c_{\gamma,m(\gamma)-1})$.
As ${\rm cf}(\lambda )=\theta>\aleph_0$, without loss of generality
$\tau_\gamma =\tau$, $n(\gamma )=n(*)$ and $m(\gamma )=m(*)$. Also for each
$\e <\theta$ there is $X_\e\in [\lambda^+_\e ]^{\lambda^+_\e}$ such that
\begin{enumerate}
\item[($*$)] $\gamma\in X_\e$ implies $\zeta_{\gamma,\ell}=\zeta_{\e,\ell}(*)
<\theta$, $c_{\gamma ,\ell}=c^*_{\e ,\ell}\in\B$.
\end{enumerate}
Without loss of generality, $\lan\zeta_{\e ,\ell}:\ell <n(*)\ran$ is
nondecreasing. We can find $Y\in [\theta ]^\theta$ such that $\lan\lan
\zeta_{\e ,\ell}(*):\ell <n\ran :\e\in Y\ran$ is a $\Delta$-system. In fact
for some $n'(*)\le n(*)$ we have 
\begin{enumerate}
\item[($*$)$_1$] $\e\in Y$ \& $\ell <n'(*)\Rightarrow\zeta_{\e,\ell}(*)=
\zeta_\ell (*)$.
\item[($*$)$_2$] $\e_1 \in Y$ \& $\e_2\in Y$ \& $\e_1<\e_2\Rightarrow
\zeta_{\e_1,n(*)-1}(*)<\zeta_{\e_2,n'(*)}(*)$.
\end{enumerate}
By renaming, without loss of generality $X_\e =[\lambda_\e ,\lambda^+_\e]$ for
$\e\in Y$. Let $w_\e =\{ \zeta_{\e ,\ell}(*):n'(*)\le\ell <n(*)\}$, so let $u$
be as required in clause (A)(c), so $u\subseteq\theta$ is finite.

Let for $\e\in u$, $\gamma_{\e ,1}<\gamma_{\e ,2}$ be members of $X_\e$.

Clearly
$$
b_{\gamma_{\e ,1}}\triangle b_{\gamma_{\e ,2}}\le\bigcup_{\ell\in [n'(*),n(*)}
a_{\zeta_{\e ,\ell}}
$$
hence
\begin{eqnarray*}
\bigcap_{\e\in u}(b_{\gamma_{\e ,2}}\triangle b_{\gamma_{\e ,2}})
&\le& \bigcap_{\e\in u}\Big(\bigcup_{\ell\in [n'(*),n(*))} a_{\zeta_{\e,\ell}}
\Big)\\
&=&\bigcup_{\e\in u}\bigcap_{\xi\in w_\e} a_\xi =\emptyset,
\end{eqnarray*}
so $\lan b_\gamma :\gamma <\lambda\ran$ is not independent.
\medskip

\n $\underline{\neg(A)\Rightarrow\neg (B)}$.

Like \cite{Sh:92} (in short: Let $\lambda=\sum_{\zeta <\theta}\lambda_\zeta$,
$(\forall\alpha <\lambda )$ $(|\alpha |^{<\kappa}<\lambda_\zeta )$,
$\lambda_\zeta =\mbox{ cf}(\lambda_\zeta )>\kappa +\theta +\sum_{\xi
<\zeta}\lambda_\xi$. Let $S_\zeta =\{\delta <\lambda_\zeta :\mbox{ cf}(\delta
)\geq\kappa\}$. Remember that by \cite{Sh:92}:
\begin{description}
\item[$\boxtimes_{\lambda_\zeta}{[\B]}$]\quad Let $\B$ be a
$\kappa$-c.c. Boolean algebra.  Then:
\begin{enumerate}
\item[($*$)] for any $\bar x=\lan x_\alpha :\alpha <\lambda_\zeta\ran$
pairwise distinct $x_\alpha\in\B$, there are $a^-<a^+$ in $\B\backslash\{
0\}$, such that: if $\lan\B_\alpha :\alpha <\lambda_\zeta\ran$ is an
increasing continuous sequence of subalgebras of $\B$ of cardinality
$<\lambda_\zeta$ satisfying $x_\alpha\in\B_{\alpha +1}$, $\{ a^-,
a^+\}\subseteq \B_0$, we have
$$
\left\{\delta\in S_\zeta:\,
\begin{array}{l}
a^-\le x_\delta\le a^+\mbox{ and}\\
(\forall y)[0<y\le a^+-a^-\;\mbox{ \& }\; y\in
B_\delta\rightarrow\\
 (x_\delta -a^-)\cap y\not= 0 \;\mbox{
\& }\; (a^+-x_\delta )\cap y\not= 0]
\end{array}\right\}
$$
is stationary.
\end{enumerate}
\end{description}
So fix $\bar x=\lan x_\gamma :\gamma <\lambda\ran$, sequence of distinct
elements of $\B$, for each $\zeta <\theta$ let $a^-_\zeta -a^+_\zeta$ be as in
($*$) (for $\bar x\restriction\lambda_\zeta$), and let $a_\zeta =a^+_\zeta
-a^-_\zeta\in\B^+$. Let $\B^\zeta_\alpha$ be the subalgebra generated by $\{
x_\gamma :\gamma <\max\{\alpha ,\bigcup_{\xi <\zeta}\lambda_\xi\}\}\cup
\{a_\xi :\xi <\theta\}$ for $\alpha <\lambda_\zeta$ and for each $\zeta
<\theta$ let $S_\zeta$ be as above.

As $\neg (A)$, necessarily there is a sequence of pairwise disjoint finite
subsets of $\theta$, say $\bar u=\lan u_\e :\e<\theta\ran$ with any finite
intersection of members $\lan\bigcup_{\zeta\in u_\e} a_\zeta :\e <\theta\ran$
is not zero.

Now we can manipulate, choosing by induction on $\e <\theta$, $\bar t^{\e
,\alpha}\in\prod_{\zeta\in u_\e}S_\zeta$ and defining
$$
a^*_{\e ,\alpha}=\bigcup_{\zeta\in u_\e}\Big(\big( a_\zeta -\bigcup_{\xi \in
u_\e\backslash (\zeta+1)} a_\xi \big)\cap x_{t_\zeta^{\e ,\alpha}}\Big).
$$

2) Similarly.{\hfill$\square_{\rm 13.1}$}

\smallskip

\n {\bf 13.2 Discussion}. 1) Note: if $\theta <\kappa$, clearly (A)$_\theta$
\& (A)$'_\theta$.

2) Note if $(\forall\alpha <\theta )$ $(|\alpha |^{<\kappa}<\theta )$, then
$\neg (A)_\theta$ \& $\neg (A)'_\theta$.

3) Note that if $\chi =\chi^{<\chi}<\chi (*)=\chi (*)^{<\chi (*)}$ then for
some $\chi^+$-c.c. $(<\chi )$-complete forcing notion of cardinality $\chi
(*)$ in $V^P$ we have $\neg (A)_\theta$ \& $\neg (A)'_\theta$ when $\theta
=\mbox{cf}(\theta )\in (\chi ,\chi (*))$.

4) It is natural to get ${\rm CON}(\kappa <\chi =\chi^{<\chi}<\theta
=\mbox{cf}(\theta )<2^\chi+ (A)_\theta$ \& $\neg (A)'_\theta)$. This is well
connected to our problems but we have not looked at it.
\medskip

\n {\bf 13.3 Claim}. In 11.3 the condition
\begin{enumerate}
\item[($*$)] $(\forall\alpha <\lambda )$ $(|\alpha |^{<\kappa} <\lambda$)
\end{enumerate}
can be replaced by the weaker one
\begin{enumerate}
\item[($*$)$^-$] for arbitrarily large regular $\lambda'<\lambda$ we have
$\boxtimes_{\lambda'}[\cBB ]$ for any $\kappa$-c.c.\ Boolean Algebra (see
13.1's proof).
\end{enumerate}

\section{Getting free caliber for regular cardinals}
Remember that $\lambda$ is a free caliber of a Boolean algebra $\B$ if for any
$X\in [\B ]^\lambda$ there is an independent $Y\in [X]^\lambda$; of course we
can replace a Boolean algebra by a locally compact topological space (which is
a slightly more general case, but the proof is not really affected).

Monk asks whether there is a $\kappa$-cc Boolean algebra $\B$ of cardinality
$\ge \lambda$ with no independent subset of cardinality $\lambda$, and $\mu$
such that
$$
\mu<\lambda<\mu^\kappa,\qquad (\forall\alpha<\mu)(|\alpha|^\kappa<\lambda).
$$
Here we deal with the case of $\lambda$ regular and give a sufficient
set-theoretic condition on $\kappa$ such that any $\kappa$-cc Boolean algebra
of cardinality $\ge \lambda$ has $\lambda$ as a free caliber, so the
consistency of a negative answer follows, but we do not directly force. So
this section is complementary to sections 12 and 11.  
\medskip

\n {\bf 14.1 Hypothesis}.
\begin{enumerate}
\item[(a)] $\lambda =\mbox{cf}(\lambda )>2^\kappa$, but for simplicity we
assume  
$$
\lambda =\mu^+,\quad \mu=\sum_{i<{\rm cf}(\mu)} \lambda_i,\quad\lambda_i=
\lambda_i^{<\kappa},\quad {\rm cf}(\mu)<\kappa.
$$
We shall use it to shorten proofs when helpful, and, later, will show what can
be done without it.
\item[(b)] $\B^*$ is a $\kappa$-cc. Boolean algebra, $a_\alpha\in\B$ for
$\alpha <\lambda$ are pairwise distinct.
\end{enumerate}
Let $\bar a\stackrel{\rm def}{=} \lan a_\alpha :\alpha <\lambda\ran$. We would
like to find $X\in [\lambda ]^\lambda$ such that $\{ a_\alpha :\alpha\in X\}$
is independent.
\medskip

\n {\bf 14.2 Definition}. {\it For $\B\subseteq\B^*$, $x\in\B^*$ let
$$
\begin{array}{l}
\mbox{Proj}^0(x,\B)\stackrel{\rm def}{=}\{ y\in\B :y\le x\}\\
\mbox{Proj}^1(x,\B)\stackrel{\rm def}{=}\{ y\in\B :y\cap x=0\}\\
\mbox{Proj}^2(x,\B)\stackrel{\rm def}{=}\left\{\begin{array}{l} y\in B:y=0\;\;
{\rm or}\\ \qquad\qquad(\forall z )(0<z\le y\; {\rm \&}\; z\in B\Rightarrow
0<z\cap x<z)\end{array}\right\}.
\end{array}
$$}
\smallskip

\n {\bf 14.3 Fact}. Let $\B\subseteq \B^*$, $x\in B^*$
\begin{enumerate}
\item[(1)] If $y_\ell\in\mbox{Proj}^\ell (x,\B )$ for $\ell <3$,
$\underline{\rm then}$ $\lan y_\ell :\ell <3\ran$ are pairwise disjoint.
\item[(2)] $\bigcup_{\ell <3}\mbox{Proj}^\ell (x,\B )$ is dense in $\B$.
\item[(3)] $\mbox{Proj}^\ell (x,\B )$ is an ideal on $\B$.
\item[(4)] ${\rm Proj}^\ell (x,\B )$ is complete inside $\B^*$ i.e.\ if in 
$\B^*$ we have $x$ is $\le$ lub of $\{x_\alpha:\alpha<\alpha^*\}$ and $\{
x_\alpha :\alpha <\alpha^*\}\subseteq {\rm Proj}^\ell (x,\B )$ and $x\in B$
$\underline{\rm then}$ $x\in {\rm Proj}^\ell(x,\B)$.
\end{enumerate} 
\smallskip

\n {\bf 14.4 Definition}. {\it
$$
\chi =\chi_{\bar a}={\rm Min}\{\|\B\| :\B \subseteq \B^*, |W_\B |=\lambda\},
$$
where
\begin{eqnarray*}
W_\B=W_{\B, \bar a} &=&\left\{ \alpha :{\rm Proj}^2(a_\alpha ,\lan\B\cup\{
a_\beta :\beta <\alpha\}\ran_{\B^*} )=\{ 0\},\right.\\
&&\qquad\quad{\rm and}\;\; {\rm Proj}^\ell (a_\alpha ,\B )\;\; {\rm
is ~dense ~in}\\
&&\qquad\quad{\rm Proj}^\ell (a_\alpha ,\lan\B\cup\{ a_\beta :\beta
<\alpha\}\ran_{\B^*})\\
&&\qquad\quad\left.{\rm for}\;\; \ell =0,1\right\}.
\end{eqnarray*}}
\smallskip

\n {\bf 14.5 Remark}.
\begin{enumerate}
\item[(1)] ${\rm Proj}^2(a_\alpha ,B)=\{0\}$ is close to saying, $a_\alpha=$
the lub in $\B^*$ of ${\rm Proj}^0(a_\alpha ,\B)$, but not the same (holds if
$\B \lcirc\B^*$).

Could have worked with a variant as indicated.

\item[(2)] Trivially $\chi \le \lambda$, use $\B=\lan a_\alpha :\alpha
<\lambda\ran_{\B^*}$.
\end{enumerate}
\smallskip

\n {\bf 14.6 Fact}. If $\chi =\lambda$, $\underline{\rm then}$ for some $X\in
[\lambda ]^\lambda$, $\lan a_\alpha :\alpha \in X\ran$ is independent.
\medskip

\n {\bf Proof}. Let $\B_\alpha\stackrel{\rm def}{=}\lan a_\beta :\beta
<\alpha\ran_{\B^*}$, so $\B_\alpha$ are increasing continuous in $\alpha$,
$\|\B_\alpha\|\le \aleph_0+|\alpha |<\lambda$. Let
$$
\begin{array}{l}
S\stackrel{\rm def}{=}\{\alpha <\lambda: {\rm Proj}^2(a_\alpha ,\B_\alpha )=\{
0\}\}\\ 
S'\stackrel{\rm def}{=}\{\alpha\in S:{\rm cf}(\alpha)\ge\kappa\}.
\end{array}
$$
Now
\begin{enumerate}
\item[($*$)] $S'$ is not stationary.
\end{enumerate}
[Why? For $\delta\in S'$, $\ell <2$ let $\II_{\delta ,\ell}\subseteq {\rm
Proj}^\ell (a_\delta ,\B_\delta )\backslash\{ 0\}$ be an antichain, maximal
under the conditions defining ${\rm Proj}^\ell$. So $|\II_{\delta,\ell}|<
\kappa$, as $\B^*\vDash\kappa$-cc. Hence for some $f(\delta )<\delta$ we have
$$
\II_{\delta ,0}\cup\II_{\delta ,1}\subseteq \B_{f(\delta )}.
$$
So if $S'$ is stationary, by Fodor lemma, for some $\alpha^*<\lambda$, $S^*=\{
\delta\in S':f(\delta )=\alpha^*\}$ is stationary.

We would like to show:
\begin{enumerate}
\item[($**$)]\qquad $\delta\in S^*\Rightarrow {\rm Proj}^2(a_\delta,\B_{
\alpha^*})=\{ 0\}$.
\end{enumerate}
If so, we have gotten that $\B_{\alpha^*}$, $S^*$ exemplify $\chi\le
\|\B_{\alpha^*}\|$, contradiction. For proving ($**$), let $\delta\in S^*$,
assume $b\in {\rm Proj}^2(a_\delta ,\B_{\alpha^*})\backslash\{ 0\}$.

So, by 14.3, (for $\B_{\alpha^*}, a_\delta$) we have ($\forall x\in
\II_{\delta ,0}\cup\II_{\delta ,1}$) $x\cap b=0$.

Now, $b\not\in {\rm Proj}^2(a_\delta ,\B_\delta )$, as the latter is $\{
0\}$. So, there is $c$ such that $\B_\delta\vDash\mbox{``}0<c\le b$ and $c\cap
a_\delta =0\vee c\leq a_\delta\mbox{''}$, that is $c\in\mbox{Proj}^0(a_\delta
,\B_\delta )\cup{\rm Proj}^1(a_\delta ,\B_\delta )$, but as $c\leq b$ we have
$$
(\forall x\in\II_{\delta ,0}\cup\II_{\delta ,1})\; (x\cap c=0).
$$
So $c$ contradicts the maximality of $\II_{\delta ,0}$ (if $c\in{\rm
Proj}^0(a_\delta ,\B_\delta )$) or of $\II_{\delta ,1}$ (if $c\in {\rm
Proj}^1(a_\delta ,\B_\delta )$).

The contradiction proves ($**$) and ($*$).]
\medskip

So $\lambda\backslash S$ is stationary. For $\delta\in\lambda\backslash S$
choose $b_\delta\in {\rm Proj}^2(a_\delta ,\B_\delta )\backslash\{ 0\}$. So by
Fodor's lemma, for some $b^*\in\bigcup_{\alpha <\lambda} \B_\alpha$ we have
$$
S^{**}\stackrel{\rm def}{=}\{\delta :\delta\in\lambda\backslash S, b_\delta
=b^*\}\;\;\mbox{is stationary}.
$$
Now we know that $\lan a_\delta :\delta\in S^*\ran$ is independent.
{\hfill$\square_{\rm 14.6}$}
\medskip

\n {\bf 14.6A Remark}. In the characteristic case, $\B^*$ is the completion of
a Boolean algebra of smaller cardinality $\B'$, so $\chi\leq\|\B'\|$.
\medskip

\n {\bf 14.7 Claim}. Now, without loss of generality
\begin{enumerate}
\item[$\boxtimes$] \qquad $\B^*=\lan\B\cup\{ a_\alpha :\alpha\in W_\B
\}\ran$ for some $\B\subseteq \B^*$,
\item[{}] $\|\B\|=\chi$, $W_\B=\lambda$.
\end{enumerate}
\smallskip

\n {\bf Proof}. $\B\subseteq\B^*$ exemplifies the value of $\chi$, let $\B^c$
 be the completion of $\B$, and we can let for $\alpha\in W_\B$
$$
a'_\alpha = {\rm lub ~in}\; \B^c\;\mbox{ of }\;{\rm Proj}^0(a_0,\B ).
$$
Now if $Y\in [W_\B ]^\lambda$, $\lan a'_\alpha :\alpha\in Y\ran$ is
independent in $\B^c$ then $\{ a^*_\alpha :\alpha\in y\}$ is independent in
$\B^*$. Alternatively use $\lan\B\cup\{a_i:\alpha\in W_\B\}\ran_{\B^*}$.
\smallskip

\n (Remember: $\B$ is not necessarily a complete subalgebra of $\B^*$.)
\medskip

\n {\bf 14.8 Definition} {\it
Let
$$
K\stackrel{\rm def}{=}\left\{
\begin{array}{l} \bar\B :\bar \B=\lan\B_i:i\le\chi\ran\;\;
\mbox{is an increasing}\\
\qquad\mbox{continuous sequence of subalgebras}\\
\qquad\mbox{of $\B^*$, $\|\B_i\|<\aleph_0+|i|^+$,}\\
\qquad\mbox{and $W_{\B_\chi}\in [\lambda ]^\lambda$,
$\B_\chi\supseteq\B$ (of $\boxtimes$ of 14.7)}.
\end{array}\right\}
$$
(so $W_{\B_\chi}$ is cobounded in $\lambda$, in fact if
$\B_\chi\subseteq\lan\B\cup\{ a_\beta :\beta <\alpha\}\ran_{\B^*}$ then
$|W_{\B_\chi}\supseteq (\alpha ,\lambda )|$).}
\medskip

\n {\bf 14.9 Fact}
\begin{enumerate}
\item[(1)] $\mbox{cf}(\chi )<\kappa$
\item[(2)] $\mbox{cov}(\chi ,\chi ,\kappa ,2)\geq \lambda$, meaning:
$$\lambda\le \min
\left\{\begin{array}{l}
 |\PP |:\PP\subseteq [\chi
]^{<\chi}\;\;\mbox{\&}\\
\qquad\quad(\forall A\in [\chi ]^{<\kappa})(\exists B\in\PP)(A\subseteq B)
\end{array}\right\}.
$$
\end{enumerate}
\smallskip

\n {\bf Proof}. \quad (1) By (2).

\n (2)\quad Assume not. Remember $\B\subseteq\B^*$, $|W_\B|=\lambda$, $\|\B\|
=\chi$.

\n For each $\alpha\in W_\B$ choose $\II_{\alpha ,\ell}\subseteq {\rm
Proj}^\ell (a_\alpha ,\B)$ for $\ell <2$ as in the proof of 14.6. Let
$\PP\subseteq [\B]^{<\chi}$, $|\PP |<\lambda$ and 
$$
(\forall A\in [\chi ]^{<\kappa})(\exists B\in\PP )(A\subseteq B).
$$
So for each $\alpha \in W_\B$, there is $A_\alpha\in\PP$ such that
$\II_{\alpha ,0}\cup\II_{\alpha ,1}\subseteq A_\alpha$.  So for some
$A^*\in\PP$
$$
W=\{\alpha\in W_\B :\II_{\alpha,0}\cup\II_{\alpha ,1}\subseteq A^*\}\in
[\lambda ]^\lambda.
$$
(exists as we divide $W_\B$ into $|\PP|$ sets, so at least one has size
$\lambda$, as $|\PP|<\lambda ={\rm cf}(\lambda )$). Now $\chi\le |\lan
A^*\ran_\B|$, contradiction, as in the proof of 14.6 (to the definition of
$\chi$.{\hfill$\square_{\rm 14.9}$}
\medskip

\n {\bf 14.10 Definition}. {\it For $\bar\B\in K$ and $\alpha\in W_{\B_\chi}$
let
$$
u (\alpha, \bar B)\stackrel{\rm def}{=}\left\{
\begin{array}{l}
\mbox{$i<\chi:$ for some $\ell <2$, ${\rm Proj}^\ell
(a_\alpha ,\B_i)$ is not a predense}\\
\qquad\qquad\mbox{subset of ${\rm Proj}^\ell (a_\alpha ,\B_{i+1})$}
\end{array}\right\}.
$$}
\medskip

\n {\bf 14.10A Discussion}. We may consider $\bar\B'=\lan\B'_i:i\le\chi\ran
\in\kappa$ when
$$
\B'_i=\lan \B_i\cup X\ran, X\;\mbox{ fixed countable $\subseteq\B^*$}.
$$
Possibly
$$
u (\alpha ,\bar\B )\not= u (\alpha ,\bar\B')
$$
or just for some $i$, ${\rm Proj}^\ell (a_\alpha ,\B_i)$ is not dense in ${\rm
Proj}^\ell (a_\alpha ,\B'_i)$.  We think of the set of such $\alpha$ as bad,
and put them all in one $\lambda$-complete ideal. But maybe $\lambda$ belongs
to it.  So we will try to find some $\bar\B$ for which this does not occur.

This will help in that we eventually try to choose $\alpha_\zeta\in W_{\bar
\B}$ for $\zeta <\lambda$ by induction on $\lambda$ such that $\lan
a_{\alpha_\zeta}:\zeta <\lambda\ran$ is independent.

So in stage $\zeta$ we consider all
$$
X\in [\{ \alpha_\xi :\xi <\zeta\} ]^{<\aleph_0}.
$$
The existence of $\bar\B$ requires some properties of $\lambda$ which
certainly hold in the main case (with $\lambda =\mu^+\ldots$).

So to ease the proof instead of every $i<\chi$, we use ``every $i<\chi$ large
enough''.
\medskip

\n {\bf 14.11 Definition}. {\it 
\begin{enumerate}
\item[(1)] We define a partial order on $K :\bar\B^1\le\bar\B^2$
$\underline{\rm if}$ for every $i$ large enough
$$
i\leq\chi\Rightarrow \B^1_i\subseteq \B^2_i.
$$
\item[(2)] We say $\bar\B^2$ is finitely generated over $\bar \B^1$
$\underline{\rm if}$ for some finite $X$
$$
\B^2_i=\lan \B^1_i\cup X\ran_{\B^*}\;\;\mbox{for $i<\chi$ large enough.}
$$
In this case we let $\bar\B^1[X]=\lan\B^1_i[X]:i\le \chi\ran$ be $\bar \B^2$.

\item[(3)] For $\bar\B^1\le\bar \B^2$ let
$$
\mbox{Bad}(\B^1,\B^2)=\left\{\begin{array}{l}
\alpha :\mbox{ if $\alpha\in W_{\B^1_\chi}\cap
W_{\B^2_\chi}$,}\\
\qquad\quad\mbox{then for arbitrarily large $i<\chi$, for some}\\
\qquad\quad\mbox{$\ell <2$, ${\rm Proj}^\ell(a_\alpha ,\B^1_i)$ is not
dense}\\
\qquad\quad\mbox{in ${\rm Proj}^\ell(a_\alpha,\B^2_i)$}\end{array}\right\}.
$$
\item[(4)] $J_{\bar\B^1}$ is the $\lambda$-complete ideal on $\lambda$
generated by all ${\rm Bad}(\bar\B^1, \bar\B^2)$, where $\bar\B^1\le \bar
\B^2$ and $\bar\B^2$ is finitely generated over $\bar\B^1$.
\end{enumerate}}
\smallskip

What do we need to carry a proof?
\medskip

\n {\bf 14.12 Lemma}. There is $\bar\B^\otimes \in K$ such that
$\lambda\not\in J_{\bar B^\otimes}$.
\medskip

\n {\bf 14.12A Remark}. We may like to have $J\supseteq J_{\bar B^\otimes}$
normal extending $I^{{\rm nst},\theta}_\lambda$ (and $\lambda\not\in J$), then
we need more work.
\medskip

\n {\bf Proof in the case $\lambda =\chi^+$}. (Enough, see 14.1(a)). Assume
there is no such $\bar\B =\bar\B^\otimes$. We choose by induction on $\zeta
<\chi$, $\bar\B^\zeta \in K$, such that $\bar\B^\zeta$ is increasing with
$\zeta$ and: for each $\zeta$, as $\lambda\in J_{\bar\B^\zeta}$ we can find
$\lan X_{\zeta ,\epsilon} :\epsilon <\epsilon_\zeta \ran$ witnessing it, i.e.\
$X_{\zeta ,\epsilon}\in [\B^*]^{<\aleph_0}$, $\epsilon_\zeta <\lambda$ (so
without loss of generality $\epsilon_\zeta\le\chi$)
$$
\lambda =\bigcup_{\epsilon<\epsilon_\zeta}{\rm Bad}(\bar\B^\zeta,\bar\B^\zeta
[X_{\zeta,\epsilon}])
$$
where
$$
\B^\zeta_i [X_{\zeta ,\epsilon}]=\lan\B^\zeta_i\cup X_{\zeta,\epsilon}
\ran_{\B^*}.
$$
Now easily $(K,\le )$ is $\chi^+$-directed, so we demand
$$
\bigwedge_{\epsilon <\epsilon_\zeta}\bar\B^\zeta \le\bar\B^\zeta
[X_{\zeta ,\epsilon}]\le \bar\B^{\zeta +1}.
$$
Also $i\in [\zeta, \lambda)$ \& $\zeta <\xi\le \chi\Rightarrow
\B^\zeta_i\subseteq \B^\xi_i$. Let $\delta^*<\lambda$ be such that
$$
\bigwedge_{\zeta <\chi}\B^\zeta_\chi\subseteq\lan \B\cup\{ a_\alpha :\alpha
<\delta^*\}\ran_{\B^*}.
$$
So for each $\zeta <\chi$ we have
$$
\delta^*\in\bigcup_{\epsilon <\epsilon_\zeta}{\rm
Bad}(\bar\B^\zeta ,\bar\B^\zeta [X_{\zeta ,\epsilon}]),
$$
hence there is $\xi (\zeta )<\epsilon_\zeta$ such that
$$
\delta^*\in {\rm Bad}(\bar\B^\zeta ,\bar\B^\zeta[X_{\zeta ,\xi (\zeta )}]).
$$
For each $\zeta <\chi$, there is $i(\zeta )<\chi$ such that $X_{\zeta ,\xi
(\zeta )}\subseteq \B^{\zeta +1}_{i(\zeta )}$, hence
$$
(\forall i)[i(\zeta )\le i\le\chi\Rightarrow \B^\zeta_i[X_{\zeta ,\zeta (\xi
)}]\subseteq\B^{\zeta +1}_i
$$
because
$$
X_{\zeta,\xi(\zeta)}\subseteq\B^{\zeta+1}_{i(\zeta)}\subseteq\B^{\zeta+1}_i.
$$
We restrict ourselves to $\xi <\kappa$. So without loss of generality
$$
\bigwedge_{\zeta_1 <\zeta_2\le \kappa}\;\;\bigwedge_{\alpha\in
[\kappa^+,\chi ]}\B^{\zeta_1}_\alpha\subseteq\B^{\zeta_2}_\alpha,
$$
and if $\zeta$ is a limit and $\alpha\in [\kappa^+,\chi]$, then
$\B^\zeta_\alpha=\bigcup_{\xi <\zeta}\B^\xi_\alpha$. As ${\rm cf}(\chi
)<\kappa$, there is $i(*)<\chi$ such that $Z=\{\zeta <\kappa:i(\zeta )\le
i(*)\}$ is unbounded (we can demand more).

Now the set $u (\delta^*,\bar\B^\kappa )$ has cardinality $<\kappa$ because
$\B^*$ satisfies the c.c.c.

Remember,
$$
u (\delta^*,\bar\B^\kappa )=\left\{
\begin{array}{rl}i<\chi :&\bigcup_{\ell =0,1}
{\rm Proj}^\ell (a_{\delta^*},\B^\kappa_i)\\
&\mbox{is not predense in $
\bigcup_{\ell =0,1}{\rm Proj}^\ell
(a_{\delta^*},\B^\kappa_{i+1})$}\end{array}\right\}.
$$
Choose for $i\in u (\delta^*,\bar\B^\kappa)\cup\{\kappa^+\}$ and $\ell =0,1$
a predense subset $\II^{\delta^*,\ell}_{\kappa,i}$ of ${\rm Proj}^\ell
(a_{\delta^*},\B^\kappa_{i+1})$ of cardinality $<\kappa$.

Now, for $i\in u (\delta^*,\bar\B^\kappa)\cup\{\kappa^+\}\backslash\kappa^+$
the sequence $\lan\B^\zeta_{i+1}:\zeta\le\kappa\ran$ is increasing
continuous. So for some $\zeta_i<\kappa$
$$
\II^{\delta^*,0}_{\kappa,i}\cup\II_{\kappa,i}^{\delta^*,1}\subseteq
\B^{\zeta_i}_{i+1}.
$$
Let
$$
\zeta (\delta^*)\stackrel{\rm def}{=}\sup_i\zeta_i<\kappa.
$$
So clearly
\begin{enumerate}
\item[($*$)]  if $i\in [\kappa^+,\chi ]$, $\ell <2$, then
$$
{\rm Proj}^\ell (a_{\delta^*},\B^{\zeta [\delta^*]}_i)={\rm
Proj}^\ell (a_{\delta^*},\B^\kappa_i)\cap\B^{\zeta
[\delta^*]}_i
$$
is a predense subset of ${\rm Proj}^\ell (a_{\delta^*}, \B^\kappa_i)$.
\end{enumerate}
[Why? By induction on $i$. If $i=\kappa^+$ directly.  If $i$ is a limit --
trivial. If $i=j+1\ge \kappa^+$, $j\not\in u (\delta^*, \bar\B^\kappa )$, then
by transitivity of being predense in. If $i=j+1$, $j\in u(\delta^*,
\bar\B^\kappa)$, using $\II_j^{\delta^*,\ell}$.]

Now, clearly 
$$
\begin{array}{l}\zeta\in [\zeta (\delta^*),\kappa )\qquad\Rightarrow\\
\bigwedge\limits_{\ell<2}\;\;\bigwedge\limits_{i\in [\kappa^+,\chi )}({\rm
Proj}^\ell(a_{\delta^*}, \B^\zeta_i)\;\mbox{ is predense in ${\rm Proj}^\ell
(a_{\delta^*}, \B^{\zeta+1}_i)$}).
\end{array}
$$
This follows from ($*$). Choose $\zeta\in Z\backslash\zeta (\delta^*)$ so we
contradict the choice of $\bar\B^{\zeta +1}$.{\hfill$\square_{\rm 14.12}$}
\medskip

\n {\bf 14.13 Convention}. We fix $\bar\B^\otimes\in K$ such
that $\lambda\not\in J_{\bar\B^\otimes}$.
\medskip

\n {\bf 14.14 Fact}. $\{\alpha <\lambda : u (\alpha ,\bar\B^\otimes)\;\mbox{
bounded in }\;\chi\}$ is bounded in $\lambda$.
\medskip

\n {\bf Proof}. By the choice of $\chi$ as minimal.  {\hfill$\square_{\rm
14.14}$}
\medskip

\n {\bf 14.15 Convention}. Let $f_\alpha$ be an increasing function from ${\rm
otp}(u (\alpha,\bar\B^\otimes ))$ onto $u (\alpha,\bar\B^\otimes )$.
\medskip

\n {\bf 14.16 Fact}. For some $j^*<\kappa$
$$
Y_0=\{\alpha <\lambda: {\rm Dom}(f_\alpha )=j^*\}\in (J_{\bar\B^\otimes})^+.
$$
So without loss of generality $(\forall \alpha)[{\rm Dom}(f_\alpha )=j^*]$.
\medskip

\n {\bf 14.17 Claim}. We can find $\lan\gamma^*_j:j<j^*\ran$, $w^*\subseteq
j^*$ such that:
\begin{enumerate}
\item[$(*)_1$] $\underline{\rm if}$ $\bar\gamma =\lan\gamma_j: j<j^*\ran$,
$\gamma_j\le \gamma^*_j$,
$$
\gamma_j=\gamma^*_j\Leftrightarrow j\in w^*,
$$
$\underline{\rm then}$ the set of $\alpha\in Y_0$ satisfying the following, is
in $(J_{\B^\otimes})^+$:
$$\begin{array}{l}
j\in w^*\Rightarrow f_\alpha (j)=\gamma^*_j\\
j\in j^*\backslash w^*\Rightarrow\gamma_j<f_\alpha
(j)<\gamma^*_j.\end{array}
$$
Also
\item[$(*)_2$] $j\in j^*\backslash w^*\Rightarrow {\rm cf}(\gamma^*_j)>
2^\kappa$ and
$$
\lambda =\max {\rm pcf}\{{\rm cf}(\gamma^*_j):j\in j^*\backslash w^*\}.
$$
\item[$(*)_3$] Moreover if we fix $\mu =\mu^{<\kappa} <\lambda$ we can
demand 
$$
j\in j^*\backslash w^*\Rightarrow {\rm cf}(\gamma^*_j)>\mu.
$$
\item[$(*)_4$] if $j^*=\sup (J^*\backslash w^*)$, and $E$ is the equivalence
relation on $j^*\backslash w$ defined by $j_1\; E\; j_1\Leftrightarrow
\gamma^*_{j_1}=\gamma^*_{j_2}$ (so the equivalence classes are convex) then
$J$ is an ideal on $j^*$ such that $J^{\rm bd}_{j^*}\subseteq J$, $w^*\in J$,
$$
A\in J\quad\Rightarrow\quad\bigcup \{ j/E:j\in A\}\in J,\qquad\mbox{ and}$$ 
\begin{description}
\item[$(\alpha)$]\quad $\prod_{j<j^*} \gamma^*_j/J$ has true cofinality
$\lambda$, so possibly shrinking $Y_0$, for $\alpha <\beta$ in $Y_0$,
$f_\alpha<_J f_\beta$.
\end{description}
\end{enumerate}
\smallskip

\n {\bf Proof}. By 7.0(0) (or \cite[6.6D]{Sh:430} or \cite[6.1]{Sh:513}), as
$j^*<\kappa$, so $2^{|j^*|}<\lambda$. ~{\hfill$\square_{\rm 14.17}$}
\medskip

\n {\bf 14.18 Observation}. $\lan \gamma^*_j:j<j^*\ran$ is non-decreasing,
with limit $\chi$, and $\gamma^*_j<\chi$ and of course, ${\rm cf}(j^*)={\rm
cf}(\chi )$.
\medskip

\n {\bf Proof}. As ${\rm Rang}(f_\alpha )\subseteq\chi$, and the fact,
$\gamma^*_j<\chi$ if $j\in w^*$, $\gamma^*_j\le \chi$ if $j\not\in w^*$, but
then
$$
{\rm cf}(\gamma^*_j)\ge 2^\kappa >\kappa >{\rm cf}(\chi ).
$$
~{\hfill$\square_{\rm 14.18}$}
\medskip

\n {\bf 14.19 Comment on the Claim}. 1) For it, possibly $\wedge_\alpha
f_\alpha =f^*$, so then we get $w^*=j^*$. Also possibly $f_\alpha (j)<\alpha$,
so $w^*=\emptyset$ and $J=\{ \phi\}$.

2) If the ideal $J_{\bar\B^\otimes}$ is normal enough, for some $X\in
(J_{\bar\B^\otimes})^+$, $\lan f_\alpha :\alpha\in X\ran$ is $<_J$-increasing.

3) If $(\forall\alpha <\lambda )(|\alpha |^{|j^*|}<\lambda )$, then
necessarily
$$
j\in j^*\backslash w^*, {\rm cf}(\gamma^*_j)=\lambda
$$
(like the $\triangle$-system lemma).  $\underline{\rm BUT}$ for the
interesting case, and in particular by our assumptions, this is not the case:
as $\gamma^*_j\le\chi <\lambda$, hence $J\supseteq [j^*]^{<\aleph_0}$.
\medskip

\n {\bf 14.20 Hypothesis}. Each $\B^\otimes_i$ is the union of $\mu$ filters
$\lan\DD_{i,\beta}:\beta <\mu\ran$, $\mu =\mu^{<\kappa}$ (we can use somewhat
less), this of course is only a consistent assumption.  
\medskip

\n {\bf 14.21 Claim}. For some
$$
\bar\iota =\lan\iota_j:j<j^*\ran\in {}^{j^*}\mu
$$
we can restrict ourselves to 
$$
Y_1=\left\{\alpha <\lambda:
\begin{array}{l}
j\in w^*\Rightarrow
f_\alpha (j)=\gamma^*_j,\\
j\in j^*\backslash w^*\Rightarrow
\gamma^{**}_j<f_\alpha (j)<\gamma^*_j \;{\rm and}\\
\displaystyle{\bigwedge_{j<j^*}({\rm
Proj}^2(a_\alpha ,\bar\B^\otimes_{\gamma^{**}_j})\cap
\DD_{\gamma^{**}_j,\iota_j}\not= \{ 0\})}\end{array}
\right\},
$$
where
$$
\gamma^{**}_j=\left\{\begin{array}{ll}
\gamma^*_j&\quad\mbox{if $j\in w^*$}\\
\cup\{\gamma^*_i:
\gamma^*_i<\gamma^*_j\}&\quad\mbox{otherwise}.\end{array}
\right.
$$
in particular $Y_1\not\in J_{\bar B^\otimes}$.
\smallskip

\n {\bf Proof}. As $\mu^{<\kappa}<\lambda$, $J_{\bar\B^\otimes}$
$\lambda$-complete and $j\in j^*\backslash w^*\Rightarrow {\rm
cf}(\gamma^*_j)>\mu$.{\hfill$\square_{\rm 14.21}$}
\medskip

\n {\bf 14.22 Claim}: For some $X\in [X]^\lambda$, the sequence $\lan a_\alpha
:\alpha \in X\ran$ is independent.
\medskip

 \n {\bf Proof}: {\bf Case 1}. $w^*$ is unbounded in $j^*$: We choose by
induction on $\beta <\lambda$,
$$
N_\beta \prec (\HH ((2^\lambda )^+), \in ,<^*_{(2^\lambda )^+})
$$
increasing continuous, $\|N_\beta\|<\lambda$, $N_\beta\cap\lambda\in\lambda$,
$\lan N_{\beta_1:\beta_1\le \beta}\ran\in N_{\beta+1}$ and $\bar\B^\otimes$,
$\B^*$, $\bar a\in N_0$. Let $\alpha_\beta =\alpha (\beta )$ be the first
$\alpha <\lambda$ such that $\alpha\in Y_1$, $\alpha\not\in \cup (J_{\bar
\B^\otimes}\cap N_\beta )$
$$
\Big(\mbox{so }\;\bigwedge_{j\in w^*}f_\alpha (j)=\gamma^*_j\Big)
$$
Clearly $\alpha_\beta\in \lambda\cap N_{\beta +1}\backslash N_\beta$,
$\lan\alpha_{\beta_1}:\beta_1<\beta\ran\in N_{\beta +1}$.  Let $n<\omega$,
$\beta_1<\cdots <\beta_n$ and we will prove that $\lan \alpha_{\alpha
(\beta_\ell)}:\ell =1,\ldots ,n\ran$ is independent.

Now 
$$
j\in w^*\Rightarrow\;\mbox{ there is }\; b_j\in\bigcap^n_{\ell =1}{\rm
Proj}^2(a_{\alpha_{\beta_\ell}}, \bar\B_{\gamma^*_j}^\otimes )\backslash\{0\}.
$$
[Why? As $\alpha_{\beta_1},\ldots ,\alpha_{\beta_n}\in Y_1$, so $\DD_{
\gamma^*_j,\iota_j}\cap {\rm Proj}^2(a_{\alpha_{\beta_\ell}},\B^\otimes_{
\gamma^*_j})\not=\emptyset$. Choose $b_{j,\ell}$ there, so $b_j=\bigcap^n_{
\ell =1}b_{j,\ell}$ is OK.]
\smallskip

Consider
$$
{\rm Bad}(\B^\otimes , \bar\B^\otimes [\{ a_{\alpha(\beta_1)},\ldots
,a_{\alpha (\beta_\ell)}\}])\in J_{\bar\B^\otimes},
$$
it belongs to $N_{\beta_{\ell +1}}$. So
$$
\alpha_{\beta_{\ell +1}}\not\in {\rm Bad}(\B^\otimes ,\B^\otimes [\{ a_{\alpha
(\beta_1)},\ldots ,a_{\alpha (\beta_\ell)}\} ]).
$$
So for each $\ell$ for some $i_\ell <\chi$, $k<2$ \& $i\in [i_\ell ,\chi
)\Rightarrow {\rm Proj}^k(a _{\alpha_{\beta_{\ell +1}}},\B^\otimes_i )$ is
predense in ${\rm Proj}^k(a_{\alpha_{\beta_{\ell +1}}},\lan\B^\otimes_i\cup\{
a_{\alpha_{(\beta_1)}},\ldots ,a_{\alpha_{(\beta_\ell)}}\}\ran )$.

So if $j\in w^*$, $\gamma^*_j>\sup_{\ell =1,\ldots ,n} i_\ell$ (exists) and
$\eta\in{}^{[1,n]}2$, we prove by induction on $\ell$ that
$$
b^\ell_j=b_j\cap\bigcap^\ell_{k=1}(a_{\alpha_{\beta_k}})^{[\eta (k)]}.
$$

For $\ell =0$ trivial.

For $\ell >0$, $b_j^{\ell -1}\in\lan \B^\otimes_{\gamma^*_j}\cup\{
a_{\alpha_{\beta_1}},\ldots ,a_{\alpha_{\beta_{\ell -1}}}\}\ran$ is $>0$, is
in
$$
{\rm Proj}^2(a_{\alpha_\beta}, \lan \B^\otimes_{\gamma ^*_j}\cup \{ a_{\alpha
(\beta_1)},\ldots a_{\alpha (\beta_{\ell -1}}\}\ran )
$$
as it is below $b_j$ and $b_j\in {\rm Proj}^2(a_{\alpha
(\beta_\ell)},\B^\otimes_{\gamma^*_j})$ by its choice and $j$ is $>i_\ell$, so
$b_j\in {\rm Proj}^2(a_{\alpha_{\beta_\ell}},\lan\bar\B^\otimes\cup\{a_{
\alpha_{\beta_1}},\ldots ,a_{\alpha_{\beta_{\ell -1}}}\}\ran )$. We use
implicitly 
\medskip

\n {\bf 14.23 Fact}. For $\alpha <\lambda$ large enough,
$$
i<\chi\Rightarrow {\rm Proj}^2(a_\alpha ,\bar\B^\otimes_i )\not= \{ 0\}.
$$
\smallskip

\n {\bf Proof}. By $\chi$'s minimality.{\hfill$\square_{\rm 14.23}$}
\medskip

\n {\bf Case 2}. Not 1, i.e.\ $w^*$ bounded in $j^*$ $\underline{\rm or}$ just
$j^*=\sup (J^*\backslash w)$. Similarly using ($*$)$_2$ of 14.17 find $j\in
j^*\setminus w^*$ such that if $j_\ell\in j/E$ for $\ell =1,\ldots ,n$ then
$f_{\alpha_{\beta_1}}(f_1)<f_{\alpha_{\beta_2}}(j_2)<\cdots<f_{\alpha_{
\beta_n}}(j_n)$. {\hfill$\square_{\rm 14.22}$}
\medskip\smallskip

\centerline{$\ast\ast\ast$}
\smallskip

\n {\bf 14.24 Conclusion}. If $\mu =\mu^{<\mu}<\theta =\theta^{<\theta}$ then
for some $\mu$-complete $\mu^+$-c.c.\ forcing notion of cardinality $\theta$,
in $V^P$:

$\underline{\rm If}$ $\B$ is a $\kappa$-c.c.\ Boolean algebra of cardinality
$\ge\lambda$, $\mu =\mu^{<\kappa}$, $\lambda ={\rm cf}(\lambda )\in (\mu,
\theta]$ $\underline{\rm then}$ $\lambda$ is a free caliber of $\B$.
\medskip

\n {\bf Proof}. By 14.1--14.23 above and \cite{Sh:80}.{\hfill$\square_{\rm
14.24}$} 
\medskip

\n {\bf 14.25 Claim}. The following implications hold: ($*$)$_1$ $\Rightarrow$
($*$)$_2$ $\Rightarrow$ ($*$)$_3$ $\Rightarrow$ ($*$)$_4$ where
\begin{enumerate}
\item[($*$)$_1$] (a) ~ $\mu^{2^{<\kappa}}=\mu <\lambda ={\rm cf}(\lambda )$.
\item[{}] (b) if a Boolean algebra $\B$ satisfies the $(2^{<\kappa})^+$-c.c.\
and $|\B |<\lambda$,

\n $\underline{\rm then}$ $\B$ is the union of $\mu$ filters.

\item[($*$)$_2$] (a) $\kappa <\lambda ={\rm cf}(\lambda )$
\item[{}] (b) if a Boolean algebra $\B$ satisfies the $\kappa$-c.c., for
$i<\lambda$, $F_i\subseteq\B\backslash\{ 0\}$ is a set of $<\kappa$ members
closed under intersection $\underline{\rm then}$ we can find $<\lambda$
filters $\DD_\alpha$ $(\alpha <\alpha^*<\lambda )$ of $\B$ such that $(\forall
i<\lambda )$ $(F_i\subseteq\DD_\alpha )$.
\item[($*$)$_3$] (a) $\kappa <\lambda = {\rm cf}(\lambda )$
\item[{}] (b) if a Boolean algebra $\B$ satisfies the $\kappa$-c.c., $\DD$ a
$\lambda$-complete uniform filter on $\lambda$, $\theta ={\rm cf}(\theta
)<\kappa$ and for $i<\lambda$, $F_i$ is a decreasing sequence of elements of
$\B\backslash\{ 0\}$ of length $\theta$ $\underline{\rm then}$ for some
$X\in\DD^+$, $\bigcup_{i\in X}F_i$ belongs to some ultrafilter on $\B$.
\item[($*$)$_4$] (a) $\kappa <\lambda ={\rm cf}(\lambda )$
\item[{}] (b) if $\B$ is a $\kappa$-c.c.\ Boolean algebra of cardinality
$\ge\lambda$ $\underline{\rm then}$ $\lambda$ is a free caliber of $\B$.
\end{enumerate}
\smallskip

\n {\bf Proof}. Should be clear from the proof in \S14.{\hfill$\square_{\rm
14.25}$}
\medskip

\section{On irr: The invariant of the ultraproduct bigger than the
ultraproduct of invariants}
We solve here some of the questions of Monk \cite{M2} on the possibility that 
$$
{\rm inv}(\prod_{\zeta<\kappa}\B_\zeta/\DD)>\prod_{\zeta<\kappa}{\rm inv}(
\B_\zeta)/\DD.
$$
In 15.1--15.10A we deal with the irredundance number ${\rm irr}$ (getting
consistency of the above and solving \cite[Problem 26]{M2}). We then prove the
existence of such examples in ZFC (improving Ros{\l}anowski Shelah
\cite{RoSh:534}) for ${\rm inv}=s,\; {\rm hd},\; {\rm hL},\; {\rm Length}$
solving \cite[Problems 46, 51, 55, 22]{M2}, respectively. See more in
\cite{Sh:641}.
\medskip

\n {\bf 15.1 Hypothesis}. $\lambda =\lambda^{<\lambda}$, $n(*)<\omega$.
\medskip

\n {\bf 15.2 Definition}. {\it $P=P_\lambda^{n(*)}$ is the set of $p=(u,\B,
\bar F)=(u^p,\B^p,\bar F^p)$ such that
\begin{enumerate}
\item[(a)] $u\in [\lambda^+]^{<\lambda}$
\item[(b)] $\B$ is a Boolean algebra generated by $\{ x_\alpha :\alpha\in u\}$
\item[(c)] $\alpha\in u\Rightarrow x_\alpha\not\in \lan\{ x_\beta :\beta\in
u\cap\alpha\}\ran_\B$
\item[(d)] in $\B$, $\{ x_\alpha :\alpha\in u\}$ is $n(*)$-independent i.e.\
any nontrivial Boolean combination of $\le n(*)$ members of $\{ x_\alpha:
\alpha\in u\}$ is not zero (in $\B$).
\item[(e)] $\bar F=\lan F_\ell :\ell\le n(*)\ran$ and $F_{\ell +1}\subseteq
F_\ell$.
\item[(f)] $F_\ell$ is a family of functions from $\{ x_\alpha :\alpha\in u\}$
to $\{ 0,1\}$ respecting the equations holding in $\B$. Call the homomorphism
(from $\B$ to $\{0,1\}$) that $f$ induces, $\hat f$.
\item[(g)] if $f\in F_{\ell +1}$, $\ell< n(*)$ and $\alpha\in u$ then for some
$f'\in F_\ell$ we have
$$
\begin{array}{c}
f'\restriction (\alpha\cap u)=f\restriction (\alpha\cap u)\\
f'(\alpha )\not= f(\alpha ).
\end{array}
$$
\item[(h)] if $f:u\to \{ 0,1\}$ and $(\forall v\in [u]^{<\aleph_0})(f
\restriction u\in F_\ell )$ then $f\in F_\ell$.
\item[(i)] if $a\in \B\backslash\{ 0\}$ then for some $f\in F_0$, we have
$\hat f(a)=1$.
\end{enumerate}
The order is $p \le q$\quad $\underline{\rm iff}$
\begin{enumerate}
\item[$\alpha$)] $u^p\subseteq u^q$
\item[$\beta$)] $\B^p$ is a subalgebra of $\B^q$
\item[$\gamma$)] $F^p_\ell =\{ f\restriction u^p:f\in F^q_\ell\}$.
\end{enumerate}
Let $\name{\B} =$ the direct limit of $\{ \B^p:p\in \Name{G}_P\}$. } 
\medskip

\n {\bf Note}. We can ignore $\B^p$ at it is reconstructible from
$F^p_0$. Also clause (d) follows from the rest.
\medskip

\n {\bf 15.3 Notation}. We let $p\restriction\alpha =(u^p\cap\alpha ,\lan
x_\beta :\beta\in u^p\cap\alpha\ran_\B$, $\lan F_\ell\restriction\alpha
:\ell\le n(*)\ran)$ where $F_\ell\restriction\alpha =\{f\restriction\alpha
:f\in F_\ell\}$.
\medskip

\n {\bf 15.4 Fact}.. $(p\restriction\alpha )\le p$ for $p\in P$.
\medskip

\n {\bf 15.5 Fact}. In $P$, every increasing sequence of length $<\lambda$ has
a $\hbox{lub}$: essentially the union.
\medskip

\n {\bf Proof}. Trivial (use compactness and clause (h) of Definition 15.2).
{\hfill$\square_{\rm 15.5}$}
\medskip

\n {\bf 15.6 Fact}. For $\alpha <\lambda$, $\{ p\in P:\alpha\in u^p\}$ is
dense open.
\medskip

\n {\bf Proof}. If $p\in P$ let us define $q=(u^q, \B^q, F^q)$, $u^q=u^p\cup\{
a\}$, $\B^q$ is $\B^p$ if $\alpha\in u^p$, and is the free extension of $\B$
by $x_\alpha$ otherwise, $F^q_\ell =\{ f\in{}^{u^q}2: f\restriction u^p\in
F^p_\ell\}$.{\hfill$\square_{\rm 15.6}$}
\medskip

\n {\bf 15.7 Fact}. 1) If $p\in P$, $p\restriction\alpha\le q$ and
$u^q\subseteq\alpha$ $\underline{\rm then}$ $p,q$ are compatible.

2) $P$ satisfies the $\lambda^+$-c.c.\ and even in $\lambda^+$-Knaster.
\medskip

\n {\bf Proof}. 1) Let us define $r=(u^r, \B^r, \bar F^r)$ by:
$$
u^r=u^p\cup u^q,\quad F^r_\ell =\{ f:f\in {}^{u^r}2\mbox{ and }f\restriction
u^p\in F^p_\ell,\ f\restriction u^q\in F^q_\ell\}.
$$
Now
\begin{enumerate}
\item[($*$)$_1$] $F^p_\ell =F^r_\ell\restriction u^p$
\end{enumerate}
[Why? if $f\in F^p_\ell$, then $f\restriction\alpha =f\restriction (\alpha\cap
u)\in F_\ell^{p\restriction\alpha}$ but $p\restriction\alpha\le q$. Hence
there is $g\in F^q_\ell$ such that $f\restriction\alpha\subseteq g$, so $f\cup
g\in F^r_\ell$, $(f\cup q)\restriction u^p=f$, so $F^p_\ell\subseteq
F^r_\ell\restriction u^p$. The other direction holds by the choice of
$F^r_\ell$.]
\begin{enumerate}
\item[($*$)$_2$] $F^q_\ell =F^r_\ell\restriction u^p$
\end{enumerate}
[Why? similarly using 15.4.]
\begin{enumerate}
\item[($*$)$_3$] $F^r_{\ell +1}\subseteq F^r_\ell$
\end{enumerate}
[Why? as $F^p_{\ell +1}\subseteq F^p_\ell$, $F^q_{\ell +1}\subseteq F^q$.]
\begin{enumerate}
\item[($*$)$_4$] if $f\in F^r_{\ell +1}$, $\beta\in u^r$ then for some $g\in
F^r_{\ell +1}$ we have $f\restriction \beta\subseteq g$, $f(\beta )\not=
g(\beta )$.
\end{enumerate}
[Why? The proof splits into two cases:
\medskip

\n {\bf Case 1}. $\beta\in u^q$.

So $f\restriction\alpha\in F^q_{\ell +1}\restriction\alpha$ but $q\in P$ so
there is $g_0\in F^q_\ell$ such that $(f\restriction\alpha
)\restriction\beta\subseteq g_0$, $(f\restriction\alpha )(\beta )\not=
g_0(\beta)$ so $g_0\in F^q_\ell=F^p_\ell\restriction\alpha$ so
$g_0\restriction (u^p\cap\alpha )\in F_\ell^{p\restriction\alpha}$ so there is
$g_1$ such that
$$
g_0\restriction (u^p\cap\alpha )\subseteq g_1\in F^p_\ell.
$$

So $g_0\cup g_1\in F^r_\ell$ is as required.
\medskip

\n {\bf Case 2}. $\beta\not\in u^q$.

So $\beta\in u^p\backslash\alpha$. Now $f\restriction u^p\in F^p_{\ell +1}$
hence there is $f'\in F^p_\ell$ such that
$$
f'\restriction (u^p\cap \beta )=f\restriction (u^p\cap\beta ),\quad f'(\beta
)\not= f(\beta ).
$$
Now $f\restriction\alpha\in F^q_{\ell +1}$ hence $f\restriction\alpha\in
F^q_\ell$ hence
$$
(f\restriction\alpha )\cup f'\in F^r_\ell\;\;\mbox{is as
required}.
$$
By $F^r_\ell$ we can define $\B^r$ and is as required.]

2) Follows from (1).{\hfill$\square_{\rm 15.7}$} 
\medskip

\n {\bf 15.8 Claim}. If $k>2n(*)+1$, $\lan\delta_\ell:\ell <k\ran$ is
increasing, $\delta_\ell <\lambda$; we stipulate $\delta_k=\lambda^+$, for
$\ell <k$, $p_\ell\in P$, $p_\ell\restriction\delta_\ell=p^*$,
$u^{p_\ell}\subseteq\delta_{\ell +1}$ and for $\ell ,m<k$, ${\rm OP}_{u^{p_m},
u^{p_\ell}}:u^{p_\ell}\rightarrow u^{p_m}$ maps $p_\ell$ to $p_m$ (the natural
meaning ${\rm otp}(u^{p_\ell})={\rm otp}(u^{p_m})$ and
$$
F^{p_\ell}_n=\{f\circ {\rm OP}_{u^{p_\ell}, u^{p_m}}:f\in F^{p_m}_n\}
$$
so ${\rm OP}_{u^{p_\ell}, u^{p_m}}$ induces an isomorphism ${\rm OP}_{p_\ell,
 p_m}$ from $\B^{p_\ell}$ onto $\B^{p_m}$), $\underline{\rm then}$ there is
 $q\in P$ such that 
\begin{enumerate}
\item[(a)] $\displaystyle{\bigwedge_{m<k}p_{m}\le q}$
\item[(b)] if $b\in\B^{p_0}$ then $\B^q\vDash {\rm ``}b=\bigcup\limits_{u
\subseteq (0,k)\atop |u|>n(*)}\Big(\bigcap\limits_{m\in u}{\rm OP}_{p_m,p_0}
(b)\Big)$''.
\end{enumerate}
\smallskip

\n {\bf Proof}. 1) Let us define $q$: put $u^q=\bigcup\limits_{m<k}u^{p_m}$
and  
$$
\begin{array}{ll}
F^q_\ell=&\{f\in {}^{(u^q)}2:\\
& n(*)-\ell\ge |\{m\in [0,k):(\exists\alpha\in u^{p_0}\backslash u^{p^*})
[f({\rm OP}_{u^{p_m},u^{p_0}}(\alpha ))\not= f(\alpha)]\}|\\
&\mbox{and }\;\; f\restriction u^{p_m}\in F^{p_m}_\ell\;\;\mbox{for}\;\;
m<k\big\}.
\end{array}
$$
Now note
\begin{enumerate}
\item[($*$)$_1$] $F^{p_m}_\ell =F^q_\ell\restriction u^{p_m}$
\end{enumerate}
[Why? if $f\in F^q_\ell$ then $f\restriction u^{p_m}\in F^{p_m}_\ell$ by the
definition of $F^q_\ell$. If $f\in F^{p_m}_\ell$, for $m_1<k$ we let
$f_{m_1}=f\circ {\rm OP}_{u^{p_{m_1}},u^{p_m}}$, so $\bigcup_{m_1<k}f_{m_1}\in
F^q_\ell$ and we are done.]
\begin{enumerate}
\item[($*$)$_2$] if $f\in F^q_{\ell +1}$, $\alpha\in u^q$ then for some $g\in
F^q_\ell$
$$
g\restriction\alpha =f\restriction\alpha,\;\; g(\alpha )\not= f(\alpha )
$$
\end{enumerate}
[Why? if $\alpha\in u^{p_0}$ we have $f\restriction u^{p_0}\in F^{p_0}_{\ell
+1}$ and there is $g_0\in F^{p_0}_\ell$, such that $g_0\restriction\alpha
=f\restriction\alpha$, $g_0(\alpha )\not= f(\alpha )$. Let $g_m={\rm
OP}_{u^{p_0},u^{p_m}}\circ g_0$. Then $g=\bigcup_{m<k}g_m$ is as required.

If not, $\alpha\in u^{p_m}\backslash u^{p^*}$ for some $m>0$, so
$\alpha\ge\delta_m$ and $f\restriction u^{p_m}\in F^{p_m}_{\ell +1}$ so there
is $g\in F^{p_m}_{\ell +1}$, $g\restriction\alpha =f\restriction\alpha$,
$g(\alpha )\not= f(\alpha )$. Now $g^*=g\cup \Big(
f\restriction\big(\bigcup_{m_1<k\atop m_1\not= m} u^{p_{m_1}})$ is as
required.] 

\n So
\begin{enumerate}
\item[($*$)$_3$] $q\in P$ and $p_m\le q$.
\end{enumerate}
So (a) of the conclusion holds. By clause (i) of Definition 15.2
and the choice of $q$ also clause (b) holds.{\hfill$\square_{\rm 15..8}$}
\medskip

\n {\bf 15.9 Conclusion}. $\Vdash_{P^{n(*)}_x}$ ``$\Name{\B}$ is a Boolean
algebra generated by $\{ x_\alpha :\alpha <\lambda^+\}$, which is $n(*)$-free
hence ${\rm irr}_{n(*)}(\Name{\B})\ge\lambda^+$ but ${\rm
irr}_{2n(*)+1}(\Name{\B} )=\lambda$''.
\medskip

\n {\bf Proof}. Putting together the claims.{\hfill$\square_{\rm 15.9}$}
\medskip

\n {\bf 15.10 Conclusion}. If $\lambda =\lambda^{<\lambda}>\aleph_0$ and the
forcing notion $P$ is $P=\prod_n P^n_\lambda$ (where $P^n_\lambda$ is from
Definition 15.2) then
\begin{enumerate}
\item[($*$)] $P$ is a $\lambda$-complete $\lambda^+$-c.c.\ forcing notion, and
in $V^P$ for some Boolean algebras $\B_n$ $(n<\omega )$ we have
\begin{enumerate} 
\item[(a)] ${\rm irr}_n(\B_n)=\lambda^+$, ${\rm irr}_{2n+1}(\B_n)=\lambda$
\item[(b)] for $\DD$ a nonprincipal ultrafilter on $\omega$, $\lambda^+\le
{\rm irr}\big(\prod_{n<\omega}\B_n/\DD \big)$, $\prod_{n<\omega}({\rm irr}\;
\B_n)/\DD =\lambda^\omega/\DD =\lambda$.
\item[(c)] So ${\rm irr}\big(\prod_{n<\omega}\B_n/\DD\big)>\prod_{n<\omega}
{\rm irr}(\B_n)/\DD$
\end{enumerate} \end{enumerate}
\medskip

\n {\bf Proof}: The $\lambda^+$-c.c.\ follows from 15.7(2). The $\B_n$ are
from 15.2. The proof that $\Vdash_P$ ``${\rm irr}_n(\name{\B}_n)=\lambda^+$
but ${\rm irr}_{2n+2}(\name{\B})=\lambda'$ is like the proof of 15.9.

\n Concerning ${\rm irr}\big(\prod_{n<\omega}\name{\B}_n/\DD )=\lambda^+$ use
$x_\alpha^*=\lan x^n_\alpha :n<\omega\ran/\DD\in
\prod_{n<\omega}\name{\B}_n/\DD$.

~{\hfill$\square_{\rm 15.10}$}
\medskip

\n {\bf 15.10A Comment}: Surely in 15.9 we can fix exactly the $n$ such that
${\rm irr}_n(\B )=\lambda^+$, ${\rm irr}_{n+1}(\B )=\lambda$. For this it
suffices to demand in 15.8 that $k=n(*)+2$. Let $a_\ell\in\B^{p_\ell}$ (for
$\ell<k$) be such that ${\rm OP}_{p_\ell,p_0}(a_0)=a_\ell$ and replace (b) in
the conclusion by
\begin{description}
\item[(b)$'$] for some Boolean term $\tau$,
$$
\B^a\models \mbox{``}a_0=\tau(a_1,\ldots,a_{n(*)+1})\mbox{''}.
$$
\end{description}
In fact, 
$$
\tau(x_1,\ldots,x_{n(*)+1})=\bigcup_{i=1}^{n(*)}[\bigcap_{m=1}^{i-1}(x_m
\triangle x_{m+1})\cap(-(x_i\triangle x_{i+1}))].
$$
In the proof we let $u^q=\bigcup\limits_{m<k} u^{p_m}$ and
$$
\begin{array}{ll}
F^q_\ell=\{f\in {}^{(u^p)}2:&\mbox{for $m<k$ we have $f_m=:f\restriction
u^{p_m}$ belongs to $F^{p_m}_\ell$ and}\\
&\mbox{for some $i\in\{1,\ldots,n(*)-\ell\}$ we have:}\\
&j\in\{1,\ldots,i-1\}\quad\Rightarrow\quad [\hat{f}_j(a_j)=0\ \Leftrightarrow\
\hat{f}_{j+1}(a_{j+1})=1]\\
&\mbox{and }\quad\hat{f}_i(a_i)=\hat{f}_{i+1}(a_{i+1})=\hat{f}_0(a_0)\ \}
\end{array}
$$
(where for $f\in F^{p_m}_\ell$, $\hat{f}$ is the homomorphism from $\B^q$ into
$\{0,1\}$ which $f$ induces). 
\medskip

\n {\bf 15.11 Claim}. Assume
\begin{enumerate}
\item[(A)] $\lambda ={\rm tcf}\Big(\prod_{i<\delta}\lambda_i/J\Big)$
\item[(B)] $\bar\lambda =\lan\lambda_i :i <\delta\ran$ is a sequence of
regular cardinals $>|\delta|$
\item[(C)] $\lambda_i >$ max pcf$\{\lambda_j:j<i\}$, so necessarily $J^{\rm
bd}_\delta\subseteq J$
\item[(D)] $\lan A_\zeta :\zeta <\kappa\ran$ is a sequence of pairwise
disjoint members of $J^+$.
\item[(E)] $\DD$ is a uniform ultrafilter on $\kappa$.
\end{enumerate}
$\underline{\rm Then}$, we can find a Boolean algebra $\B_\zeta$ for $\zeta
<\kappa$ such that for ${\rm inv}\in \{s,hd,hL\}$ (see Monk \cite{M2})
\begin{enumerate}
\item[(a)] ${\rm inv}^+(\B_\zeta )\le \lambda$ so $\lambda =\chi^+\Rightarrow
{\rm inv} (\B_\zeta )\le \chi$ (moreover ${\rm inv}^+_2(\B_3)\le \lambda$; see
\cite{RoSh:534}.) 
\item[(b)] ${\rm inv}^+\Big(\prod_{\zeta <\kappa}\B_\zeta /\DD\Big)>\lambda$
(so if $\lambda=\chi^+$ then ${\rm inv}\Big(\prod_{\zeta <\kappa} \B_\zeta
/\DD \Big) \ge\lambda$)
\end{enumerate}
\smallskip

\n {\bf Proof}. Let $\bar\eta =\lan\eta_\alpha :\alpha <\lambda\ran$ be a
$<_J$-increasing cofinal sequence of members of $\prod_{i<\delta}\lambda_i$
such that 
$$
\zeta<\kappa\quad\Rightarrow\quad\lambda_\zeta>|\{\eta_\alpha\restriction
\zeta:\alpha<\lambda\}|
$$
(such $\bar{\eta}$ exists by \cite[II 3.5]{Sh:g}). We define $(\B_{\zeta
,i}^*\bar x^-_{\zeta,i},\bar x^+_{\zeta ,i})$ for $\zeta <\kappa$, $i<\delta$
as follows. Let $I_i=J^{\rm bd}_{\lambda_i}$, so $\bar x^-_{\zeta ,i}=\lan
x^-_{\zeta,i,\alpha}:\alpha<\lambda_i\ran$, $\bar x^+_{\zeta,i}=\lan
x^+_{\zeta,i,\alpha}:\alpha <\lambda_i\ran$.
\medskip

\n {\bf Case 1}. $i\not\in\cup\{ A_\e :\e\in [\zeta ,\kappa )\}$.

Let $\B_{\zeta ,i}$ be the Boolean algebra generated by $\{ x^-_{\zeta
,i,\alpha}, x^+_{\zeta ,i,\alpha}:\alpha <\lambda_i\}$ freely except that
$x^-_{\zeta ,i,\alpha}\le x^+_{\zeta ,i,\alpha}$, and $(x^+_{\zeta
,i,\alpha}-x^-_{\zeta ,i,\alpha})\cap (x^+_{\zeta ,i,\beta}-x^-_{\zeta
,i,\beta})=0$ when $\alpha <\beta <\lambda_i$.
\medskip

\n {\bf Case 2}. $i\in\cup\{ A_\e :\e\in [\zeta ,\kappa )\}$.

Let $\B_{\zeta ,i}$ be the Boolean algebra generated by $\{ x^-_{\zeta
,i,\alpha}, x^+_{\zeta ,i,\alpha}:\alpha <\lambda_i\}$ freely except that
$$
\alpha <\beta\Rightarrow x^-_{\zeta ,i,\alpha}\le x^+_{\zeta ,i,\alpha}\le
x^-_{\zeta ,i,\beta }\le x^+_{\zeta ,i,\beta}
$$
(e.g.\ $\B_{\zeta ,i}\subseteq\PP (\lambda_i)$, $x^-_{\zeta,i,\alpha}=
[0,4\alpha +1)$, $x^+_{\zeta ,i,\alpha}=[0,4\alpha +2))$.

Let $\B_\zeta$ be constructed as in 12.1, 12.3 from $\bar \lambda$, $\lan
I_i:i<\delta\ran$, $(\B_{\zeta ,i},\bar x^-_{\zeta ,i}, \bar x^+_{\zeta ,i})$
for $i<\delta$, and let $y^\zeta_\alpha$, $y^\zeta_\eta$ be as there.

Now ${\rm inv}^+\Big(\prod_{\zeta <\kappa} \B_\zeta /\DD \Big) >\lambda$ is
exemplified by $\lan y^*_\alpha :\alpha <\lambda\ran$ where $y^*_\alpha =\lan
y^\zeta_\alpha :\zeta <\kappa\ran /\DD$, because for $\alpha <\lambda$,
$u\subseteq\lambda\setminus\{\alpha\}$ finite, for some $\zeta^*<\kappa$, we
have $\beta\in u\ \Rightarrow\ {\ell g}(\eta_\alpha\cap\eta_\beta)\in\delta
\backslash\bigcup\limits_{\e\in [\zeta ,\kappa )}A_\e$, hence 
$$
\zeta\in [\zeta_{\alpha,\beta},\kappa)\quad\Rightarrow\quad \B_\zeta\vDash
y^\zeta_\alpha-\bigcap_{\beta\in u} y^\zeta_\beta>0.
$$
Hence $\{\zeta <\kappa :\B_\zeta \vDash y^\zeta_\alpha\cap y^\zeta_\beta =0\}
\supseteq [\zeta_{\alpha ,\beta}, \kappa )\in\DD$ and therefore
$\prod\limits_{\zeta<\kappa}\B_\zeta /\DD\vDash \mbox{``}y^*_\alpha
-\bigcap\limits_{\beta\in u} y^*_\beta >0\mbox{''}$.
\smallskip

Lastly ${\rm inv}^+_{(2)}(\B_\zeta)\le\lambda$ follows by 12.17(2) for
$\tau(x_0,x_1, x_2)=(x_1-x_0\cup x_2)$ with the variables permuted according
to the particular inv. {\hfill$\square_{\rm 15.11}$}
\medskip

\n {\bf 15.12 Claim}. Claim 15.11 holds for Length too.
\medskip

\n {\bf Proof}. We repeat the proof of 15.11, but in the definition of
$\B_{\zeta ,i}$ just interchange the two cases.
\medskip

\n {\bf Case 1}. $i\not\in\cup\{ A_\e :\e\in [\zeta ,\kappa ]\}$.

Let $\B_{\zeta ,i}$ be as $\B_{\zeta ,i}$ in case 2 in the proof of 15.11.
\medskip

\n {\bf Case 2}. $i\in \cup\{ A_\e :\e\in [\zeta ,\kappa )\}$.
\medskip

As in Case 1 in the proof of 15.11 or just let $\B_{i,\zeta}$ be generated by
$\{x^-_\alpha , x^+_\alpha :\alpha <\lambda_i\}$, $\{ x^-_{\zeta ,i,\alpha},
x^+_{\zeta , i,\alpha}:\alpha <\lambda_i\}$ freely except $x^-_{\zeta
,i,\alpha}\le x^+_{\zeta ,i,\alpha}$.

Now for $\alpha <\beta <\lambda$, letting $i(\alpha ,\beta )={\rm Min}\{
i:\eta_\alpha (i)\not=\eta_\beta (i)\}$ and $\zeta_{\alpha ,\beta}={\rm
Min}\{\zeta :$ $i(\alpha ,\beta )\not\in\cup\{A_\e :\e\in [\zeta ,\kappa )\}$
we have
$$
\zeta\in [\zeta_{\alpha ,\beta},\kappa )\Rightarrow\B_{\zeta ,i}\vDash
\mbox{``} y^\zeta_\alpha \le y^\zeta_\beta\;\; {\rm or}\;\; y^\zeta_\beta <
y^\zeta_\alpha\mbox{''}.
$$
hence
$$
\prod_{\zeta <\kappa} \B_\zeta\vDash \mbox{``}y^*_\alpha <y^*_\beta\;\; {\rm
or}\;\; y^*_\beta <y^*_\alpha",
$$
where $y^*_\alpha =\lan y^\zeta_\alpha :\zeta <\kappa\ran /\DD$.
\medskip

As for ${\rm Length}^+(\B_{\zeta})\le \lambda$, it is by 12.7(1).
{\hfill$\square_{\rm 15.12}$}
\medskip

\n {\bf 15.13 Conclusion}: 1) If $\DD$ is a uniform ultrafilter on $\kappa$,
then for a class of cardinals $\chi =\chi^\kappa$ and Boolean algebras $\B_i$
for $i<\kappa$ such that, for ${\rm inv}\in \{s,hL,hd\}$ we have:
\begin{enumerate}
\item[(a)] ${\rm inv}(\B_i)\le \chi$ hence $\prod_{i<\kappa}{\rm
inv}(\B_i)\le\chi$

\n or

 \item[(b)] ${\rm inv}\Big(\prod_{i<\kappa}\B_i/\DD\Big)=\chi^+$
\end{enumerate}

2) Similarly with inv $=$ Length.
\medskip

\n {\bf Proof}: Let $\chi$ be any strong limit singular cardinal of cofinality
$>\kappa$. So by \cite[VIII, \S1]{Sh:g} we can find $\lan\lambda_i:i<{\rm
cf}(\chi)\ran$, strictly increasing sequence of regular cardinals $<\chi$
with
${\rm tcf}\Big(\prod_{i<{\rm cf}(\chi)}\lambda_i/J^{\rm bd}_{{\rm cf}(\chi)}
\Big)=\chi^+$. Without loss of generality $\prod_{i<j}\lambda_i<\lambda_j$ and
let for $i<\kappa$, $A_i\{\alpha\kappa +i:\alpha<{\rm cf}(\chi )\}$. So we can
apply 15.11 (for part (1)) or 15.12 (for part (2)).{\hfill$\square_{\rm
15.13}$} 

\n {\bf Remark}: For cellularity similar results hold (in ZFC), i.e.~$c(\B_n)
\leq \lambda$, $c(\prod_{n<\omega}\B_n)>\lambda$, see on it in Monk
[\cite{M2} p.\ 61--62]; by
\cite[III 4.11, p.181, 4.12]{Sh:g} so this applies to $\lambda=\mu^+$ for
$\lambda>\aleph_1$ by \cite[II, 4.1]{Sh:g}, \cite{Sh:572}, to $\lambda$
inaccessible not Mahlo by \cite[III 4.8(2) p.177]{Sh:g} and for many Mahlo
cardinals (see \cite[III 4.10A p.178]{Sh:g}. For incomparability number ({\rm
Inc}) similar results are proved ``almost in ZFC'', see \cite{Sh:462}.

\bibliographystyle{literal-plain}
\bibliography{listb,lista,listc,liste,listx}

\end{document}